\documentclass[11pt,reqno]{amsart}
\usepackage[foot]{amsaddr}
\usepackage{amsmath,dsfont}
\usepackage{verbatim}
\usepackage{amssymb}
\usepackage{amsthm}
\usepackage{cite}
\usepackage{pstricks,pst-node,pst-text,pst-3d}
\usepackage{color}
\usepackage{eucal}
\usepackage{tikz}
\usetikzlibrary{matrix,arrows,decorations.markings}
\usepackage{amscd}
\usepackage{bm}
\usepackage{yfonts}
\usepackage{url}
\usepackage[outline]{contour}
\usetikzlibrary{decorations.pathmorphing}
\usetikzlibrary{shapes.geometric,tqft}
\usepackage{tqft}
\tikzset{tqft/use   nodes=false}
\usepackage{datetime2}

\DeclareSymbolFont{rsfscript}{OMS}{rsfs}{m}{n}
\DeclareSymbolFontAlphabet{\mathrsfs}{rsfscript}
\def\malcev{\mathop{\hbox{$\bigcirc$\kern-9pt\raise1pt\hbox{\scriptsize$m$}\kern1.5pt}}}
\def\malcevk{\mathop{\hbox{$\bigcirc$\kern-8.5pt\raise-0.5pt\hbox{$k$}\kern1.5pt}}}

\DeclareMathOperator{\Id}{Eq}
 
\def\softd{{\leavevmode\setbox1=\hbox{d}\hbox
            to 1.15\wd1{d\kern-0.2ex{\char039}\hss}}}    
\def\softl{l\kern-0.3ex\raise0.1ex\hbox{'}\kern-0.3ex}   
\numberwithin{equation}{section}
\newtheorem{Thm}{Theorem}[section]
\newtheorem{Prop}[Thm]{Proposition}
\newtheorem{Lemma}[Thm]{Lemma}
\newtheorem{Res}[Thm]{Result}
\newtheorem{Cor}[Thm]{Corollary}

\DeclareSymbolFont{rsfscript}{OMS}{rsfs}{m}{n}
\DeclareSymbolFontAlphabet{\mathrsfs}{rsfscript}

\theoremstyle{remark}
\newtheorem{Rmk}{Remark}
\newtheorem*{Problem}{Problem}

\newtheorem{Exm}{Example}
\newtheorem{Def}{Definition}[section]

\DeclareMathOperator{\var}{\mathrm{var}}
\DeclareMathOperator{\occ}{occ}
\DeclareMathOperator{\con}{con}
\newcommand{\inc}{\operatorname{inc}}
\newcommand{\rk}{\operatorname{rk}}

\newcommand{\Ob}{\operatorname{Ob}}
\newcommand{\iv}{\operatorname{iv}}
\newcommand{\ov}{\operatorname{ov}}
\newcommand{\ve}{\operatorname{v}}
\newcommand{\g}{\operatorname{inc}}
\newcommand{\cir}{\operatorname{cir}}
\newcommand{\ub}{\operatorname{ub}}
\newcommand{\lb}{\operatorname{lb}}

\newcommand{\Rms}{Rees matrix semigroup}

\def\til#1{\widetilde {#1}}
\def\ol#1{\overline{#1}}

\def\mS{\mathcal S}

\def\al{\alpha}

\def\vp{\varepsilon}

\begin{document}

\title[Combinatorial skeletons]{Combinatorial skeletons\\ of $2$-cobordism and annular categories \\ with  applications to equational logic}
\author{K.~Auinger and M.~V.~Volkov}
\address{Fakult\"at f\"ur Mathematik, Universit\"at Wien,
Oskar-Morgenstern-Platz 1, A-1090 Wien, Austria}
\email{karl.auinger@univie.ac.at}
\address{620075 Yekaterinburg, Russia}
\email{m.v.volkov@urfu.ru}

\date{}
\begin{abstract}
We introduce a complete set of combinatorial data that encode the category $2\mathfrak{Cob}$ of all $2$-cobordisms. As an application, we show that the local monoids of $2\mathfrak{Cob}$ do not have finitely axiomatizable equational theories. As yet another application, we construct a von-Neumann-regular extension of this category. Similar results are provided for the topological annular category and various quotients of the latter, like the affine Temperley--Lieb category.
\end{abstract}
\maketitle

\tableofcontents

\section{Introduction}
The generic name `\emph{diagram monoids}' refers to a broad class of algebraic-combinatorial objects arising in various branches of mathematics and statistical physics. The very first diagram monoids were the ones introduced by  Brauer~\cite{Br37} back in 1937 in his studies of the linear representations of orthogonal and symplectic groups. The next series of diagram monoids was related to linear algebras introduced by Temperley and Lieb~\cite{TL71} in 1971 in connection to certain graph-theoretic issues in statistical mechanics; in 1990, the same species of monoids arose in Kauffman's work on knots \cite{kauffmanIsotopy}.

A very general family of diagram monoids (of which the ones mentioned above are specializations) consists of so-called partition monoids. The latter serve as bases of certain semigroup algebras relevant in statistical mechanics and representation theory, the so-called \emph{partition algebras}. Partition algebras were discovered and studied in depth by Martin \cite{Martin91,martin,Martin96,Martin00} and, independently, by Jones \cite{Jones94} in the context of  statistical mechanics; their remarkable role in representation theory is nicely presented in the introduction of~\cite{halverson}.

It is natural to consider partition monoids as the local monoids of partition categories, also widely studied and having various applications; see, e.g., \cite{BrundanVargas,comes,kujawatharp,penneys,HeTub2025} and references therein.
Here we present a novel use of partition categories: we analyze various categories with a topological flavor, seeking their characterization with discrete data, and partition categories turn out to form a crucial ingredient of these data. A typical object of our study is the category $2\mathfrak{Cob}$ of all $2$-cobordisms.

Whilst we hope that our characterization results might be useful for the areas where topological categories come from (such as low-dimensional topology and topological quantum field theory), the applications we selected for the present paper are of algebraic kind. For algebraists, diagram monoids have evolved into an appealing field of research lately since, first, they are `fresh' and natural, and second, their study has revealed that their structural properties are both neat and quite nontrivial. All sorts of standard objects of monoid theory have been investigated for diagram monoids, including idempotents~\cite{DE17,dolinkaetal,Dea19,Ea19a,EF12}, ideals~\cite{east_gray,LF06}, Green's relations~\cite{DE18}, congruences~\cite{congruencelattice,ER22a,ER22b,ER22c}, endomorphisms~\cite{maz:endomorphisms}, presentations~\cite{Ea11,Ea18,Ea19b,KM06,CarEasFre2025}, etc. Since the posting to arXiv in February 2020 of a preliminary version \cite{{av2020}} of the present article many more papers on the subject have been posted to arXiv or published in journals. It is virtually impossible to list them all; for a more complete overview on the literature, the reader is referred to the paper~\cite{FreSteTub2025} and the references cited therein. 

An interesting feature of diagram monoids is that even infinite monoids in this class admit nontrivial identities. This makes them rather exceptional since, as a rule, infinite noncommutative monoids that naturally arise in mathematics --- monoids of transformations of an infinite set, monoids of relations on an infinite domain, monoids of matrices over an infinite ring, etc. --- all are `too big' to satisfy any identity apart from the associativity law and its consequences. In~\cite{auingeretal}, it has been shown that identities holding in the monoids introduced in~\cite{kauffmanIsotopy} and several other, similarly defined diagram monoids are rich and complicated so that they do not admit any finite axiomatization. The present paper continues this line of research. Here we establish similar facts for infinite partition monoids, and our characterization results allow us to further extend these facts to local monoids of categories such as $2\mathfrak{Cob}$, topological annular category and affine Temperley--Lieb category.

As yet another application, we construct von-Neumann-regular extensions for the aforementioned categories. We believe that these regular extensions are of independent interest in that they contribute to the understanding of the algebraic structure of the original categories and, at the same time, they maintain certain key properties of the latter. In particular, this applies to the properties of identities of local monoids on which we focus in this paper. 

The paper is organized as follows. In Sections \ref{partition categories} and \ref{sec:annular type} we discuss the category of $2$-cobordisms and various versions of the annular category and present them in a way suitable for our purpose. We also introduce and discuss von-Neumann-regular extensions of these categories. In Section \ref{preliminaries} we collect prerequisites concerning identities of algebraic structures and semigroup theoretic constructions. Finally, in Section \ref{finite basis problem} we show that the identities of the local monoids of all aforementioned categories do not admit finite axiomatization (with a few exceptions of small degree in which cases we provide an explicit finite list of axioms). In addition, it is shown that analogous results are true if one considers identities involving along with multiplication certain naturally defined unary operations. The appendix contains a table summarizing the results of Section~\ref{finite basis problem}.

A fair effort has been made to keep the paper reasonably self-contained. It uses no results of category theory or topological quantum field theory (except for the Ti{\textsl k}Z package \texttt{tqft}), and all the necessary background from equational logic and semigroup theory is provided in Section~\ref{preliminaries}.

\section{Partition categories and their regular extensions}\label{partition categories} In this section we start with recalling the central notions concerning the partition and $2$-cobordism categories and presenting them in a notation which will be convenient for our later purpose. Then we define von-Neumann-regular extensions of these categories.

A short remark on terminology and conventions seems to be appropriate. All categories studied in this paper are \textsl{small} and are treated like algebraic structures rather than as a tool for classifying mathematical objects. (That is why no acquaintance with results and constructions of category theory is necessary to understand this paper.) \emph{Morphisms} are usually structure preserving maps between (carriers of) structures (including categories!). In contrast, morphisms between objects of a category are referred to as \emph{arrows}. Morphisms as well as arrows are composed from left to right; an argument of a morphism is usually written on the left. Arguments of other types of functions (for example labeling functions) are usually written on the right.

A \emph{small category} $\mathfrak{C}$ consists of a set of objects $\Ob\mathfrak{C}$ and a set of arrows, the latter being the disjoint union of sets of arrows $a\to b$ between objects $a,b\in\Ob\mathfrak{C}$ over all
pairs $(a,b)$ of objects. The set of all arrows $a\to b$ for two objects $a,b\in\Ob\mathfrak{C}$ is denoted by $\mathfrak{C}(a,b)$. For all $a,b,c\in\Ob\mathfrak{C}$, every arrow $\alpha\in\mathfrak{C}(a,b)$ and every arrow $\beta\in\mathfrak{C}(b,c)$ can be \emph{composed}, producing the arrow $\alpha\beta\in\mathfrak{C}(a,c)$. The composition is subject to the following two axioms.
\begin{description}
  \item[Associativity] for all arrows $\alpha\in\mathfrak{C}(a,b)$, $\beta\in\mathfrak{C}(b,c)$, $\gamma\in\mathfrak{C}(c,d)$,
  \[
  (\alpha\beta)\gamma=\alpha(\beta\gamma).
  \]
  \item[Existence of identities] for every $a\in\Ob\mathfrak{C}$, there exists an arrow $\varepsilon_a\in\mathfrak{C}(a,a)$ (called the \emph{local identity} at object $a$) such that
  \[
  \beta\varepsilon_a=\beta\ \text{ and } \varepsilon_a\alpha=\alpha
  \]
  for all $\beta\in\mathfrak{C}(b,a)$ and $\alpha\in\mathfrak{C}(a,b)$.
\end{description}
The two axioms imply that for every $a\in\Ob\mathfrak{C}$, the set $\mathfrak{C}(a,a)$ constitutes a monoid under composition, in which $\varepsilon_a$ serves as the neutral element. This monoid is called the \emph{local monoid} of $\mathfrak{C}$ at object $a$.

\subsection{Undeformed and deformed partition categories}\label{partition category}
Throughout, we let $\mathbb{N}_0$ be the set of all non-negative integers. For any $n\in \mathbb{N}_0$ we set $[n]:=\{1,2,\dots, n\}$; in particular, $[0]=\varnothing$.

We shall be widely concerned with \emph{disjoint unions} $X\sqcup Y$ and $X_1\sqcup\cdots\sqcup X_n$ of not necessarily disjoint sets $X$, $Y$ and $X_1,\dots, X_n$. In order to avoid too cumbersome notation this will be understood as follows. Let $X_1,\dots, X_n$ be sets; then $$U=X_1\sqcup \cdots \sqcup X_n$$ is the union of pairwise disjoint copies of the sets $X_i$ such that, in addition, the blocks $X_i$ of $U$ are linearly ordered as $X_1>X_2>\cdots >X_n$. (This could be made precise in terms of relational structures or by means of category theory, however we would prefer to avoid technical details.) This implies that for $X\ne Y$, $X\sqcup Y$ {\sl is not the same as} $Y\sqcup X$ (in contrast to $X\cup Y=Y\cup X$). The block $X_i$ then is the $i$th \emph{component} of $U$. There is a unique first (leftmost) component of $U$ which will be denoted $\lambda U$ and likewise a unique last (rightmost) component, to be denoted $U\!\varrho$. Note that some of the sets $X_i$ may be equal. In particular, if $X=Y$ then in $X\sqcup X$ we have a first component $X$ and a second component $X$. There is a standard way to make the sets $X_i$ pairwise disjoint, namely by replacing each element $x\in X_i$ by the ordered pair $(x,i)$. However, we would like to leave the method of distinguishing equal elements of distinct blocks $X_i$ unspecified. In particular, also for $n> 2$  we would like to consider $\lambda U\sqcup U\!\varrho=X_1\sqcup X_n$ as a subset of $U=X_1\sqcup\cdots \sqcup X_n$.

Let $m,n\in \mathbb{N}_0$; an \emph{$(m,n)$-partition} $\alpha$, or a \emph{partition} $\alpha\colon[m]\leadsto  [n]$ is a set partition of $[m]\sqcup [n]$. Throughout, any set partition is identified with its corresponding equivalence relation whenever this seems to be convenient. For $m,n\in \mathbb{N}_0$, let $\mathfrak{P}(m,n)$ be the set of all partitions $[m]\leadsto  [n]$. In this context, the elements of $[m]$ (first component of $[m]\sqcup [n]$) will usually be referred to as \emph{incoming vertices}, while those of $[n]$ as \emph{outgoing vertices} of $\alpha$. More generally, we sometimes will refer to the elements of a disjoint union $[n_0]\sqcup\cdots\sqcup [n_k]$ as its \emph{vertices}.

We shall define a partial composition of set partitions of this form.
Let $k\ge 1$, $n_0,\dots, n_k\in \mathbb{N}_0$ and $\alpha_i\in \mathfrak{P}(n_{i-1},n_i)$, $i=1,\dots,k$. Let
$$U=U(n_0,\dots,n_k)=[n_0]\sqcup[n_1]\sqcup \cdots\sqcup[n_k].$$
Every $\alpha_i$ is a symmetric and transitive relation on $U$; the union $\alpha_1\cup\cdots \cup\alpha_k$ is a reflexive and symmetric relation on $U$; the transitive closure
\begin{equation}\label{eq:partition composition}
(\alpha_1\cup\cdots\cup \alpha_k)^t    
\end{equation}
is an equivalence relation on $U$. We set
\begin{equation*}
\alpha_1\cdots\alpha_k:=(\alpha_1\cup\cdots\cup \alpha_k)^t\vert_{{\lambda U}\sqcup {U\!\varrho}},
\end{equation*}
considered as a partition $\lambda U\leadsto  U\!\varrho$. It is well known that for $k,\ell,m,n\in \mathbb{N}_0$ and $\alpha\colon[k]\leadsto  [\ell]$, $\beta\colon[\ell]\leadsto  [m]$, and $\gamma\colon[m]\leadsto  [n]$ we get
$$(\alpha\beta)\gamma=\alpha\beta\gamma=\alpha(\beta\gamma)\colon[k]\leadsto  [n].$$

In the following illustrative example, two partitions $\alpha\colon[6]\leadsto [8]$ (the incoming vertices are drawn at the top, the outgoing at the bottom, and partition blocks are encircled with closed curves)
\begin{center}
\begin{tikzpicture}
[scale=1.3]
\foreach \x in {1,2,...,6}\draw (\x,0) node {$\x$};
\foreach \x in {1,2,...,8}\filldraw (\x,-1) node {$\x$};
\draw[rounded corners=5pt](0.9,0)--(0.9,0.15)--(2.1,0.15)-- (2.1,-0.15)
 .. controls (1.7,-0.5) and (1.3,-0.5).. (0.9,-0.15)--(0.9,0);
\draw[rounded corners=5pt] (3,0.15) -- (3.15,0.15) -- (3.15,-1.15) -- (0.9,-1.15)-- (0.9,-0.85)..controls (2.85,-0.85) and (2.85,0.15) ..(3,0.15);
\draw[rounded corners=5pt] (3.9,0)--(3.9,0.15)--(6.1,0.15)--(6.1,-0.15) .. controls (5.5,-0.5) and (4.5,-0.5) .. (3.9,-0.15)--(3.9,0);
\draw[rounded corners=5pt] (4.9,-1)--(4.9,-1.15)--(7.1,-1.15)--(7.1,-0.85) .. controls (6.5,-0.5) and (5.5,-0.5) .. (4.9,-0.85)--(4.9,-1);
\draw[rounded corners=5pt] (3.85,-1)--(3.85,-0.85)--(4.15,-0.85)--(4.15,-1.15)--(3.85,-1.15)--(3.85,-1);
\draw[rounded corners=5pt] (7.85,-1)--(7.85,-0.85)--(8.15,-0.85)--(8.15,-1.15)--(7.85,-1.15)--(7.85,-1);
\end{tikzpicture}
\end{center}
{\parindent=0pt and $\beta\colon[8]\leadsto [9]$}\newline
\begin{center}
\begin{tikzpicture}
 [scale=1.3]
\foreach \x in {1,2,...,8}\draw (\x,0) node {$\x$};
\foreach \x in {1,2,...,9}\draw (\x,-1) node {$\x$};
\draw[rounded corners=5pt] (1.5,0.15)--(2.1,0.15)--(2.1,-0.15) .. controls (1.5,-0.5) .. (1.1,-1.15) --(0.9,-1.15)--(0.9,0.15) --(1.5,0.15);
\draw[rounded corners=5pt] (3.1,-0.15)--(3.1,-1.15)--(1.9,-1.15)--(1.9,-0.85).. controls (2.5,-0.5) .. (2.9,0.15)--(3.1,0.15)--(3.1,-0.15);
\draw[rounded corners=5pt] (4.1,0.15)--(7.1,0.15)--(7.1,-0.15).. controls (6,-0.5) and (5,-0.5) .. (3.9,-0.15)--(3.9,0.15)--(4.1,0.15);
\draw[rounded corners=5pt] (4.1,-1.15)--(6.1,-1.15)--(6.1,-0.85).. controls (5.5,-0.5) and (4.5,-0.5)..  (3.9,-0.85)--(3.9,-1.15)--(4.1,-1.15);
\draw[rounded corners=5pt](7.1,-1.15)--(9.1,-1.15)--(9.1,-0.85)..controls (8.5,-0.5) .. (8.1,0.15)--(7.9,0.15).. controls (7.5,-0.5) ..(6.9,-0.85)--(6.9,-1.15)--(7.1,-1.15);
\end{tikzpicture}\end{center}
\vskip0.2cm{\parindent=0pt
are being composed. Forming the relation $\alpha\cup\beta$ on $[6]\sqcup[8]\sqcup[9]$, we get}\newline
\begin{center}\begin{tikzpicture}
 [scale=1.3]
\foreach \x in {1,2,...,6}\draw (\x,0) node {$\x$};
\foreach \x in {1,2,...,8}\draw (\x,-1) node {$\x$};
\draw[rounded corners=5pt](0.9,0)--(0.9,0.15)--(2.1,0.15)-- (2.1,-0.15) .. controls (1.7,-0.5) and (1.3,-0.5).. (0.9,-0.15)--(0.9,0);
\draw[rounded corners=5pt] (3,0.15) -- (3.15,0.15) -- (3.15,-1.15) -- (0.9,-1.15)-- (0.9,-0.85)..controls (2.85,-0.85) and (2.85,0.15) ..(3,0.15);
\draw[rounded corners=5pt] (3.9,0)--(3.9,0.15)--(6.1,0.15)--(6.1,-0.15) .. controls (5.5,-0.5) and (4.5,-0.5) .. (3.9,-0.15)--(3.9,0);
\draw[rounded corners=5pt] (4.9,-1)--(4.9,-1.15)--(7.1,-1.15)--(7.1,-0.85) .. controls (6.5,-0.5) and (5.5,-0.5) .. (4.9,-0.85)--(4.9,-1);
\draw[rounded corners=5pt] (3.85,-1)--(3.85,-0.85)--(4.15,-0.85)--(4.15,-1.15)--(3.85,-1.15)--(3.85,-1);
\draw[rounded corners=5pt] (7.85,-1)--(7.85,-0.85)--(8.15,-0.85)--(8.15,-1.15)--(7.85,-1.15)--(7.85,-1);
\foreach \x in {1,2,...,9}\draw (\x,-2) node {$\x$};
\draw[rounded corners=5pt] (1.5,-0.85)--(2.1,-0.85)--(2.1,-1.15) .. controls (1.5,-1.5) .. (1.1,-2.15) --(0.9,-2.15)--(0.9,-0.85) --(1.5,-0.85);
\draw[rounded corners=5pt] (3.1,-1.15)--(3.1,-2.15)--(1.9,-2.15)--(1.9,-1.85).. controls (2.5,-1.5) .. (2.9,-0.85)--(3.1,-0.85)--(3.1,-1.15);
\draw[rounded corners=5pt] (4.1,-0.85)--(7.1,-0.85)--(7.1,-1.15).. controls (6,-1.5) and (5,-1.5) .. (3.9,-1.15)--(3.9,-0.85)--(4.1,-0.85);
\draw[rounded corners=5pt] (4.1,-2.15)--(6.1,-2.15)--(6.1,-1.85).. controls (5.5,-1.5) and (4.5,-1.5)..  (3.9,-1.85)--(3.9,-2.15)--(4.1,-2.15);
\draw[rounded corners=5pt](7.1,-2.15)--(9.1,-2.15)--(9.1,-1.85)..controls (8.5,-1.5) .. (8.1,-0.85)--(7.9,-0.85).. controls (7.5,-1.5) ..(6.9,-1.85)--(6.9,-2.15)--(7.1,-2.15);
\end{tikzpicture}.\end{center}
\vskip0.2cm{\parindent=0pt Passing to the transitive closure $(\alpha\cup\beta)^t$ gives}\newline
\begin{center}\begin{tikzpicture}
 [scale=1.3]
\foreach \x in {1,2,...,6}\draw (\x,0) node {$\x$};
\foreach \x in {1,2,...,8}\draw (\x,-1) node {$\x$};
\draw[rounded corners=5pt](0.9,0)--(0.9,0.15)--(2.1,0.15)-- (2.1,-0.15) .. controls (1.7,-0.5) and (1.3,-0.5).. (0.9,-0.15)--(0.9,0);
\draw[rounded corners=5pt] (3,0.15) -- (3.15,0.15) -- (3.15,-2.15) -- (0.9,-2.15) -- (0.9,-1.15)-- (0.9,-0.85)..controls (2.85,-0.85) and (2.85,0.15) ..(3,0.15);
\draw[rounded corners=5pt] (3.9,0)--(3.9,0.15)--(6.1,0.15)--(6.1,-0.15) .. controls (5.5,-0.5) and (4.5,-0.5) .. (3.9,-0.15)--(3.9,0);
\foreach \x in {1,2,...,9}\draw (\x,-2) node {$\x$};
\draw[rounded corners=5pt] (4.1,-0.85)--(7.1,-0.85)--(7.1,-1.15) -- (3.9,-1.15)--(3.9,-0.85)--(4.1,-0.85);
\draw[rounded corners=5pt] (4.1,-2.15)--(6.1,-2.15)--(6.1,-1.85).. controls (5.5,-1.5) and (4.5,-1.5)..  (3.9,-1.85)--(3.9,-2.15)--(4.1,-2.15);
\draw[rounded corners=5pt](7.1,-2.15)--(9.1,-2.15)--(9.1,-1.85)..controls (8.5,-1.5) .. (8.1,-0.85)--(7.9,-0.85).. controls (7.5,-1.5) ..(6.9,-1.85)--(6.9,-2.15)--(7.1,-2.15);
\end{tikzpicture},\end{center}
\vskip0.2cm{\parindent=0pt and restricting to $[6]\sqcup[9]$ finally gives $\alpha\beta\colon[6]\leadsto [9]$}
\newline
\begin{center}\begin{tikzpicture}
 [scale=1.3]
\foreach \x in {1,2,...,6}\draw (\x,0) node {$\x$};
\foreach \x in {1,2,...,9}\draw (\x,-1) node {$\x$};
\draw[rounded corners=5pt](0.9,0)--(0.9,0.15)--(2.1,0.15)-- (2.1,-0.15) .. controls (1.7,-0.5) and (1.3,-0.5).. (0.9,-0.15)--(0.9,0);
\draw[rounded corners=5pt] (3,0.15) -- (3.15,0.15) -- (3.15,-1.15) -- (0.9,-1.15)-- (0.9,-0.85)..controls (2.85,-0.85) and (2.85,0.15) ..(3,0.15);
\draw[rounded corners=5pt] (3.9,0)--(3.9,0.15)--(6.1,0.15)--(6.1,-0.15) .. controls (5.5,-0.5) and (4.5,-0.5) .. (3.9,-0.15)--(3.9,0);
\draw[rounded corners=5pt] (4.1,-1.15)--(6.1,-1.15)--(6.1,-0.85).. controls (5.5,-0.5) and (4.5,-0.5)..  (3.9,-0.85)--(3.9,-1.15)--(4.1,-1.15);
\draw[rounded corners=5pt](7.1,-1.15)--(9.1,-1.15)--(9.1,-0.85)..controls (8.5,-0.5) and (7.5,-0.5)..(6.9,-0.85)--(6.9,-1.15)--(7.1,-1.15);
\end{tikzpicture}.\end{center}
\vskip0.3cm
We have thus already defined the (undeformed) \emph{partition category} $\mathfrak{P}$: it is the (small) category whose set of objects $\Ob\mathfrak{P}$ is the set of non-negative integers $\mathbb{N}_0$, and, for $m,n\in\Ob\mathfrak{P}$ the set of arrows $m\to n$ is $\mathfrak{P}(m,n)$,  endowed with the composition $\mathfrak{P}(\ell,m)\times \mathfrak{P}(m,n)\to \mathfrak{P}(\ell,n)$, as defined above. The local identity at object $n$ is the partition of  $[n]\sqcup[n]$ into $n$ two-element blocks of the form $\{i$th incoming vertex, $i$th outgoing vertex$\}$.

In the formation of the product $\alpha\beta\colon[k]\leadsto  [m]$ for $\alpha\colon[k]\leadsto  [\ell]$ and $\beta\colon[\ell]\leadsto  [m]$, it may happen that some block of the partition $(\alpha\cup\beta)^t$ of $[k]\sqcup[\ell]\sqcup [m]$ contains only vertices of the second component $[\ell]$, that is, only \emph{intermediary} vertices. (This situation occurs in the above example.) Let $\mathrm{b}(\alpha,\beta)$ be the number of such \emph{dead blocks}. By definition, such dead blocks are ignored in the category $\mathfrak{P}$. However, they are essential for the \emph{deformed} partition category $\mathfrak{P}^{\mathrm{d}}$, which we define next.

The set of objects $\Ob\mathfrak{P}^{\mathrm{d}}$ is again the set of non-negative integers $\mathbb{N}_0$. For $m,n\in \mathbb{N}_0$ the set of arrows $m\to n$ is the set of all pairs $(\alpha,s)$ where $\alpha\in \mathfrak{P}(m,n)$ and $s\in \mathbb{N}_0$. The composition of two composable arrows is defined by
\begin{equation}\label{composition deformed partition category}
(\alpha,s)(\beta,t)=(\alpha\beta,s+t+\mathrm{b}(\alpha,\beta)).
\end{equation}
It is known \cite{martin} that for composable partitions $\alpha,\beta,\gamma$ the equality
\begin{equation}\label{associativity deformed}
\mathrm{b}(\alpha,\beta)+\mathrm{b}(\alpha\beta,\gamma)=\mathrm{b}(\alpha,\beta\gamma)+\mathrm{b}(\beta,\gamma)
\end{equation}
holds from which associativity of (\ref{composition deformed partition category}) follows.

The local monoids at object $n$ of the categories $\mathfrak{P}$ and $\mathfrak{P}^{\mathrm{d}}$ will be denoted $\mathfrak{P}_n$ and $\mathfrak{P}^{\mathrm{d}_n}$, respectively. Observe that the \emph{partition monoid} $\mathfrak{P}_n$ is finite for each $n\in\mathbb{N}_0$ while the \emph{deformed} (or \emph{twisted}) \emph{partition monoid} $\mathfrak{P}^{\mathrm{d}_n}$ is infinite for each $n>0$.

\begin{Rmk}
We mentioned in the introduction that the `early' diagram monoids considered by Brauer~\cite{Br37} are specializations of the partition monoids. Now we can make this precise: for each $n>0$, the \emph{Brauer monoid} $\mathfrak{B}_n$ is the submonoid of $\mathfrak{P}_n$ consisting of all partitions with two-element blocks.
\end{Rmk}

\subsection{The category $2\mathfrak{Cob}$ of $2$-cobordisms and its discrete skeleton}\label{sec:2Cob} For the category $2\mathfrak{Cob}$ of $2$-cobordisms, the reader is referred to Kock \cite{kock}, Carqueville and Runkel \cite{carq} and Lauda and Pfeiffer \cite{lauda}. We shall present here an \textsl{intuitive view} of this category since we are mainly interested in its underlying combinatorial structure. The precise definition involves delicate concepts from differential geometry and category theory, for the details of which the reader may consult \cite{kock,carq,lauda}.

As for the categories $\mathfrak{P}$ and $\mathfrak{P}^{\mathrm{d}}$, the set of objects of $2\mathfrak{Cob}$ is $\mathbb{N}_0$. For $m,n\in \mathbb{N}_0$, an $(m,n)$-\emph{2-cobordism} $\bm\alpha$ is a compact, oriented surface having a boundary which consists of $m+n$ connected components. Every such component is a circle (a compact $1$-dimensional manifold without boundary). This collection of $m+n$ pairwise disjoint circles is grouped into two parts: the first of which, consisting of $m$ circles, is denoted $\mathbf{m}$ and is called the \emph{incoming boundary} (its circles are the \emph{incoming circles}), while the second part, consisting of $n$ circles, is denoted $\mathbf n$ and is called the \emph{outgoing boundary} (its circles are the \emph{outgoing circles}). The entire boundary then is viewed as the disjoint union $\mathbf{m}\sqcup\mathbf{n}$ and $\bm\alpha$ is a cobordism $\mathbf{m}\mathrel{\contour{black}{${\rightsquigarrow}$}}\mathbf{n}$. In addition, it is assumed that a linear order $C_1>\cdots>C_m$ is imposed on the circles of $\mathbf m$, and likewise on the circles of $\mathbf{n}$. In other words, the circles of $\mathbf{m}$ (respectively $\mathbf{n}$) are attached with unique markers from $[m]$ (respectively $[n]$) so that we have a canonical map $\wp\colon\mathbf{m}\sqcup\mathbf{n}\to[m]\sqcup[n]$ mapping the $i$th circle of $\mathbf{m}$ (respectively $\mathbf{n}$) to $i$. Throughout, every such surface is considered only up to diffeomorphisms fixing the boundary $\mathbf{m}\sqcup\mathbf{n}$.

The next illustration shows a cobordism $\mathbf{1}\mathrel{\contour{black}{${\rightsquigarrow}$}}\mathbf{2}$ known as a \emph{pair of pants} (the left circle is incoming and the right circles are outgoing):
\begin{center}
  \begin{tikzpicture}[
  tqft,
  cobordism edge/.style={draw},
  view from=incoming,
  cobordism height=2cm,
]
\begin{scope}[every node/.style={rotate=90}]
\pic[tqft/pair of pants,
every incoming boundary component/.style={draw, rotate=90},
every outgoing lower boundary component/.style={draw,dotted},
every outgoing upper boundary component/.style={draw},
at={(-5,0)},
];
\end{scope}
  \end{tikzpicture}
\end{center}

An important issue of a component of a cobordism will be its \emph{genus}: this is the maximal number of cuts along (pairwise non-intersecting) closed simple curves (not going through the boundary) leaving the surface connected. This number is essentially given by the number of the ``holes'' of the surface (see the example below, composed by two pairs of pants). Every such hole provides the possibility for one cut along a circle running either inside around the hole or outside around one of the two ``pipes'' leading past the hole.
\begin{center}
	\begin{tikzpicture}[
		tqft,
		cobordism edge/.style={draw},
		view from=incoming,
		cobordism height=2cm,
		]
		\begin{scope}[every node/.style={rotate=90}]
				\pic[tqft/pair of pants,
			every incoming boundary component/.style={draw, rotate=90},
			every outgoing upper boundary component/.style={},
			at={(-5,0)},
			];
			\pic[tqft/reverse pair of pants,
			at = {(-3,-1)},
			every outgoing lower boundary component/.style={draw,dotted},
			every outgoing upper boundary component/.style={draw},
			];
		\end{scope};
		\end{tikzpicture}
\end{center}

Let $\ell,m,n\in\mathbb{N}_0$ and $\bm\alpha\colon\contour{black}{$\ell$}\mathrel{\contour{black}{${\rightsquigarrow}$}}\mathbf{m}$ and $\bm\beta\colon\mathbf{m}\mathrel{\contour{black}{${\rightsquigarrow}$}}\mathbf{n}$ be $2$-cobordisms; we assume that $\bm\alpha\cap\bm\beta=\mathbf{m}$: the outgoing boundary of $\bm\alpha$ coincides with the incoming boundary of $\bm\beta$, and otherwise $\bm\alpha$ and $\bm\beta$ are disjoint. We define $\bm\alpha\bm\beta\colon\contour{black}{$\ell$}\mathrel{\contour{black}{${\rightsquigarrow}$}}\mathbf{n}$ by \emph{gluing $\bm\alpha$ with $\bm\beta$ along the common boundary $\mathbf m$}:
\begin{equation}\label{cobordcomp}
\bm\alpha\bm\beta=\bm\alpha\mathrel{\underset{\mathbf{m}}{\sqcup}}\bm\beta:=\bm\alpha\cup\bm\beta
\end{equation}
(the latter under the proviso that $\bm\alpha\cap\bm\beta=\mathbf{m}$ and that the merged surface is smooth along the gluing circles of $\mathbf{m}$). That is, the incoming boundary of $\bm\alpha$ inside $\bm\alpha\cup\bm\beta$ is declared to be the incoming boundary of $\bm\alpha\bm\beta$, and the outgoing boundary of $\bm\beta$ inside $\bm\alpha\cup\bm\beta$ to be the outgoing boundary of $\bm\alpha\bm\beta$. This provides an associative (partial) composition of $2$-cobordisms. We present an example: the composition of the $2$-cobordism $\bm\alpha\colon\mathbf{6}\mathrel{\contour{black}{${\rightsquigarrow}$}}\mathbf{5}$ given by (upper circles are incoming, lower are outgoing)
\newline\vskip0.01cm  \begin{center}
\begin{tikzpicture}[
  tqft,
  every incoming boundary component/.style={draw},
  every outgoing lower boundary component/.style={draw,dotted},
  cobordism/.style={draw},
  cobordism edge/.style={draw},
  view from=incoming,
  cobordism height=2cm,
]
\begin{scope}[every node/.style={rotate=0}]
\pic[name=a,
  tqft,
  incoming boundary components=3,
  skip incoming boundary components=2,
  outgoing boundary components=0
  ];
\pic[name=b,
  tqft,
  incoming boundary components=3,
  skip incoming boundary components=2,
  outgoing boundary components=3,
  skip outgoing boundary components=2,
  offset=1,
  anchor=incoming boundary 1,
  at=(a-incoming boundary 2)
];
\pic[name=c,
  tqft,
  incoming boundary components=2,
  outgoing boundary components=0,
  anchor={(0,0)},
  at=(b-incoming boundary 3)
];
\pic[name=d,
  tqft,
  incoming boundary components=0,
  outgoing boundary components=4,
  skip outgoing boundary components=3,
  at=(a-incoming boundary 1)
];
\end{scope}
\end{tikzpicture}\end{center}
\vskip0.01cm{\parindent=0pt
with the $2$-cobordism $\bm\beta\colon\mathbf{5}\mathrel{\contour{black}{${\rightsquigarrow}$}}\mathbf{2}$ given by}\newline\vskip0.01cm \begin{center}
\begin{tikzpicture}[
  tqft,
  every incoming boundary component/.style={draw},
  every outgoing lower boundary component/.style={draw,dotted},
  cobordism/.style={draw},
  cobordism edge/.style={draw},
  view from=incoming,
  cobordism height=2cm,
]
\begin{scope}[every node/.style={rotate=0}]
\pic[name=a,
  tqft,
  incoming boundary components=3,
  skip incoming boundary components=2,
  outgoing boundary components=0
  ];
\pic[name=b,
  tqft,
  incoming boundary components=3,
  skip incoming boundary components=2,
  outgoing boundary components=3,
  skip outgoing boundary components=2,
  offset=0.5,
  anchor=incoming boundary 1,
  at=(a-incoming boundary 2)
];
\pic[name=c,
  tqft,
  incoming boundary components=1,
  outgoing boundary components=0,
  anchor={(0,0)},
  at=(b-incoming boundary 3)
];
\end{scope}
\end{tikzpicture}\end{center}
\vskip0.01cm
{\parindent=0pt yields $\bm\alpha\bm\beta\colon\mathbf{6}\mathrel{\contour{black}{${\rightsquigarrow}$}}\mathbf{2}$}
\newline\vskip0.01cm
\begin{center}\begin{tikzpicture}[
  tqft,
  cobordism edge/.style={draw},
  view from=incoming,
  cobordism height=2cm,
]
\begin{scope}[every node/.style={rotate=0}]
\pic[name=a,
tqft,
every incoming boundary component/.style={draw},
incoming boundary components=3,
skip incoming boundary components=2,
outgoing boundary components=0
];
\pic[name=b,
tqft,
every incoming boundary component/.style={draw},
 every outgoing upper boundary component/.style={draw,dotted},
incoming boundary components=3,
skip incoming boundary components=2,
outgoing boundary components=3,
skip outgoing boundary components=2,
offset=1,
anchor=incoming boundary 1,
at=(a-incoming boundary 2)
];
\pic[name=c,
  tqft,
  every incoming boundary component/.style={draw},
  incoming boundary components=2,
  outgoing boundary components=0,
  anchor={(0,0)},
  at=(b-incoming boundary 3)
];
\pic[name=d,
  tqft,
  every outgoing upper boundary component/.style={draw,dotted},
  incoming boundary components=0,
  outgoing boundary components=4,
  skip outgoing boundary components=3,
  at=(a-incoming boundary 1)
];
\pic[name=aa,
  tqft,
  incoming boundary components=3,
  skip incoming boundary components=2,
  outgoing boundary components=0,
  anchor=incoming boundary 1,
  at={(0,-2)}
  ];
\pic[name=bb,
  tqft,
  every outgoing lower boundary component/.style={draw,dotted},
  every outgoing upper boundary component/.style={draw},
  incoming boundary components=3,
  skip incoming boundary components=2,
  outgoing boundary components=3,
  skip outgoing boundary components=2,
  offset=0.5,
  anchor=incoming boundary 1,
  at={(2,-2)}
  ];
\pic[name=cc,
  tqft,
  incoming boundary components=1,
  outgoing boundary components=0,
  anchor=incoming boundary 1,
  at={(8,-2)}
];
\end{scope}
\end{tikzpicture}\end{center}
\vskip0.01cm
{\parindent=0pt which (up to diffeomorphism) is the same as }\newline\vskip0.01cm
\begin{center}\begin{tikzpicture}[
  tqft,
  every incoming boundary component/.style={draw},
  every outgoing lower boundary component/.style={draw,dotted},
  cobordism/.style={draw},
  cobordism edge/.style={draw},
  view from=incoming,
  cobordism height=2cm,
]
\begin{scope}[every node/.style={rotate=0}]
\pic[name=a,
  tqft,
  incoming boundary components=3,
  skip incoming boundary components=2,
  outgoing boundary components=0
  ];
\pic[name=b,
  tqft,
  incoming boundary components=3,
  skip incoming boundary components=2,
  outgoing boundary components=3,
  skip outgoing boundary components=2,
  genus=1,
  offset=0.5,
  anchor=incoming boundary 1,
  at=(a-incoming boundary 2)
];
\pic[name=c,
  tqft,
  incoming boundary components=2,
  outgoing boundary components=0,
  anchor={(0,0)},
  at=(b-incoming boundary 3)
];
\end{scope}
\end{tikzpicture}
\end{center}\vskip0.01cm
{\parindent=0pt
Observe that all components in $\bm\alpha$ and $\bm\beta$ have genus $0$, but in $\bm\alpha\bm\beta$ there is a component having genus $1$. This behavior will be discussed and will be important later.}

Every $2$-cobordism $\bm\alpha\colon\mathbf{m}\mathrel{\contour{black}{${\rightsquigarrow}$}}\mathbf{n}$ may contain, apart from the connected components containing the boundary circles of $\mathbf{m}\sqcup\mathbf{n}$ also connected components without boundary. In the process of forming the composition $\bm\alpha\bm\beta$ of two $2$-cobordisms $\bm\alpha$ and $\bm\beta$, the components of $\bm\alpha$ and $\bm\beta$ without boundary are retained, but there may be created new components without boundary (this is analogous to the formation of dead blocks in the composition of partitions).

The (small) category $2\mathfrak{Cob}$ of all $2$-cobordisms then is defined as follows: the set of objects $\Ob(2\mathfrak{Cob})$ is $\mathbb{N}_0$, and for $m,n\in \mathbb{N}_0$ the set of all arrows $m\to n$ is the set of all $2$-cobordisms $\mathbf{m}\mathrel{\contour{black}{${\rightsquigarrow}$}}\mathbf{n}$, endowed with the composition defined above. The local identity at object $n$ is the $2$-cobordism $\mathbf{n}\mathrel{\contour{black}{${\rightsquigarrow}$}}\mathbf{n}$ whose components are $n$ `parallel' cylinders as illustrated below:
\begin{center}
\begin{tikzpicture}[
  tqft,
  cobordism edge/.style={draw},
  view from=incoming,
  cobordism height=2cm,
]
\begin{scope}[every node/.style={rotate=0}]
\pic[tqft/cylinder,  at={(1,0)},
 every incoming boundary component/.style={draw},
 every outgoing lower boundary component/.style={draw,dotted},
 every outgoing upper boundary component/.style={draw},
 ];
\pic[tqft/cylinder,  at={(2.5,0)},
 every incoming boundary component/.style={draw},
 every outgoing lower boundary component/.style={draw,dotted},
 every outgoing upper boundary component/.style={draw},
 ];
\pic[tqft/cylinder,  at={(5.5,0)},
 every incoming boundary component/.style={draw},
 every outgoing lower boundary component/.style={draw,dotted},
 every outgoing upper boundary component/.style={draw},
 ];
\draw (1,0) node {$1$};
\draw (1,-2) node {$1$};
\draw (2.5,0) node {$2$};
\draw (2.5,-2) node {$2$};
\draw (4,0) node {$\dots$};
\draw (4,-2) node {$\dots$};
\draw (5.5,0) node {$n$};
\draw (5.5,-2) node {$n$};
\end{scope}
\end{tikzpicture}
\end{center}

In the following we shall describe every $2$-cobordism by discrete and `finitary' data and thus shall relate the $2$-cobordism category $2\mathfrak{Cob}$ with the partition categories $\mathfrak{P}^{\mathrm{d}}$ and $\mathfrak{P}$. Let $\bm\alpha\colon\mathbf{m}\mathrel{\contour{black}{${\rightsquigarrow}$}}\mathbf{n}$ be a $2$-cobordism; in general, $\bm\alpha$ has components with and without boundary. We define
\begin{equation*} \label{def:components without boundary}
s_{\bm\alpha}:=(s_0,s_1,\dots)\in \bigoplus_{i\in\mathbb{N}_0}\mathbb{N}_0
\end{equation*}
where $s_i$ is the number of components of $\bm\alpha$ without boundary which have genus $i$. All but finitely many  entries of the sequence $s_{\bm\alpha}$ are $0$. Next we recall that there is a canonical map $\wp\colon\mathbf{m}\sqcup \mathbf{n}\twoheadrightarrow [m]\sqcup[n]$ which maps the $i$th circle of $\mathbf{m}$ (respectively $\mathbf{n}$) onto the $i$th vertex of $[m]$ (respectively $[n]$). Then $\bm\alpha$ induces a partition $\pi_{\bm\alpha}\colon[m]\leadsto [n]$ such that vertices $u,v\in [m]\sqcup[n]$ lie in the same block of $\pi_{\bm\alpha}$ if and only if the circles $\wp^{-1}(u)$ and $\wp^{-1}(v)$ belong to the same connected component of $\bm\alpha$. In order to determine $\bm\alpha$ completely (up to diffeomorphisms), we need to know, apart from the induced partition $\pi_{\bm\alpha}$ and the sequence $s_{\bm\alpha}$,  only the genera of its components with non-empty boundary. If we consider the induced partition $\pi_{\bm\alpha}$ as the set of its blocks $\pi_{\bm\alpha}=\{B_1,\dots,B_k\}$ then the respective genera of the components of $\bm\alpha$ with boundary can be expressed as a function $g_{\bm\alpha}\colon\pi_{\bm\alpha}=\{B_1,\dots,B_k\}\to \mathbb{N}_0$. Altogether, the $2$-cobordism $\bm\alpha$ is completely and uniquely determined by the triple
\begin{equation}\label{def:discrete data for cobs}
(\pi_{\bm\alpha},g_{\bm\alpha},s_{\bm\alpha}).
\end{equation}
Conversely, every triple $(\alpha,g,s)$, where
\begin{enumerate}
\item $\alpha$ is a partition $[m]\leadsto [n]$,
\item $g$ is a function from (the set of blocks of) $\alpha$ to $\mathbb{N}_0$,
\item $s$ is a sequence of non-negative integers with only finitely many positive entries,
\end{enumerate}
defines a unique $2$-cobordism $\mathbf{m}\mathrel{\contour{black}{${\rightsquigarrow}$}} \mathbf{n}$. Thus, we can coordinatize the arrows of the category $2\mathfrak{Cob}$ with the set of all triples satisfying items (1)--(3). In order to complete the picture we need to explain composition of $2$-cobordisms in terms of their triple representation. In other words, given triples $(\eta,g,s)$ and $(\kappa,h,t)$ describing composable $2$-cobordisms $\bm\alpha$ and $\bm\beta$, respectively, how can we obtain the data $(\pi_{\bm\alpha\bm\beta},g_{\bm\alpha\bm\beta},s_{\bm\alpha\bm\beta})$ from the triples $(\eta,g,s)$ and $(\kappa,h,t)$?

The first component is the most obvious one: $\pi_{\bm\alpha\bm\beta}=\pi_{\bm\alpha}\pi_{\bm\beta}=\eta\kappa$. In order to describe $g_{\bm\alpha\bm\beta}$ and $s_{\bm\alpha\bm\beta}$ we have to determine the genus of every component of $\bm\alpha\bm\beta$. So, let $\bm\alpha\colon\contour{black}{$\ell$}\mathrel{\contour{black}{${\rightsquigarrow}$}}\mathbf{m}$ and $\bm\beta\colon \mathbf{m}\mathrel{\contour{black}{${\rightsquigarrow}$}} \mathbf{n}$ be two composable $2$-cobordisms and let $C$ be a component of $\bm\alpha\bm\beta$. If $C$ coincides with a component of $\bm\alpha$ (a component without boundary or having only incoming circles as boundary components) then the genus of $C$ as a component of $\bm\alpha\bm\beta$ coincides with the genus of $C$ as a component of $\bm\alpha$; likewise, if $C$ coincides with a component of $\bm\beta$ (a component without boundary or having only outgoing circles as boundary components) then the genus of $C$  as a component of $\bm\alpha\bm\beta$ coincides with the genus of $C$ as a component of $\bm\beta$. Otherwise $C$ is the union of some components $A_1,\dots,A_p$ of $\bm\alpha$ and $B_1,\dots,B_q$ of $\bm\beta$:
\begin{equation}\label{union of blocks}
C=\bigcup_{i=1}^p A_i\cup \bigcup_{j=1}^q B_j.
\end{equation}
The genus $g(C)$ of $C$ is at least the sum of the genera of the components $A_i$ and $B_j$:
$$g(C)\ge \sum_{i=1}^pg(A_i)+\sum_{j=1}^q g(B_j)$$ and we are interested in computing the \emph{increment}
\begin{equation*}
\iota(C;\bm\alpha,\bm\beta):= g(C)-\sum_{i=1}^pg(A_i)-\sum_{j=1}^qg(B_j).
\end{equation*}
This is most easily done by looking  at the \emph{Euler characteristic}. Recall that for a compact oriented surface $M$, the Euler characteristic $\chi(M)$ is the quantity
\begin{equation*}
\chi(M):=V-E+F
\end{equation*}
were $V,E,F$ are the numbers of vertices, edges, faces, respectively, of a triangularization of $M$ (and the quantity is independent of the chosen triangularization). We recall two fundamental properties of the Euler characteristic, see \cite[p. 63]{kock}.

The first one is that for such a surface $M$ with boundary the equality
\begin{equation}\label{chi(M)}
\chi(M)=2-2g-k
\end{equation}
holds where $g$ is the genus of $M$ and $k$ the number of its boundary components. The second property is the formula
\begin{equation*}
\chi(M\cup N)=\chi(M)+\chi(N)-\chi(M\cap N)
\end{equation*} for two surfaces $M$ and $N$. Using the fact that $\chi(S)=0$ for every circle $S$,  from (\ref{union of blocks}) we get
\begin{equation}\label{genus of C}
\chi(C)=\sum_{i=1}^p\chi(A_i)+\sum_{j=1}^q\chi(B_j).
\end{equation}
Let us denote:
\begin{itemize}
    \item $g$ the genus of $C$;
    \item $g_i, h_j$ the genera of $A_i, B_j$, respectively; 
    \item $\ell, r$ the number of incoming and outgoing circles of $C$, respectively;
    \item  $\ell_i, r_j$ the number of incoming circles of $A_i$ and outgoing circles of $B_j$, respectively;
    \item $s_i,t_j$ the number of outgoing circles of $A_i$ and incoming circles of $B_j$, respectively.
\end{itemize}
Then (\ref{chi(M)}) and (\ref{genus of C}) imply
\begin{equation*}
-\chi(C)=2g+\ell+r-2=\sum_{i=1}^p (2g_i+\ell_i+s_i-2)+\sum_{j=1}^q(2h_j+t_j+r_j-2).
\end{equation*}
Since $\ell=\sum_{i=1}^p \ell_i$ and $r=\sum_{j=1}^q r_j$ we get
\begin{equation}\label{last equality}
2g-2=\sum_{i=1}^p(2g_i+s_i)-2p+\sum_{j=1}^q(2h_j+t_j)-2q.
\end{equation}
Now we have:
$\sum_{i=1}^p s_i=\sum_{j=1}^q t_j$ and this coincides with the number of `inner' circles of $\mathbf{m}$ involved in the formation of the component $C$; observe also that $p+q$ is the number of components of $\bm\alpha$ and $\bm \beta$ involved in the formation of the component $C$. Since
\begin{equation}\label{definition of increment}
\iota(C;\bm\alpha,\bm\beta)=g-\sum_{i=1}^p g_i-\sum_{j=1}^q h_j
\end{equation}
we can use (\ref{last equality}) to get
\begin{equation*}
2\iota(C;\bm\alpha,\bm\beta)=\sum_{i=1}^p s_i+\sum_{j=1}^q t_j -2(p+q)+2
\end{equation*}
which eventually implies:
\begin{Prop}
\begin{equation}\label{determination of increment}
\iota(C;\bm\alpha,\bm\beta)=\sharp(\mbox{involved inner circles})-\sharp(\mbox{involved components})+1.
\end{equation}
\end{Prop}

The above analysis suggests defining a discrete counterpart $2\mathfrak{DisCob}$ of the category $2\mathfrak{Cob}$ as follows. The set of objects is again $\mathbb{N}_0$. For $m,n\in \mathbb{N}_0$ we let $2\mathfrak{DisCob}(m,n)$ be the set of all triples $(\alpha,g,s)$ where
\begin{itemize}
\item $\alpha$ is a partition $[m]\leadsto  [n]$,
\item $g$ is a function $\alpha\to\mathbb{N}_0$ ($\alpha$ is considered as the set of its blocks),
\item $s\in \bigoplus_{i\in\mathbb{N}_0}\mathbb{N}_0$ has only finitely many positive entries.
\end{itemize}
For $\ell,m,n\in \mathbb{N}_0$ and $(\alpha,g,s)\in 2\mathfrak{DisCob}(\ell,m)$, $(\beta,h,t)\in 2\mathfrak{DisCob}(m,n)$, we define
\begin{equation}\label{def:composition in 2Cob discrete}
(\alpha,g,s)(\beta,h,t):=\left(\alpha\beta,\mathbf{g}_{\alpha,g;\beta,h},s+t+\mathbf{b}_{\alpha,g;\beta,h}\right)
\end{equation}
where $\mathbf{g}$ and $\mathbf{b}$ are defined below (Definitions~\ref{def:label} and \ref{def:sequence}). First, we introduce some further notation: for a block $B$ of a partition $\alpha$ we denote by $\iv(B),\ov(B)$ and $\ve(B)$ the number of, respectively, incoming vertices, outgoing vertices and vertices of $B$. Moreover $B$ is an \emph{upper} block if $B$ contains only incoming vertices (that is, no outgoing vertices), a \emph{lower} block if it contains only outgoing vertices (no incoming vertices), and, otherwise, it is a \emph{transversal} block.

\begin{Def}
\label{def:label}
For a block $C$ of $\alpha\beta$, we define the number $\mathbf{g}_{\alpha,g;\beta,h}(C)\in\mathbb{N}_0$ which we refer to as the \emph{label} of $C$. If $C$ is an upper block of $\alpha$, then $\mathbf{g}_{\alpha,g;\beta,h}(C)=g(C)$; if it is a lower block of $\beta$, then $\mathbf{g}_{\alpha,g;\beta,h}(C)=h(C)$; otherwise $C$ is composed from blocks of $\alpha$ and $\beta$. In this case let $C'$ be the corresponding block of $(\alpha\cup\beta)^t$ (the partition of $[\ell]\sqcup[m]\sqcup[n]$ generated by $\alpha\cup \beta$). Suppose that $C'=\bigcup_{i=1}^a A_i\cup \bigcup_{j=1}^b B_j$ where the $A_i$ are the involved blocks of $\alpha$ and the $B_j$ are the involved blocks of $\beta$. Let $v$ denote the number of elements of $[m]$ involved in $C'$. Note that $v=\sum_{i=1}^a\ov(A_i)=\sum_{j=1}^b\iv(B_j)$. Motivated by \eqref{definition of increment} and \eqref{determination of increment}, $\mathbf{g}_{\alpha,g;\beta,h}$ is defined as
\begin{equation}\label{def:g}
\mathbf{g}_{\alpha,g;\beta,h}(C)=\sum_{i=1}^a g(A_i) + \sum_{j=1}^b h(B_j)+v-(a+b)+1.
\end{equation}
In this context, we call the quantity
\begin{equation}\label{def:increment(discrete)}
\g(C;\alpha,\beta):=v-(a+b)+1
\end{equation}
the \emph{increment of the label} of $C$ (with respect to $\alpha$ and $\beta$), where $\mathbf{g}_{\alpha,g;\beta,h}(C)$  is defined by \eqref{def:g}. Notice that the formula \eqref{def:g} serves for all blocks of $\alpha\beta$ (dead or not) that involve blocks of both $\alpha$ and $\beta$; in fact, it also agrees with the above conventions about upper blocks of $\alpha$ and lower blocks of $\beta$ that are inherited by $\alpha\beta$. Say, for any upper block $C$ of $\alpha$ considered as a block of $\alpha\beta$, then we have $v=b=0$ since $C$ has no vertices from $[m]$ and involves no blocks of $\beta$. Further, $a=1$ and the only block of $\alpha$ involved in $C$ is $C$ itself. Therefore, the right hand side of \eqref{def:g} reduces to just $g(C)$. Similarly, for any lower block $C$ of $\beta$ considered as a block of $\alpha\beta$, the right hand side of \eqref{def:g} reduces to $h(C)$. Thus, \eqref{def:g} defines the function $\mathbf{g}_{\alpha,g;\beta,h}$ for all blocks of $\alpha\beta$.
\end{Def}

\begin{Def}
\label{def:sequence}
We define $\mathbf{b}_{\alpha,g;\beta,h}:=(b_0,b_1,\dots)$ where for every $i=0,1,\dotsc$, the number $b_i$ shows how many dead blocks $C$ created in the composition of $\alpha$ and $\beta$ have $\mathbf{g}_{\alpha,g;\beta,h}(C)=i$.
\end{Def}

Clearly, the definition of the category $2\mathfrak{DisCob}$ has been designed such that the following holds:
\begin{Thm}
\label{thm:skeleton of 2Cob}
The category $2\mathfrak{Cob}$ is isomorphic to the category $2\mathfrak{DisCob}$ under the identity map on objects and the map $\bm\alpha\mapsto (\pi_{\bm\alpha},g_{\bm\alpha},s_{\bm\alpha})$ on arrows.
\end{Thm}

Due to Theorem~\ref{thm:skeleton of 2Cob} we may (and will) identify the category $2\mathfrak{Cob}$ with its `discrete skeleton'  $2\mathfrak{DisCob}$.

\begin{Rmk}
Roughly speaking, Theorem~\ref{thm:skeleton of 2Cob} means that the category $2\mathfrak{Cob}$ --- up to the components without boundary --- is just the partition category whose arrows are decorated with extra labels, namely, genera of the corresponding cobordisms. Attaching labels to partitions is a frequently used idea. Say, partition \textsl{algebras} whose underlying diagrams carry labels from a finite group $G$ have been used by Bloss \cite{bloss} for his study of the representation theory of the wreath product $G\wr\mathrm{S}_n$ (where $\mathrm{S}_n$ denotes the symmetric group of degree $n$). In the Brauer and Temperley--Lieb world (that is, partition categories/algebras/monoids for which all blocks of the underlying partitions have size two), versions using labeled diagrams are abundant; see \cite{pseudovarbrauer,ernst,ernst2,green1,green2,kujawatharp,blob,towers}, to mention just a few.
\end{Rmk}

For every $n\in\mathbb{N}_0$, the local monoid at object $n$ of the category $2\mathfrak{Cob}$  is denoted $2\mathfrak{Cob}_n$. The maps
\begin{equation}\label{category morphisms}
(\alpha,g,s)\mapsto (\alpha,\sum s_i)\mapsto \alpha
\end{equation}
provide quotient morphisms $2\mathfrak{Cob}\twoheadrightarrow \mathfrak{P}^{\mathrm d}\twoheadrightarrow \mathfrak{P}$ of categories,  which restrict to quotient morphisms $2\mathfrak{Cob}_n\twoheadrightarrow \mathfrak{P}^{\mathrm d}_n\twoheadrightarrow \mathfrak{P}_n$ of the corresponding monoids.

\begin{Rmk}\label{rmk:monoidal} Without going into the details we mention that, as a \emph{monoidal category}, $2\mathfrak{Cob}$ can be nicely presented by a small number of generators and a transparent set of defining relations. (This presentation was found by Abrams~\cite{Abrams}; see \cite[Section 1.4]{kock} for a neat and accessible explanation.) To form the quotient $2\mathfrak{Cob}\twoheadrightarrow \mathfrak{P}^{\mathrm{d}}$ one needs to add just one more defining relation among $1\to 1$ arrows, namely
\begin{center}
\begin{tikzpicture}[
  tqft,
  cobordism edge/.style={draw},
  view from=incoming,
  cobordism height=2cm,
]
\begin{scope}[every node/.style={rotate=90}]
\pic[tqft/cylinder,  at={(1,0)},
 every incoming boundary component/.style={draw, rotate=90},
 every outgoing lower boundary component/.style={draw,dotted},
 every outgoing upper boundary component/.style={draw},
 ];
\pic[tqft/pair of pants,
every incoming boundary component/.style={draw, rotate=90},
every outgoing upper boundary component/.style={draw,dotted},
at={(-5,0)},
];
\pic[tqft/reverse pair of pants,
at = {(-3,-1)},
every outgoing lower boundary component/.style={draw,dotted},
every outgoing upper boundary component/.style={draw},
];
\end{scope};
\draw (0,0) node {$=$};
\end{tikzpicture}
\end{center}
and
to get the quotient $2\mathfrak{Cob}\twoheadrightarrow \mathfrak{P}$ one needs to add one further relation among $0\to 0$ arrows, namely
\begin{center}
\begin{tikzpicture}[
  tqft,
  cobordism edge/.style={draw},
  view from=incoming,
  cobordism height=2cm,
]
\begin{scope}[every node/.style={rotate=90}]
\pic[tqft/cap,
every incoming boundary component/.style={draw},
every outgoing upper boundary component/.style={draw,dotted},
at={(-3.5,0)},
];
\pic[tqft/cup,
at = {(-1.5,0)},
every outgoing lower boundary component/.style={draw,dotted},
every outgoing upper boundary component/.style={draw},
];
\end{scope};
\draw (0,0) node {$=$};
\draw (1,0) node {$\mathrm{id}_\mathbf{0}$};

\end{tikzpicture}
\end{center}
where $\mathrm{id}_\mathbf{0}$ is the identity at object $\mathbf{0}$ which essentially is the empty cobordism.
Comes~\cite{comes} used the quotient map $2\mathfrak{Cob}\twoheadrightarrow \mathfrak{P}$ and the presentation of $2\mathfrak{Cob}$ from \cite{kock} for his construction of a presentation of the linearized version of the partition category $\mathfrak{P}$.
\end{Rmk}

\subsection{Types of dead blocks in partition categories}\label{sec: types}
Mazorchuk \cite[Section 3]{maz:endomorphisms} observed that in the \emph{partial Brauer category} (that is, in the subcategory $\mathfrak{PB}$ of $\mathfrak{P}$ consisting of all partitions all of whose blocks have size at most $2$), two topologically different types of dead blocks can occur when two partitions are composed: dead blocks that involve a singleton block of one of the partitions and dead blocks that involve no singleton blocks. In \cite{MaMa14} the former dead blocks are called \emph{open strings} and the latter are called \emph{loops}. In the next illustration, we represent partitions $\alpha\colon[6]\leadsto [8]$ and $\beta\colon[8]\leadsto [7]$ from $\mathfrak{PB}$ by graphs with vertex sets $[6]\sqcup[8]$ and $[8]\sqcup[7]$, respectively; the edges of the graphs represent two-element blocks while the isolated vertices stand for singleton blocks.
\begin{center}
\begin{tikzpicture}
\node[] at (0.4,-1) {$\alpha$};
\foreach \x in {1,2,...,6} \draw (\x,0) node[above=1pt] {$\x$};
\foreach \x in {1,2,...,6} \filldraw (\x,0) circle (2pt);
\foreach \x in {1,2,...,8} \draw (\x,-2) node[below=4pt] {$\x$};
\foreach \x in {1,2,...,8} \filldraw (\x,-2) circle (2pt);
\draw[thick] (1,0)--(3,-2); \draw[thick](4,0)--(5,-2);
\draw[thick] (2,-2)..controls (3,-1.3)..(4,-2); \draw[thick] (3,0)..controls (4,-0.7)..(5,0);\draw[thick] (6,-2)..controls (6.5,-1.5)..(7,-2);
\node[] at (0.4,-3.8) {$\beta$};
\foreach \x in {1,2,...,8} \filldraw (\x,-2.8) circle (2pt);
\foreach \x in {1,2,...,7} \filldraw (\x,-4.8) circle (2pt);
\foreach \x in {1,2,...,7} \draw (\x,-4.8) node[below=4pt] {$\x$};
\draw[thick] (2,-4.8)--(3,-2.8); \draw[thick] (6,-4.8)--(5,-2.8);
\draw[thick] (1,-2.8).. controls (1.5,-3.3)..(2,-2.8);
\draw[thick] (6,-2.8).. controls (6.5,-3.3)..(7,-2.8);
\draw[thick] (1,-4.8).. controls (2,-4.1)..(3,-4.8);
\draw[thick] (5,-4.8).. controls (6,-4.1)..(7,-4.8);
\end{tikzpicture}
\end{center}
\noindent In the process of computing the product $\alpha\beta$, three dead blocks arise: $\{1,2,4\}$, $\{6,7\}$, and $\{8\}$, of which the first and the last are open strings and the second is a loop.

Mazorchuk \cite[Section 3]{maz:endomorphisms} asked if there exists an analogue of this classification of dead blocks in the context of the full partition category $\mathfrak{P}$. Let us call open strings \emph{type $0$ blocks} and loops \emph{type $1$ blocks}. We believe that the correct generalization to dead blocks of arbitrary partitions comes from the above analysis of composition in $2\mathfrak{Cob}$ and is given by formula \eqref{def:increment(discrete)}: if a dead block $C$ in the composition $\alpha\beta$ of two composable partitions $\alpha$ and $\beta$ is formed of, say, $a$ blocks from $\alpha$ and $b$ blocks from $\beta$ and contains $v$ intermediary vertices, then its \emph{type} should be the quantity $v-(a+b)+1$.

To see that the suggested notion of the type naturally generalizes the two types of dead blocks occurring in $\mathfrak{PB}$, represent partitions as cobordisms rather than graphs. For this, we `inflate' each edge to a cylindric pipe, each incoming isolated vertex to a cup, and each outgoing isolated vertex to a cap; see the next picture that demonstrates such inflation for the partitions $\alpha\colon[6]\leadsto [8]$ and $\beta\colon[8]\leadsto [7]$ from the previous illustration.

\medskip

\begin{center}
\begin{tikzpicture}
[
  tqft,
  every incoming boundary component/.style={draw},
  every outgoing lower boundary component/.style={draw,dotted},
  cobordism/.style={draw},
  cobordism edge/.style={draw},
  view from=incoming,
  cobordism height=2cm,
]
\foreach \x in {1,2,...,6} \draw (\x,0) node {$\x$};
\foreach \x in {1,2,...,8} \draw (\x,-2) node {$\x$};
\node[] at (0.4,-1) {$\alpha$};
\begin{scope}[every node/.style={rotate=0}]
\pic[name=a,
  tqft,
  incoming boundary components=2,
  outgoing boundary components=0,
  at={(3,0)}
];
\pic[name=b,
  tqft,
  incoming boundary components=1,
  outgoing boundary components=1,
  offset=.5,
  anchor=incoming boundary 1,
  at={(4,0)}
];
\pic[name=bb,
  tqft,
  incoming boundary components=1,
  outgoing boundary components=1,
  offset=1,
  anchor=incoming boundary 1,
  at={(1,0)}
];
\pic[name=c,
  tqft,
  incoming boundary components=0,
  boundary separation=28,
  outgoing boundary components=2,
  anchor=outgoing boundary 1,
  at={(6,-2)}
];
\pic[name=cc,
  tqft,
  incoming boundary components=0,
  outgoing boundary components=1,
  anchor=outgoing boundary 1,
  at={(8,-2)}
];
\pic[name=dd,
  tqft,
  incoming boundary components=1,
  outgoing boundary components=0,
  anchor=outgoing boundary 1,
  at={(2,0)}
];
\pic[name=ee,
  tqft,
  incoming boundary components=1,
  outgoing boundary components=0,
  anchor=outgoing boundary 1,
  at={(6,0)}
];
\pic[name=ff,
  tqft,
  incoming boundary components=0,
  outgoing boundary components=1,
  anchor=outgoing boundary 1,
  at={(1,-2)}
];
\pic[name=pp,
  tqft,
  incoming boundary components=0,
  outgoing boundary components=2,
  at={(2,0)}
];
\end{scope}
\end{tikzpicture}

\medskip

\begin{tikzpicture}
[
  tqft,
  every incoming boundary component/.style={draw},
  every outgoing lower boundary component/.style={draw,dotted},
  cobordism/.style={draw},
  cobordism edge/.style={draw},
  view from=incoming,
  cobordism height=2cm,
]
\foreach \x in {1,2,...,8} \draw (\x,0) node {$\x$};
\foreach \x in {1,2,...,7} \draw (\x,-2) node {$\x$};
\node[] at (0.4,-1) {$\beta$};
\begin{scope}[every node/.style={rotate=0}]
\pic[name=a12,
  tqft,
  incoming boundary components=2,
  boundary separation=28,
  outgoing boundary components=0,
  at={(1,0)}
];
\pic[name=b3,
  tqft,
  incoming boundary components=1,
  outgoing boundary components=1,
  offset=-.5,
  anchor=incoming boundary 1,
  at={(3,0)}
];
\pic[name=e4,
  tqft,
  incoming boundary components=1,
  outgoing boundary components=0,
  anchor=outgoing boundary 1,
  at={(4,0)}
];
\pic[name=d5,
  tqft,
  incoming boundary components=1,
  outgoing boundary components=1,
  offset=.5,
  anchor=incoming boundary 1,
  at={(5,0)}
];
\pic[name=a67,
  tqft,
  incoming boundary components=2,
  boundary separation=28,
  outgoing boundary components=0,
  at={(6,0)}
];
\pic[name=e8,
  tqft,
  incoming boundary components=1,
  outgoing boundary components=0,
  anchor=outgoing boundary 1,
  at={(8,0)}
];
\pic[name=p1,
  tqft,
  incoming boundary components=0,
  outgoing boundary components=2,
  at={(1,0)}
];
\pic[name=p5,
  tqft,
  incoming boundary components=0,
  outgoing boundary components=2,
  at={(5,0)}
];
\pic[name=f4,
  tqft,
  incoming boundary components=0,
  outgoing boundary components=1,
  at={(4,0)}
];
\end{scope}
\end{tikzpicture}
\end{center}

\medskip

When we compose two partitions in $\mathfrak{PB}$, we get pipes (and/or cups/caps) with no holes. However, some of the components without boundary arising in the course of composition may have holes, and it is the phenomenon that causes the distinction observed by Mazorchuk: components corresponding to open strings have genus 0 while components corresponding to loops have genus 1. Continuing with the above example, the product cobordism has the following three components without boundary, of which the first and the last have genus 0 and the second has genus 1.
\begin{center}
\begin{tikzpicture}
[
  tqft,
  every incoming boundary component/.style={draw,dotted},
  every outgoing lower boundary component/.style={draw,dotted},
  cobordism/.style={draw},
  cobordism edge/.style={draw},
  view from=incoming,
  cobordism height=2cm,
]
\foreach \x in {1,2,4,6,7,8} \draw (\x,0) node {$\x$};
\begin{scope}[every node/.style={rotate=0}]
\pic[name=a12,
  tqft,
  incoming boundary components=2,
  boundary separation=28,
  outgoing boundary components=0,
  at={(1,0)}
];
\pic[name=ff,
  tqft,
  incoming boundary components=0,
  outgoing boundary components=1,
  at={(1,2)}
];
\pic[name=e4,
  tqft,
  incoming boundary components=1,
  outgoing boundary components=0,
  anchor=outgoing boundary 1,
  at={(4,0)}
];
\pic[name=pp,
  tqft,
  incoming boundary components=0,
  outgoing boundary components=2,
  at={(2,2)}
];
\pic[name=a67,
  tqft,
  incoming boundary components=2,
  boundary separation=28,
  outgoing boundary components=0,
  at={(6,0)}
];
\pic[name=b67,
  tqft,
  incoming boundary components=0,
  boundary separation=28,
  outgoing boundary components=2,
  at={(6,2)}
];
\pic[name=e8,
  tqft,
  incoming boundary components=1,
  outgoing boundary components=0,
  at={(8,0)}
];
\pic[name=e8,
  tqft,
  incoming boundary components=0,
  outgoing boundary components=1,
  at={(8,2)}
];
\end{scope}
\end{tikzpicture}
\end{center}

Similarly, when composing partitions $\alpha$ and $\beta$ in the general partition category $\mathfrak{P}$, we can `inflate' $\alpha$ and $\beta$ to cobordisms without holes. Dead blocks arising in the course of composition correspond to components without boundary that may have holes, and it is quite natural to use the genus of the component corresponding to a dead block $C$ as the type of $C$. Formula \eqref{def:increment(discrete)} gives the expression $v-(a+b)+1$ for this genus provided that $C$ involves $a$ blocks from $\alpha$ and $b$ blocks from $\beta$ and contains $v$ intermediary vertices.

We mention in passing that the expression $v-(a+b)+1$ can be given yet another meaning if one represents partitions as forests whose trees represent blocks. (Such a graphical representation is widely used in the literature; see, e.g., \cite{dolinkaetal,Dea19,Ea11,Ea18,Ea19b,EF12,ER22a,ER22b,ER22c,halverson}.) Indeed, let a dead block $C$ arising when partitions $\alpha$ and $\beta$ are composed involve $a$ blocks from $\alpha$ and $b$ blocks from $\beta$ and contain $v$ intermediary vertices. If the involved blocks from $\alpha$ and $\beta$ are represented as trees, then the cyclomatic number of the (multi)graph obtained as the union of these trees is exactly $v-(a+b)+1$.
This is not an accidental coincidence. Indeed, the topological invariant of cobordisms expressed by the number $v-(a+b)+1$: namely  the \textsl{increment} of the genus of a merged component $C$ of a product $\bm\alpha\bm\beta$ of two cobordisms $\bm\alpha$ and $\bm\beta$  (with $v$ the number of intermediary boundary components and $a$ and $b$ the number of involved connected components of $\bm\alpha$ and $\bm\beta$) is \textsl{preserved} in the transition $2\mathfrak{Cob}\twoheadrightarrow\mathfrak{P}$, and is expressed as a topological invariant of graphs that encode the resulting partitions $\alpha$ and $\beta$. The transition $2\mathfrak{Cob}\twoheadrightarrow\mathfrak{P}$ can be realized by ignoring the genus of all involved surfaces and letting them (considered as surfaces of genus $0$) shrink to trees (with the boundary components shrinking to vertices). This way, the partitions are represented as forests.

Now let a (not necessarily dead) block $C$ arising when partitions $\alpha$ and $\beta$ are composed involve $a$ blocks and $k+v$ vertices from $\alpha$ and $b$ blocks and $\ell+ v$ vertices from $\beta$ with $v$ the number of the (common) intermediary vertices. Then the cyclomatic number of the (multi)graph obtained as the union of these trees is exactly $v-(a+b)+1$. Indeed, any forest on $k+v$ vertices having $a$  connected components has $k+v-a$ edges; likewise, any forest on $\ell+v$ vertices having $b$ connected components has $\ell+v-b$ edges. The connected (multi)graph formed by merging the two forests along their common vertices has $k+\ell+2v-(a+b)$ edges and $k+\ell+v$ vertices; its cyclomatic number therefore is $v-(a+b)+1$. In other words, the quantity $v-(a+b)+1$ measures the increment of the cyclomatic number of the graph, obtained as the result of merging partition blocks, over the respective cyclomatic numbers of the involved components.

Why should one care about any classification of dead blocks in the category $\mathfrak{P}$ if the composition in $\mathfrak{P}$ eventually ignores such blocks? Answering this natural question requires a digression in the role of diagram monoids in representation theory.

In Brauer's pioneering work~\cite{Br37}, Brauer monoids serve as linear bases of what is called \emph{Brauer's centralizer algebras} nowadays. Let $F$ be a field of characteristic 0 and $\theta$ an element in $F\setminus\{0\}$. Brauer's centralizer algebra $B_n(\theta)$ is the $F$-algebra spanned by the Brauer monoid $\mathfrak{B}_n$ as a linear space, with the multiplication of the basis vectors being defined by the rule
\begin{equation}\label{eq:multrule}
\pi_1\pi_2=\theta^{\mathrm{b}(\pi_1,\pi_2)}\pi,
\end{equation}
where $\mathrm{b}(\pi_1,\pi_2)$ is the number of dead blocks and $\pi$ is the product that arises when the partitions $\pi_1$ and $\pi_2$ are multiplied in $\mathfrak{B}_n$. In~\cite{Br37}, Brauer's centralizer algebra $B_n(m)$, where $m$ is a positive integer, was used to study the natural representation of the orthogonal group $\mathrm{O}_m$  on the $n$-th tensor power $(F^m)^{\otimes n}$ of the space $F^m$. (The algebra $B_n(m)$ is exactly the centralizer of the diagonal action of $\mathrm{O}_m$ on $(F^m)^{\otimes n}$, hence the name.) When $m$ is an even positive integer, Brauer's centralizer algebra $B_n(-m)$ allowed for a similar study of the representation of the symplectic group $\mathrm{Sp}_m$ on $(F^m)^{\otimes n}$. It was also present in~\cite{Br37}, albeit implicitly; see~\cite{HW89} for a detailed analysis. For the case where the parameter $\theta$ is arbitrary, the ring-theoretic structure of the algebra $B_n(\theta)$ has been determined by Wenzl~\cite{Wen88} (in particular, he has proved that if $\theta$ is not an integer, then the algebra $B_n(\theta)$ is semisimple).

Formula \eqref{eq:multrule} can also be used to define $F$-algebras $P_n(\theta)$ whose bases are the partition monoids $\mathfrak{P}_n$, and these $F$-algebras are nothing but Martin and Jones's partition algebras mentioned in the introduction. In representation theory, the partition algebra $P_n(m)$, where $m$ is a positive integer, appears as the centralizer of the diagonal action of the symmetric group $\mathrm{S}_m$ on $(F^m)^{\otimes n}$, where the action of $\mathrm{S}_m$ on $F^m$ comes from the standard representation of $\mathrm{S}_m$ by permutation matrices; see \cite{Jones94,Martin91,Martin00,Martin98}.

Mazorchuk~\cite[Section 3]{maz:endomorphisms} suggested utilizing two types of dead blocks occurring in the composition of arrows in the partial Brauer category $\mathfrak{PB}$ to define a two-parameter family of $F$-algebras. Let  $\mathfrak{PB}_n$ stand for the \emph{partial Brauer monoid}, that is, the local monoid of $\mathfrak{PB}$ at object $n$. One fixes two elements $\theta_0,\theta_1\in F\setminus\{0\}$ and introduces the algebra $PB_n(\theta_0,\theta_1)$ as the $F$-algebra spanned by $\mathfrak{PB}_n$ as a linear space, with the multiplication of the basis vectors being defined by the rule
\begin{equation*}
\pi_1\pi_2=\theta_0^{\mathrm{b}_0(\pi_1,\pi_2)}\theta_1^{\mathrm{b}_1(\pi_1,\pi_2)}\pi,
\end{equation*}
where $\mathrm{b}_i(\pi_1,\pi_2)$ is the number of dead blocks of type $i=0,1$ and $\pi$ is the product that arises  when the partitions $\pi_1$ and $\pi_2$ are multiplied in $\mathfrak{PB}_n$. The algebras $PB_n(\theta_0,\theta_1)$ and their representations were studied in depths by Martin and Mazorchuk~\cite{MaMa14}.

Now it should be clear that the suggested classification of dead blocks in the category $\mathfrak{P}$ allows one to define a multi-parameter version of partition algebras. Namely, fixing any sequence
$\Theta:=\{\theta_i\}_{i=0,1,\dots}$ of elements from $F\setminus\{0\}$, one can introduce the algebra $P_n(\Theta)$ as the $F$-algebra spanned by the partition monoid $\mathfrak{P}_n$ as a linear space, with the multiplication of the basis vectors being defined by the rule
\begin{equation*}
\pi_1\pi_2=\prod \theta_i^{\mathrm{b}_i(\pi_1,\pi_2)}\pi,
\end{equation*}
where $\mathrm{b}_i(\pi_1,\pi_2)$ is the number of dead blocks of type $i=0,1,\dots$ and $\pi$ is the product that arises when the partitions $\pi_1$ and $\pi_2$ are multiplied in $\mathfrak{P}_n$. (The expression $\prod\theta_i^{\mathrm{b}_i(\pi_1,\pi_2)}$ is well-defined because only finitely many among the numbers $\mathrm{b}_i(\pi_1,\pi_2)$ can be non-zero.) Even though the study of $P_n(\Theta)$ goes out of the scope of the present paper, this multi-parameter partition algebra appears to us to constitute a worthwhile object relevant to certain issues of representation theory\footnote{The proposed multi-parameter partition algebras are different from  algebras arising from multi-parameter colored partition categories recently introduced and studied by Mazorchuk and Srivastava \cite{MaSr22}.}.

\subsection{Regular versions of $\mathfrak{P}^{\mathrm{d}}$ and $2\mathfrak{Cob}$}\label{regular version}
We say that a category $\mathfrak{C}$ is (\emph{von-Neumann})-\emph{regular} if for every pair of objects $a,b\in\Ob\mathfrak{C}$ and every arrow $\alpha\in\mathfrak{C}(a,b)$, there exists an arrow $\gamma\in\mathfrak{C}(b,a)$ satisfying $\alpha\gamma\alpha=\alpha$. It is known that in this case there exists an arrow $\beta\in \mathfrak{C}(b,a)$ satisfying $\alpha\beta\alpha=\alpha$ as well as $\beta\alpha\beta=\beta$. In the latter case, $\alpha$ and $\beta$ are \emph{mutually inverse arrows}. A mapping $^*\colon \alpha\to \alpha^*$ on the set of arrows of $\mathfrak{C}$ such that $\alpha^*$ is an inverse of $\alpha$ for every arrow $\alpha$ is an \emph{inverse unary operation} of $\mathfrak{C}$. By the axiom of choice, every regular category admits an inverse unary operation $\alpha\mapsto \alpha^*$. In the following we shall be concerned with  unary operations satisfying 
the identities
\begin{equation}\label{weak involution}
 \alpha\alpha^*\alpha=\alpha\mbox{ and }  (\alpha^*)^*=\alpha
\end{equation}
in which case we call $^*$ a \emph{weak regular involution}. We note that the identity $\alpha^*\alpha\alpha^*=\alpha^*$ is a consequence of~\eqref{weak involution}. In particular, any weak regular involution is an inverse unary operation. If, in addition, the equality
\begin{equation}\label{involution law}
(\alpha\beta)^*=\beta^*\alpha^*,
\end{equation}
holds for all composable arrows $\alpha$ and $\beta$ of $\mathfrak{C}$, we refer to the unary operation $\alpha\mapsto\alpha^*$ as a \emph{regular involution} and say that $\mathfrak{C}$ is a \emph{regular $*$-category} (in analogy with \emph{regular $*$-semigroups} as defined in~\cite{nordal_scheiblich}, see also \cite{dudoea}).

For a partition $\alpha\colon[m]\leadsto [n]$, we define $\alpha^*\colon[n]\leadsto [m]$ to be `the same partition' as $\alpha$, but with incoming and outgoing vertices swapped. Then it is well known (and easy to verify) that
$\alpha\mapsto\alpha^*$ is a regular involution, so that the partition category $\mathfrak{P}$ is a regular $*$-category. In contrast, the deformed category $\mathfrak{P}^{\mathrm{d}}$ is not regular because of the second component of its arrows. Indeed, for any pair $(\alpha,s)$ with $\alpha\in \mathfrak{P}(m,n)$ and $s>0$ and any pair $(\beta,t)$ with $\beta\in \mathfrak{P}(n,m)$, the second component of the product $(\alpha,s)(\beta,t)(\alpha,s)$ is $2s+t+\mathrm{b}(\alpha,\beta)+\mathrm{b}(\alpha\beta,\alpha)>s$, whence the product is not equal to $(\alpha,s)$.

However, the category $\mathfrak{P}^{\mathrm{d}}$ can be embedded into a regular category $\overline{\mathfrak{P}^{\mathrm{d}}}$ by allowing negative second entries in the pairs $(\alpha,s)$. More precisely, the object set of $\overline{\mathfrak{P}^{\mathrm{d}}}$ is again $\mathbb{N}_0$, while the set of arrows $\overline{\mathfrak{P}^{\mathrm{d}}}(m,n)$ consists of all pairs $(\alpha,s)$ with $\alpha\in\mathfrak{P}(m,n)$ and $s\in \mathbb{Z}$. The composition of composable arrows is defined as per \eqref{composition deformed partition category}, and it is clear that the composition remains associative when extended to $\overline{\mathfrak{P}^{\mathrm{d}}}$. Now, for each $(\alpha,s)\in \overline{\mathfrak{P}^{\mathrm{d}}}(m,n)$ set
\begin{equation}\label{inverse unary operation deformed}
(\alpha,s)^*:=(\alpha^*,-s-\mathrm{b}(\alpha,\alpha^*)-\mathrm{b}(\alpha^*,\alpha)).
\end{equation}
Notice that $\mathrm{b}(\alpha,\alpha^*)$ is just the number $\lb(\alpha)$ of lower blocks of $\alpha$ while $\mathrm{b}(\alpha^*,\alpha)$ is the number $\ub(\alpha)$ of upper blocks of $\alpha$.
\begin{Prop}
The map $(\alpha,s)\mapsto (\alpha,s)^*$ defined by \eqref{inverse unary operation deformed} is a weak regular involution on the category $\overline{\mathfrak{P}^{\mathrm{d}}}$.
\end{Prop}

\begin{proof}
We verify the equalities \eqref{weak involution} by a direct computation. For the last of these equalities, this is pretty straightforward:
\begin{align*}
((\alpha,s)^*)^*
&=(\alpha^*,-s-\mathrm{b}(\alpha,\alpha^*)-\mathrm{b}(\alpha^*,\alpha))^*&&\text{by \eqref{inverse unary operation deformed}}\\
&=\big((\alpha^*)^*,\,-(-s-\mathrm{b}(\alpha,\alpha^*)-\mathrm{b}(\alpha^*,\alpha))\\ 
&\phantom{=(((\alpha^*)^*,}
-\mathrm{b}(\alpha^*,(\alpha^*)^*)
-\mathrm{b}((\alpha^*)^*,\alpha^*)\big)&&\text{by \eqref{inverse unary operation deformed}}\\
&=(\alpha,s)&&\text{as $(\alpha^*)^*=\alpha$}.
\end{align*}

In order to show that $(\alpha,s)(\alpha,s)^*(\alpha,s)=(\alpha,s)$, first compute the product $(\alpha,s)(\alpha,s)^*$:
\begin{align*}
(\alpha,s)(\alpha,s)^*&=(\alpha,s)(\alpha^*,-s-\mathrm{b}(\alpha,\alpha^*)-\mathrm{b}(\alpha^*,\alpha))&&\text{by \eqref{inverse unary operation deformed}}\\
&=(\alpha\alpha^*,s-s-\mathrm{b}(\alpha,\alpha^*)-\mathrm{b}(\alpha^*,\alpha)+\mathrm{b}(\alpha,\alpha^*))&&\text{by \eqref{composition deformed partition category}}\\
&=(\alpha\alpha^*,-\mathrm{b}(\alpha^*,\alpha)).
\end{align*}
Using this and \eqref{composition deformed partition category}, we then get
\begin{align*}
(\alpha,s)(\alpha,s)^*(\alpha,s)&=(\alpha\alpha^*,-\mathrm{b}(\alpha^*,\alpha))(\alpha,s)\\ 
&=(\alpha\alpha^*\alpha,-\mathrm{b}(\alpha^*,\alpha)+s+\mathrm{b}(\alpha\alpha^*,\alpha)).
\end{align*}
Since $\alpha\alpha^*\alpha=\alpha$, it remains to observe that $\mathrm{b}(\alpha^*,\alpha)=\mathrm{b}(\alpha\alpha^*,\alpha)$. Indeed, the lower blocks of the partition $\alpha^*$ are in one-to-one correspondence with the upper blocks of $\alpha$ and each lower block of $\alpha^*$ is inherited by the product $\alpha\alpha^*$. From this we see that $\mathrm{b}(\alpha^*,\alpha)=\ub(\alpha)=\mathrm{b}(\alpha\alpha^*,\alpha)$.
\end{proof}

\begin{Rmk}
The map $(\alpha,s)\mapsto (\alpha,s)^*$ is not a regular involution on $\ol{\mathfrak{P}^\mathrm{d}}$. Indeed, one can check that for composable arrows $(\alpha,s),(\beta,t)$ in $\mathfrak{P}^{\mathrm{d}}$ the equality
\begin{equation*}\left((\alpha,s)(\beta,t)\right)^*=(\beta,t)^*(\alpha,s)^*
\end{equation*}
holds if and only if $\lb(\alpha)+\ub(\beta)=2\mathrm{b}(\alpha,\beta)$ which can only happen if $\lb(\alpha)=\ub(\beta)$ and every lower block of $\alpha$ is merged with exactly one upper block of $\beta$ (and conversely).
\end{Rmk}

Similarly as the deformed partition category $\mathfrak{P}^{\mathrm d}$, also the cobordism category $2\mathfrak{Cob}$ can be extended to a regular category $\overline{2\mathfrak{Cob}}$ by allowing negative entries. More precisely, the object set of  $\overline{2\mathfrak{Cob}}$ is $\mathbb{N}_0$, and for $m,n\in\mathbb{N}_0$, the set of arrows $\overline{2\mathfrak{Cob}}(m,n)$ consists of all triples $(\alpha,g,s)$ where $\alpha\colon[m]\leadsto[n]$ is a partition, $g$ is a function from (the set of blocks of) $\alpha$ to $\mathbb{Z}$, and $s\in \bigoplus_{i\in\mathbb{Z}} \mathbb{Z}$ is a bi-infinite sequence $(\dots,s_{-1},s_0,s_1,\dots)$ of integers with only finitely many non-zero entries. The composition of composable arrows $(\alpha,g,s)$ and $(\beta,h,t)$ is defined as per \eqref{def:composition in 2Cob discrete}, where the function $\mathbf{g}_{\alpha,g;\beta,h}$ is still defined by \eqref{def:g} while $\mathbf{b}_{\alpha,g;\beta,h}$ now stands for the bi-infinite sequence $(\dots,b_{-1},b_0,b_1,\dots)$, where for every $i\in\mathbb{Z}$, the number $b_i$ shows how many dead blocks $C$ created in the composition of $\alpha$ and $\beta$ have $\mathbf{g}_{\alpha,g;\beta,h}(C)=i$. We need to show that this composition is associative.

\begin{Prop}
The composition defined by \eqref{def:composition in 2Cob discrete} is associative in $\overline{2\mathfrak{Cob}}$.
\end{Prop}

\begin{proof}
Let $\alpha\colon[k]\leadsto [\ell]$, $\beta\colon[\ell]\leadsto [m]$, $\gamma\colon[m]\leadsto [n]$ be partitions. Let $K$ be a block of $\alpha\beta\gamma$; we need to show that the accumulated increment of its label (over the sum of the labels of its blocks from $\alpha,\beta,\gamma$, respectively) is independent of whether we calculate this increment with respect to the factorization $(\alpha\beta)\gamma$ or $\alpha(\beta\gamma)$.

So, let $K$ be a block of $\alpha\beta\gamma$ and let $K'$ be the block of the partition $(\alpha\beta\cup\gamma)^t$ of $[k]\sqcup[m]\sqcup[n]$ corresponding to $K$. Then $K'=\bigcup_{i=1}^d D_i\cup\bigcup_{j=1}^cC_j$ where  $D_i$ and $C_j$ are, respectively,  the involved blocks of $\alpha\beta$ and  $\gamma$. Then
\begin{equation}\label{increment ab,c}
\g(K;\alpha\beta,\gamma)=\sum_{i=1}^d\ov(D_i)-(d+c)+1.
\end{equation}

Let us consider the partition $(\alpha\cup\beta)^t$ of $[k]\sqcup[\ell]\sqcup[m]$; for every $i$ let $D'_i$ be the block of $(\alpha\cup\beta)^t$ corresponding to $D_i$. Now there exist blocks $A_{i1},\dots,A_{ip_i}$ of $\alpha$ and $B_{i1},\dots B_{iq_i}$ of $\beta$ such that
$$D'_i=\bigcup_{r=1}^{p_i}A_{ir}\cup \bigcup_{t=1}^{q_i}B_{it}$$
and
\begin{equation}\label{increment D_i,a,b}
\g(D_i;\alpha,\beta)=\sum_{t=1}^{q_i}\iv(B_{it})-(p_i+q_i)+1.
\end{equation}
The accumulated increment of the label of $K$ with respect to the factorization $(\alpha\beta)\gamma$ then is
$$\g(K;\alpha\beta,\gamma)+ \sum_{i=1}^d\g(D_i;\alpha,\beta)$$
which by use of (\ref{increment ab,c}) and (\ref{increment D_i,a,b}) and the equality $\ov(D_i)=\sum_{t=1}^{q_i}\ov(B_{it})$ can be expanded as
\begin{align}
&\sum_{i=1}^d\sum_{t=1}^{q_i}\ov(B_{it})-(d+c)+1 +\sum_{i=1}^d((\sum_{t=1}^{q_i} \iv(B_{it }))-(p_i+q_i)+1)\nonumber\\
=&\sum_{i=1}^d\sum_{t=1}^{q_i}\ve(B_{it})-(d+c)+1-\sum_{i=1}^d(p_i+q_i) +d\nonumber\\
=&\sum_{i=1}^d\sum_{t=1}^{q_i}\ve(B_{it}) -(a+b+c)+1 \label{total increment}
\end{align}
were $\sum_{i=1}^d\sum_{t=1}^{q_i}\ve(B_{it})$ is the number of intermediary vertices from $[\ell]\sqcup[m]$ involved in the formation of block $K''$ within $(\alpha\cup\beta\cup\gamma)^t$ corresponding to $K$, $a=\sum_{i=1}^d p_i$ is the number of involved blocks from $\alpha$ and $b=\sum_{i=1}^d q_i$ is the number of involved blocks from $\beta$; since $c$ is the number of involved blocks from $\gamma$ the quantity in (\ref{total increment}) only depends on $\alpha\beta\gamma$ but not on the factorization $(\alpha\beta)\gamma$. That is, the calculation of the accumulated increment with respect to the factorization $\alpha(\beta\gamma)$ leads to the same result.

Altogether, the composition (\ref{def:composition in 2Cob discrete}) is associative with respect to the second component; in combination with (\ref{associativity deformed}) we get associativity also with respect to the third component.
\end{proof}

We define a unary operation $^*\colon\overline{2\mathfrak{Cob}}(m,n)\to \overline{2\mathfrak{Cob}}(n,m)$ for which the equalities \eqref{weak involution} hold, letting
\begin{equation*}
(\alpha,g,s)^*:=(\alpha^*,g^*,s^*)
\end{equation*}
where the blocks of the partition $\alpha^*$ are obtained from those of the partition $\alpha$ by swapping the incoming and outgoing vertices, the function $g^*\colon\alpha^*\to \mathbb{Z}$ is defined by
\begin{equation}\label{def:gstar}
g^*(C^*):=-g(C)-\ve(C)+2
 \end{equation}
for each block $C^*$ of $\alpha^*$, where $C$ stands for $C^*$ with the incoming and outgoing vertices swapped, and the sequence $s^*$ is given by
\begin{equation}\label{def:sstar}
s^*:=(\dots,-s_{-1},-s_0,-s_1-\ub(\alpha)-\lb(\alpha),-s_2,\dots).
\end{equation}

\begin{Prop}
\label{prop:2Cob_regular}
The map $(\alpha,g,s)\mapsto (\alpha^*,g^*,s^*)$ is a weak regular involution on the category $\overline{2\mathfrak{Cob}}$.
\end{Prop}

\begin{proof}
Again, we verify the equalities \eqref{weak involution} by a direct computation. First we show that $((\alpha,g,s)^*)^*=(\alpha^*,g^*,s^*)^*=((\alpha^*)^*,(g^*)^*,(s^*)^*)$ coincides with $(\alpha,g,s)$. The equality  $(\alpha^*)^*=\alpha$ is clear. Further, if $C$ is any block of $\alpha$ and $C^*$ is the corresponding block of $\alpha^*$, then we have from \eqref{def:gstar}
\[
(g^*)^*(C)=-g^*(C^*)-\ve(C^*)+2=-(-g(C)-\ve(C)+2)-\ve(C^*)+2=g(C)
\]
since $\ve(C^*)=\ve(C)$. Similarly, \eqref{def:sstar} yields
\begin{multline*}
(s^*)^*=\\
(\dots,-(-s_{-1}),-(-s_0),-(-s_1-\ub(\alpha)-\lb(\alpha))-\ub(\alpha^*)-\lb(\alpha^*),-(-s_2),\dots)\\
=(\dots,s_{-1},s_0,s_1,s_2,\dots)=s
\end{multline*}
since $\ub(\alpha)=\lb(\alpha^*)$ and $\lb(\alpha)=\ub(\alpha^*)$.

In order to show that $(\alpha,g,s)(\alpha^*,g^*,s^*)(\alpha,g,s)=(\alpha,g,s)$, first compute the product $(\alpha,g,s)(\alpha^*,g^*,s^*)$. By the formula \eqref{def:composition in 2Cob discrete} this product is equal to
$\left(\alpha\alpha^*,\mathbf{g}_{\alpha,g;\alpha^*,g^*},s+s^*+\mathbf{b}_{\alpha,g;\alpha^*,g^*}\right)$. Let us analyze the second and the third components of this triple.

The function $\mathbf{g}_{\alpha,g;\alpha^*,g^*}$ assigns a certain integer value to each block of the partition $\alpha\alpha^*$. Each block of $\alpha\alpha^*$ belongs to one of the following three sorts
:
\begin{itemize}
  \item the blocks of $\alpha$ having no outgoing vertices;
  \item the blocks of $\alpha^*$ having no incoming vertices;
  \item the blocks arising from unions $C\cup C^*$ where $C$ is a block of $\alpha$ having an outgoing vertex and $C^*$ is the corresponding block of $\alpha^*$.
\end{itemize}
We apply the definition of $\mathbf{g}_{\alpha,g;\beta,h}$ (see Definition~\ref{def:label}) to compute the value $\mathbf{g}_{\alpha,g;\alpha^*,g^*}(D)$ for blocks $D$ of each of these three sorts.

If $D=C$, where $C$ is a block of $\alpha$ having no outgoing vertices, then
\[
\mathbf{g}_{\alpha,g;\alpha^*,g^*}(D)=g(C).
\]

If $D=C^*$, where $C^*$ is a block of $\alpha^*$ having no incoming vertices, and $C$ is the block of $\alpha$ corresponding to $C^*$, then by \eqref{def:gstar}
\[
\mathbf{g}_{\alpha,g;\alpha^*,g^*}(D)=g^*(C^*)=-g(C)-\ve(C)+2.
\]

If $D$ is obtained by omitting the intermediary vertices from the union $C\cup C^*$ where $C$ is a block of $\alpha$ containing an outgoing vertex and $C^*$ is the corresponding block of $\alpha^*$, then computing the increment via the formula \eqref{def:increment(discrete)} yields
\[
\g(D;\alpha,\alpha^*)=v-(a+b)+1=\ov(C)-(1+1)+1=\ov(C)-1.
\]
Indeed, $C\cup C^*$ involves exactly one block of $\alpha$ and exactly one block of $\alpha^*$ so that the parameters $a$ and $b$ are both equal to 1, while the intermediary vertices of $C\cup C^*$ are nothing but the outgoing vertices of $C$ so that the parameter $v$ is equal to $\ov(C)$. Therefore, we have
\begin{align*}
\mathbf{g}_{\alpha,g;\alpha^*,g^*}(D)&=g(C)+g^*(C^*)+\g(D;\alpha,\alpha^*)&&\text{by~\eqref{def:g}}\\
                                     &=g(C)-g(C)-\ve(C)+2+\ov(C)-1&&\text{by \eqref{def:gstar}}\\
                                     &=\ov(C)-\ve(C)+1.&&
\end{align*}
An important consequence of the formula
\begin{equation}\label{eq:intermediate g}
\mathbf{g}_{\alpha,g;\alpha^*,g^*}(D)=\ov(C)-\ve(C)+1
\end{equation}
is that $\mathbf{g}_{\alpha,g;\alpha^*,g^*}(D)=1$ whenever $D$ is a dead block. Indeed, $C\cup C^*$ gives rise to a dead block if and only if $C$ has no incoming vertices (which implies that $C^*$ has no outgoing vertices), that is, $\ov(C)=\ve(C)$. In this case, $\ov(C)$ and $\ve(C)$ in the right hand side of \eqref{eq:intermediate g} cancel. Since the number of dead blocks created in the composition of $\alpha$ and $\alpha^*$ is equal to the number $\lb{\alpha}$ of lower blocks of $\alpha$, we conclude that all terms of the sequence $\mathbf{b}_{\alpha,g;\alpha^*,g^*}=(\dots,b_{-1},b_0,b_1,\dots)$ are equal to zero, except for $b_1=\lb(\alpha)$. Therefore, the third component of the product $(\alpha,g,s)(\alpha^*,g^*,s^*)$ is
\begin{align*}
s+s^*+\mathbf{b}_{\alpha,g;\alpha^*,g^*}&=(\dots,s_{-1},s_0,s_1,s_2,\dots)+{}\\
&\rule{14pt}{0pt}(\dots,-s_{-1},-s_0,-s_1-\ub(\alpha)-\lb(\alpha),-s_2,\dots)+{}\\
&\rule{14pt}{0pt}(\dots,\underset{-1\rule{0pt}{7.8pt}}{0},\underset{0\rule{0pt}{7.8pt}}{0},\underset{1}{\lb(\alpha)},\underset{2\rule{0pt}{7.8pt}}{0},\dots)\\
&=(\dots,\underset{-1\rule{0pt}{7.8pt}}{0},\underset{0\rule{0pt}{7.8pt}}{0},\underset{1}{-\ub(\alpha)},\underset{2\rule{0pt}{7.8pt}}{0},\dots).
\end{align*}

Now we are ready to compute $(\alpha,g,s)(\alpha^*,g^*,s^*)(\alpha,g,s)$. To avoid multistory subscripts, re-denote the function $\mathbf{g}_{\alpha,g;\alpha^*,g^*}$ by $f$. Then  the formula \eqref{def:composition in 2Cob discrete} yields
\begin{multline}\label{eq:expression for aa*a}
(\alpha,g,s)(\alpha^*,g^*,s^*)(\alpha,g,s)={}\\
\left(\alpha\alpha^*\alpha,\mathbf{g}_{\alpha\alpha^*,f;\alpha,g},(\dots,\underset{-1\rule{0pt}{7.8pt}}{0},\underset{0\rule{0pt}{7.8pt}}{0},\underset{1}{-\ub(\alpha)},\underset{2\rule{0pt}{7.8pt}}{0},\dots)+s+\mathbf{b}_{\alpha\alpha^*,f;\alpha,g}\right).
\end{multline}
We aim to show that the triple on the right hand side is equal to $(\alpha,g,s)$. As we know that $\alpha\alpha^*\alpha=\alpha$, it remains to analyze the second and the third components. For the second component, we have to verify the equality $\mathbf{g}_{\alpha\alpha^*,f;\alpha,g}(C)=g(C)$ for each block $C$ of $\alpha$. This holds by the definition of the function $\mathbf{g}_{\alpha\alpha^*,f;\alpha,g}$ for the lower blocks of $\alpha$. If $C$ is an upper block of $\alpha$, then the definition of $\mathbf{g}_{\alpha\alpha^*,f;\alpha,g}$ tells us that $\mathbf{g}_{\alpha\alpha^*,f;\alpha,g}(C)=f(C)=\mathbf{g}_{\alpha,g;\alpha^*,g^*}(C)$, but in the above analysis of the function $\mathbf{g}_{\alpha,g;\alpha^*,g^*}$, it was shown that $\mathbf{g}_{\alpha,g;\alpha^*,g^*}(C)=g(C)$ for upper blocks of $\alpha$. Thus, assume that $C$ is a transversal block of $\alpha$. It is easy to see that as a block of the product of $\alpha\alpha^*$ and $\alpha$, it is obtained by omitting the intermediary vertices from the union $D\cup C$ where, in turn, the block $D$ of $\alpha\alpha^*$ is obtained by omitting the intermediary vertices from the union $C\cup C^*$. The formula \eqref{def:g} gives
\begin{equation}\label{eq:compute 2nd comp}
\mathbf{g}_{\alpha\alpha^*,f;\alpha,g}(C)=f(D)+g(C)+\g(C;\alpha\alpha^*,\alpha).
\end{equation}

Computing the increment $\g(C;\alpha\alpha^*,\alpha)$ via the formula \eqref{def:increment(discrete)} yields
\[
\g(C;\alpha\alpha^*,\alpha)=v-(a+b)+1=\iv(C)-(1+1)+1=\iv(C)-1.
\]
Indeed, $D\cup C$ involves exactly one block of $\alpha\alpha^*$ and exactly one block of $\alpha$ so that the parameters $a$ and $b$ are both equal to 1, while the intermediary vertices of $D\cup C$ are nothing but the incoming vertices of $C$ so that the parameter $v$ is equal to $\iv(C)$.

Substituting the obtained expression for the increment in~\eqref{eq:compute 2nd comp}, we get
\begin{align*}
\mathbf{g}_{\alpha\alpha^*,f;\alpha,g}(C)&=f(D)+g(C)+\iv(C)-1&&\\
                                         &=\ov(C)-\ve(C)+1+g(C)+\iv(C)-1&&\text{by \eqref{eq:intermediate g}}\\
                                         &=g(C)&&
\end{align*}
since $\ov(C)+\iv(C)=\ve(C)$. The equality $\mathbf{g}_{\alpha\alpha^*,f;\alpha,g}=g$ is thus verified.

If $C$ is an upper block of $\alpha$, then $C^*$ is a lower block of $\alpha^*$, and this block is inherited by $\alpha\alpha^*$. Merging $C^*$ and $C$ creates a dead block, and this is the only way a dead block can arise when $\alpha\alpha^*$ and $\alpha$ are composed. Let $B$ be the dead block obtained by omitting the intermediary vertices from the union $C^*\cup C$. Computing its label, we have
\begin{align}
\mathbf{g}_{\alpha\alpha^*,f;\alpha,g}(B)&=f(C^*)+g(C)+\g(B;\alpha\alpha^*,\alpha)&&\notag\\
                                         &=g^*(C^*)+g(C)+\g(B;\alpha\alpha^*,\alpha)&&\notag\\
                                         &=-g(C)-\ve(C)+2+g(C)+\g(B;\alpha\alpha^*,\alpha)&&\text{by \eqref{def:gstar}}\notag\\
                                         &=-\ve(C)+2+\g(B;\alpha\alpha^*,\alpha).\label{eq:dead block label}
\end{align}
Computing the increment $\g(B;\alpha\alpha^*,\alpha)$ is completely analogous to the above computation of the transversal block increment and leads to the same result:
\[
\g(C;\alpha\alpha^*,\alpha)=\iv(C)-1.
\]
Substituting this in the expression \eqref{eq:dead block label} yields
\[
\mathbf{g}_{\alpha\alpha^*,f;\alpha,g}(B)=-\ve(C)+2+\iv(C)-1=1
\]
since $\ve(C)=\iv(C)$ for $C$ because an upper block of $\alpha$ has only incoming vertices. Thus, all dead blocks arising when $\alpha\alpha^*$ and $\alpha$ are composed have label 1. The number of dead blocks created in the composition of $\alpha\alpha^*$ and $\alpha$ is equal to the number $\ub(\alpha)$ of upper blocks of $\alpha$. Hence 
\[\mathbf{b}_{\alpha\alpha^*,f;\alpha,g}=(\dots,\underset{-1\rule{0pt}{7.8pt}}{0},\underset{0\rule{0pt}{7.8pt}}{0},\underset{1}{\ub(\alpha)},\underset{2\rule{0pt}{7.8pt}}{0},\dots).\] 
Therefore in the expression
\[
(\dots,\underset{-1\rule{0pt}{7.8pt}}{0},\underset{0\rule{0pt}{7.8pt}}{0},\underset{1}{-\ub(\alpha)},\underset{2\rule{0pt}{7.8pt}}{0},\dots)+s+\mathbf{b}_{\alpha\alpha^*,f;\alpha,g}
\]
for the third component of the product $(\alpha,g,s)(\alpha^*,g^*,s^*)(\alpha,g,s)$ in \eqref{eq:expression for aa*a}, the first and the last summands cancel, leaving the expression equal to $s$, as required.
\end{proof}

\begin{Rmk}
There exist other unary operations $(\alpha,g,s)\mapsto(\alpha,g,s)^*$ on arrows of the category $\overline{2\mathfrak{Cob}}$ that satisfy the equalities \eqref{weak involution} and extend the regular involution $\alpha\mapsto\alpha^*$ of the category $\mathfrak{P}$. For instance, one can modify the definition of $g^*\colon\alpha^*\to\mathbb{Z}$ as follows: if $C$ is a block of $\alpha$, then for the corresponding block $C^*$ of $\alpha^*$,
\[
g^*(C^*):=\begin{cases}
             -g(C)-\ve(C)+2 & \hbox{if $C$ is transversal;} \\
             -g(C) & \hbox{otherwise.}
           \end{cases}
\]
Then, letting $f:=\mathbf{g}_{\alpha,g;\alpha^*,g^*}$, it is easy to check that $\mathbf{g}_{\alpha\alpha^*,f;\alpha,g}(C)=g(C)$ for each block $C$ of $\alpha\alpha^*\alpha=\alpha$. Indeed, for the transversal blocks, the argument from the proof of Proposition~\ref{prop:2Cob_regular} works, and for the non-transversal blocks, the equality follows from the definition of the function $\mathbf{g}_{\alpha,g;\beta,h}$; see Definition~\ref{def:label}.

When $\alpha$ and $\alpha^*$ or $\alpha\alpha^*$ and $\alpha$ are composed, each dead block of the resulting partition arises from the union $C\cup C^*$ with $C$ a lower block of $\alpha$ or, respectively, $C^*\cup C$ with $C$ an upper block of $\alpha$. The increment of such dead blocks was computed in the course of the proof of Proposition~\ref{prop:2Cob_regular}: it is equal to $\ve(C)-1$ in either case. Since our modified definition implies $g^*(C^*)+g(C)=0$ for each non-transversal block $C$ of $\alpha$, the labels of dead blocks arising when $\alpha$ and $\alpha^*$ or $\alpha\alpha^*$ and $\alpha$ are composed are of the form $\ve(C)-1$ where $C$ runs over the set of lower, respectively, upper blocks of $\alpha$. For $i=1,2\dotsc$, let $\ub_i(\alpha)$ and $\lb_i(\alpha)$ stand for the number of upper and, respectively, lower blocks of $\alpha$ with $i$ vertices. The above analysis shows that
\begin{align*}
\mathbf{b}_{\alpha,g;\alpha^*,g^*}&=(\dots,\underset{-1\rule{0pt}{7.8pt}}{0},\underset{0\rule{0pt}{7.8pt}}{\lb_1(\alpha)},\underset{1}{\lb_2(\alpha)},\underset{2\rule{0pt}{7.8pt}}{\lb_3(\alpha)},\dots),\\
\mathbf{b}_{\alpha\alpha^*,f;\alpha,g}&=(\dots,\underset{-1\rule{0pt}{7.8pt}}{0},\underset{0\rule{0pt}{7.8pt}}{\ub_1(\alpha)},\underset{1}{\ub_2(\alpha)},\underset{2\rule{0pt}{7.8pt}}{\ub_3(\alpha)},\dots).
\end{align*}
Now letting
\begin{multline*}
s^*:=\\
(...,-s_{-1},-s_0-\ub_1(\alpha)-\lb_1(\alpha),-s_1-\ub_2(\alpha)-\lb_2(\alpha),-s_2-\ub_3(\alpha)-\lb_3(\alpha),...),
\end{multline*}
we see that
\[
(s+s^*+\mathbf{b}_{\alpha,g;\alpha^*,g^*})+s+\mathbf{b}_{\alpha\alpha^*,f;\alpha,g}=s.
\]

It should be clear now that it is possible to build further unary operations $(\alpha,g,s)\mapsto(\alpha,g,s)^*$, retaining the equalities \eqref{weak involution} as follows: one can vary the values of the function $g^*\colon\alpha^*\to\mathbb{Z}$ at non-transversal blocks of $\alpha$, compensating changes in the sequences of labels of dead blocks by a suitable modification of the map $s\mapsto s^*$. On the other hand, it seems that for the transversal blocks, the rule $g^*(C^*)=-g(C)-\ve(C)+2$ is unavoidable if one wants to ensure \eqref{weak involution}.
\end{Rmk}

As in the classical (non-regular) case, the maps \eqref{category morphisms} provide quotient morphisms $\overline{2\mathfrak{Cob}}\twoheadrightarrow \overline{\mathfrak{P}^{\mathrm{d}}}\twoheadrightarrow \mathfrak{P}$ which respect the unary operation $^*$ and restrict to quotient morphisms on the local monoids.

A presentation for $\overline{2\mathfrak{Cob}}$ as a monoidal category can be obtained from the one of $2\mathfrak{Cob}$, mentioned in Remark~\ref{rmk:monoidal} at the end of Subsection~\ref{sec:2Cob}, by adding firstly one $1\to 1$ arrow (which may be thought of an \emph{antihole}):
\begin{equation}\label{eq:antihole}
   \begin{tikzpicture}[
		tqft,
		cobordism edge/.style={draw},
		view from=incoming,
		cobordism height=2cm
		]
		\begin{scope}[every node/.style={rotate=90}]
						\pic[tqft/cylinder,  at={(-1,0)},
			every incoming boundary component/.style={draw, rotate=90},
			every outgoing lower boundary component/.style={draw,dotted},
			every outgoing upper boundary component/.style={draw},
			]; 
        \end{scope}
        \filldraw[gray](0,0)circle(5pt);
        \draw(-1.8,0)node{$\mathbf{ah}:=$};
    \end{tikzpicture}
\end{equation} and the relations:
\begin{center}
	\begin{tikzpicture}[
		tqft,
		cobordism edge/.style={draw},
		view from=incoming,
		cobordism height=2cm
		]
		\begin{scope}[every node/.style={rotate=90}]
			\pic[tqft/cylinder,  at={(1,0)},
			every incoming boundary component/.style={draw, rotate=90},
			every outgoing lower boundary component/.style={draw,dotted},
			every outgoing upper boundary component/.style={draw},
			];
			\pic[tqft/cylinder,  at={(-3,0)},
			every incoming boundary component/.style={draw, dotted, rotate=90},
			every outgoing lower boundary component/.style={draw,dotted},
			every outgoing upper boundary component/.style={draw},
			];
			\pic[tqft/pair of pants,
			every incoming boundary component/.style={draw, rotate=90},
			every outgoing upper boundary component/.style={draw,dotted},
			at={(-7,0)},
			];
			\pic[tqft/reverse pair of pants,
			at = {(-5,-1)},
			every outgoing lower boundary component/.style={draw,dotted},
			every outgoing upper boundary component/.style={draw, dotted},
			];
		\end{scope};
		\filldraw[gray](-2,0)circle(5pt);
		\draw (0,0) node {$=$};
        \draw(3.4,0)node{,};
       	\end{tikzpicture}
\end{center}
\begin{center}
	\begin{tikzpicture}[
		tqft,
		cobordism edge/.style={draw},
		view from=incoming,
		cobordism height=2cm,
		]
		\begin{scope}[every node/.style={rotate=90}]
			\pic[tqft/cylinder,  at={(1,0)},
			every incoming boundary component/.style={draw, rotate=90},
			every outgoing lower boundary component/.style={draw,dotted},
			every outgoing upper boundary component/.style={draw},
			];
			\pic[tqft/cylinder,  at={(-7,0)},
			every incoming boundary component/.style={draw, rotate=90},
			every outgoing lower boundary component/.style={draw,dotted},
			every outgoing upper boundary component/.style={draw,dotted},
			];
			
			\pic[tqft/pair of pants,
			every incoming boundary component/.style={draw, dotted, rotate=90},
			every outgoing upper boundary component/.style={draw,dotted},
			at={(-5,0)},
			];
			\pic[tqft/reverse pair of pants,
			at = {(-3,-1)},
			every outgoing lower boundary component/.style={draw,dotted},
			every outgoing upper boundary component/.style={draw},
			];
		\end{scope};
		\draw (0,0) node {$=$};
		\filldraw[gray](-6,0)circle(5pt);
         \draw(3.4,0)node{.};
	\end{tikzpicture}
\end{center}
In the triple representation this $\mathbf{1}\to \mathbf{1}$ (generalized) cobordism can be written as $\mathbf{ah}:=(\mathrm{id}_{[1]},-1,\underline{0})$ where $\mathrm{id}_{[1]}$ is the identity partition $[1]\leadsto [1]$, that is, the partition that identifies both elements of $[1]\sqcup[1]$, $-1$ in the second entry labels the unique block of $\mathrm{id}_1$ by $-1$, and $\underline{0}$ is the bi-infinite sequence consisting entirely of zeros: $\underline{0}=(0)_{i\in \mathbb{Z}}$. One can check that the subcategory of $\overline{2\mathfrak{Cob}}$ generated by $2\mathfrak{Cob}\cup\{\mathbf{ah}\}$ consists of all elements $(\al,g,s)$ for which $s$ has only non-negative entries. 

The category $\langle 2\mathfrak{Cob},\mathbf{ah}\rangle$ admits a geometric model. The objects are the same as for $2\mathfrak{Cob}$, namely $\{\mathbf{n}\colon n\in \mathbb{N}_0\}$. The arrows $\mathbf{m}\to\mathbf{n}$ are similar to those of $2\mathfrak{Cob}$, except that some of the involved surfaces of genus $0$ may carry finitely many labels of the form \begin{tikzpicture}[baseline = {(0, -0.1)}]
\filldraw[gray](0,0)circle(5pt);    
\end{tikzpicture} (antiholes). Altogether, the involved surfaces may carry holes or antiholes (but not both), or non of those.

For every integer $k$ let $\underline{k}$ denote the bi-infinite sequence $(s_i)_{i\in \mathbb{Z}}$ with 
\begin{equation*}
    s_i=\begin{cases}
        1\mbox{ if }i=k\\
        0\mbox{ otherwise}
    \end{cases}
\end{equation*}
and let $\malcevk:=(\mathrm{id}_{[0]},\varnothing,\underline{k})$.  Note that $\malcevk$ is an arrow $\mathbf{0}\to\mathbf{0}$.

In this geometric model, $\malcevk$ is a  compact surface without boundary  
having $k$ holes if $k>0$ and otherwise it is a sphere having $-k$ labels  of the form
\begin{tikzpicture}[baseline = {(0, -0.1)}]
\filldraw[gray](0,0)circle(5pt);    
\end{tikzpicture} (antiholes). In any case, $\malcevk$ is an arrow $\mathbf{0}\to\mathbf{0}$.

It is clear how $\malcevk$ can be expressed as a product of generators. In order to generate all of $\overline{2\mathfrak{Cob}}$ we need,  for each such $\malcevk$, a new generator $\overline{\malcevk}\colon \mathbf{0}\to \mathbf{0}$ representing the triple $(\mathrm{id}_{[0]}, \varnothing,-\underline{k})$ that is able to erase the element $\malcevk$: the corresponding defining relations are:
\begin{equation}\label{eq:boundaryless items}
\malcevk\overline{\malcevk}=\overline{\malcevk}\malcevk=\mathrm{id}_\mathbf{0}\mbox{ for all }k\in\mathbb{Z}.
\end{equation}
It follows that the geometric models of $2\mathfrak{Cob}$ and $\langle 2\mathfrak{Cob},\mathbf{ah}\rangle$ can be extended to all of $\overline{2\mathfrak{Cob}}$: for every $\malcevk$ add just a copy of it labeled by ``$-$'' and let this signed version represent $\overline{\malcevk}$.

Finally, we mention another sort of intermediary categories  between $2\mathfrak{Cob}$ and $\mathfrak{P}$ and $\overline{2\mathfrak{Cob}}$ and $\mathfrak{P}$, respectively. These are obtained by ignoring the third component in (\ref{def:discrete data for cobs}) and (\ref{def:composition in 2Cob discrete}) and we denote them by $2\mathfrak{Cob}^\circ$ and $\overline{2\mathfrak{Cob}^\circ}$, respectively; for cobordisms this means that we ignore  all existing components without boundary and neglect the components without boundary which are created in the composition of two cobordisms. The arrows of these categories are thus partitions the blocks of which carry labels from $\mathbb{N}_0$, respectively $\mathbb{Z}$. The labeling of the blocks of two composed partitions is obtained as the sum of the labels of the blocks involved, `deformed' by the increment (\ref{def:increment(discrete)}). The maps
\begin{equation}\label{category morphisms 2}
  (\alpha,g,s)\mapsto (\alpha,g)\mapsto \alpha
\end{equation}
then provide quotient morphisms $\overline{2\mathfrak{Cob}}\twoheadrightarrow \overline{2\mathfrak{Cob}^\circ}\twoheadrightarrow \mathfrak{P}$ which restrict to quotient morphisms   $2\mathfrak{Cob}\twoheadrightarrow 2\mathfrak{Cob}^\circ\twoheadrightarrow\mathfrak{P}$.
Finally, we show that the operation $^*$ becomes a regular involution, when restricted to $\overline{2\mathfrak{Cob}^\circ}$.

\begin{Prop} The category $\overline{2\mathfrak{Cob}^\circ}$ is a regular $*$-category with respect to the involution $^*$.
\end{Prop}
\begin{proof}
Let $(\alpha,g)$ and $(\beta,h)$ be composable arrows of $\overline{2\mathfrak{Cob}^\circ}$.
{We define the labeling functions $g^*, h^*$ by \eqref{def:gstar}, $g_1$ and $g_2$ by
\[ (\alpha,g)(\beta,h)=:(\alpha\beta,g_1)\mbox{ and }(\beta^*,h^*)(\alpha^*,g^*)=:(\beta^*\alpha^*,g_2)\] and $g_1^*$ again by \eqref{def:gstar}, that is, \[(\alpha\beta,g_1)^*=:((\alpha\beta)^*,g_1^*)\]  and we need to verify that $g_1^*=g_2$.}

Let $C$ be a block of $\alpha\beta$ and $C^*$ be the corresponding block of $(\alpha\beta)^*=\beta^*\alpha^*$. If $C$ is an upper block of $\alpha$ or a lower block of $\beta$ then $g_1^*(C^*)=g_2(C^*)$ is obviously true. So assume that $C$ is composed of blocks of $\alpha$ as well as $\beta$. Let $C_1$ be the corresponding block of $(\alpha\cup \beta)^t$; then
\begin{equation*}
  C_1=\bigcup_{i=1}^a A_i\cup \bigcup_{j=1}^b B_j
\end{equation*}
for some blocks $A_i$ of $\alpha$ and $B_j$ of $\beta$. By \eqref{def:g} and \eqref{def:increment(discrete)}
\begin{equation*}
  g_1(C)=\sum_{i=1}^{a}g(A_i)+\sum_{j=1}^{b}h(B_j)+\g(C;\alpha,\beta)
\end{equation*}
and hence by (\ref{def:gstar})
\begin{equation}\label{g1star}
  g_1^*(C^*)=-\sum_{i=1}^{a}g(A_i)-\sum_{j=1}^{b}h(B_j)-\g(C;\alpha,\beta)-\ve(C^*)+2
\end{equation}
and we note that
\begin{equation}\label{degree}
  \ve(C^*)=\ve(C)=\sum_{i=1}^{a}\iv(A_i)+\sum_{j=1}^{b}\ov(B_j).
\end{equation}
On the other hand,
\begin{equation*}
  C_1^*=\bigcup_{j=1}^bB_j^*\cup \bigcup_{i=1}^aA_i^*
\end{equation*}
where $C_1^*$ is the block of $(\beta^*\cup \alpha^*)^t$ corresponding to $C^*$ and $A_i^*$ and $B_j^*$ have the obvious meaning.
From
\begin{equation*}
  g^*(A_i^*)=-g(A_i)-\ve(A_i)+2
\end{equation*}
and
\begin{equation*}
  h^*(B_j^*)=-h(B_j)-\ve(B_j)+2
\end{equation*}
we obtain
\begin{multline*}
  g_2(C^*)=\\-\sum_{i=1}^{a}g(A_i)-\sum_{j=1}^{b}h(B_j)-\sum_{i=1}^{a}\ve(A_i)+2a-\sum_{j=1}^{b}\ve(B_j)+2b+\g(C^*;\beta^*,\alpha^*).
\end{multline*}
Hence, in view of \eqref{g1star} and \eqref{degree} we need to verify that
\begin{equation}
  \label{crucial equation}
  \begin{split}
& -\g(C;\alpha,\beta)-\sum_{i=1}^{a}\iv(A_i)-\sum_{j=1}^{b}\ov(B_j)+2 \\
 & =  -\sum_{i=1}^{a}\ve(A_i)+2a-\sum_{j=1}^{b}\ve(B_j)+2b+\g(C^*;\beta^*,\alpha^*) 
\end{split}
\end{equation}

Starting with the left hand side of (\ref{crucial equation}) we get
\begin{align*}
  &-\g(C;\alpha,\beta)-\sum_{i=1}^{a}\iv(A_i)-\sum_{j=1}^{b}\ov(B_j)+2 \\
 = &-\sum_{i=1}^{a}\ov(A_i)+(a+b)-1-\sum_{i=1}^{a}\iv(A_i)-\sum_{j=1}^{b}\ov(B_j) +2\\
 = &-\sum_{i=1}^{a}\ve(A_i)-\sum_{j=1}^{b}\ve(B_j)+\sum_{j=1}^{b}\iv(B_j) +(a+b)+1\\
 = &-\sum_{i=1}^{a}\ve(A_i)-\sum_{j=1}^{b}\ve(B_j)+(2a+2b)+\big(\sum_{j=1}^{b}\ov(B_j^*)-(a+b)+1\big)\\
 = &-\sum_{i=1}^{a}\ve(A_i)-\sum_{j=1}^{b}\ve(B_j) +(2a+2b) +\g(C^*;\beta^*,\alpha^*)
\end{align*}
we end up with the right hand side of (\ref{crucial equation}), as required.
\end{proof}

Again, a presentation of $\overline{2\mathfrak{Cob}^\circ}$ as a monoidal category can be obtained by the one of $2\mathfrak{Cob}$ mentioned in Remark~\ref{rmk:monoidal} (end of Section~\ref{sec:2Cob}) by adding the new generator $\mathbf{ah}$~\eqref{eq:antihole} together with its two natural relations. Erasing the items without boundary  is now easier since they do not exist any more: for this, no further generators are needed and instead of the relations~\eqref{eq:boundaryless items} we just add the relations
\[\malcevk=\mathrm{id}_\mathbf{0}\mbox{ for all }k\in\mathbb{Z}.\]

\begin{Rmk} At the first glance it may seem to be unnatural and unmotivated to consider the aforementioned regular extension $\overline{\mathfrak{C}}$ of the category $\mathfrak{C}$ in question. Indeed, the geometric intuition gets lost. For example, in case of $\overline{2\mathfrak{Cob}}$ a geometric realization would lead to surfaces having a negative genus which is absurd\footnote{But who knows? Many `negative' notions once considered absurd---such as negative numbers and square roots of negative numbers, or, more recently, negative probabilities (see, e.g., \cite{Feynman1987,Khrennikov2009})---have been proved to be useful.}. And indeed, the ``geometric model'' of $\overline{2\mathfrak{Cob}}$ we presented has some very artificial flavor:  surfaces with negative genus have to be represented as surfaces carrying ``antiholes'', that is,  as surfaces of genus $0$  having attached some labels \begin{tikzpicture}[baseline = {(0, -0.1)}]
\filldraw[gray](0,0)circle(5pt);    
\end{tikzpicture} ``coming from outside''.  However, from a semigroup theoretic point of view such a regular extension may be helpful to get a better understanding of the algebraic structure of the category, respectively  monoid, in question. For example the starred Green relations $\mathrsfs{R}^*$ and $\mathrsfs{L}^*$ \cite{gould} seem to be `what one expects'. What is more, the embedding of a category $\mathfrak{C}$ into a regular category $\overline{\mathfrak{C}}$ has been already successfully employed in some proofs of \cite{kitovvolkov,kitovvolkov23,east_gray_azeef_nik2026,east_gray_azeef_nik}.
\end{Rmk}

\section{Categories of annular type}\label{sec:annular type}
\subsection{Topological annular category} We are going to define the \emph{topological annular category}. This has been mentioned, for example, in   \cite{towers}; the local monoids of some subcategory appear in \cite{dosenpetric}, \cite{ernst2} and elsewhere. An important quotient is the  affine Temperley--Lieb category of \cite{grahamlehrer}.

Let $a>b>0$ be (positive) real numbers and $A(a,b)$ be the annulus in the complex plane (around $0$) with radii $a$ and $b$:
\begin{equation*}
A(a,b):=\{z\in \mathbb{C}\mid b\le \vert z\vert \le a\}.
\end{equation*}
Let further $m,n\in \mathbb{N}_0$ and assume that $[m]\sqcup [n]$ is embedded in $A(a,b)$ via
\begin{equation}\label{eq:m cup n in A(a,b)}
k\mapsto ae^{\frac{2\pi i(2k-1)}{2m}}\ (k\in [m])\mbox{ and }
\ell \mapsto be^{\frac{2\pi i(2\ell-1)}{2n}}\ (\ell\in [n]).
\end{equation}
Throughout the present context we shall identify sets of the form $[m]\sqcup[n]$ with their realization in $A(a,b)$ (unless mentioned otherwise).
{
\begin{Def}\label{def:A-diagram}
An \emph{$A(a,b)$-diagram} $\bm\alpha\colon[m]\mathrel{\contour{black}{${\rightsquigarrow}$}}[n]$ is a compact $1$-manifold  with boundary, homeomorphically embedded in $A(a,b)$ such that
\begin{enumerate}
\item the boundary of $\bm\alpha$ coincides with $[m]\sqcup [n]$ (that is, its realization in $A(a,b)$)
\item with the exception of its boundary, $\bm\alpha$ is contained in the open annulus $A^\circ(a,b)=\{z\in \mathbb{C}\mid b<\vert z\vert <a\}$.
\end{enumerate}
\end{Def}}
Moreover, an $A$-diagram is an $A(a,b)$-diagram for some $a>b>0$.
In this context, the elements of the first component $[m]$ of $[m]\sqcup [n]$ are called the \emph{incoming points} while those of the second component $[n]$ are the \emph{outgoing points}.
Every connected component of such a $1$-manifold is homeomorphic either with the unit interval {$I=[0,1]$}, in which case this component has two boundary points, or  with the unit circle {$C=\{c\in \mathbb{C}\mid |z|=1\}$}. The former components are called \emph{strings} the latter \emph{circles}; a string which connects an incoming point with an outgoing point is a \emph{transversal string}, a string connecting two incoming [outgoing] points is an \emph{in-string} [\emph{out-string}]. Every $A$-diagram $\bm\alpha\colon[m]\mathrel{\contour{black}{${\rightsquigarrow}$}}[n]$ has $\frac{m+n}{2}$ strings; in particular, $m\equiv n\ (\bmod\ 2)$. For later use we define the \emph{rank} of an $A$-diagram as the number of its transversal strings. Note that in case of rank $0$ there may exist circles which, inside the annulus, wrap around the inner boundary circle.

Instead of the annulus $A(a,b)$ one could consider the cylinder
\begin{equation*}
Z(a,b)=C\times [b,a],
\end{equation*} embed $[m]\sqcup[n]$ in $Z(a,b)$ via $[m]\to C\times\{a\}$, $k\mapsto (e^{\frac{2\pi i(2k-1)}{2m}},a)$ and $[n]\to  C\times \{b\}$, $\ell\mapsto(e^{\frac{2\pi i(2\ell-1)}{2n}},b)$, then could define, in a similar way, a $Z(a,b)$-diagram and then proceed in a way analogous as described below.

Now let $a>b$ and $c>d$ be positive real numbers and $\bm\alpha_0\colon[m]\mathrel{\contour{black}{${\rightsquigarrow}$}}[n]$ and $\bm\alpha_1\colon[m]\mathrel{\contour{black}{${\rightsquigarrow}$}}[n]$ be an $A(a,b)$- and an $A(c,d)$-diagram, respectively. Then $\bm\alpha_0$ and $\bm\alpha_1$ are \emph{equivalent} if
{they are related by an (ambient) isotopy relative to $[m]\sqcup[n]$. According to Definition \ref{def:A-diagram} there is a compact $1$-manifold $M= \bigsqcup_{i=1}^{\frac{m+n}{2}}I_i\sqcup\bigsqcup_{j=1}^kC_j$ consisting of $\frac{m+n}{2}$ copies $I_i$ of the unit interval and, say, $k$ copies $C_j$ of the unit circle such that $\bm\alpha_0$ is determined by a homeomorphic embedding $\alpha_0\colon M\to A(a,b)$ and $\bm\alpha_1$ is determined by a homeomorphic embedding $\alpha_1\colon M\to A(c,d)$. Then $\bm\alpha_0$ and $\bm\alpha_1$ are isotopic (relative to $[m]\sqcup [n]$) if the embedding $\alpha_0$ can be continuously deformed to the embedding $\alpha_1$ while leaving fixed the boundary points $[m]\sqcup [n]$. More precisely, for every $t\in I=[0,1]$ consider the intermediary annulus $A_t:=A((1-t)a+tc,(1-t)c+td)$ (then $A_0=A(a,b)$ and $A_1=A(c,d)$).  Then $\alpha_0$ and $\alpha_1$ are \emph{isotopic relative to $m\sqcup n$} if for every $t\in (0,1)$ there is  a homeomorphic embedding $\alpha_t\colon M\to A_t$ satisfying the conditions of Definition \ref{def:A-diagram} such that the map $H\colon M\times I\to \mathbb{C}\times I$ defined by  $H(x,t):=(\alpha_t(x),t)$ ($x\in M,\ t\in I$) is a homeomorphic  embedding, and such that, for every boundary point $u$ of $M$, $\alpha_t(u)$ is the same element of $[m]\sqcup [n]$ for every $t\in I$ (here $[m]\sqcup [n]$ is interpreted canonically as subset of every $A_t$ via \eqref{eq:m cup n in A(a,b)}).}

 A \emph{topological $A$-diagram} $[m]\mathrel{\contour{black}{${\rightsquigarrow}$}}[n]$ (t$A$-diagram for short) is an equivalence class of  $A$-diagrams $[m]\mathrel{\contour{black}{${\rightsquigarrow}$}}[n]$. The \emph{rank} of a t$A$-diagram is the rank of any of its representatives.
Less formally  speaking, a t$A$-diagram $\bm\alpha\colon[m]\mathrel{\contour{black}{${\rightsquigarrow}$}}[n]$ consists of a collection of $\frac{m+n}{2}$ strings drawn in an annulus $A(a,b)$ linking pairwise the elements of $[m]\sqcup [n]$, together with a collection of finitely many circles, also drawn in $A(a,b)$, such that any two distinct lines (strings or circles) do not intersect; the relative positions of the circles to each other and with respect to the strings do matter, but otherwise lines and circles may be moved and deformed {continuously, but must not intersect with each other during this deformation process}.

The composition of two composable t$A$-diagrams is very similar to  the composition of $2$-cobordisms (\ref{cobordcomp}). Let $\ell,m,n\in \mathbb{N}_0$ and $\bm\alpha\colon[\ell]\mathrel{\contour{black}{${\rightsquigarrow}$}}[m]$ and $\bm\beta\colon[m]\mathrel{\contour{black}{${\rightsquigarrow}$}}[n]$ be two t$A$-diagrams. Choose real numbers $a>b>c>0$ and let $\bm\alpha$ be realized as an $A(a,b)$-diagram $\bm\alpha_{a,b}$ and $\bm\beta$ be realized as an $A(b,c)$-diagram $\bm\beta_{b,c}$. Then $\bm\alpha_{a,b}\cup\bm\beta_{b,c}$ is an $A(a,c)$-diagram $[\ell]\mathrel{\contour{black}{${\rightsquigarrow}$}}[n]$  and by definition, $\bm\alpha\bm\beta$ is the equivalence class of $\bm\alpha_{a,b}\cup\bm\beta_{b,c}$.

Let us look at the example of two t$A$-diagrams $\bm\alpha, \bm\beta\colon [8]\mathrel{\contour{black}{${\rightsquigarrow}$}} [8]$ depicted in Figure~\ref{tikz:2diagrams}.
\begin{figure}[ht]
\begin{center}
\begin{tikzpicture}
[scale=0.25,rotate=22.5]
\centering
\draw(0,0)node{$\bm\alpha$};
\draw[thick,dotted] (0,0) circle [radius=6cm];
\draw[thick,dotted] (0,0) circle [radius=9cm];
\foreach \x in {0,45,...,315} \foreach \y in {6,9} \filldraw (\x:\y) circle (3pt);
\draw[blue,thick,domain=0:225,smooth, variable=\t]
plot (-\t-45:6+3*\t/225);
\draw[blue,thick,domain=0:225,smooth, variable=\t]
plot (-\t+180:6+3*\t/225);
\draw[red,thick,domain=0:225,smooth, variable=\t]
plot (-\t:6+3*\t/225);
\draw[red,thick,domain=0:225,smooth, variable=\t]
plot (-\t+45:6+3*\t/225);
\draw[red,thick] (90:6)..controls (112.5:6.7)..(135:6);
\draw[thick] (-90:6) ..controls (-112.5:6.7)..(-135:6);
\draw[thick] (0:9) ..controls (22.4:8.8)..(45:9);
\draw[thick] (-90:9) ..controls (-112.5:8.8)..(-135:9);
\end{tikzpicture}\quad
\begin{tikzpicture}
[scale=0.25,rotate=22.5]
\draw(0,0)node{$\bm\beta$};
\draw (0,-9) circle (0.1pt);
\draw[thick,dotted] (0,0) circle [radius=3cm];
\draw[thick,dotted] (0,0) circle [radius=6cm];
\foreach \x in {0,45,...,315} \foreach \y in {3,6} \filldraw (\x:\y) circle (3pt);
\draw[blue,thick,domain=0:135, smooth, variable=\t]
plot (-\t-45:3+3*\t/135);
\draw[blue,thick,domain=0:135, smooth, variable=\t]
plot (-\t+90:3+3*\t/135);
\draw[red,thick,domain=0:135, smooth, variable=\t]
plot (-\t-90:3+3*\t/135);
\draw[red,thick,domain=0:135, smooth, variable=\t]
plot (-\t-135:3+3*\t/135);
\draw[red,thick] (0:6)..controls (22.5:5)..(45:6);
\draw[thick] (-90:6)..controls (-112.5:5)..(-135:6);
\draw[thick] (45:3) ..controls (22.5:3.6)..(0:3);
\draw[thick] (135:3)..controls (157.5:3.6)..(180:3);
\end{tikzpicture}
\end{center}\caption{Two t$A$-diagrams $\bm\alpha$  and $\bm\beta$ }\label{tikz:2diagrams}
\end{figure}
Composition of these two diagrams is shown in Figure \ref{tikz:smoothing} (left hand side), after some `smoothening' of the strings we end up with the diagram depicted in Figure \ref{tikz:smoothing} (right hand side).
\begin{figure}[ht]
\begin{center}
\begin{tikzpicture}
[scale=0.25, rotate=22.5]
\centering
\draw[thick,dotted] (0,0) circle [radius=9cm];
\foreach \x in {0,45,...,315} \foreach \y in {6,9} \filldraw (\x:\y) circle (3pt);
\draw[blue,thick,domain=0:225,smooth, variable=\t]
plot (-\t-45:6+3*\t/225);
\draw[blue,thick,domain=0:225,smooth, variable=\t]
plot (-\t+180:6+3*\t/225);
\draw[red,thick,domain=0:225,smooth, variable=\t]
plot (-\t:6+3*\t/225);
\draw[red,thick,domain=0:225,smooth, variable=\t]
plot (-\t+45:6+3*\t/225);
\draw[red,thick] (90:6)..controls (112.5:6.7)..(135:6);
\draw[thick] (-90:6) ..controls (-112.5:6.7)..(-135:6);
\draw[thick] (0:9) ..controls (22.4:8.8)..(45:9);
\draw[thick] (-90:9) ..controls (-112.5:8.8)..(-135:9);
\draw (0,-9) circle (0.1pt);
\draw[thick,dotted] (0,0) circle [radius=3cm];
\draw[very thin,dotted] (0,0) circle [radius=6cm];
\foreach \x in {0,45,...,315} \foreach \y in {3,6} \filldraw (\x:\y) circle (3pt);
\draw[blue,thick,domain=0:135, smooth, variable=\t]
plot (-\t-45:3+3*\t/135);
\draw[blue,thick,domain=0:135, smooth, variable=\t]
plot (-\t+90:3+3*\t/135);
\draw[red,thick,domain=0:135, smooth, variable=\t]
plot (-\t-90:3+3*\t/135);
\draw[red,thick,domain=0:135, smooth, variable=\t]
plot (-\t-135:3+3*\t/135);
\draw[red,thick] (0:6)..controls (22.5:5)..(45:6);
\draw[thick] (-90:6)..controls (-112.5:5)..(-135:6);
\draw[thick] (45:3) ..controls (22.5:3.6)..(0:3);
\draw[thick] (135:3)..controls (157.5:3.6)..(180:3);
\end{tikzpicture}
\begin{tikzpicture}
\draw(0,0)node{};
    \draw(0,2.5)node{$=$};
\end{tikzpicture}
\begin{tikzpicture}
[scale=0.25,rotate=22.5]
\draw[thick,dotted] (0,0) circle [radius=3cm];
\foreach \x in {0,45,...,315}  \filldraw (\x:3) circle (3pt);
\draw[blue,thick,domain=0:135, smooth, variable=\t]
plot (-\t-45:3+3*\t/135);
\draw[blue,thick,domain=0:135, smooth, variable=\t]
plot (-\t+90:3+3*\t/135);
\draw[red,thick] (-90:3)..controls (-112.5:3.6)..(-135:3);
\draw[thick] (45:3) ..controls (22.5:3.6)..(0:3);
\draw[thick] (135:3)..controls (157.5:3.6)..(180:3);
\draw[thick,dotted] (0,0) circle [radius=9cm];
\foreach \x in {0,45,...,315}  \filldraw (\x:9) circle (3pt);
\draw[blue,thick,domain=0:225,smooth, variable=\t]
plot (-\t-45:6+3*\t/225);
\draw[blue,thick,domain=0:225,smooth, variable=\t]
plot (-\t+180:6+3*\t/225);
\draw[red,thick] (135:9)..controls (157.5:8.8)..(180:9);
\draw[thick] (0:9) ..controls (22.4:8.8)..(45:9);
\draw[thick] (-90:9) ..controls (-112.5:8.8)..(-135:9);
\draw[thick] (-2,-5.3) circle [radius=0.75cm];
\end{tikzpicture}\caption{The t$A$-diagram $\bm{\alpha\beta}$}\label{tikz:smoothing}
\end{center}
\end{figure}

This way we get the \emph{topological annular category} $\mathfrak{tAnn}$: its set of objects is again $\mathbb{N}_0$, and for $m,n\in \mathbb{N}_0$, the set of arrows $m\to n$ is the set of t$A$-diagrams $[m]\mathrel{\contour{black}{${\rightsquigarrow}$}}[n]$, subject to composition of arrows as defined above. As already mentioned, if $m\not\equiv n\ (\bmod\ 2)$ then there is no such diagram (the set of arrows $[m]\mathrel{\contour{black}{${\rightsquigarrow}$}}[n]$ is empty in this case). Hence $\mathfrak{tAnn}$ is decomposed into two connected components, the even part $\mathfrak{tAnn}^\mathrm{even}$ and the odd part $\mathfrak{tAnn}^\mathrm{odd}$ which could as well be  treated separately. As earlier, denote the local monoid at object $n$ by $\mathfrak{tAnn}_n$; then $\mathfrak{tAnn}_n=\mathfrak{tAnn}^\mathrm{odd}_n$ or $\mathfrak{tAnn}_n=\mathfrak{tAnn}^\mathrm{even}_n$ depending on whether $n$ is odd or even.

For later use, yet another view of (topological) annular diagrams is useful \cite{grahamlehrer}. This is by `lifting a diagram to the universal covering space of $A(a,b)$'. Let $a>b$ be again positive real numbers, let
\begin{equation*}
R(a,b)=\mathbb{R}\times [b,a]
\end{equation*} and denote the horizontal shift $R(a,b)\to R(a,b)$, $(x,y)\mapsto(x+1,y)$ by $\Theta$.
For $m,n\in \mathbb{N}_0$ let $\mathbb{Z}\times [m]$ and $\mathbb{Z}\times [n]$ be endowed with the lexicographic order and let $\left.\mathbb{Z}\times [m]\right.\sqcup\left.\mathbb{Z}\times [n]\right.$ be embedded in $R(a,b)$ via
\begin{align}\label{embedding points in R(a,b)}
\mathbb{Z}\times [m]\to \mathbb{R}\times\{a\},&\ (t,k)\mapsto (t+\frac{2k-1}{2m},a),\\\label{embedding points in R(a,b) II}
\mathbb{Z}\times [n]\to \mathbb{R}\times\{b\},&\ (t,\ell)\mapsto (t+\frac{2\ell-1}{2n},b).
\end{align}
By an \emph{$R(a,b)$-diagram} $\bm\alpha\colon[m]\mathrel{\contour{black}{${\rightsquigarrow}$}}[n]$ we understand a $1$-manifold inside $R(a,b)$ consisting of (at most) countably many circles in the interior $R^\circ(a,b)$ of $R(a,b)$, of countably many `strings' (being isomorphic with the unit interval $[0,1]$)  whose boundary points coincide with (the realization inside $R(a,b)$ of) $\left.\mathbb{Z}\times [m]\right.\sqcup\left.\mathbb{Z}\times [n]\right.$, and, at most finitely many components homeomorphic with $\mathbb{R}$ (these come from circles wrapping around in the annulus and are realized as straight lines parallel to the $x$-axis and can occur only if the rank of $\bm\alpha$ is $0$), and such that $\Theta(\bm\alpha)=\bm\alpha$ and  $\bm\alpha\cap [-1,1]\times [b,a]$ is compact (which essentially means that inside every bounded region, $\bm\alpha$ contains only finitely many circles). An $R$-diagram is an $R(a,b)$-diagram for some $a>b>0$. In an obvious way we define equivalence of $R$-diagrams  and then define a topological annular diagram as an equivalence class of $R$-diagrams; in this case equivalence can be also described by means of homeomorphisms $R(a,b)\to R(c,d)$ which pointwise `fix' the boundary, that is, map $(x,a)\mapsto (x,c)$ and $(x,b)\mapsto (x,d)$ for every $x\in \mathbb{R}$.

\subsection{Quotients and a subcategory of $\mathfrak{tAnn}$} We start with the category $\mathfrak{L}$, a decorated version of which occurs in \cite{green2} and elsewhere and the local monoids of which have been introduced and studied in \cite{dosenpetric}. Intuitively speaking, instead of annuli  we consider bounded rectangles in the plane with $[m]$ and $[n]$ embedded on opposite sides of the rectangle; otherwise the definition of suitable (topological) diagrams is analogous. However, we can realize this category as a subcategory of $\mathfrak{tAnn}$. Let us call an $A(a,b)$-diagram {$[m+1]\mathrel{\contour{black}{${\rightsquigarrow}$}}[n+1]$ \emph{rectangular} if it contains a string connecting the incoming point $m+1$ with the outgoing point $n+1$. Obviously, to be rectangular is a property of the topological diagram rather than of some representative.} The category $\mathfrak{L}$ then by definition has $\mathbb{N}_0$ as its object set, and, for $m,n\in \mathbb{N}_0$, $\mathfrak{L}(m,n)$ consists of all rectangular t$A$-diagrams $[m+1]\mathrel{\contour{black}{${\rightsquigarrow}$}}[n+1]$.

The local monoid $\mathfrak{L}_n$ is exactly the monoid $\mathcal{L}_n$ of \cite{dosenpetric}. The category $\mathfrak{L}$ may be considered as a Temperley--Lieb category containing circles in which the relative positions of the circles to each other and to the strings \textsl{do matter}, as opposed to the classical Temperley--Lieb category where only the \textsl{number of circles}, but not their position matters.

As already mentioned, a t$A$-diagram $[m]\mathrel{\contour{black}{${\rightsquigarrow}$}}[n]$ without a transversal string (in which case $m$ and $n$ have to be even) may admit two kinds of circles: circles which are $0$-homotopic (on the annulus) and circles which are not (which wrap around the inner circle). We call the former $0$-circles and the latter $\omega$-circles. For a t$A$-diagram $\bm\alpha$ denote by
$\underline{\bm\alpha}$ the diagram obtained from $\bm\alpha$ by removing all circles and call $\underline{\bm\alpha}$ the \emph{skeleton} of $\bm\alpha$. Moreover denote by $\cir_0(\bm\alpha)$ the number of all $0$-circle components contained in $\bm\alpha$ and by $\cir_\omega(\bm\alpha)$ the number of all $\omega$-circle components of $\bm\alpha$, and set $\cir(\bm\alpha):=\cir_\omega(\bm\alpha)+\cir_0(\bm\alpha)$. Finally, for a t$A$-diagram $\bm\alpha\colon[m]\mathrel{\contour{black}{${\rightsquigarrow}$}}[n]$ denote by $\pi_{\bm\alpha}$ the partition $[m]\leadsto [n]$ induced by $\bm\alpha$. We shall define five equivalence relations ${\equiv}_i$ ($i=1,\dots,5$) on $\mathfrak{tAnn}$ by: for $\bm\alpha,\bm\beta\in \mathfrak{tAnn}(m,n)$ let
\begin{align*}
\bm\alpha\mathrel{{\equiv}_1}\bm\beta &\Longleftrightarrow \underline{\bm\alpha}=\underline{\bm\beta},\ \cir_\omega(\bm\alpha)=\cir_\omega(\bm\beta),\ \cir_0(\bm\alpha)=\cir_0(\bm\beta);\\
\bm\alpha\mathrel{{\equiv}_2}\bm\beta &\Longleftrightarrow \underline{\bm\alpha}=\underline{\bm\beta},\ \cir_\omega(\bm\alpha)=\cir_\omega(\bm\beta);\\
\bm\alpha\mathrel{{\equiv}_3}\bm\beta &\Longleftrightarrow \underline{\bm\alpha}=\underline{\bm\beta};\\
\bm\alpha\mathrel{{\equiv}_4}\bm\beta &\Longleftrightarrow \pi_{\bm\alpha}=\pi_{\bm\beta};\\
\bm\alpha\mathrel{{\equiv}_5}\bm\beta &\Longleftrightarrow \pi_{\bm\alpha}=\pi_{\bm\beta},\ \cir(\bm\alpha)=\cir(\bm\beta).\end{align*}
It is easy to verify that these relations are category congruences. Under inclusion, they form a pentagon poset as shown:
\begin{center}
 \begin{tikzpicture}[scale=1.5,
  every node/.style={circle, inner sep=2pt}]

\node (e1) at (1,3) {$\equiv_4$};
\node (e2) at (0,2) {$\equiv_3$};
\node (e3) at (0,1) {$\equiv_2$};
\node (e4) at (1,0) {$\equiv_1$};
\node (e5) at (2,1.5) {$\equiv_5$};

\draw (e1) -- (e2);
\draw (e2) -- (e3);
\draw (e3) -- (e4);
\draw (e4) -- (e5);
\draw (e5) -- (e1);

\end{tikzpicture}   
\end{center}

We set:
\begin{equation}\label{annular}
\mathfrak{Ann}:=\mathfrak{tAnn}/{\equiv}_4
\end{equation}
and call it the \emph{annular} category; it is the category introduced by Jones \cite{Jonesannular}, see also \cite{penneys};
\begin{equation}\label{affine annular}
\mathfrak{aTL}:=\mathfrak{tAnn}/{\equiv}_2
\end{equation}
and note that this is the \emph{affine Temperley--Lieb category} introduced by Graham--Lehrer \cite{grahamlehrer}. Moreover, we call
\begin{equation}\label{deformed affine, essential affine}
\mathfrak{aTL}^{\mathrm d}:=\mathfrak{tAnn}/{\equiv}_1\mbox{ and }\mathfrak{aTL}^{\mathrm e}:=\mathfrak{tAnn}/{\equiv}_3
\end{equation}
the \emph{deformed} and the \emph{essential} affine Temperley--Lieb category, respectively. Finally, we call
\begin{equation}\label{deformed annular}
\mathfrak{Ann}^{\mathrm d}:=\mathfrak{tAnn}/{\equiv}_5
\end{equation} the \emph{deformed} annular category. We have quotient morphisms
\begin{equation*}
\mathfrak{tAnn}\twoheadrightarrow  \mathfrak{aTL}^{\mathrm d} \twoheadrightarrow \mathfrak{aTL} \twoheadrightarrow \mathfrak{aTL}^{\mathrm e} \twoheadrightarrow \mathfrak{Ann}
\end{equation*}
and
\begin{equation*}
\mathfrak{tAnn}\twoheadrightarrow \mathfrak{aTL}^{\mathrm{d}}\twoheadrightarrow\mathfrak{Ann}^{\mathrm d}\twoheadrightarrow \mathfrak{Ann}
\end{equation*} which induce quotient morphisms among the respective local monoids
\begin{equation*}
\mathfrak{tAnn}_n\twoheadrightarrow  \mathfrak{aTL}^{\mathrm d}_n \twoheadrightarrow \mathfrak{aTL}_n \twoheadrightarrow \mathfrak{aTL}^{\mathrm e}_n \twoheadrightarrow \mathfrak{Ann}_n
\end{equation*} and
\begin{equation*}
\mathfrak{tAnn}_n\twoheadrightarrow \mathfrak{aTL}_n^{\mathrm{d}}\twoheadrightarrow\mathfrak{Ann}^{\mathrm d}_n\twoheadrightarrow \mathfrak{Ann}_n.
\end{equation*}
We note that, for every $n$, $\mathfrak{Ann}_n$ is finite while $\mathfrak{aTL}^\mathrm{e}_n$ is infinite for all $n\ge 1$ and all other monoids are infinite for all $n$.

\subsection{Combinatorial description of affine Temperley--Lieb categories and the morphism $\mathfrak{aLT}^{\mathrm e}\twoheadrightarrow \mathfrak{Ann}$}\label{subsection:combinatorial description} The arrows of $\mathfrak{aTL}^{\mathrm e}$ can be viewed as t$A$-diagrams without circles; their composition in $\mathfrak{aTL}^{\mathrm e}$ is `the same' as in $\mathfrak{tAnn}$ except that all possible circles which are created through the composition process are dismissed. This leads to a description of the arrows of $\mathfrak{aTL}^{\mathrm{e}}$ which is entirely combinatorial. Let $a>b>0$ be real numbers and $\left.\mathbb{Z}\times[m]\right.\sqcup\left.\mathbb{Z}\times[n]\right.$ be embedded in $R(a,b)$ according to (\ref{embedding points in R(a,b)}) and (\ref{embedding points in R(a,b) II}). Let $\bm\alpha\colon[m]\mathrel{\contour{black}{${\rightsquigarrow}$}}[n]$ be a t$A$-diagram without any circles (that is, it has neither $\omega$- nor $0$-circles). Realized as an $R(a,b)$-diagram it has only connected components homeomorphic with the unit interval $I$, each of these is connecting exactly two endpoints of $\left.\mathbb{Z}\times[m]\right.\sqcup\left.\mathbb{Z}\times[n]\right.$. It follows that $\bm\alpha$ is uniquely determined by the induced partition $\mathbb{Z}\times[m]\leadsto  \mathbb{Z}\times[n]$ --- two elements of $\left.\mathbb{Z}\times[m]\right.\sqcup\left.\mathbb{Z}\times[n]\right.$ are related if they are connected by a string (a connected component of $\bm\alpha$). This induced partition is such that every block is of size exactly $2$. On the other hand we have \cite{grahamlehrer}:
\begin{Prop}\label{affine diagrams}
 A partition $\bm\pi\colon \mathbb{Z}\times [m]\leadsto  \mathbb{Z}\times [n]$ all of whose blocks have size $2$ determines a unique element of $\mathfrak{aTL}^{\mathrm e}(m,n)$ if and only if
\begin{enumerate}
\item $(s,k)\mathrel{\bm\pi}(t,\ell)\Longleftrightarrow (s+1,k)\mathrel{\bm\pi}(t+1,\ell)$ for all $s,t\in \mathbb{Z}, k,\ell\in [m]\sqcup[n]$,
\item the  blocks of $\bm\pi$ can be drawn in the strip $R(a,b)$ without intersection.
\end{enumerate}
\end{Prop}
Condition (2) of Proposition \ref{affine diagrams} may still be reformulated purely combinatorially. Denote by $\le$ the  lexicographic order on $\mathbb{Z}\times [m]$ as well as on $\mathbb{Z}\times [n]$ and define a linear order $\preceq$ on $\left.\mathbb{Z}\times [m]\right.\sqcup\left.\mathbb{Z}\times[n]\right.$ by taking the lexicographic order on $\mathbb{Z}\times [m]$, the dual of the lexicographic order on $\mathbb{Z}\times [n]$, and letting every element of $\mathbb{Z}\times[m]$ be above every element of $\mathbb{Z}\times[n]$; more formally
\begin{enumerate}
\item if $x,y\in \mathbb{Z}\times[m]$ then $x\preceq y$ if and only if $x\le y$,
\item if $y\in \mathbb{Z}\times [m]$ and $x\in \mathbb{Z}\times [n]$ then $x\preceq y$,
\item if $x,y\in \mathbb{Z}\times [n]$ then $x\preceq y$ if and only if $x\ge y$.
\end{enumerate}
The mentioned combinatorial characterization may be formulated as follows  \cite{grahamlehrer}.
\begin{Prop}\label{affine diagrams 2} Let $\bm\pi\colon\mathbb{Z}\times [m]\leadsto \mathbb{Z}\times[n]$ be a partition all of whose blocks have size $2$. Then the blocks of $\bm\pi$ can be drawn without intersection in the strip $R(a,b)$ if and only if, for any two blocks $\{x,x'\}$ and $\{y,y'\}$ of $\bm\pi$ the following holds:
$$x\preceq y\preceq x'\mbox{ if and only if }x\preceq y'\preceq x'.$$
\end{Prop}

Altogether, t$A$-diagrams $[m]\mathrel{\contour{black}{${\rightsquigarrow}$}}[n]$ without circle components may as well be represented by partitions $\mathbb{Z}\times[m]\leadsto \mathbb{Z}\times[n]$ all of whose blocks have size $2$ and which satisfy the conditions mentioned in Propositions \ref{affine diagrams} and \ref{affine diagrams 2}. In the following we shall use this description widely. However, we shall usually retain the vocabulary of t$A$-diagrams: we still have `incoming points', `outgoing points', `strings' instead of `vertices' and `blocks'.
Fix $n\in \mathbb{N}$ and consider the partition $\bm\zeta_n\colon\mathbb{Z}\times[n]\leadsto \mathbb{Z}\times [n]$ defined by:
\begin{enumerate}
\item for all $t\in \mathbb{Z}$ and $k=1,\dots,n-1$, $\bm\zeta_n$ connects the incoming point $(t,k)$ with the outgoing point $(t,k+1)$,
\item for all $t\in \mathbb{Z}$, $\bm\zeta_n$ connects the incoming point $(t,n)$ with the outgoing point $(t+1,1)$.
\end{enumerate}
Note that all strings of $\bm\zeta_n$ are transversal and $\bm\zeta_n$ is a generator of the unit group of the monoid $\mathfrak{aLT}^{\mathrm e}_n$. Set  $\bm\lambda_n\colon=\bm\zeta_n^n$; then $\bm\lambda_n$ connects every incoming point $(t,k)$ with the outgoing point $(t+1,k)$, and, more generally, for every $s\in\mathbb{Z}$, $\bm\lambda_n^s$ connects every incoming point $(t,k)$ with the outgoing point $(t+s,k)$.

Next we describe how elements of the form $\bm\lambda_n^s$ act by left and right multiplication. For a partition $\bm\alpha\colon\mathbb{Z}\times[m]\leadsto \mathbb{Z}\times[n]$ a block consisting of the incoming point $(s,k)$ and the outgoing point $(t,\ell)$ (that is, a transversal string connecting the incoming point $(s+\frac{2k-1}{m},a)$ with the outgoing point $(t+\frac{2\ell-1}{n},b)$ in the $R(a,b)$ model) will be denoted as $(s,k)\underset{\bm\alpha}{\!\!\!-\!\!-\!\!\!}(t,\ell)$.
\begin{Lemma}\label{lem:action of lambda} Let $m,n\in\mathbb{N}$, $r\in\mathbb{Z}$ and $\bm\alpha\in \mathfrak{aTL}^{\mathrm e}(m,n)$; then the following hold.
\begin{enumerate}
\item $\bm\lambda_m^r\bm\alpha$, $\bm\alpha$ and $\bm\alpha\bm\lambda_n^r$ have the same in-strings and out-strings.
\item for all $s,t\in \mathbb{Z},\ k\in [m],\ \ell\in [n]$:
$$(s,k)\underset{\bm\lambda_m^r\bm\alpha}{\!\!-\!\!\!-\!\!\!\!-\!\!\!-\!\!\!-\!\!}(t+r,\ell) \Longleftrightarrow (s,k)\underset{\bm\alpha}{\!\!\!-\!\!-\!\!\!}(t,\ell) \Longleftrightarrow
(s,k)\underset{\bm\alpha\bm\lambda_n^r}{\!\!-\!\!\!-\!\!\!\!-\!\!\!-\!\!\!-\!\!}(t+r,\ell),$$
\item $\bm\lambda_m^r\bm\alpha=\bm\alpha\bm\lambda_n^r$.
\end{enumerate}
\end{Lemma}
\begin{proof}
It is clear that $\bm\lambda_m^r\bm\alpha$ and $\bm\alpha$ have the same out-strings. Let $\{(s,k),(s',k')\}$ be an in-string of $\bm\alpha$ (note that $s'\in\{s-1,s,s+1\}$). Since $(t-r,\ell)
\underset{\bm\lambda_m^r}{\!\!-\!\!\!\!-\!\!\!-\!\!\!-\!\!}(t,\ell)$ for all $t\in \mathbb{Z}$ and $\ell\in [m]$ we get
$$(s-r,k)\mathrel{\bm\lambda_m^r}(s,k)\mathrel{\bm\alpha}(s',k')
\mathrel{\bm\lambda_m^r}(s'-r,k')$$
(here $(s-r,k)$ and $(s'-r,k')$ are incoming points while $(s,k)$ and $(s',k')$ are outgoing points of $\bm\lambda_m^r$)
hence $\{(s-r,k), (s'-r,k')\}$  is an in-string of $\bm\lambda_m^r\bm\alpha$ and so is $\{(s,k),(s',k')\}$ by condition (1) of Proposition \ref{affine diagrams}. For the reverse inclusion we note that $\bm\alpha=\bm\lambda_m^{-r}(\bm\lambda_m^r\bm\alpha)$ so that we may employ exactly the same argument. The statement for $\bm\alpha\bm\lambda_n^r$ is proven by an entirely analogous argument. Altogether statement (1) is proved.

Now suppose that $(s,k)\underset{\bm\alpha}{\!\!\!-\!\!-\!\!\!}(t,\ell)$; since $(s-r,k)
\underset{\bm\lambda_m^r}{\!\!-\!\!\!\!-\!\!\!-\!\!\!-\!\!}(s,k)$ we have that
$$(s-r,k)\underset{\bm\lambda_m^r\bm\alpha}{\!\!-\!\!\!-\!\!\!\!-\!\!\!-\!\!\!-\!\!}(t,\ell)$$ and, again by condition (1) of Proposition \ref{affine diagrams},
$$(s,k)\underset{\bm\lambda_m^r\bm\alpha}{\!\!-\!\!\!-\!\!\!\!-\!\!\!-\!\!\!-\!\!}(t+r,\ell).$$
 For the reverse inclusion we may again employ that $\bm\alpha=\bm\lambda_m^{-r}(\bm\lambda_m^r\bm\alpha)$. The statement concerning $\bm\alpha\bm\lambda_n^r$ is proved analogously. Finally, statement (3) is an immediate consequence of (1) and (2).
\end{proof}
The invertible elements $\bm\lambda_m^r$ allow us to describe the kernel of the canonical morphism $\Psi\colon\mathfrak{aTL}^{\mathrm e}\twoheadrightarrow \mathfrak{Ann}$ (at least to a good extent).
\begin{Lemma}\label{lem:map atle->ann} Let $m,n\in\mathbb{N}$, $\bm \alpha,\bm\beta\in\mathfrak{aTL}^{\mathrm e}(m,n)$ and suppose that $\bm\alpha$ and $\bm\beta$ have non-zero rank. Then $\bm\alpha\Psi=\bm\beta\Psi$ if and only if there exists $q\in \mathbb{Z}$ such that $\bm\lambda_m^q\bm\alpha=\bm\beta$.
\end{Lemma}
\begin{proof}
We note that the rank of $\bm\alpha$ coincides with the rank of $\bm\alpha\Psi$ (where, in this context, $\bm\alpha\Psi$ is identified with the partition it induces on $[m]\sqcup[n]$) hence $\bm\alpha$ and $\bm\beta$ have  the same rank if $\bm\alpha\Psi=\bm\beta\Psi$. Moreover $\bm\lambda_m^q$ induces the identity partition $[m]\leadsto [m]$ hence $(\bm\lambda_m^q\bm\alpha)\Psi=\bm\alpha\Psi$ which is the `if' part of the statement. So we are left with the `only if' part.

Let $\bm\alpha,\bm\beta$ be such that $\bm\alpha\Psi=\bm\beta\Psi$ and both of them have rank $r>0$. Let $1\le a_1<a_2<\cdots<a_r\le m$ be such that, for some $t_i\in \mathbb{Z}$ and $b_i\in [n]$,
\begin{equation}\label{trans strings}
(0,a_i)\underset{\bm\alpha}{\!\!\!-\!\!-\!\!\!}(t_i,b_i)\mbox{ for }i=1,\dots,r.
\end{equation}
By periodicity we then also have 
\begin{equation}\label{next string}
(1,a_1)\underset{\bm\alpha}{\!\!\!-\!\!-\!\!\!}(t_1+1,b_1).
\end{equation}
Since the strings~\eqref{trans strings} and~\eqref{next string} do not cross we have
$$t_1\le t_2\le\cdots\le t_r\le t_1+1.$$ It follows that either
\begin{equation}\label{eq:tis_1}
  t_1=t_2=\cdots=t_r  
\end{equation}
or
\begin{equation}\label{eq:tis_2}
 t_1=\cdots =t_i<t_{i+1}=\cdots =t_r=t_1+1\mbox{ for some }i.   
\end{equation}
In case~\eqref{eq:tis_1} we have
\begin{equation}\label{eq:bis_1}
b_1<b_2<\cdots<b_r  
\end{equation}
while in case~\eqref{eq:tis_2},
\begin{equation}\label{eq:bis_2}
 b_{i+1}<\cdots <b_r<b_1<\cdots<b_i.   
\end{equation}

Since $\bm\alpha$ and $\bm\beta$ induce the same partition $[m]\leadsto [n]$, there are transversal strings of the form
$(0,a_i)\underset{\bm\beta}{\!\!\!-\!\!-\!\!\!}(s_i,b_i)$ for some $s_i\in \mathbb{Z}$, $i=1,\dots,r$ and for the same $a_i,b_i$ as for $\bm\alpha$. In case~\eqref{eq:tis_1}, because of~\eqref{eq:bis_1} we have
\[s_1=s_2=\cdots=s_r\] while is case of \eqref{eq:tis_2} we have
\[s_{1}=\cdots=s_i<s_{i+1}=\cdots s_r=s_1+1\] because of~\eqref{eq:bis_2}.
It follows that $q:=s_j-t_j$ does not depend on the index $j$. From Lemma \ref{lem:action of lambda} it follows that $\bm\lambda_m^q\bm\alpha$ and $\bm\beta$ have the same transversal strings.

We are left with showing that $\bm\lambda_m^q\bm\alpha$ and $\bm\beta$ have the same in-strings and out-strings, altogether $\bm\lambda_m^q\bm\alpha=\bm\beta$. We start with in-strings. It suffices to consider in-strings, one endpoint of which is of the form $(0,a)$ (the claim for all other in-strings then follows from shift-invariance, that is, condition (1) of Proposition \ref{affine diagrams}). So, let $(0,a)$ be the endpoint of an in-string $S$, say, of $\bm\alpha$.  If, for some $i$, $a_i<a<a_{i+1}$ then the other endpoint of $S$ must be $(0,b)$ for some $b$ with $a_i<b<a_{i+1}$. Without loss of generality assume that $a<b$. Since $\bm\beta$ induces the same partition $[m]\leadsto [n]$ as $\bm\alpha$, $\bm\beta$ also has an in-string $T$, say, one endpoint of which is $(0,a)$; then there are only two choices for the other endpoint of $T$, namely $(-1,b)$ and $(0,b)$; but since $\bm\beta$ admits a transversal string starting at incoming point $(0,a_i)$ the first case is impossible. Consequently, $T$ connects $(0,a)$ and $(0,b)$, hence coincides with $S$. Now suppose that $1\le a<a_1$; then the second incoming point of $S$ is, say, $(-1,b)$ for some $b>a_r$; by the same reasoning as before, $\bm\beta$ contains an in-string with the same endpoints. Altogether, $\bm\beta$ contains all in-strings of $\bm \alpha$ from which equality follows. By an analogous reasoning one obtains that $\bm\alpha$ and $\bm\beta$ also contain the same out-strings. Indeed, there are $1\le c_1<c_2<\cdots<c_r\le n$ such that the points $(0,c_1)<(0,c_2)<\cdots<(0,c_r)$ are the outgoing points of transversal strings of $\bm\alpha$ as well as $\bm\beta$. (The sequence $c_1,c_2,\dots,c_r$ is just a cyclic shift of the sequence $b_1,b_2,\dots,b_r$ encountered above.) One can then proceed in the same way as for in-strings. Altogether, $\bm\lambda_m^q\bm\alpha$ and $\bm\beta$ admit the same strings and hence are equal.
\end{proof}

From the definitions of the congruences ${\equiv}_1$ and ${\equiv}_2$ it follows immediately that the arrows of $\mathfrak{aTL}$ are uniquely parameterized by pairs $(\bm\alpha,k)$ where $\bm\alpha$ is an arrow in $\mathfrak{aTL}^{\mathrm e}$ and $k\in \mathbb{N}_0$ such that $k=0$ whenever $\bm\alpha$ has non-zero rank. Likewise, arrows of $\mathfrak{aTL}^{\mathrm d}$ can be parameterized by triples $(\bm\alpha,k,k')$ where $k,k'\in \mathbb{N}_0$ and $k=0$ whenever $\bm\alpha$ has non-zero rank. For composable t$A$-diagrams $\bm\alpha$ and $\bm\beta$ let $\mathrm{b}_0(\bm\alpha,\bm\beta)$ be the number of $0$-circles formed in the composition $(\bm\alpha,\bm\beta)\mapsto \bm\alpha\bm\beta$, and likewise, let $\mathrm{b}_\omega(\bm\alpha,\bm\beta)$ be the number of $\omega$-circles formed in this composition. Then composition of arrows in  $\mathfrak{aTL}$ respectively $\mathfrak{aTL}^{\mathrm d}$ is given by the following formulas, and is analogous to (\ref{composition deformed partition category}):
\begin{align}
(\bm\alpha,k)(\bm\beta,\ell)&=(\bm\alpha\bm\beta,k+\ell+\mathrm{b}_\omega(\bm\alpha,\bm\beta))\label{TL,discrete},\\
(\bm\alpha,k,k')(\bm\beta,\ell,\ell')&=(\bm\alpha\bm\beta,k+\ell+\mathrm{b}_\omega(\bm\alpha,\bm\beta),k'+\ell'+\mathrm{b}_0(\bm\alpha,\bm\beta))\label{deformedTL,discrete}.
\end{align}

\subsection{Regular versions of the affine Temperley--Lieb categories}
The annular category $\mathfrak{Ann}$ is a subcategory of the partition category $\mathfrak{P}$ closed under the involution $^*$, whence $\mathfrak{Ann}$ is a regular $*$-category.

As already mentioned, the arrows of the essential affine Temperley--Lieb category $\mathfrak{aTL}^{\mathrm{e}}$ can be viewed to be comprised of all t$A$-diagrams without circle components. The composition of two such diagrams $\bm\alpha$ and $\bm\beta$ within $\mathfrak{aTL}^{\mathrm{e}}$ then is `the same as in $\mathfrak{tAnn}$ except that all circle components which are produced during the composition are dismissed'. In view of this it is not hard to see that $\mathfrak{aTL}^{\mathrm{e}}$ is a regular $*$-category. Indeed,
for an $A(a,b)$-diagram $\alpha\colon[m]\mathrel{\contour{black}{${\rightsquigarrow}$}} [n]$ let $\alpha^*\colon[n]\mathrel{\contour{black}{${\rightsquigarrow}$}}[m]$ be the $A(\frac{1}{b},\frac{1}{a})$-diagram obtained as the image of $\alpha$ under the homeomorphism $\mathbb{C}^*\to\mathbb{C}^*$ given by $z\mapsto \frac{1}{\bar z}$ ($\alpha^*$ is, in a sense, the mirror image of $\alpha$). Finally for a t$A$-diagram $\bm\alpha\colon[m]\mathrel{\contour{black}{${\rightsquigarrow}$}}[n]$ define $\bm\alpha^*\colon[n]\mathrel{\contour{black}{${\rightsquigarrow}$}}[m]$ by choosing a representative $A(a,b)$-diagram $\alpha$ of $\bm\alpha$ and letting $\bm\alpha^*$ be the equivalence class of $\alpha^*$. This provides a unary operation on the set of arrows of $\mathfrak{aTL}^{\mathrm{e}}$ satisfying~\eqref{weak involution} and~\eqref{involution law}, thus $\mathfrak{aTL}^{\mathrm{e}}$ is also a regular $*$-category. Moreover, the canonical morphism $\mathfrak{aTL}^{\mathrm{e}}\twoheadrightarrow \mathfrak{Ann}$ respects the operation $^*$. The operation $^*$ on $\mathfrak{aTL}^{\mathrm e}$ is most easily understood if we interpret the diagram $\bm\alpha$ in terms of an $R(a,b)$-diagram:  reflection on the horizontal line $\mathbb{R}\times\{\frac{a+b}{2}\}$ of the realization of $\bm\alpha$ as an $R(a,b)$-diagram yields $\bm\alpha^*$ as an $R(a,b)$-diagram.

None of the categories $\mathfrak{tAnn}$, $\mathfrak{L}$, $\mathfrak{aTL}^{\mathrm{d}}$, $\mathfrak{aTL}$ and $\mathfrak{Ann}^{\mathrm{d}}$ is regular. However, from the pair, respectively, triple representation of the arrows of $\mathfrak{aTL}$ and $\mathfrak{aTL}^{\mathrm d}$ it is clear that if we allow negative entries in the second, respectively, second and third components, then, just as in (\ref{inverse unary operation deformed}) for $\overline{\mathfrak{P}^{\mathrm d}}$ we get extended categories $\overline{\mathfrak{aTL}^{\mathrm{d}}}$ and $\overline{\mathfrak{aTL}}$ which admit a weak regular involution $^*$ as in the category $\overline{\mathfrak{P}^{\mathrm d}}$. 
More precisely, from~\eqref{TL,discrete} and~\eqref{deformedTL,discrete} we see that for arrows $(\bm\alpha,k)$ in $\overline{\mathfrak{aTL}}$ and $(\bm\alpha,k,\ell)$ in  $\overline{\mathfrak{aTL}^\mathrm{d}}$, respectively,  weak regular involutions can be defined as follows: 
\begin{align*}
(\bm\alpha,k)^*:=&(\bm\alpha^*,-k-\mathrm{b}_\omega(\bm\alpha,\bm\alpha^*)-\mathrm{b}_\omega(\bm\alpha^*,\bm\alpha)) \\
(\bm\alpha,k,\ell)^*:=&(\bm\alpha^*,-k-\mathrm{b}_\omega(\bm\alpha,\bm\alpha^*)-\mathrm{b}_\omega(\bm\alpha^*,\bm\alpha),-\ell-\mathrm{b}_0(\bm\alpha,\bm\alpha^*)-\mathrm{b}_0(\bm\alpha^*,\bm\alpha)).
\end{align*}
The same applies to $\mathfrak{Ann}^{\mathrm{d}}$ since its arrows can be parameterized by pairs $(\alpha,k)$ with $\alpha$ an arrow in $\mathfrak{Ann}$ and $k\in \mathbb{N}_0$. The cases $\mathfrak{tAnn}$ and $\mathfrak{L}$ will be treated in the next subsection.

\subsection{Combinatorial description of $\mathfrak{tAnn}$ and $\mathfrak{L}$, and their regular versions}\label{combinatorial tAnn, L} We are going to look more carefully at the structure of $\mathfrak{tAnn}$. This will be also useful for the treatment of the finite basis problem in Section \ref{finite basis problem}. We need some preparation.

\subsubsection{{Basics of the combinatorial description}}\label{sec:combinatorial description basics} A \emph{unary monoid} is a monoid with an additional unary operation. Let $\mathbb{M}$ be the \emph{free commutative unary monoid generated by the empty set}. We denote the neutral element of $\mathbb{M}$ by $0$, the binary operation by $\underline{\phantom x}+\underline{\phantom x}$ and the unary operation by $(\underline{\phantom{x}})$. The elements of $\mathbb{M}$ are $0$, $(0)$, $(0)+(0)$, $((0))$, $(0)+((0))=((0))+(0)$, $((0)+(0))$, etc. Every such element can be uniquely identified with a finite collection of pairwise nonintersecting circles in the plane, see \cite{dosenpetric} and \cite{kauffmanLogic} (assuming that two such collections of circles are `equal' if they can be mapped onto each other by a homeomorphism of the plane). Let $\bm\alpha$ be a t$A$-diagram $[m]\mathrel{\contour{black}{${\rightsquigarrow}$}}[n]$ and suppose first that $\bm\alpha$ has no $\omega$-circles. We choose some realization of $\bm\alpha$ on an annulus $A$ and observe that the strings of $\bm\alpha$ decompose $A$ into \emph{regions} (we just take into account the strings and ignore the circles!): there are \emph{in-regions} (which are bordered only by in-strings and segments of the outer circle of $A$), \emph{out-regions} (bordered only by out-strings and segments of the inner circle of $A$) and \emph{transversal regions} (whose border contains segments of the inner circle as well as of the outer circle of $A$). In the `$R$-strip model', every region of $\bm\alpha$ uniquely determines a block of a partition $\mathbb{Z}\times([m]+\frac{1}{2})\rightsquigarrow\mathbb{Z}\times([n]+\frac{1}{2})$ and conversely, every region is determined by `its' block --- this will be used in Figures~\ref{regions},~\ref{fig: compatible composition},~\ref{fig:omega-circle},~\ref{fig:increments} and in Section~\ref{sec: tann inverses}. Denote by $\mathrm{reg}(\bm\alpha)$ the set of all regions of $\bm\alpha$ (its representation on $A$, respectively) and note that $\mathrm{reg}(\bm\alpha)=\mathrm{reg}(\underline{\bm\alpha})$. Then, $\bm\alpha$ is uniquely determined by its {skeleton} $\underline{\bm\alpha}$, and the collection of circles sitting on each of its regions. In other words, $\bm\alpha$ is uniquely determined by $\underline{\bm\alpha}$ and a (labeling) function $\mathbf{c}\colon\mathrm{reg}(\underline{\bm\alpha})\to M$. (In case of rank $0$ t$A$-diagrams one cannot use the $R$-strip model in this context. There would be one infinite transversal region possibly carrying infinitely many circles which cannot be encoded by an element of $\mathbb{M}$.)

For later use we note that $\mathbb{M}$ is freely generated as a commutative monoid by its \emph{indecomposable} elements (those of the form $(a)$ for $a\in \mathbb{M}$). {Just as the free commutative unary monoid $\mathbb{M}$ generated by the empty set we may consider the free unary abelian  group $\mathbb{FA}$ generated by the empty set; likewise, $\mathbb{FA}$ is freely generated as an abelian group by its indecomposable elements and $\mathbb{M}$ naturally embeds into $\mathbb{FA}$. Note that $\mathbb{FA}$ is not the abelian group freely generated by $\mathbb{M}$ (the \emph{Grothendieck group} of $\mathbb{M}$), the latter is a proper subgroup of $\mathbb{FA}$.}

In case the t$A$-diagram $\bm\alpha$ admits $\omega$-circles the latter create new regions which may again contain collections of circles. However, this can happen only in case the t$A$-diagram $\bm\alpha$ in question has no transversal string. In that case, $\underline{\bm\alpha}$ has a unique transversal region and the entire collection of  $\omega$-circles and $0$-circles sitting in that  region of $\underline{\bm\alpha}$ may as well be interpreted as a label, this time not from $\mathbb{M}$ but from a more complicated monoid: let $\mathbb{O}$ be the free cyclic monoid generated by a single $\omega$-circle \tikz[scale=0.1]{ \draw circle[radius=1cm];}. We write the composition in the monoid $\mathbb{O}$ as
\tikz[scale=0.1]{\draw circle[radius=1cm];}+\tikz[scale=0.1]{ \draw circle[radius=1cm];} = \tikz[scale=0.1,baseline={(0,-0.3)}]{ \draw (0,-2) circle[radius=2.4cm];\draw(0,-2) circle[radius=1cm];}.

The collection of all circles in the (unique) transversal region of $\bm\alpha$ can be uniquely described by a single element of the free product {$\mathbb{O}\mathop{*}\mathbb{M}$} of the monoids $\mathbb{O}$ and $\mathbb{M}$ (the free product being taken within the class of all monoids). A typical element in $\mathbb{O}\mathop{*}\mathbb{M}$ can be written as
\begin{equation*}\label{element of O*C}
c_0+\tikz[scale=0.1]{\draw circle[radius=1cm];}+c_1+\tikz[scale=0.1]{\draw circle[radius=1cm];}+c_2\cdots +c_{k-1}+\tikz[scale=0.1]{\draw circle[radius=1cm];}+ c_k
\end{equation*}
with $c_0,\dots,c_k\in \mathbb{M}$ and the operation $+$  here is no longer commutative. We shall represent this element graphically as in Figure~\ref{fig:rank0label}; recall that the elements $c_i$ are finite collections of circles, they may be moved  inside their respective annulus region but not beyond. The dotted lines indicate the inner and outer circles of the annulus.
\begin{figure}
\begin{tikzpicture}[scale=0.6]
\draw[dotted] circle [radius=0.25cm];
\draw[thick] circle [radius=1cm];
\draw[thick] circle [radius=2.4cm];
\draw[thick] circle [radius=5cm];
\draw[thick] circle [radius=4.3cm];
\draw (-0.6,0) node {$c_k$};
\draw(-1.7,0) node {$c_{k-1}$};
\draw(-3.4,0) node {$\dots$};
\draw(-4.6,0) node {$c_{1}$};
\draw(-5.4,0) node {$c_0$};
\draw[dotted] circle [radius=6cm];
\end{tikzpicture}
\caption{Graphical representation of the element \\ \centering{$c_0+\bigcirc+c_1+\cdots+c_{k-1}+\bigcirc+c_k$.}}\label{fig:rank0label}
\end{figure}

Altogether, in order to complete the above characterization: every t$A$-diagram $\bm\alpha$ is uniquely determined by
\begin{enumerate}
\item its skeleton $\underline{\bm\alpha}$ which is a member of $\mathfrak{aTL}^{\mathrm{e}}$
\item a labeling function $\mathbf{c}\colon\mathrm{reg}(\underline{\bm\alpha})\to \mathbb{O}\mathop{*}\mathbb{M}$ such that $\mathbf{c}(R)\notin \mathbb{M}$ can occur only if $\bm\alpha$ (equivalently $\underline{\bm\alpha}$) has no transversal string and $R$ is the (then unique) transversal region of $\underline{\bm\alpha}$, and $\mathbf{c}(R)\in \mathbb{M}$ in all other cases.
\end{enumerate}
Conversely every pair $(\bm\alpha,\mathbf{c})$ where $\bm\alpha\in \mathfrak{aTL}^{\mathrm{e}}$ and $\mathbf c$ satisfies condition (2) above determines a unique element of $\mathfrak{tAnn}$. Similarly as for cobordisms (\ref{def:composition in 2Cob discrete}) one can express the composition in $\mathfrak{tAnn}$ in terms of its pair representation. One then can define a regular category $\overline{\mathfrak{tAnn}}$ by use of the pair representation $(\bm\alpha,\mathbf{c})$ and allowing $\mathbf{c}$ to have values in the group free product $\mathbb{Z}*\mathbb{FA}$ rather than the monoid $\mathbb{O}\mathop{*}\mathbb{M}$ (again subject to a properly modified condition (2) above)  in which $\mathfrak{tAnn}$ is naturally embedded. From this, a regular extension  $\overline{\mathfrak{L}}$ of $\mathfrak{L}$ follows immediately. In the latter case, all labels are from $\mathbb{FA}$. More details to this construction are given in Section~\ref{sec: tann inverses}.

\subsubsection{{Composition of arrows}}\label{sec:description of arrow composition} {Here we intend to roughly discuss how the composition of arrows in $\mathfrak{tAnn}$ and $\mathfrak{L}$ works.}
{For the following discussion it is convenient to view elements $\bm\al\in \mathfrak{aTL}^{\mathrm e}$ in their  $R$-strip interpretation.}  We need to analyze how the labels of the regions of a product $\bm\alpha\bm\beta$ of two t$A$-diagrams are formed. Every region $C$ in $\bm\alpha\bm\beta$ is the union of certain regions $A_i$ of $\bm\alpha$ and $B_j$ of $\bm\beta$.  \textsl{Formally} rather similar to the situation of cobordisms, the label of the region $C$ then is the sum of the labels of $A_i$ and of $B_j$ plus, perhaps, a certain `increment' which is formed where the regions $A_i$ and $B_j$ have a common border.  Figure \ref{a+b+c} shows an example of how a product of three members $\bm\alpha$, $\bm\beta$ and $\bm\gamma$ of $\mathfrak{tAnn}_{21}$ is formed; {more precisely, the picture shows the \emph{product graph} (in the sense of \cite{pseudovarbrauer}) of $\bm\al\circ\bm\beta\circ\bm\gamma:=(\bm\alpha\cup \bm\beta\cup \bm\gamma)^t$ (according to~\eqref{eq:partition composition}). The result $\bm\al\bm\beta\bm\gamma$ of the concatenation is depicted  in Figure \ref{abc}.}
\begin{figure}[ht]
\begin{tikzpicture}[scale=0.5]
\foreach \x in {0,...,21} \foreach \y in {0,4,8,12} \filldraw (\x,\y) circle (2.5pt);
\draw[dotted](-1,0)--(22,0);
\draw[dotted](-1,4)--(22,4);
\draw[dotted](-1,8)--(22,8);
\draw[dotted](-1,12)--(22,12);
\draw (0,12.5) node {$1$};
\draw (9,12.5) node {$10$};
\draw (20,12.5) node {$21$};
\draw (21,12.5) node {$1$};
\draw (15,-.5) node {$16$};
\draw (-1,10) node {$\bm\alpha$};
\draw (-1,6) node {$\bm\beta$};
\draw (-1,2) node {$\bm\gamma$};
\draw[very thick] (9,12)--(15,8);
\draw (10,12)--(16,8);
\draw (0,12)--(0,8);
\draw (21,12)--(21,8);

\draw (1,12)  to [out=315,in=225] (2,12);
\draw (3,12)  to [out=315,in=225] (4,12);
\draw (5,12)  to [out=315,in=225] (8,12);
\draw (6,12)  to [out=315,in=225] (7,12);
\draw (11,12)  to [out=315,in=225] (20,12);
\draw (12,12)  to [out=315,in=225] (19,12);
\draw (13,12)  to [out=315,in=225] (14,12);
\draw (15,12)  to [out=315,in=225] (18,12);
\draw (16,12)  to [out=315,in=225] (17,12);

\draw(4.5,10) node{$f$};
\draw((18.5,8.3) node {$g$};

\draw (1,8)  to [out=90,in=90] (4,8);
\draw(2.5,8.5) node{$a$};
\draw (2,8)  to [out=45,in=135] (3,8);
\draw[very thick] (5,8)  to [out=45,in=135] (12,8);
\draw (6,8)  to [out=45,in=135] (11,8);
\draw (7,8)  to [out=45,in=135] (10,8);
\draw [very thick] (8,8)  to [out=45,in=135] (9,8);
\draw (13,8)  to [out=45,in=135] (14,8);
\draw (17,8) .. controls (17,10) and (20,10).. (20,8);
\draw[very thick] (18,8)  .. controls (18,9) and (19,9) ..  (19,8);

\draw [very thick](9,8)--(15,4);
\draw (10,8)--(16,4);
\draw (0,8)--(0,4);
\draw (21,8)--(21,4);

\draw (1,8)  to [out=315,in=225] (2,8);
\draw (3,8)  to [out=315,in=225] (4,8);
\draw [very thick](5,8)  to [out=315,in=225] (8,8);
\draw (6,8)  to [out=315,in=225] (7,8);
\draw (11,8)  to [out=315,in=225] (20,8);
\draw[very thick] (12,8)  to [out=315,in=225] (19,8);
\draw (13,8)  to [out=315,in=225] (14,8);
\draw [very thick](15,8)  to [out=315,in=225] (18,8);
\draw (16,8)  to [out=315,in=225] (17,8);

\draw(4.5,6) node{$h$};
\draw((14.5,3.4) node {$k$};of which is

\draw (1,4)  to [out=45,in=135] (4,4);
\draw (2,4)  to [out=45,in=135] (3,4);
\draw [very thick](5,4)  to [out=45,in=135] (12,4);
\draw (6,4)  to [out=45,in=135] (11,4);
\draw (7,4)  to [out=45,in=135] (10,4);
\draw [very thick](8,4)  to [out=45,in=135] (9,4);
\draw (13,4)  to [out=45,in=135] (14,4);
\draw (17,4)  to [out=45,in=135] (20,4);
\draw [very thick](18,4)  to [out=45,in=135] (19,4);

\draw [very thick](9,4)--(15,0);
\draw (10,4)--(16,0);
\draw (0,4)--(0,0);
\draw (21,4)--(21,0);

\draw (1,4) .. controls (1,3) and (2,3) .. (2,4);
\draw (3,4)  .. controls (3,3) and (4,3) .. (4,4);
\draw (1.5,3.6) node {$u$};
\draw(3.5,3.6) node {$v$};
\draw [very thick](5,4)  to [out=315,in=225] (8,4);
\draw (6,4)  to [out=315,in=225] (7,4);
\draw (11,4)  to [out=315,in=225] (20,4);
\draw [very thick](12,4)  to [out=315,in=225] (19,4);
\draw (13,4)  to [out=315,in=225] (14,4);
\draw [very thick](15,4)  to [out=315,in=225] (18,4);
\draw (16,4)  to [out=315,in=225] (17,4);

\draw(4.5,2) node{$\ell$};

\draw (1,0)  to [out=45,in=135] (4,0);
\draw (2,0)  to [out=45,in=135] (3,0);
\draw (5,0)  to [out=45,in=135] (12,0);
\draw (6,0)  to [out=45,in=135] (11,0);
\draw (7,0)  to [out=45,in=135] (10,0);
\draw (8,0)  to [out=45,in=135] (9,0);
\draw (13,0)  to [out=45,in=135] (14,0);
\draw (17,0)  to [out=45,in=135] (20,0);
\draw (18,0)  to [out=45,in=135] (19,0);
\end{tikzpicture}
\caption{$\bm\alpha\circ\bm\beta\circ\bm\gamma$}\label{a+b+c}
\end{figure}

In Figure \ref{abc} we look at the transversal region shown on the left from the transversal string starting at incoming point $10$ and ending at outgoing point~$16$. Its label in $\bm\alpha\bm\beta\bm\gamma$ is
\[f+(0)+g+(a)+h+k+(0)+(u+v)+\ell\]
where the sum over the labels of the involved regions is
\[f+g+h+k+\ell,\]
while the `increment' generated in the composition $\bm\alpha\circ\bm\beta$ is $(a)+(0)$ and that resulting from the composition $\bm\beta\circ\bm\gamma$ is $(u+v)+(0)$. Note that the example is chosen such that $\underline{\bm\alpha}=\underline{\bm\beta}=\underline{\bm\gamma}$ and is an idempotent if composed in $\mathfrak{aTL}^{\mathrm e}$.

\begin{figure} [ht]
\begin{tikzpicture}[scale=0.5]
\foreach \x in {0,...,21} \foreach \y in {0,4} \filldraw (\x,\y) circle (2.5pt);
\draw[dotted](-1,0)--(22,0);
\draw[dotted](-1,4)--(22,4);
\draw (0,4.5) node {$1$};
\draw (9,4.5) node {$10$};
\draw (20,4.5) node {$21$};
\draw (21,4.5) node {$1$};
\draw (15,-.5) node {$16$};
\draw [very thick](9,4)--(15,0);
\draw (10,4)--(16,0);
\draw (0,4)--(0,0);
\draw (21,4)--(21,0);
\draw (1,4) to [out=315, in=225] (2,4);
\draw (3,4) to [out=315, in =225] (4,4);
\draw (5,4)  to [out=315,in=225] (8,4);
\draw (6,4)  to [out=315,in=225] (7,4);
\draw (11,4)  to [out=315,in=225] (20,4);
\draw (12,4)  to [out=315,in=225] (19,4);
\draw (13,4)  to [out=315,in=225] (14,4);
\draw (15,4)  to [out=315,in=225] (18,4);
\draw (16,4)  to [out=315,in=225] (17,4);
\draw(4,2.9) node{$f+(0)+g+(a)+h$};
\draw(7,1.9) node{$+k+(0)+(u+v)+\ell$};
\draw (1,0)  to [out=45,in=135] (4,0);
\draw (2,0)  to [out=45,in=135] (3,0);
\draw (5,0)  to [out=45,in=135] (12,0);
\draw (6,0)  to [out=45,in=135] (11,0);
\draw (7,0)  to [out=45,in=135] (10,0);
\draw (8,0)  to [out=45,in=135] (9,0);
\draw (13,0)  to [out=45,in=135] (14,0);
\draw (17,0)  to [out=45,in=135] (20,0);
\draw (18,0)  to [out=45,in=135] (19,0);
\end{tikzpicture}

\caption{$\bm\alpha\bm\beta\bm\gamma$}\label{abc}
\end{figure}
{This composition is reminiscent of the composition of the arrows of $2\mathfrak{Cob}$ and $\overline{2\mathfrak{Cob}}$ interpreted as partitions labeled by elements of $\mathbb{N}_0$ and $\mathbb{Z}$, respectively. We shall discuss this aspect further.}

{The collection of \emph{regions} of a t$A$-diagram $\bm\alpha$ is obtained as the complement
of $\bm\alpha$ in $R(a,b)$ which forms a set of simply connected subsets of $R(a,b)$. By doing so for the t${A}$-diagram $\underline{\bm\alpha}$ in Figure \ref{a+b+c} we see in Figure \ref{regions} that the set of regions $\mathrm{reg}(\underline{\bm\alpha})$ gives rise to a  partition $\mathbb{Z}\times ([21]+\frac{1}{2})\leadsto \mathbb{Z}\times ([21]+\frac{1}{2})$ which is conveniently also denoted by $\mathrm{reg}(\underline{\bm\alpha})$. This partition is invariant under the shift map $(s,k)\mapsto (s+1,k)$ and is uniquely determined by and uniquely determines the skeleton $\underline{\bm\alpha}$ of the original t$A$-diagram $\bm\alpha$. It follows that the original t$A$-diagram $\bm\alpha$ is uniquely determined by the pair $(\mathrm{reg}(\underline{\bm\alpha}),\mathrm{c}_{\bm\alpha})$: the partition
$\mathrm{reg}(\underline{\bm\alpha})\colon \mathbb{Z}\times ([21]+\frac{1}{2})\leadsto \mathbb{Z}\times ([21]+\frac{1}{2})$ together with a labeling function $\mathrm{c}_{\bm\alpha}\colon \mathrm{reg}(\underline{\bm\alpha})\to \mathbb{O}\mathop{*}\mathbb{M}$ where the domain of $\mathrm{c}_{\bm\alpha}$ here is understood as the set of blocks of the partition $\mathrm{reg}(\underline{\bm\alpha})$. Recall that only in case of a rank $0$ diagram, the value $c_{\bm\alpha}(C)$ of the unique transversal region $C$ is in  $\mathbb{O}\mathop{*}\mathbb{M}$, all other values of $c_{\bm\alpha}$ are in $\mathbb{M}$. In this sense it is justified to say that the category $\mathfrak{tAnn}$ can be described entirely by combinatorial data. Throughout the following, we shall identify every set $[m]+\frac{1}{2}$ with $[m]$ via the mapping $x+\frac{1}{2}\mapsto x$.
}
\begin{figure} [ht]
\begin{tikzpicture}[scale=0.5]
\draw[thick](-1,0)--(22,0);
\draw[thick](-1,4)--(22,4);
\filldraw[gray] (-1,0)--(22,0)--(22,4)--(-1,4);
\draw (0.5,4.5) node {$1$};
\draw (20.5,4.5) node {$21$};
\foreach \x in {-1,...,22} \foreach \y in {0,4} \filldraw[white] (\x,\y) circle (2pt);
\foreach \x in {-1,...,21} \foreach \y in {0,4} \filldraw (\x+.5,\y) circle (3pt);
\draw[ultra thick, white] (9,4)--(15,0);
\draw [very thick, white] (10,4)--(16,0);
\draw[ultra thick, white] (0,4)--(0,0);
\draw[ultra thick, white] (21,4)--(21,0);
\draw[ultra thick, white] (1,4) to [out=315, in=225] (2,4);
\draw[ultra thick, white] (3,4) to [out=315, in =225] (4,4);
\draw[ultra thick, white] (5,4)  to [out=315,in=225] (8,4);
\draw[ultra thick, white] (6,4)  to [out=315,in=225] (7,4);
\draw[ultra thick, white] (11,4)  to [out=315,in=225] (20,4);
\draw[ultra thick, white] (12,4)  to [out=315,in=225] (19,4);
\draw[ultra thick, white] (13,4)  to [out=315,in=225] (14,4);
\draw[ultra thick, white] (15,4)  to [out=315,in=225] (18,4);
\draw[ultra thick, white] (16,4)  to [out=315,in=225] (17,4);
\draw[ultra thick, white] (1,0)  to [out=45,in=135] (4,0);
\draw[ultra thick, white] (2,0)  to [out=45,in=135] (3,0);
\draw[ultra thick, white] (5,0)  to [out=45,in=135] (12,0);
\draw[ultra thick, white] (6,0)  to [out=45,in=135] (11,0);
\draw[ultra thick, white] (7,0)  to [out=45,in=135] (10,0);
\draw[ultra thick, white] (8,0)  to [out=45,in=135] (9,0);
\draw[ultra thick, white] (13,0)  to [out=45,in=135] (14,0);
\draw[ultra thick, white] (17,0)  to [out=45,in=135] (20,0);
\draw[ultra thick, white] (18,0)  to [out=45,in=135] (19,0);
\end{tikzpicture}
\caption{Regions}\label{regions}
\end{figure}

{ In this way one can in general encode any t$A$-diagram $\bm\alpha\colon [m]\leadsto [n]$ by means of a suitable (periodic) partition $\mathrm{reg}(\underline{\alpha})\colon\mathbb{Z}\times[m]\leadsto \mathbb{Z}\times [n]$ and a labeling function $\mathrm{c}_{\bm\alpha}\colon\mathrm{reg}(\underline{\bm\alpha})\to \mathbb{O}\mathop{*}\mathbb{M}$. The correspondence $\bm\alpha\mapsto (\mathrm{reg}(\underline{\bm\alpha}),\mathrm{c}_{\bm\alpha})$  can be made to an isomorphism of categories by defining an appropriate composition of arrows of the form $(\mathrm{reg}(\underline{\bm\alpha}),\mathrm{c}_{\bm\alpha})$. Indeed, the collection of all such pairs $(\mathrm{reg}(\underline{\bm\alpha}),\mathrm{c}_{\bm\alpha})$ forms the arrows of a suitable category where the composition of two composable arrows of this form, say $(\alpha,\mathrm{c})$ and $(\beta,\mathrm{d})$ is given by $(\alpha,\mathrm{c})(\beta,\mathrm{d})=(\alpha\circ\beta,\mathrm{e})$ where, in a first approximation,  $\alpha\circ\beta$ is the composition in the appropriate category of partitions, and the labeling function $\mathrm{e}\colon \alpha\circ\beta\to \mathbb{O}\mathop{*}\mathbb{M}$ is obtained similarly as in the case of cobordisms: for every block $C$ of $\alpha\circ\beta$, $\mathrm{e}(C)$ is obtained as the sum of all $\mathrm{c}$- and $\mathrm{d}$-values of the blocks of $\alpha$ and $\beta$, respectively, involved in the formation of the new block $C$ plus a certain increment which is formed through the merging process of $C$ and which we shall discuss later.}

{When the product $\bm\alpha\bm\beta$ of two arrows $\bm\alpha$ and $\bm\beta$ does not lead to the formation of $\omega$-circles then the map $\bm\alpha\mapsto(\mathrm{reg}(\underline{\bm\al}),\mathrm{c}_{\bm\alpha})$ is compatible with the composition, that is, $\mathrm{reg}(\underline{\bm\alpha\bm\beta})=\mathrm{reg}(\underline{\bm\alpha})\circ\mathrm{reg}(\underline{\bm\beta})$; an example of arrows $[6]\to [6]$ is presented in Figure \ref{fig: compatible composition}.
}
\begin{figure}[ht]
\begin{tikzpicture}[scale=0.45]

\foreach \x in {1,...,12} \foreach \y in {-1,1,2,4,6} \filldraw (\x,\y) circle (2.5pt);
\draw[dotted](1,2)--(12,2);
\draw[dotted](1,4)--(12,4);
\draw[dotted](1,6)--(12,6);
\draw(4,6.5) node{$1$};
\draw(9,6.5) node{$6$};
\draw(10,6.5) node{$1$};

\draw (1,6)--(0.5,5.66);
\draw (2,6) to [out=315,in=225] (3,6);
\draw(5,6) to [out=315,in=225] (6,6);
\draw(8,6) to [out=315,in=225] (9,6);
\draw(11,6) to [out=315,in=225] (12,6);

\draw (1,4) to [out=45,in=135] (2,4);
\draw(2,4) to [out=315,in=225] (3,4);
\draw(3,4) to (4,6);
\draw (4,4) to  (7,6);
\draw(4,4) to [out=315,in=225] (5,4);
\draw(5,4) to [out=45,in=135] (6,4);
\draw(6,4) to (6,2);
\draw(7,4)to (7,2);
\draw (7,4) to [out=45,in=135] (8,4);
\draw (8,4) to [out=315,in=225] (9,4);
\draw(9,4) to (10,6);
\draw(11,4) to [out=45,in=135] (12,4);
\draw(10,4) to [out=315,in=225] (11,4);
\draw(10,4) to (12.5,5.66);
\draw(1,2) to (1,4);
\draw(2,2)to[out=45,in=135](3,2);
\draw(4,2) to [out=45,in=135] (5,2);
\draw(8,2) to [out=45,in=135] (9,2);
\draw(10,2) to [out=45, in=135] (11,2);
\draw(12,2) to (12,4);

\draw[dotted](1,-1)--(12,-1);
\draw[dotted](1,1)--(12,1);

\draw (2,1) to [out=315,in=225] (3,1);
\draw (5,1) to [out=315,in=225] (6,1);
\draw (8,1) to [out=315,in=225] (9,1);
\draw (11,1) to [out=315,in=225] (12,1);

\draw (2,-1) to [out=45,in=135] (3,-1);
\draw(4,-1) to [out=45,in=135] (5,-1);
\draw(8,-1) to [out=45,in=135] (9,-1);
\draw(10,-1) to [out=45,in=135] (11,-1);
\draw(1,-1) to ((4,1);
\draw(6,-1) to (7,1);
\draw (7,-1) to (10,1);
\draw (1,1) to (0.5,0.66);
\draw (12,-1) to (12.5,-0.66);
\end{tikzpicture}
\begin{tikzpicture}[scale=0.45]
\draw(-2,1)node{};
\filldraw[gray] (0.625,2)--(12.5,2)--(12.5,6)--(0.625,6);
\filldraw[gray] (0.625,-1)--(12.5,-1)--(12.5,1)--(0.625,1);
\draw[thick](0.625,2)--(12.375,2);
\draw[thick](0.625,4)--(12.375,4);
\draw[thick](0.625,6)--(12.375,6);
\draw[thick](0.625,-1)--(12.375,-1);
\draw[thick](0.625,1)--(12.375,1);
\foreach \x in {1,...,12} \foreach \y in {-1,1,2,4,6} \filldraw (\x+.5,\y) circle (2.5pt);
\foreach \x in {1,...,12} \foreach \y in {-1,1,2,4,6} \filldraw[white] (\x,\y) circle (2.5pt);
\draw(4,6.5) node{$1$};
\draw(9,6.5) node{$6$};
\draw(10,6.5) node{$1$};
\draw[ultra thick,white] (1,6)--(0.5,5.66);
\draw[ultra thick, white] (2,6) to [out=315,in=225] (3,6);
\draw[ultra thick, white](5,6) to [out=315,in=225] (6,6);
\draw[ultra thick, white](8,6) to [out=315,in=225] (9,6);
\draw[ultra thick, white] (11,6) to [out=315,in=225] (12,6);
\draw[ultra thick, white] (1,4) to [out=45,in=135] (2,4);
\draw[ultra thick, white](2,4) to [out=315,in=225] (3,4);
\draw[ultra thick, white](3,4) to (4,6);
\draw[ultra thick, white] (4,4) to  (7,6);
\draw[ultra thick, white](4,4) to [out=315,in=225] (5,4);
\draw[ultra thick, white](5,4) to [out=45,in=135] (6,4);
\draw[ultra thick, white](6,4) to (6,2);
\draw[ultra thick, white](7,4) to (7,2);
\draw[ultra thick, white] (7,4) to [out=45,in=135] (8,4);
\draw[ultra thick, white] (8,4) to [out=315,in=225] (9,4);
\draw[ultra thick, white](9,4) to (10,6);
\draw[ultra thick, white](11,4) to [out=45,in=135] (12,4);
\draw[ultra thick, white](10,4) to [out=315,in=225] (11,4);
\draw[ultra thick, white](10,4) to (12.5,5.66);
\draw[ultra thick, white](1,2) to (1,4);
\draw[ultra thick, white](2,2) to[out=45,in=135] (3,2);
\draw[ultra thick, white](4,2) to [out=45,in=135] (5,2);
\draw[ultra thick, white](8,2) to [out=45,in=135] (9,2);
\draw[ultra thick, white](10,2) to [out=45, in=135] (11,2);
\draw[ultra thick, white](12,2) to (12,4);
\draw[ultra thick, white] (2,1) to [out=315,in=225] (3,1);
\draw[ultra thick, white] (5,1) to [out=315,in=225] (6,1);
\draw[ultra thick, white] (8,1) to [out=315,in=225] (9,1);
\draw[ultra thick, white] (11,1) to [out=315,in=225] (12,1);
\draw[ultra thick, white] (2,-1) to [out=45,in=135] (3,-1);
\draw[ultra thick, white](4,-1) to [out=45,in=135] (5,-1);
\draw[ultra thick, white](8,-1) to [out=45,in=135] (9,-1);
\draw[ultra thick, white](10,-1) to [out=45,in=135] (11,-1);
\draw[ultra thick, white](1,-1) to ((4,1);
\draw[ultra thick, white](6,-1) to (7,1);
\draw[ultra thick, white] (7,-1) to (10,1);
\draw[ultra thick, white] (1,1) to (0.5,0.66);
\draw[ultra thick, white] (12,-1) to (12.5,-0.66);
\end{tikzpicture}
\caption{Composition in $\mathfrak{tAnn}$ and in the corresponding partition category}\label{fig: compatible composition}
\end{figure}

{Problems which obstruct the map $\bm\alpha\mapsto(\mathrm{reg}(\underline{\bm\al}),\mathrm{c}_{\bm\alpha})$ from being compatible with composition of partitions in the first component only occur in products $\bm\alpha\bm\beta$ which produce $\omega$-circles. In particular these problems do not occur in the subcategories $\mathfrak{tAnn}^{\mathrm{odd}}$ and $\mathfrak{L}$. Figure \ref{fig:omega-circle} shows an example where the composition of two arrows $\bm\alpha$ and $\bm\beta$ produces an $\omega$-circle (left hand side). This is reflected by the resulting partition (right hand side) in that it has no longer a transversal block but instead two infinite blocks: one incoming and one outgoing\footnote{In the cylinder or the annulus model, the `infinite' incoming and outgoing block are not infinite, of course; in these models, `infinite' has to be replaced by block being maxima in the ${\sqsubseteq}_{in}$ respectively ${\sqsubseteq}_{out}$-orders, to be defined in Section~\ref{sec:rank0inverses}}. In order to match with the t$A$-diagram $\bm\alpha\bm\beta$ one needs to merge the two infinite blocks to one infinite transversal block on the partition side. On this side, the neglected $\omega$-circle has to be taken into account as increment in the corresponding labeling function $\mathrm{c}_{\bm\alpha\bm\beta}$. Hence, in order that the map $\bm\alpha\mapsto (\mathrm{reg}(\underline{\bm\alpha}),\mathrm{c}_{\bm\alpha})$ is compatible with composition, in such cases the partition $\mathrm{reg}(\underline{\bm\alpha})\circ\mathrm{reg}(\underline{\bm\beta})$ has to be replaced with the one obtained by merging its two infinite blocks to one infinite transversal block. That block then has to be labeled by the appropriate element of $\mathbb{O}\mathop{*}\mathbb{M}$ indicating the $\omega$-circles (and the circles between them) obtained in the concatenation process of the diagrams $\bm\alpha$ and $\bm\beta$.
}
\begin{figure}[ht]
\begin{tikzpicture}[scale=0.45]
\draw[dotted](0.625,2)--(12.375,2);
\draw[dotted](0.625,4)--(12.375,4);
\draw[dotted](0.625,6)--(12.375,6);
\draw[dotted](0.625,-1)--(12.375,-1);
\draw[dotted](0.625,1)--(12.375,1);
\foreach \x in {1,...,12} \foreach \y in {-1,1,2,4,6} \filldraw (\x,\y) circle (2.5pt);
\draw(4,6.5) node{$1$};
\draw(9,6.5) node{$6$};
\draw(10,6.5) node{$1$};
\draw (1,6) to [out=315,in=225] (2,6);
\draw(3,6) to [out=315,in=225] (4,6);
\draw(7,6) to [out=315,in=225] (8,6);
\draw (9,6) to [out=315,in=225] (10,6);
\draw (1,4) to [out=45,in=135] (4,4);
\draw(2,4) to [out=45,in=135] (3,4);
\draw(5,4) to (5,6);
\draw(6,4) to  (6,6);
\draw (11,4) to (11,6);
\draw (12,4) to (12,6);
\draw(4,4) to [out=315,in=225] (7,4);
\draw(5,4) to [out=315,in=225] (6,4);
\draw(1,4) to (0.5,3.66);
\draw (7,4) to [out=45,in=135] (10,4);
\draw (8,4) to [out=45,in=135] (9,4);
\draw(11,4) to [out=315,in=225] (12,4);
\draw(10,4) to [out=315,in=225] (12.5,3.66);
\draw (1,2) to (0.5,2.33);
\draw(12,2) to  (12.5,2.33);
\draw(2,2) to (2,4);
\draw(3,2) to (3,4);
\draw(4,2) to [out=45,in=135] (5,2);
\draw (8,2) to (8,4);
\draw (9,2) to (9,4);
\draw(6,2) to [out=45,in=135] (7,2);
\draw(10,2) to [out=45, in=135] (11,2);
\draw(1,1) to [out=315,in=225] (2,1);
\draw (3,1) to [out=315,in=225] (4,1);
\draw (5,1) to [out=315,in=225] (6,1);
\draw (7,1) to [out=315,in=225] (8,1);
\draw (9,1) to [out=315,in=225] (10,1);
\draw (11,1) to [out=315,in=225] (12,1);
\draw (1,-1) to   (0.5,-0.66);
\draw (12,-1) to  (12.5,-0.66);
\draw(4,-1) to [out=45,in=135] (5,-1);
\draw(2,-1) to [out=45,in=135] (3,-1);
\draw(8,-1) to [out=45, in=135] (9,-1);
\draw(6,-1) to [out=45,in=135] (7,-1);
\draw(10,-1) to [out=45, in=135] (11,-1);
\end{tikzpicture}
\begin{tikzpicture}[scale=0.45]
\draw(-2,1)node{};
\filldraw[gray] (0.625,2)--(12.5,2)--(12.5,6)--(0.625,6);
\filldraw[gray] (0.625,-1)--(12.5,-1)--(12.5,1)--(0.625,1);
\draw[thick](0.625,2)--(12.375,2);
\draw[thick](0.625,4)--(12.375,4);
\draw[thick](0.625,6)--(12.375,6);
\draw[thick](0.625,-1)--(12.375,-1);
\draw[thick](0.625,1)--(12.375,1);
\foreach \x in {1,...,12} \foreach \y in {-1,1,2,4,6} \filldraw (\x+.5,\y) circle (2.5pt);
\foreach \x in {1,...,12} \foreach \y in {-1,1,2,4,6} \filldraw[white] (\x,\y) circle (2.5pt);
\draw(4,6.5) node{$1$};
\draw(9,6.5) node{$6$};
\draw(10,6.5) node{$1$};
\draw[ultra thick, white] (1,6) to [out=315,in=225] (2,6);
\draw[ultra thick, white](3,6) to [out=315,in=225] (4,6);
\draw[ultra thick, white](7,6) to [out=315,in=225] (8,6);
\draw[ultra thick, white] (9,6) to [out=315,in=225] (10,6);
\draw[ultra thick, white] (1,4) to [out=45,in=135] (4,4);
\draw[ultra thick, white](2,4) to [out=45,in=135] (3,4);
\draw[ultra thick, white](5,4) to (5,6);
\draw[ultra thick, white] (6,4) to  (6,6);
\draw[ultra thick, white] (11,4) to (11,6);
\draw[ultra thick, white] (12,4) to (12,6);
\draw[ultra thick, white](4,4) to [out=315,in=225] (7,4);
\draw[ultra thick, white](5,4) to [out=315,in=225] (6,4);
\draw[ultra thick, white](1,4) to [out=225,in=240] (0.5,3.66);
\draw[ultra thick, white] (7,4) to [out=45,in=135] (10,4);
\draw[ultra thick, white] (8,4) to [out=45,in=135] (9,4);
\draw[ultra thick, white](11,4) to [out=315,in=225] (12,4);
\draw[ultra thick, white](10,4) to [out=315,in=225] (12.5,3.66);
\draw[ultra thick, white] (1,2) to [out=135,in=170]  (0.5,2.33);
\draw[ultra thick, white] (12,2) to [out=45, in=10] (12.5,2.33);
\draw[ultra thick, white](2,2) to (2,4);
\draw[ultra thick, white](3,2) to (3,4);
\draw[ultra thick, white](4,2) to [out=45,in=135] (5,2);
\draw[ultra thick, white] (8,2) to (8,4);
\draw[ultra thick, white] (9,2) to (9,4);
\draw[ultra thick, white](6,2) to [out=45,in=135] (7,2);
\draw[ultra thick, white](10,2) to [out=45, in=135] (11,2);
\draw[ultra thick, white] (1,1) to [out=315,in=225] (2,1);
\draw[ultra thick, white] (3,1) to [out=315,in=225] (4,1);
\draw[ultra thick, white] (5,1) to [out=315,in=225] (6,1);
\draw[ultra thick, white] (7,1) to [out=315,in=225] (8,1);
\draw[ultra thick, white] (9,1) to [out=315,in=225] (10,1);
\draw[ultra thick, white] (11,1) to [out=315,in=225] (12,1);
\draw[ultra thick, white] (1,-1) to [out=135,in=170]  (0.5,-0.66);
\draw[ultra thick, white] (12,-1) to [out=45, in=10] (12.5,-0.66);
\draw[ultra thick, white](4,-1) to [out=45,in=135] (5,-1);
\draw[ultra thick, white](2,-1) to [out=45,in=135] (3,-1);
\draw[ultra thick, white](8,-1) to [out=45, in=135] (9,-1);
\draw[ultra thick, white](6,-1) to [out=45,in=135] (7,-1);
\draw[ultra thick, white](10,-1) to [out=45, in=135] (11,-1);
\draw[ultra thick, white](0.5,0)--(12.5,0);
\end{tikzpicture}
\caption{Composition producing an $\omega$-circle}\label{fig:omega-circle}
\end{figure}
\subsubsection{Formation of increments}\label{sec: increments} {Let us again consider an example. The result of the composition of two t$A$-diagrams of $\mathfrak{tAnn}_{12}$ is a t$A$-diagram of rank $0$ having one transversal region whose increment is $(((0)+(0))+(0))$ (the left hand side of Figure \ref{fig:increments}).}
\begin{figure}[ht]
\begin{tikzpicture}[scale=0.45]
\draw[dotted](0.625,4)--(12.375,4);
\draw[dotted](0.625,7)--(12.375,7);
\draw[dotted](0.625,1)--(12.375,1);
\foreach \x in {1,...,12} \foreach \y in {1,4,7} \filldraw (\x,\y) circle (2.5pt);
\draw (1,7) to [out=315,in=225] (2,7);
\draw(3,7) to [out=315,in=225] (4,7);
\draw(5,7) to [out=315, in =225] (6,7);
\draw(7,7) to [out=315,in=225] (8,7);
\draw (9,7) to [out=315,in=225] (10,7);
\draw(11,7) to [out=315,in=225] (12,7);
\draw (1,4) to[out=45,in=135] (12,4);
\draw (2,4) to [out=45,in=135] (3,4);
\draw (4,4) to [out=45, in=135] (9,4);
\draw (10,4) to [out=45,in=135] (11,4);
\draw (5,4) to [out=45,in=135] (6,4);
\draw (7,4) to [out=45,in=135] (8,4);
\draw(1,4) to [out=315,in=225] (2,4);
\draw(3,4) to [out=315,in=225] (12,4);
\draw(4,4) to[out=315,in=225] (9,4);
\draw(5,4) to [out=315,in=225] (6,4);
\draw(7,4) to [out=315,in=225] (8,4);
\draw (10,4) to [out=315,in=225] (11,4);
\draw (1,1) to [out=45,in=135] (2,1);
\draw(3,1) to [out=45,in=135] (4,1);
\draw(5,1) to [out=45, in =135] (6,1);
\draw(7,1) to [out=45,in=135] (8,1);
\draw (9,1) to [out=45,in=135] (10,1);
\draw(11,1) to [out=45,in=135] (12,1);
\end{tikzpicture}
\begin{tikzpicture}[scale=0.45]
\draw (-2,1) node{};
\filldraw[gray] (0.5,1)--(12.5,1)--(12.5,7)--(0.5,7);
\draw[dotted](0.625,4)--(12.375,4);
\draw[dotted](0.625,7)--(12.375,7);
\draw[dotted](0.625,1)--(12.375,1);
\foreach \x in {1,...,12} \foreach \y in {1,4,7} \filldraw[white] (\x,\y) circle (2pt);
\foreach \x in {1,...,12} \foreach \y in {1,4,7} \filldraw (\x+0.5,\y) circle (2.5pt);
\draw[ultra thick, white]  (1,7) to [out=315,in=225] (2,7);
\draw[ultra thick, white](3,7) to [out=315,in=225] (4,7);
\draw[ultra thick, white](5,7) to [out=315, in =225] (6,7);
\draw[ultra thick, white](7,7) to [out=315,in=225] (8,7);
\draw[ultra thick, white] (9,7) to [out=315,in=225] (10,7);
\draw[ultra thick, white](11,7) to [out=315,in=225] (12,7);
\draw[ultra thick, white] (1,4) to[out=45,in=135] (12,4);
\draw[ultra thick, white] (2,4) to [out=45,in=135] (3,4);
\draw[ultra thick, white] (4,4) to [out=45, in=135] (9,4);
\draw [ultra thick, white](10,4) to [out=45,in=135] (11,4);
\draw[ultra thick, white] (5,4) to [out=45,in=135] (6,4);
\draw[ultra thick, white] (7,4) to [out=45,in=135] (8,4);
\draw[ultra thick, white](1,4) to [out=315,in=225] (2,4);
\draw[ultra thick, white](3,4) to [out=315,in=225] (12,4);
\draw[ultra thick, white](4,4) to[out=315,in=225] (9,4);
\draw[ultra thick, white](5,4) to [out=315,in=225] (6,4);
\draw[ultra thick, white](7,4) to [out=315,in=225] (8,4);
\draw[ultra thick, white] (10,4) to [out=315,in=225] (11,4);
\draw[ultra thick, white] (1,1) to [out=45,in=135] (2,1);
\draw[ultra thick, white](3,1) to [out=45,in=135] (4,1);
\draw[ultra thick, white](5,1) to [out=45, in =135] (6,1);
\draw[ultra thick, white](7,1) to [out=45,in=135] (8,1);
\draw[ultra thick, white] (9,1) to [out=45,in=135] (10,1);
\draw[ultra thick, white](11,1) to [out=45,in=135] (12,1);
\end{tikzpicture}
\caption{Formation of increment}\label{fig:increments}
\end{figure}
{How is the increment encoded in the formation of the product of the corresponding partitions? For this we have to look at the dead blocks formed in the composition of these partitions. Let us denote the upper partition $\alpha$, the lower one $\beta$; the vertices involved in dead blocks are $1,3,4,5,6,7,8,9,10,11$, the corresponding $\alpha$-blocks are
\[\{1,3,9,11\},\{4,6,8\},\{5\},\{7\},\{10\},\]
the $\beta$-blocks are
\[\{1\},\{3,9,11\},\{4,6,8\},\{5\},\{7\},\{10\}\] (to be seen on the right hand side of Figure \ref{fig:increments}). In order to derive from this information about the above mentioned increment we proceed as follows: draw the vertices in the plane on a line ordered as \[1,3,4,5,6,7,8,9,10,11,\]
connect $\alpha$-related vertices by edges in the upper half plane, $\beta$-related vertices by 
\begin{figure}[ht]
\begin{tikzpicture}[scale=0.5]
\draw[dotted](0.625,4)--(12.375,4);
\foreach \x in {3,...,11} \foreach \y in {4} \filldraw (\x,\y) circle (2.5pt);
\filldraw (1,4) circle (2.5pt);
\draw[above] (12,5) node{$\alpha$};
\draw[thick] (1,4) to [out=45,in=135] (3,4);
\draw[thick] (3,4) to [out=90,in=90] (9,4);
\draw[thick] (9,4) to [out=60, in=120] (11,4);
\draw[thick] (4,4) to [out=60,in=120] (6,4);
\draw[thick] (6,4) to [out=60, in=120] (8,4);
\draw[below] (12,3) node {$\beta$};
\draw[thick] (3,4) to [out=270,in=270] (9,4);
\draw[thick] (9,4) to [out=300, in=240] (11,4);
\draw[thick] (4,4) to [out=300,in=240] (6,4);
\draw[thick] (6,4) to [out=300, in=240] (8,4);
\draw[thick, dashed, red] (5,4) circle (9pt);
\draw[thick, dashed, red] (7,4) circle (9pt);
\draw[thick,dashed, red] (10,4) circle (9pt);
\draw[thick, dashed, red] plot [smooth cycle] coordinates {(4,3)(6,2.5) (8,3) (8,5)(6,5.5) (4,5)};
\draw[thick,dashed, red] plot [smooth cycle] coordinates {(1.1,2.5)(6,1.7) (11,2.5) (11,5.5)(6,6.3) (1.1,5.5)};
\draw (15,4) node {};
\end{tikzpicture}
\begin{tikzpicture}[scale=0.75]
\node {$\{1,3,9,11\}$}
	child {node {$\{4,6,8\}$}
		child {node {$\{5\}$}}
		child {node {$\{7\}$}}}
	child {node {$\{10\}$}};	
\end{tikzpicture}
\caption{Increment encoded in dead blocks of $\alpha\beta$}\label{fig:increments2}
\end{figure} edges in the lower half plane as indicated in the left hand side of Figure \ref{fig:increments2} (every $\alpha$- and $\beta$-block should be represented by a circle-free graph) and consider the resulting graph as a subset of the plane. The circles of the increment then are in bijective correspondence with the connected components of this graph in the plane and have to be drawn closely around each component in the surrounding region (in the picture drawn dashed). The structure of this circle collection can be also described in terms of a partial order the dead $\alpha\beta$-blocks can be imposed on, the Hasse diagram of which, if considered as a graph, turns out to have the structure of a forest. This partial order can be defined as follows: a block $B$ is above another block $C$ if there are vertices $b_1,b_2\in B$ such that for every vertex $c\in C$: $b_1<c<b_2$ (here $<$ denotes the linear order on the appropriate set $\mathbb{Z}\times [n]$, cf.~the ${\sqsubseteq}_{in}$ and ${\sqsubseteq}_{out}$ in Section~\ref{sec:rank0inverses}). The circle structure then can be read directly from the Hasse diagram of this poset: every node determines a circle and every circle contains the circles of all nodes below. For the present example the Hasse diagram can be seen in the right hand side of Figure \ref{fig:increments2}.}

In general, increments will be formed from merged regions carrying already labels (from the appropriate monoid or group). The process of forming the increment is exactly the same except that in the result the original labels have to be taken into account. In the left hand side of Figure~\ref{fig:increments with labels},  we took the same initial data as in the previous example, except that several regions carry labels (in the example denoted $a,b,c,d,f,g$). An application of the prescribed procedure is depicted in the right hand side of Figure~\ref{fig:increments with labels}. It is obtained  from the picture of the  left hand side of Figure~\ref{fig:increments2}  by deleting everything except the dashed circles and inscribing the labels appropriately.
\begin{figure}[ht]
\begin{tikzpicture}[scale=0.45]
	\draw[dotted](0.625,4)--(12.375,4);
	\draw[dotted](0.625,7)--(12.375,7);
	\draw[dotted](0.625,1)--(12.375,1);
	\foreach \x in {1,...,12} \foreach \y in {1,4,7} \filldraw (\x,\y) circle (2.5pt);
	\draw (1,7) to [out=315,in=225] (2,7);
	\draw(3,7) to [out=315,in=225] (4,7);
	\draw(5,7) to [out=315, in =225] (6,7);
	\draw(7,7) to [out=315,in=225] (8,7);
	\draw (9,7) to [out=315,in=225] (10,7);
	\draw(11,7) to [out=315,in=225] (12,7);
	\draw (1,4) to[out=45,in=135] (12,4);
	\draw (2,4) to [out=45,in=135] (3,4);
	\draw (4,4) to [out=80, in=100] (9,4);
	\draw (10,4) to [out=45,in=135] (11,4);
\draw(6.5,5)node{$b$};
\draw(6.5,3.2)node{$d$};
\draw(8.5,2.5)node{$f$};
\draw(3.5,4.5)node{$g$};
\draw (5,4) .. controls (5,5) and (6,5) .. (6,4);
\draw (5.5,4.4)node{$a$};
	\draw (7,4) to [out=45,in=135] (8,4);
	\draw(1,4) to [out=315,in=225] (2,4);
	\draw(3,4) to [out=290,in=250] (12,4);
	\draw(4,4) to[out=290,in=250] (9,4);
	\draw(5,4) to [out=315,in=225] (6,4);
\draw (7,4)..controls (7,3) and (8,3).. (8,4);
\draw (7.5,3.6) node{$c$};
	\draw (10,4) to [out=315,in=225] (11,4);
	\draw (1,1) to [out=45,in=135] (2,1);
	\draw(3,1) to [out=45,in=135] (4,1);
	\draw(5,1) to [out=45, in =135] (6,1);
	\draw(7,1) to [out=45,in=135] (8,1);
	\draw (9,1) to [out=45,in=135] (10,1);
	\draw(11,1) to [out=45,in=135] (12,1);
\end{tikzpicture}  \quad\quad
\begin{tikzpicture}[scale=0.5]
	\draw[thick] (5,4) circle (14pt);
	\draw(5,4)node{$a$};
	\draw[thick] (7,4) circle (14pt);
	\draw(7,4)node{$c$};
	\draw[thick] (10,4) circle (14pt);
	\draw[thick] plot [smooth cycle] coordinates {(4,3)(6,2.5) (8,3) (8,5.5)(6,6) (4,5.5)};
	\draw(6,5.3)node{$b+d$};
	\draw[thick] plot [smooth cycle] coordinates {(1.1,1.7)(6,0.9) (11,1.7) (11,6)(6,6.8) (1.1,6)};
	\draw (6,1.7) node {$f+g$};
\end{tikzpicture}
\caption{Increment from regions carrying labels}\label{fig:increments with labels}
\end{figure}

Finally, the most general case is as follows: every region $R$ in a product $\bm\alpha\bm\beta$ is the result of merging a set of regions $A_i$ from $\bm\alpha$ and $B_j$ from $\bm\beta$. The label of $R$ then is the sum of the labels of the regions $A_i$ and $B_j$ plus the sum of the increments. Every such summand then is of the form of the one on the right hand side of Figure~\ref{fig:increments with labels} and is produced by the procedure described before.

\subsubsection{Forming inverses}\label{sec: tann inverses} 
From now on we assume that the labels of regions are from $\mathbb{FA}$, the free unary abelian group generated by the empty set, with the exception of the unique transversal region in case of t$A$-diagrams of rank $0$: the resulting  regions  carry labels from the free product $\mathbb{Z}*\mathbb{FA}$.
\paragraph{\ref{sec: tann inverses}.1. \emph{Inverses of t$A$-diagrams of positive rank}}\label{sec:inverses for positive rank} 
In this subsection all t$A$-diagrams have positive rank, hence all labels of regions are  from $\mathbb{FA}$.
\begin{figure}[ht]
\begin{tikzpicture}[scale=0.58]
	
	\foreach \x in {0,...,21} \foreach \y in {0,4,8,12} \filldraw (\x,\y) circle (2.5pt);
	\draw[dotted](0,0)--(21,0);
	\draw[dotted](0,4)--(21,4);
	\draw[dotted](0,8)--(21,8);
	\draw[dotted](0,12)--(21,12);
	\draw (0,12.5) node {$1$};
	\draw (9,12.5) node {$10$};
	\draw(20,12.5) node {$21$};
	\draw (21,12.5) node {$1$};
    \draw (15,-.5) node {$16$};
    
	\draw[thick] (9,12)--(15,8);
	\draw[thick](15,8)--(9,4);
	\draw (10,12)--(16,8);
	\draw(16,8)--(10,4);
	\draw[thick] (0,12)--(0,8);
	\draw (21,12)--(21,8);
	
	\draw (1,12)  to [out=315,in=225] (2,12);
	
	\draw (3,12)  to [out=315,in=225] (4,12);
	
	\draw (5,12)  to [out=315,in=225] (8,12);
	\draw (6,12)  to [out=315,in=225] (7,12);
	\draw (11,12)  to [out=315,in=225] (20,12);
	\draw (12,12)  to [out=315,in=225] (19,12);
	\draw (13,12)  to [out=315,in=225] (14,12);
	\draw (15,12)  to [out=315,in=225] (18,12);
	\draw (16,12)  to [out=315,in=225] (17,12);
	
	\draw(4.5,10) node{$f$};
	\draw((18.5,8.3) node {$g$};

     \draw (5.8,6) node{{$-f-(a+(0))-((((0))))-3(0) -((0))$}};
    \filldraw[gray, nearly transparent](0.05,5.7)--(11.55,5.7)--(11.55,6.45)--(0.05,6.45);
    \draw (17.5,6.5) node{{$-((g))-(((0)+((0))))$}};
	\filldraw[gray, nearly transparent](14.2,6.1)--(20.8,6.1)--(20.8,6.85)--(14.2,6.85);
	\draw (1,8)  to [out=90,in=90] (4,8);
	\draw(1,8) to [out=310, in =230] (4,8);
	\draw(2.5,8.5) node{$a$};
	\draw (2,8)  to [out=45,in=135] (3,8);
	\draw(2,8) to [out=310, in=230] (3,8);
	\draw (5,8)  to [out=45,in=135] (12,8);
	\draw (5,8)  to [out=310,in=230] (12,8);
	\draw (6,8)  to [out=45,in=135] (11,8);
	\draw (6,8)  to [out=310,in=230] (11,8);
	\draw (7,8)  to [out=45,in=135] (10,8);
	\draw (7,8)  to [out=310,in=230] (10,8);
	\draw  (8,8)  to [out=45,in=135] (9,8);
	\draw  (8,8)  to [out=310,in=230] (9,8);
	\draw (13,8)  to [out=45,in=135] (14,8);
	\draw (13,8)  to [out=310,in=230] (14,8);
	\draw (17,8) .. controls (17,10) and (20,10).. (20,8);
	\draw (17,8) to [out=310,in=230] (20,8); 
	\draw (18,8)  .. controls (18,9) and (19,9) ..  (19,8);
	\draw (18,8) to [out=310,in=230] (19,8);

	\draw [thick](0,8)--(0,4);
	\draw (21,8)--(21,4);

	\draw [thick](9,4)--(15,0);
	\draw (10,4)--(16,0);
	\draw [thick](0,4)--(0,0);
	\draw (21,4)--(21,0);
	
	\draw (1,4)  to [out=310, in=230] (2,4);
	\draw(1,4) to [out=50,in=130] (2,4);
	\draw (3,4) to [out=310, in=230](4,4);
	\draw(3,4) to [out=50, in=130] (4,4);

	\draw(5,4)  to [out=310,in=230] (8,4);
	\draw(5,4)  to [out=50,in=130] (8,4);
	\draw (6,4)  to [out=310,in=230] (7,4);
	\draw (6,4)  to [out=50,in=130] (7,4);
	\draw (11,4)  to [out=310,in=230] (20,4);
	\draw (11,4)  to [out=50,in=130] (20,4);
	\draw (12,4)  to [out=310,in=230] (19,4);
	\draw (12,4)  to [out=50,in=130] (19,4);
	\draw (13,4)  to [out=310,in=230] (14,4);
	\draw (13,4)  to [out=50,in=130] (14,4);
	\draw (15,4)  to [out=310,in=230] (18,4);
	\draw (15,4)  to [out=50,in=130] (18,4);
	\draw (16,4)  to [out=310,in=230] (17,4); 
	\draw (16,4)  to [out=50,in=130] (17,4);
	
	\draw (1,0)  to [out=90,in=90] (4,0);
	\draw(2.5,0.5) node{$a$};
	\draw (2,0)  to [out=50,in=130] (3,0);
	
	\draw(4.5,2) node{$f$};
	\draw((18.5,0.3) node {$g$};
	
	\draw (17,0) .. controls (17,2) and (20,2).. (20,0);
	\draw (18,0)  .. controls (18,1) and (19,1) ..  (19,0);

	\draw (5,0)  to [out=45,in=135] (12,0);
	\draw (6,0)  to [out=45,in=135] (11,0);
	\draw (7,0)  to [out=45,in=135] (10,0);
	\draw (8,0)  to [out=45,in=135] (9,0);
	\draw (13,0)  to [out=45,in=135] (14,0);
\end{tikzpicture}
\caption{The t$A$-diagrams $\bm\varepsilon:=\bm\alpha \underline{\bm\alpha}^*\bm\alpha$ and ${\bm\beta}$}\label{fig:inverse for positive rank}
\end{figure}
Let us look at the example depicted in Figure~\ref{fig:inverse for positive rank}, but first ignore the shaded entries. For the t$A$-diagram $\bm\alpha$ of Figure~\ref{a+b+c} we calculated the t$A$-diagram $\bm\varepsilon:=\bm\alpha\underline{\bm\al}^*\bm\alpha$. All of the involved t$A$-diagrams have three transversal regions (the `left', the `middle' and the `right' one, respectively). The label in $\bm\alpha$ of the left region is $f$ while the labels in $\bm\alpha$ of the other two are empty. The label in $\bm\varepsilon$ of the left region consists of 
\begin{align*}
    &2f &&\mbox{(coming from the two factors $\bm\alpha$)}\\
     &+(a+(0))+((((0))))+(0) &&\mbox{(increment resulting from $\bm\alpha\underline{\bm\alpha}^*$)}\\
   & +2(0)+((0)) &&\mbox{(increment resulting from $\underline{\bm\alpha}^*\bm\al$)}\\
\end{align*}
Similarly, the label  in $\bm\varepsilon$ of the right region is 
\[r:=((g))+(((0)+((0))))\] (consisting entirely of the increment)
while the label of the middle region is empty.

Now we define the t$A$-diagram $\bm\beta$ by defining its skeleton $\underline{\bm\beta}:=\underline{\bm\alpha}^*$ and the labeling function $c$ on its regions: the label of the left transversal region is 
\[-f-(a+(0))-((((0))))-3(0)-((0))\] while the one of the right transversal region should be
\[-((g))-(((0)+((0))))\]
and all other regions being equipped with the empty label. The t$A$-diagram $\bm\beta$ then is an inverse of $\bm\alpha$ and the procedure gives rise to the definition of weak regular involution $\bm\alpha\mapsto\bm\alpha^*$. The procedure is as follows: every transversal region of $\underline{\bm\alpha}^*$ is the mirror image $C^*$ of the corresponding transversal region $C$ of $\bm\alpha$. The region $C$ in ${\bm\alpha}\underline{\bm\alpha}^*{\bm\alpha}$ is obtained by merging $C\cup C^*\cup C$. The label of $C$ in ${\bm\alpha}\underline{\bm\alpha}^*{\bm\alpha}$ is twice the label $f$ of $C$ in $\bm\alpha$ plus the increments of merging $C\cup C^*$ and $C^*\cup C$. Hence if we endow $C^*$ with the label $-f$ minus the increment resulting from merging $C\cup C^*\cup C$ (and putting the empty label on all involved in- and out regions) and do so for every transversal region $C$ of $\bm\alpha$, then we obtain an inverse diagram $\bm\beta$. It is easy to see that application of the same procedure to $\bm\beta$ leads back to $\bm\alpha$. Altogether the prescribed procedure results in a weak regular involution $\bm\alpha\mapsto \bm\alpha^*$ (defined on arrows of positive rank).
\begin{figure}[ht]
\begin{tikzpicture}[scale=0.60]
	
	\foreach \x in {0,...,21} \foreach \y in {0,4,8,12} \filldraw (\x,\y) circle (2.5pt);
	\draw[dotted](0,0)--(21,0);
	\draw[dotted](0,4)--(21,4);
	\draw[dotted](0,8)--(21,8);
	\draw[dotted](0,12)--(21,12);
	\draw (0,12.5) node {$1$};
	\draw (9,12.5) node {$10$};
	\draw(20,12.5) node {$21$};
	\draw (21,12.5) node {$1$};
	\draw (15,-.5) node {$16$};
    
	\draw[thick] (9,12)--(15,8);
	\draw[thick](15,8)--(9,4);
	\draw (10,12)--(16,8);
	\draw(16,8)--(10,4);
	\draw[thick] (0,12)--(0,8);
	\draw (21,12)--(21,8);
	
	\draw (1,12)  to [out=315,in=225] (2,12);
	
	\draw (3,12)  to [out=315,in=225] (4,12);
	
	\draw (5,12)  to [out=315,in=225] (8,12);
	\draw (6,12)  to [out=315,in=225] (7,12);
	\draw (11,12)  to [out=315,in=225] (20,12);
	\draw (12,12)  to [out=315,in=225] (19,12);
	\draw (13,12)  to [out=315,in=225] (14,12);
	\draw (15,12)  to [out=315,in=225] (18,12);
	\draw (16,12)  to [out=315,in=225] (17,12);
	
	\draw(4.5,10) node{$f$};
	\draw(18.5,8.3) node {$g$};
    \draw(18.5,7.7) node{$-g$};
    \draw(18.5,6.9)node{$-\circ$};
    \draw(18.5,6.3)node{$-\circ$};

   \draw (4.5,6) node{{$-f$}};
	
	\draw (1,8)  to [out=90,in=90] (4,8);
	\draw(1,8) ..controls (1,6.5) and (4,6.5).. (4,8);
    \draw(2.5,6.6)node{$-\circ$};
	\draw(2.5,8.5) node{$a$};
       \draw(2.5,7.5) node{$-a-\circ$};
	\draw (2,8)  to [out=45,in=135] (3,8);
	\draw(2,8) to [out=310, in=230] (3,8);
	\draw (5,8)  to [out=45,in=135] (12,8);
	\draw (5,8)  to [out=310,in=230] (12,8);
     \draw(8.5,6.1)node{$-\circ$};
	\draw (6,8)  to [out=45,in=135] (11,8);
    \draw(8.5,6.6)node{$-\circ$};
	\draw (6,8)  to [out=310,in=230] (11,8);
	\draw (7,8)  to [out=45,in=135] (10,8);
	\draw (7,8)  to [out=310,in=230] (10,8);
   \draw(8.5,7.1)node{$-\circ$}; 
    \draw  (8,8)  to [out=45,in=135] (9,8);
  	\draw  (8,8)  to [out=310,in=230] (9,8);
      \draw(8.5,7.6)node{$-\circ$};
	\draw (13,8)  to [out=45,in=135] (14,8);
      \draw(13.5,7.6)node{$-\circ$};
	\draw (13,8)  to [out=310,in=230] (14,8);
      
	\draw (17,8) .. controls (17,10) and (20,10).. (20,8);
    
	\draw (17,8) ..controls (17,6) and (20,6).. (20,8); 
	\draw (18,8)  .. controls (18,9) and (19,9) ..  (19,8);
	\draw (18,8) .. controls (18,7) and (19,7) .. (19,8);

	\draw [thick](0,8)--(0,4);
	\draw (21,8)--(21,4);

	\draw [thick](9,4)--(15,0);
	\draw (10,4)--(16,0);
	\draw [thick](0,4)--(0,0);
	\draw (21,4)--(21,0);
	
	\draw(1.5,4.5)node{$-\circ$};
    \draw (1,4)  to [out=310, in=230] (2,4);
	\draw(1,4) to [out=50,in=130] (2,4);
	\draw(3.5,4.5)node{$-\circ$};
    \draw (3,4) to [out=310, in=230](4,4);
	\draw(3,4) to [out=50, in=130] (4,4);

	\draw(5,4)  to [out=310,in=230] (8,4);
    \draw (6.5,5.2)node{$-\circ$};
	\draw(5,4)  to [out=90,in=90] (8,4);
    \draw(6.5,4.5)node{$-\circ$};
	\draw (6,4)  to [out=310,in=230] (7,4);
	\draw (6,4)  to [out=50,in=130] (7,4);
	\draw (11,4)  to [out=310,in=230] (20,4);
     \draw(15.5,6.2)node{$-\circ$};
	\draw (11,4)  to [out=50,in=130] (20,4);
	\draw (12,4)  to [out=310,in=230] (19,4);
    \draw(15.5,5.8)node{$-\circ$};
	\draw (12,4)  to [out=50,in=130] (19,4);
    \draw(13.5,4.5)node{$-\circ$};
	\draw (13,4)  to [out=310,in=230] (14,4);
	\draw (13,4)  to [out=50,in=130] (14,4);
	\draw (15,4)  to [out=310,in=230] (18,4);
    \draw(16.5,4.9)node{$-\circ$};
	\draw (15,4)  to [out=50,in=130] (18,4);
    \draw(16.5,4.4)node{$-\circ$};
	\draw (16,4)  to [out=310,in=230] (17,4); 
	\draw (16,4)  to [out=50,in=130] (17,4);
	
	\draw (1,0)  to [out=90,in=90] (4,0);
	\draw(2.5,0.5) node{$a$};
	\draw (2,0)  to [out=50,in=130] (3,0);
	
	\draw(4.5,2) node{$f$};
	\draw((18.5,0.3) node {$g$};
	
	\draw (17,0) .. controls (17,2) and (20,2).. (20,0);
	\draw (18,0)  .. controls (18,1) and (19,1) ..  (19,0);
	
	\draw (5,0)  to [out=45,in=135] (12,0);
	\draw (6,0)  to [out=45,in=135] (11,0);
	\draw (7,0)  to [out=45,in=135] (10,0);
	\draw (8,0)  to [out=45,in=135] (9,0);
	\draw (13,0)  to [out=45,in=135] (14,0);
\end{tikzpicture}
\caption{The t$A$-diagram $\bm\alpha \underline{\bm\alpha}^*\bm\alpha$}\label{fig:inverse,subtle}
\end{figure}

However, the prescribed procedure amounts to be a bit of brute force. A more subtle procedure is indicated in Figure~\ref{fig:inverse,subtle} (in the picture, $\circ$ stands for $(0)$). It also results in a weak regular involution. We shall describe its principle. We first assume that the diagram $\bm\alpha$ is given in its $R$-strip representation (and the regions considered as partitions $\mathbb{Z}\times[m]\leadsto \mathbb{Z}\times [n]$). We define a partial order $\sqsubseteq$ on the set of all regions as follows. For two regions $C,D$, we set $D\sqsubseteq C$ if and only if there are $u,v\in C$ such that for all $x\in D$, $u\le x\le v$ where ${\le}={\le}_{in}\cup{\le}_{out}$ and ${\le}_{in}$ denotes the lexicographic order on the incoming vertices $\mathbb{Z}\times[m]$ while ${\le}_{out}$ the lexicographic order on the outgoing vertices $\mathbb{Z}\times[n]$. The definition (is meant as and) implies that \textsl{both} elements $u$ and $v$ are either incoming or outgoing vertices. It follows that
\begin{itemize}
    \item the transversal regions are the maximal elements in this order,
    \item every in- and out-region has a unique cover,
    \item no in-region is comparable with any out-region (and conversely).
\end{itemize}
Let now $\bm\alpha$ be a t$A$-diagram of non-zero rank. We define the t$A$-diagram $\bm\alpha^*$ as follows: the skeleton of $\bm\alpha^*$ is $\underline{\bm\alpha}^*$, the mirror image of the skeleton $\underline{\bm\alpha}$ of $\bm\alpha$; each region of $\underline{\bm\alpha}^*$ is the mirror image $C^*$ of the corresponding region $C$ in $\bm\alpha$. We define the label $c_{{\bm\alpha}^*}(C^*)$ by
\[c_{{\bm\alpha}^*}(C^*):=-c_{{\bm\alpha}}(C)-\sharp C(0)\]
where $\sharp C$ is the number of regions covered by $C$ in $\bm\alpha$ (equivalently, the number of regions covered by $C^*$ in $\underline{\bm\alpha}^*$ (with respect to the aforementioned partial order on the set of regions of $\bm\alpha$). Altogether, the prescribed procedure provides another weak regular involution. The discussion of inverses of t$A$-diagrams of positive rank, so far provides  weak regular involutions on the categories $\overline{\mathfrak{tAnn}}^{\mathrm{odd}}$ and $\overline{\mathfrak L}$.

\paragraph{\ref{sec: tann inverses}.2. \emph{Inverses of t$A$-diagrams of rank} $0$}\label{sec:rank0inverses} In order to get regular involutions on all of $\overline{\mathfrak{tAnn}}$ we need to discuss the rank $0$ case. Let $\bm\alpha$ be a rank $0$ t$A$-diagram with unique transversal region $T$ and denote the label $c(T)$ of $T$ by $\ell$ (note that $\ell\in \mathbb{Z}*\mathbb{FA}$ and recall the the binary operation in this group is denoted $+$, the inverse of an element $x$ is denoted $-x$ despite the fact that this group is not abelian). We define partial orders: ${\sqsubseteq}_{in}$ on the set of in-regions and ${\sqsubseteq}_{out}$ on the set of out-regions in the same way as the partial order $\sqsubseteq$ for regions in the previous section and define the transversal region $T$ to be ${\sqsubseteq}_{in}$-above all in-regions and ${\sqsubseteq}_{out}$-above all out regions. We are ready to defined an inverse $\bm\alpha^*$ by providing the skeleton $\underline{\bm\alpha^*}$ and the labeling $c_{\bm\alpha^*}$ of its regions. The skeleton is, of course,  the mirror image of the skeleton of $\bm\alpha$: $\underline{\bm\alpha^*}:=\underline{\bm\alpha}^*$. Every region of $\underline{\bm\alpha}^*$ is the mirror image $R^*$ of its corresponding region $R$ of $\bm\alpha$. If $R$ is an in- or an out-region of $\bm\alpha$ define the label $c_{\bm\alpha^*}$ of $R^*$ by
\[c_{\bm\alpha^*}(R^*):=-c_{\bm\alpha}^*(R)-\sharp R(0)\]
where $\sharp R$ is the number of in- or out-regions covered by $R$ (with respect to ${\sqsubseteq}_{in}$ or ${\sqsubseteq}_{out}$, depending on whether $R$ is an in- or out-region). Finally, the label of the transversal region $T$ is as follows: 
\[c_{\bm\alpha^*}(T^*):=-\sharp_{out}T(0)-\ell-\sharp_{in}T(0)\]
where $\sharp_{out}T$ is the number of out-regions covered by $T$ in $\bm\alpha$ and $\sharp_{in}T$ the number of in-regions covered by $T$ in $\bm\alpha$ (with respect to  ${\sqsubseteq}_{out}$ respectively ${\sqsubseteq}_{in}$). 

Recall that composition of annular diagrams is `from out to in' in the sense that in the composition of two diagrams $\bm\alpha\bm\beta$, the out-vertices of $\bm\alpha$ are merged with the in-vertices of $\bm \beta$; for elements $a,b,c$ in the non-abelian group $\mathbb{Z}*\mathbb{FA}$, $-a-b-c$ stands for $(-a)+(-b)+(-c)$. The definition of $\bm\alpha^*$ obtained this way corresponds to the second version of inverse provided in the previous section for t$A$-diagrams of positive rank. One could as well define  an inverse corresponding to the first version of inverses for t$A$-diagrams of positive rank. The procedure essentially is as follows: denote the increment in $\bm\alpha\underline{\bm\alpha}^*$ formed by merging the outgoing regions of $\bm\alpha$ with the ingoing regions of $\underline{\bm\alpha}^*$ by $\ell_{in}$ and the one in $\underline{\bm\alpha}^*\bm\al$ formed by merging the outgoing regions of $\underline{\bm\alpha}^*$ with the ingoing regions of $\bm\alpha$ by $\ell_{out}$. Setting
\[c_{\bm\alpha^*}(T^*):=-\ell_{in}-\ell-\ell_{out}\] for the transversal region $T^*$ and having labeled every other region by the empty symbol $0$ we are led to an inverse $\bm\alpha^*$ that corresponds to the first version of inverses discussed in Section~\ref{sec:inverses for positive rank}.
In any case, combination with the previous section provides us with two examples of weak regular involutions on all of $\overline{\mathfrak{tAnn}}$.

\subsection{Another congruence ${\equiv}_0$ on $\mathfrak{tAnn}$}\label{sec:factored tAnn} Endowed with addition $+$ and the successor function $x\mapsto x^+$ of Peano arithmetic, the set $\mathbb{N}_0$ of all non-negative integers is a commutative unary monoid generated by the empty set $\varnothing$. Denote the structure $\langle\mathbb{N}_0;+,{}^+,0\rangle$ by $\mathbb{P}$. There is a unique surjective morphism $\mathbb{M}\twoheadrightarrow \mathbb{P}$, which essentially maps a collection of circles to the number of its circles. This gives raise to a another congruence ${\equiv}_0$ on $\mathfrak{tAnn}$, contained in ${\equiv}_1$, that is, for $\mathfrak{t^\flat Ann}:=\mathfrak{tAnn}/{\equiv}_0$ we have 
\[\mathfrak{tAnn}\twoheadrightarrow \mathfrak{t^\flat Ann}\twoheadrightarrow \mathfrak{aTL}^{\mathrm{d}}.\]
Intuitively, the category $\mathfrak{t^\flat Ann}$ retains the regions induced by the strings and the $\omega$-circles, but ignores the \textsl{structure} of the circles sitting on the regions and retains only the number of circles in each region. In other words, the data for the arrows of $\mathfrak{t^\flat Ann}$ are almost the same as those for $\mathfrak{tAnn}$, except that the labeling function $\mathbf{c}$ is composed with the morphism $\mathbb{M}\twoheadrightarrow \mathbb{P}$. (In particular, the unique transversal region in a rank $0$ arrow may carry a label from $\mathbb{O}\mathop{*}\mathbb{P}$.) Altogether, every arrow $\bm\alpha$ of $\mathfrak{t^\flat Ann}$ is described by its skeleton $\underline{\bm\alpha}$ (an arrow of $\mathfrak{aTL}^\mathrm{e}$) whose regions are labeled by non-negative integers (the unique transversal region of a rank $0$ arrow is labeled by an element of $\mathbb{O}\mathop{*}\mathbb{P}$). It should be clear how a regular extension $\overline{\mathfrak{t^\flat Ann}}$ can be constructed: the labels have to be taken from $\mathbb{Z}$ (or from $\mathbb{Z}*\mathbb{Z}$ for the transversal region in the rank $0$ case). We then have quotient morphisms
\[\overline{\mathfrak{tAnn}}\twoheadrightarrow \overline{\mathfrak{t^\flat Ann}}\twoheadrightarrow \overline{\mathfrak{aTL}^{\mathrm{d}}}\] and for every $n$ the local monoids are denoted $\mathfrak{t^\flat Ann}_n$ and 
$\overline{\mathfrak{t^\flat Ann}_n}$.

We may restrict the congruence ${\equiv}_0$ to the category $\mathfrak{L}$ to get the quotient category $\mathfrak{L}^\flat$. Again, $\mathfrak{L}^\flat$ retains the regions but ignores the structure of the circles sitting on the regions but retains only the number of circles sitting on each region. Hence there are quotient morphisms
\[\mathfrak{L}\twoheadrightarrow \mathfrak{L}^\flat\twoheadrightarrow\mathfrak{TL}\]
where $\mathfrak{TL}$ denotes the classical Temperley--Lieb category (the regions are ignored and only the total number of circles is retained). If, in the context of $\mathfrak{L}^\flat$, we allow labels from all of 
$\mathbb{Z}$ (instead of just non-negative integers) we get regular extensions and quotient morphisms
\[\overline{\mathfrak{L}}\twoheadrightarrow \overline{\mathfrak{L}^\flat}\twoheadrightarrow \overline{\mathfrak{TL}}.\]
These quotient morphisms restrict to the corresponding local monoids. We note that the local monoids $\mathfrak{TL}_n$ of $\mathfrak{TL}$ are the \emph{Kauffman monoids} and will be denoted $\mathfrak{K}_n$ as in \cite{auingeretal, bodope, kauffmanIsotopy}. Correspondingly, the local monoids $\overline{\mathfrak{TL}_n}$ of $\overline{\mathfrak{TL}}$ are denoted $\overline{\mathfrak{K}_n}$. 

\section{Prerequisites on the Finite Basis Problem}\label{preliminaries} 
As mentioned in the introduction, our combinatorial characterizations of categories with a topological flavor allow for a comprehensive analysis of the finite axiomatization problem (aka the Finite Basis Problem) of the local monoids in these categories. Here we collect all necessary prerequisites.  First we recall the required essentials of equational logic, afterwards a result from \cite{auingeretal} which will be an essential ingredient of the proofs in Section \ref{finite basis problem}, and some constructions from semigroup theory.

\subsection{Equational logic and varieties of universal algebras} The concepts of identity and  identity basis are intuitively clear. Nevertheless, any precise reasoning about these concepts requires a formal framework, especially when one aims at `negative' results as we do in this paper. Such a framework, provided by equational logic, is concisely presented, e.g., in \cite[Chapter~II]{BuSa81}. For the reader's convenience, we briefly overview the basic vocabulary of equational logic in a form adapted to the use in this paper.

A non-empty set $A$ endowed with operations
$f_1\colon\underbrace{A\times\dots\times A}_{n_1}\to
A$,\linebreak $f_2\colon\underbrace{A\times\dots\times
A}_{n_2}\to A$, \dots\ is called an \emph{algebraic structure of
type $(n_1,n_2,\dots)$ with carrier $A$} or simply an \emph{algebra} of type $(n_1,n_2,\dots)$. In the present paper we shall be concerned only with algebras of types $(2)$ and $(2,1)$ and the binary operation is always assumed to satisfy the associative law. That is, the algebraic structures in this paper will be either \emph{semigroups} $\mathcal{S}=\left<S,\cdot\right>$ or \emph{unary} semigroups $\mathcal{T}=\left<T,\cdot,^*\right>$ where a unary semigroup is simply a semigroup endowed with an additional unary operation $^*$. Most semigroups considered will actually admit a neutral element, so we mostly refer to them as \emph{monoids}. We shall, however, usually not assume that the corresponding nullary operation be contained in the signature of the respective structure (with the exception of the unary monoids $\mathbb{M}$ and $\mathbb{P}$ introduced in Subsections \ref{combinatorial tAnn, L} and~\ref{sec:factored tAnn}).

Let us fix a countably infinite set $X=\{x_1,x_2,\dots\}$ of \emph{letters} or \emph{variables} and let $X^+$ be the \emph{free semigroup} on $X$, that is, the set of all non-empty (finite) words formed by the letters of $X$, endowed with the binary operation of concatenation of words. Let further $\mathcal{S}=\left<S,\cdot\right>$ be a semigroup; a map $\vp\colon X\to \mathcal{S}$ is a \emph{substitution}, and each such substitution admits a unique extension to a morphism $X^+\to \mathcal{S}$, denoted also by $\vp\colon X^+\to \mathcal{S}$. A \emph{semigroup identity}, or simply an \emph{identity} is a pair $(u,v)\in X^+\times X^+$, written as $u\bumpeq v$. We say that a semigroup $\mathcal{S}$ \emph{satisfies} the identity, or the identity $u\bumpeq v$ \emph{holds in} $\mathcal{S}$ if $u\vp=v\vp$ for every substitution $\vp\colon X\to \mathcal{S}$. The identity $u\bumpeq v$ is \emph{nontrivial} if $u$ and $v$ are distinct elements of $X^+$.

A unary operation $^*$ on a semigroup $\mS$ is an \emph{involution} if, for any $a,b\in \mathcal{S}$,
\begin{equation}\label{involution}(ab)^*=b^*a^*\mbox{ and }(a^*)^*=a.
\end{equation}
A unary semigroup $\mS=\left<S,\cdot,^*\right>$ is an \emph{involutory semigroup} or a \emph{semigroup with involution} if its unary operation $^*$ is an involution. In the context of involutory semigroups it is necessary to replace the free  semigroup $X^+$ by the \emph{free involutory semigroup} $\mathcal{I}(X)$ which can be realized as follows. Let $X^*:=\{x_1^*,x_2^*,\dots\}$ be a disjoint copy of $X$, with $x\mapsto x^*$ being a bijection $X\to X^*$, form $X\cup X^*$, extend the bijection $x\mapsto x^*$ to a bijection $X\cup X^*\to X\cup X^*$ by letting $(x^*)^*:=x$ for every $x\in X$, form the free semigroup $(X\cup X^*)^+$ on $X\cup X^*$, and let $\mathcal{I}(X):=\left<(X\cup X^*),\cdot, ^*\right>$ where $\cdot$ is concatenation of words and $^*$ is the operation $y_1\cdots y_n\mapsto y_n^*\cdots y_1^*$ with $y_i\in X\cup X^*$. Then $\mathcal{I}(X)$ is the \emph{free involutory semigroup} on $X$, every map $\varphi\colon X\to \mS$ where $\mS$ is any involutory semigroup can be uniquely extended to a morphism of involutory semigroups $\mathcal{I}(X)\to \mathcal{S}$. Again, the concept of identity (of involutory semigroups) and of satisfaction of an identity by an involutory semigroup are defined accordingly.

The following concepts and results are entirely parallel for semigroups and involutory semigroups.  Hence we shall formulate them only for the `plain' semigroup case. Given a semigroup $\mathcal{S}$, we denote by $\Id\mathcal{S}$ the set of all semigroup identities $u\bumpeq v$ satisfied by $\mathcal{S}$ (with $u,v\in X^+$).

Given any collection $\Sigma$ of semigroup identities, we say that an identity $u\bumpeq v$ \emph{is a consequence of} $\Sigma$ or that $\Sigma$ \emph{implies} $u\bumpeq v$ if every semigroup satisfying all
identities of $\Sigma$ satisfies the identity $u\bumpeq v$, as well. Birkhoff's completeness theorem of equational logic (see \cite[Theorem 14.17]{BuSa81}) shows that this notion (which we have given a semantic definition) can be captured by a very transparent set of inference rules. These rules in fact formalize the most natural things one does with identities: substitution of a term for a variable, application of operations to identities (such as, say, multiplying both sides of an identity  on the right by the same term) and using symmetry and transitivity of equality. We need not go into more detail here because the completeness theorem is not utilized in this paper.

Given a semigroup $\mathcal{S}$, an \emph{identity basis} (or \emph{equational basis}) for $\mathcal{S}$ is any set $\Sigma\subseteq\Id\mathcal{S}$ such that every identity of $\Id\mathcal{S}$ is a consequence $\Sigma$. A  semigroup $\mathcal{S}$ is said to be \emph{finitely based} if it possesses a finite identity basis; otherwise $\mS$ is called \emph{non-finitely based}.

The class $[\Sigma]$ of all semigroups satisfying all identities from a given set $\Sigma$ of semigroup identities is called the \emph{variety defined by $\Sigma$}. A variety $\mathsf{V}$ is said to be \emph{finitely based} if  $\mathsf{V}=[\Sigma]$ for some finite set $\Sigma$ of identities; otherwise it is called \emph{non-finitely based}.

It is easy to see that the satisfaction of an identity is inherited by forming direct products, taking subsemigroups and morphic images so that each variety is closed under these operators. In fact, varieties can be characterized by this closure property (the HSP-theorem, see \cite[Theorem 11.9]{BuSa81}). Given a semigroup $\mathcal{S}$, then the variety defined by $\Id\mS$ is the \emph{variety generated by $\mS$} (it is the smallest variety of semigroups which contains $\mathcal{S}$); we denote this variety by $\var\mS$. Observe that a semigroup and the variety it generates are simultaneously finitely or non-finitely based. From the HSP-theorem it follows that every member of $\var\mS$ is a morphic image of a subsemigroup of a direct product of several copies of $\mS$. For a class $\mathsf{C}$ of semigroups, $\var\mathsf{C}$ is the smallest variety of semigroups containing $\mathsf{C}$, and, again it is the class obtained from $\mathsf{C}$ by closing under forming morphic images of subsemigroups of direct products of semigroups in $\mathsf{C}$.

Recall that a semigroup $\mathcal{S}$  is \emph{locally finite} if every finitely generated subsemigroup of $\mathcal{S}$ is finite, while a semigroup is \emph{periodic} if the latter holds for every one-generator subsemigroup. Hence every locally finite semigroup is periodic but the converse is not true in general \cite{sapir}. A variety $\mathsf{V}$ is said to be \emph{locally finite} if every  member of $\mathsf{V}$ is locally finite, and likewise, a variety is \emph{periodic} if every member is periodic. Every variety generated by a finite semigroup  is locally finite (see \cite[Theorem II.10.16]{BuSa81}).

A sufficient condition for a non-periodic variety $\mathsf{V}$ of semigroups to contain only locally finite periodic members has been found by Sapir \cite{sapir}. In order to formulate it, we define the sequence of \emph{Zimin words} inductively by $Z_1:=x_1$ and $Z_{n+1}:=Z_nx_{n+1}Z_n$.
\begin{Thm}[{\cite[Lemma 3.3]{sapir}}]\label{sapirlocallyfinite} If a non-periodic variety $\mathsf{V}$ of semigroups satisfies some nontrivial identity of the form $Z_n\bumpeq T$ then all periodic members of $\mathsf{V}$ are locally finite.
\end{Thm}

\subsection{Sufficient condition for non-finite basedness}
Let $\mathsf{A}$ and $\mathsf{B}$ be two subclasses of a fixed class $\mathsf{C}$ of algebras. The \emph{Mal'cev product} $\mathsf{A}\malcev\mathsf{B}$ of $\mathsf{A}$ and $\mathsf{B}$ (within $\mathsf{C}$) is the class of all algebras $\mathcal{C}\in\mathsf{C}$ for which there exists a congruence $\theta$ such that the quotient algebra $\mathcal{C}/\theta$ lies in $\mathsf{B}$ while all $\theta$-classes that are subalgebras in $\mathcal{C}$ belong to $\mathsf{A}$. Note that for a congruence $\theta$ on a semigroup $\mathcal{S}$, a congruence class $s\theta$ forms a subsemigroup of $\mathcal{S}$ if and only if the element $s\theta$ is an idempotent of the quotient $\mathcal{S}/\theta$.

For two  semigroup varieties $\mathsf{V}_1$ and $\mathsf{V}_2$, their Mal'cev product $\mathsf{V}_1\malcev\mathsf{V}_2$ within the class of all semigroups may fail to be a variety, but it is always closed under forming subsemigroups and direct products, see \cite[Theorems~1 and~2]{Ma67}. Therefore the variety $\var(\mathsf{V}_1\malcev\mathsf{V}_2)$ generated by $\mathsf{V}_1\malcev\mathsf{V}_2$ is comprised of all
morphic images of the members of $\mathsf{V}_1\malcev\mathsf{V}_2$.

We are now ready to formulate a key result  which is a combination of Theorem 3.4 and Remark 3.5 in \cite{auingeretal} with Theorem~\ref{sapirlocallyfinite}. Recall that a semigroup variety $\mathsf{V}$ is non-periodic if and only if it contains the infinite cyclic semigroup if and only if it contains the variety $\mathsf{Com}$ of all commutative semigroups if and only if every identity $u\bumpeq v$ satisfied by $\mathsf{V}$ is \emph{balanced}, which means that every letter occurs the same number of times in $u$ and in $v$.
\begin{Thm}\label{thm:sufficient condition for nfb} Let $\mathsf{V}$ be a semigroup variety which does not satisfy any nontrivial identity $Z_k\bumpeq W$. If $\mathsf{V}\subseteq \var(\mathsf{U}\malcev\mathsf{L})$ for some locally finite variety $\mathsf{L}$ and a variety $\mathsf{U}$ which satisfies some nontrivial balanced identity $Z_n\bumpeq T$ then $\mathsf{V}$ is non-finitely based.
\end{Thm}
\begin{proof} According to Theorem 3.4 and Remark 3.5 in \cite{auingeretal} in order that $\mathsf V$ be non-finitely based we require that $\mathsf U$ is such that all periodic members are locally finite; by Theorem \ref{sapirlocallyfinite}, the latter is guaranteed if it satisfies a nontrivial balanced identity of the form $Z_n\bumpeq T$.
\end{proof}
There is also an involutory version of this theorem.
\begin{Thm}\label{thm:sufficient involutory version} Let $\mathsf{V}$ be a variety of involutory semigroups which does not satisfy any nontrivial (involutory semigroup) identity $Z_k\bumpeq W$. If the class of all semigroup reducts of $\mathsf{V}$ is contained in $\var(\mathsf{U}\malcev\mathsf{L})$ for some locally finite variety of semigroups $\mathsf{L}$ and a semigroup variety $\mathsf{U}$ which satisfies some nontrivial balanced identity $Z_n\bumpeq T$  then $\mathsf{V}$ is non-finitely based as a variety of involutory semigroups.
\end{Thm}

In Section \ref{finite basis problem} we shall be concerned with local monoids $\mathcal{M}$ of certain categories. These will be usually non-periodic, that is, they  admit elements of infinite order. We intend to prove that they do not admit finite equational bases. The strategy to do so will be to consider suitable epimorphisms $\Phi\colon\mathcal{M}\twoheadrightarrow \mathcal{F}$ onto finite quotients $\mathcal{F}$ and study the subsemigroups $e\Phi^{-1}$ of $\mathcal{M}$ which are inverse images of idempotents $e$ of $\mathcal{F}$ under $\Phi$. We shall identify some nontrivial identities satisfied by all such semigroups $e\Phi^{-1}$ and thereby shall see that all such satisfy some nontrivial balanced identity of the form $Z_n\bumpeq T$. From this it follows that all semigroups $e\Phi^{-1}$ belong to the variety $\mathsf{U}=[Z_n\bumpeq T]$. This means that $\mathcal{M}$ belongs to $\mathsf{U}\malcev\mathsf{L}$ for $\mathsf{L}:=\var\mathcal{F}$ which clearly is locally finite. Altogether $\var\mathcal{M}\subseteq \var(\mathsf{U}\malcev\mathsf{L})$. If, in addition, $\mathcal{M}$ itself and therefore also $\var\mathcal{M}$ does not satisfy any nontrivial identity of the form $Z_k\bumpeq W$ then, according to Theorem \ref{thm:sufficient condition for nfb}, $\var\mathcal{M}$ and hence also $\mathcal{M}$ itself is non-finitely based.

\subsection{Rees quotients and Rees matrix semigroups}\label{Rees matrix}
We shall use two constructions, both due to David Rees \cite{Rees}, which we briefly recall here.

\subsubsection{Rees quotients} A non-empty subset $A$ of a semigroup $\mathcal{S}$ is called an \emph{ideal} of $\mathcal{S}$ if $as,sa\in A$ for all $a\in A$ and $s\in\mathcal{S}$. Given an ideal $A$ of $\mathcal{S}$, take a fresh symbol $\mathbf{0}\notin\mathcal{S}$ and define a multiplication on the set $\mathcal{S}/A:=(\mathcal{S}\setminus A)\cup\{\mathbf{0}\}$ by keeping products within $\mathcal{S}\setminus A$ as they were in $\mathcal{S}$ and letting all other products be $\mathbf{0}$:
\[
s\cdot t=\begin{cases}
st &\text{if } s,t,st\in\mathcal{S}\setminus A,\\
\mathbf{0} &\text{otherwise}.
\end{cases} 
\]
Under this multiplication, $\mathcal{S}/A$ is a semigroup called the \emph{Rees quotient of $\mathcal{S}$ over $A$}. Clearly, $\mathbf{0}$ is a zero of $\mathcal{S}/A$. In this context, the semigroup $\mathcal S$ is sometimes called an \emph{ideal extension of $A$ by $\mathcal{S}/A$}. 


\subsubsection{Rees matrix semigroups and their identities} Let $\mathcal{S}$ be a semigroup, and $I,\Lambda$ non-empty sets. Given a $\Lambda\times I$-matrix $P=(p_{\lambda i})$ with entries from $\mathcal{S}$, we define a multiplication on the set $I\times\mathcal{S}\times\Lambda$ by the following rule:
\[
(i,s,\lambda)\cdot(j,t,\mu):=(i,sp_{\lambda j}t,\mu)\ \text{for all $i,j\in I$, $\lambda,\mu\in\Lambda$, $s,t\in\mathcal{S}$}.
\]
The multiplication is easily seen to be associative so that $I\times\mathcal{S}\times\Lambda$ becomes a semigroup. We denote it by $\mathcal{M}(I,\mathcal{S},\Lambda;P)$ and call it the \emph{\Rms{} over $\mathcal{S}$ with sandwich matrix $P$}.

If $\mathcal{S}$ has a zero $z$, then the set $Z$ of all triples with the middle coordinate $z$ is an ideal of $\mathcal{M}(I,\mathcal{S},\Lambda;P)$. The Rees quotient $\mathcal{M}(I,\mathcal{S},\Lambda;P)/Z$ is denoted by
$\mathcal{M}^0(I,\mathcal{S},\Lambda;P)$. Thus, $\mathcal{M}^0(I,\mathcal{S},\Lambda;P)$ is the set $I\times(\mathcal{S}\setminus\{z\})\times\Lambda\cup\{\mathbf{0}\}$ with the multiplication 
\begin{gather*}
a\cdot\mathbf{0}=\mathbf{0}\cdot a:=\mathbf{0}\ \text{ for all } a\in\mathcal{M}^0(I,\mathcal{S},\Lambda;P),\\
(i,s,\lambda)\cdot(j,t,\mu):=\begin{cases}
(i,sp_{\lambda j}t,\mu)&\ \text{if } sp_{\lambda j}t\ne z,\\
\mathbf{0} &\text{otherwise}.
\end{cases}
\end{gather*}
The name `\Rms' applies also to semigroups $\mathcal{M}^0(I,\mathcal{S},\Lambda;P)$.

\Rms{}s play a distinguished role in the structure theory of semigroups, but we do not need this theory. Rather, we shall use a few special instances of \Rms{}s, some of which we present here.

A \emph{rectangular band} is a \Rms{} over the one-element semigroup. In this case, the middle coordinate can be discarded as it is the same for all triples, and thus, a rectangular band can be identified with the set $I\times\Lambda$ whose pairs are multiplied by the rule $(i,\lambda)\cdot(j,\mu):=(i,\mu)$. One can readily check that every rectangular band satisfies the identity $xyx\bumpeq x$.

The 5-element semigroups $\mathcal{A}_2$ and $\mathcal{B}_2$ are the \Rms{}s $\mathcal{M}^0(I,\mathcal{S},\Lambda;P)$ in which
\begin{itemize}
  \item $I=\Lambda=\{1,2\}$,
  \item $\mathcal{S}=\{0,1\}$ with the usual multiplication of integers,
  \item the sandwich matrix $P$ is $\begin{pmatrix}1&1\\0&1\end{pmatrix}$ for $\mathcal{A}_2$ and $\begin{pmatrix}1&0\\0&1\end{pmatrix}$ for $\mathcal{B}_2$.
\end{itemize}
The semigroup $\mathcal{B}_2$ is known as the 5-element \emph{Brandt semigroup}.

The non-zero elements of each of the semigroups $\mathcal{A}_2$ and $\mathcal{B}_2$ are the triples with the same middle coordinate 1, so it is safe to use the pair $(i,\lambda)$ instead of the triple $(i,1,\lambda)$. In this notation, the carrier of both $\mathcal{A}_2$ and $\mathcal{B}_2$ is  
\[
\{\mathbf{0},(1,1),(1,2),(2,1),(2,2)\}.
\]
The non-zero elements are multiplied in $\mathcal{A}_2$ according to the rule 
\[
(i,\lambda)\cdot(j,\mu)=\begin{cases}
\mathbf{0}&\ \text{if } \lambda=2,\, j=1,\\
(i,\mu)&\text{otherwise},
\end{cases}
\]
and in $\mathcal{B}_2$ according to the rule 
\[
(i,\lambda)\cdot(j,\mu)=\begin{cases}
\mathbf{0}&\ \text{if } \lambda\ne j,\\
(i,\mu)&\text{otherwise}.
\end{cases}
\]

For a semigroup $\mathcal{S}$, the least monoid containing $\mathcal{S}$ is denoted by $\mathcal{S}^1$. Thus $\mathcal{S}^1=\mathcal{S}$ if $\mathcal{S}$ has a neutral element and  $\mathcal{S}^1:=\mathcal{S}\cup\{1\}$ otherwise; in the latter case the multiplication on $\mathcal{S}$ is extended to $\mathcal{S}^1$ such that the fresh symbol $1$ becomes the neutral element of $\mathcal{S}^1$.

\begin{Lemma}[{\cite[Lemma 3.7]{sapir}}]\label{lem:sapir}
 Neither of the monoids $\mathcal{A}_2^1$ and $\mathcal{B}_2^1$ satisfies any nontrivial identity $Z_k\bumpeq W$.
\end{Lemma}

A semigroup $\mathcal{S}$ \emph{divides} a semigroup $\mathcal{T}$ if $\mathcal{S}$ is a morphic image of a subsemigroup of $T$.

\begin{Cor}\label{cor:a2b2}
If either of the semigroups $\mathcal{A}_2$ and $\mathcal{B}_2$ divides a monoid $\mathcal{T}$, then $\mathcal{T}$ does not satisfy any nontrivial identity $Z_k\bumpeq W$.
\end{Cor}

\begin{proof}
Let  $\mathcal{S}$ be any of the semigroups $\mathcal{A}_2$ and $\mathcal{B}_2$. If $\mathcal{S}$ divides $\mathcal{T}$, there is a subsemigroup $\mathcal{T}'$ of $\mathcal{T}$ such that $\mathcal{S}$ is a morphic image of $\mathcal{T}'$. The neutral element 1 of the monoid $\mathcal{T}$ cannot lie in the subsemigroup $\mathcal{T}'$ as otherwise its image would be a neutral element of $\mathcal{S}$ and $\mathcal{S}$ has no neutral  element. The morphism $\mathcal{T}'\twoheadrightarrow\mathcal{S}$ extends in an obvious way to a morphism from the submonoid $\mathcal{T}'\cup\{1\}$ onto the monoid $\mathcal{S}^1$.  Since identities are inherited by submonoids and morphic images, Lemma~\ref{lem:sapir} implies that $\mathcal{T}$ cannot satisfy any nontrivial identity $Z_k\bumpeq W$.
\end{proof}

\begin{Def}\label{def:RMAb}
 Denote by $\mathrm{RMAb}$ the set of semigroup identities $u\bumpeq v$ satisfying the following properties:
\begin{enumerate}
\item the first letter of $u$ is the same as the first letter of $v$,
\item the last letter of $u$ is the same as the last letter of $v$,
\item every $2$-letter factor occurs the same number of times in $u$ and $v$.
\end{enumerate}
\end{Def}
The set  $\mathrm{RMAb}$ is known (see \cite[Theorem 9]{KR79}) to be precisely the set of all identities satisfied by every \Rms{} over an abelian group; hence the notation. We note that every identity in $\mathrm{RMAb}$ is balanced  \cite[Lemma 1]{kitovvolkov}. The variety $[\mathrm{RMAb}]$ is finitely based, but we do not need this fact. For our purposes the following observation well suffices.

\begin{Lemma}\label{lem:z3}
 The set $\mathrm{RMAb}$ contains the identity 
 \begin{equation}\label{eq:Z_3}
 Z_3\bumpeq x_1x_3x_1x_2x_1x_2x_1.
 \end{equation}
 \end{Lemma}

\begin{proof}
By the construction of Zimin words, $Z_3=x_1x_2x_1x_3x_1x_2x_1$. The claim then follows by direct inspection.  
\end{proof}

\section{The finite basis problem}\label{finite basis problem}
Our first goal is to show that the monoids $\overline{2\mathfrak{Cob}_n}$, $2\mathfrak{Cob}_n$ and related monoids do not have a finite basis for their identities.

Let $n\in \mathbb{N}_0$; in order to apply the results from \cite{auingeretal}, we shall study the canonical morphism $\Phi\colon\overline{2\mathfrak{Cob}_n}\twoheadrightarrow \mathfrak{P}_n$. More precisely, our purpose is to determine the (respectively, some) identities satisfied by the inverse image under $\Phi$ of an arbitrary idempotent of $\mathfrak{P}_n$. Our second goal is to study the same question for various monoids of annular type. The method is essentially the same: we study the canonical morphism $\Phi\colon \overline{\mathfrak{tAnn}_n}\twoheadrightarrow \mathfrak{Ann}_n$ and also for $\overline{\mathfrak{tAnn}_n}$ replaced with the other types of annular monoids. We first recall the characterization of the idempotents of $\mathfrak{P}_n$ \cite{dolinkaetal}.

\subsection{Idempotents of $\mathfrak{P}_n$} In this context it is convenient of consider partitions $\alpha\colon U\leadsto U$ where $U$ is a subset of $[n]$. Let $U\subseteq [n]$; the definition of a partition $\alpha\colon U\leadsto U$, and the composition $\alpha\beta$ of two such partitions $U\leadsto U$ is completely analogous to the case $U=[n]$, introduced in Subsection \ref{partition category}. Thereby we obtain the partition monoid $\mathfrak{P}_U$ (which can be interpreted as a submonoid of $\mathfrak{P}_n$). Let
$$[n]=U_1\cup\cdots\cup U_k$$
be a partition of $[n]$. Every $k$-tuple $(\alpha_1,\dots,\alpha_k)$ with $\alpha_i\colon U_i\leadsto U_i$ gives rise to the partition $\alpha_1\cup\cdots\cup\alpha_k\colon [n]\leadsto  [n]$, and the direct product $\mathfrak{P}_{U_1}\times\cdots\times\mathfrak{P}_{U_k}$ thereby embeds in $\mathfrak{P}_n$.

Every partition  $\alpha\colon U\leadsto U$ induces two partitions $\overleftarrow{\alpha}$ and $\overrightarrow{\alpha}$ of $U$ by letting $\overleftarrow{\alpha}$ be the restriction of $\alpha$ to the first component of $U\sqcup U$ and $\overrightarrow{\alpha}$ be the restriction to the second component. According to \cite{dolinkaetal} we call $\alpha$ \emph{irreducible} if $\overleftarrow{\alpha}\vee \overrightarrow{\alpha}$ is the universal relation on $U$ (to be denoted $\omega_U$ in the following). Recall that the \emph{rank} of a partition is the number of its transversal blocks.
\begin{Res}[{\cite[Theorem 5]{dolinkaetal}}]\label{idempotents in Pn} A partition $\varepsilon\colon[n]\leadsto [n]$ is idempotent if and only if there exists a partition $[n]=U_1\cup\cdots\cup U_k$ such that
\begin{enumerate}
\item $\varepsilon=\varepsilon_1\cup\cdots \cup\varepsilon_k$ where $\varepsilon_i=\varepsilon\vert_{U_i\sqcup U_i}$,
\item every $\varepsilon_i$ is an irreducible partition $U_i\leadsto U_i$ of rank at most $1$.
\end{enumerate}
\end{Res}

\subsection{The finite basis problem  for $2\mathfrak{Cob}_n$ and $\overline{2\mathfrak{Cob}_n}$}\label{subsec:nfb-cobord} Result \ref{idempotents in Pn} allows us to prove a crucial theorem.
\begin{Thm}\label{Thm:ids in e phi^{-1}} Let $\Phi\colon \overline{2\mathfrak{Cob}_n}\twoheadrightarrow\mathfrak{P}_n$ be the canonical morphism $$(\alpha,g,s)\mapsto \alpha$$ and let $\varepsilon$ be an idempotent of $\mathfrak{P}_n$. Then the semigroup $\varepsilon\Phi^{-1}$ satisfies the identity \eqref{eq:Z_3}.
\end{Thm}
\begin{proof} Let $k\ge 1$ and for $i=1,\dots,k$ let $\varepsilon_i=(\varepsilon,g_i,s_i)\in \overline{2\mathfrak{Cob}_n}$ with $\varepsilon$ an idempotent in $\mathfrak{P}_n$ and
\begin{equation}\label{product of eps}
\varepsilon_1\cdots\varepsilon_k=(\varepsilon,g,s).
\end{equation}
According to Result \ref{idempotents in Pn} it is sufficient to prove the claim for $\varepsilon$ an irreducible idempotent and every such idempotent has rank $0$ or $1$. Assume first that the rank of $\varepsilon$ is $1$. Since $\overleftarrow{\varepsilon}\vee\overrightarrow{\varepsilon}=\omega_{[n]}$, in the composition of $\varepsilon$ with itself no dead blocks are created, that is $\mathrm{b}(\varepsilon,\varepsilon)=0$. From this it follows that
\begin{equation}\label{total s}
s=\sum_{\ell=1}^k s_{\ell}.
\end{equation}
We need to calculate the function $g\colon\varepsilon\to \mathbb{Z}$. Let $L_1,\dots, L_p, T, R_1,\dots R_q$ be the upper blocks, the transversal block, and the lower blocks of $\varepsilon$, respectively. From the definition (\ref{def:composition in 2Cob discrete}) it follows that $g(L_i)=g_1(L_i)$ and $g(R_j)=g_k(R_j)$ for all $i=1,\dots,p$ and $j=1,\dots,q$. It remains to calculate $g(T)$. In the composition $\underbrace{\varepsilon\cdots\varepsilon}_k$ the block $T$ is formed by merging:
\begin{enumerate}
\item all lower blocks of the first $\varepsilon$,
\item all upper blocks of the last $\varepsilon$,
\item the transversal block of every $\varepsilon$,
\item all one-sided blocks of all but the first and the last $\varepsilon$.
\end{enumerate}
The label $g(T)$ then is obtained by adding all labels of these blocks and the increments obtained in the compositions $\varepsilon_\ell\varepsilon_{\ell+1}$ for $\ell=1,\dots,k-1$. All these increments coincide and are equal to $n-(p+1+q+1)+1=n-p-q-1$. Altogether we obtain
\begin{equation}\label{label of T}
g(T)=\sum_{\ell=1}^kg_{\ell}(T)+\sum_{\ell=1}^{k-1}\sum_{j=1}^qg_\ell(R_j)+\sum_{\ell=2}^{k}\sum_{i=1}^p g_\ell(L_i)+(k-1)(n-p-q-1).
\end{equation}
The equalities (\ref{total s}) and (\ref{label of T}) imply that $s$ and $g$ only depend on the first element $\varepsilon_1$, the last  $\varepsilon_k$ of the sequence $\varepsilon_1,\dots,\varepsilon_k$, the number $k$, and the number how often every element $\varepsilon_i$ occurs in this sequence.

Now consider the case when the rank of $\varepsilon$ is $0$. In this case $\varepsilon$ has upper blocks and lower blocks $L_1,\dots, L_p, R_1\dots, R_q$, say (and no transversal block). In this case again, $g(L_i)=g_1(L_i)$ and $g(R_j)=g_k(R_j)$ for all $i$ and $j$. It follows that $\varepsilon$ and $g$ only depend on the first factor $\varepsilon_1$ and the last factor $\varepsilon_k$. Once more using the fact that $\overleftarrow{\varepsilon}\vee\overrightarrow{\varepsilon}=\omega_{[n]}$ we see that the composition of $\varepsilon$ with itself creates exactly one dead block, that is, $\mathrm{b}(\varepsilon,\varepsilon)=1$. In the composition $\varepsilon_\ell\varepsilon_{\ell+1}$ the label of the thereby created dead block by (\ref{def:g}) is
\begin{equation}\label{genus of new component}
\sum_{j=1}^qg_\ell(R_j)+\sum_{i=1}^pg_{\ell+1}(L_i) +n-(p+q)+1.
\end{equation}
It follows that the third entry $s$ in $(\varepsilon,g,s)$ in (\ref{product of eps})  is determined (only) by the number how often every $\varepsilon_i$ occurs in the sequence $\varepsilon_1,\dots,\varepsilon_k$ (via $s_1,\dots,s_k$) and on how  often every pair  $(\varepsilon_\ell,\varepsilon_{\ell+1})$ occurs in the sequence $(\varepsilon_1,\varepsilon_2),\dots(\varepsilon_{k-1},\varepsilon_k)$, the latter by (\ref{genus of new component}).

The above analysis shows that whenever words $u$ and $v$ satisfy conditions (1)--(3) of Definition~\ref{def:RMAb}, the images of $u$ and $v$ under every substitution $X^+\to\varepsilon\Phi^{-1}$ coincide, that is, the semigroup satisfies every identity in the set $\mathrm{RMAb}$, including~\eqref{eq:Z_3}. 
\end{proof}

\begin{Prop}\label{2Cob1} The monoid $2\mathfrak{Cob}_1$ \textup(and therefore also $2\mathfrak{Cob}_n$ for every $n\ge 1$\textup) does not satisfy any nontrivial identity $Z_k\bumpeq W$.
\end{Prop}
\begin{proof}
By Corollary~\ref{cor:a2b2}, it suffices to show that the semigroup $\mathcal{A}_2$ divides the monoid $2\mathfrak{Cob}_1$.

Let $\mathcal{S}$ be the subsemigroup of $2\mathfrak{Cob}_1$ consisting of all $(\alpha,g,s)$ where
\begin{itemize}
\item $\alpha=\mathbf{0}_{[1]\sqcup[1]}$,
\item $g$ takes only the values $0$ and $1$
\end{itemize}
and $\mathbf{0}_{[1]\sqcup[1]}$ denotes the identity relation on $[1]\sqcup[1]$ (every equivalence class is a singleton), interpreted as a partition $[1]\leadsto [1]$.
All possible choices of $(\mathbf{0}_{[1]\sqcup[1]},g)$ can be realized by the set of pairs $$\{(0,0),(0,1),(1,0),(1,1)\}$$ depending on whether $g$ takes the value $0$ or $1$ on the incoming respectively outgoing vertex. Hence $\mathcal{S}$ can be realized as the set of all triples $(i,s,j)$ with $i,j\in \{0,1\}$ and sequences $s=(s_i)_{i\ge 0}\in\oplus_{i\in \mathbb{N}_0}\mathbb{N}_0$. The composition can be described by the formula
\begin{equation*}
(i,s,j)(k,t,l)=(i,s+b(j,k)+t,\ell)
\end{equation*}
where
\begin{equation*}
b(j,k)=(\delta_{n,j+k})_{n\ge 0}
\end{equation*}
and $\delta$ is the Kronecker symbol. In other words, $b(j,k)$ adds $1$ to the $(j+k)$th entry of $s+t$ and leaves the other entries unchanged.

We define a map $\mathcal{S}\to \mathcal{A}_2$ by
\begin{equation}\label{A21 dividing 2Cob1}
(i,s,j)\mapsto\begin{cases} (i,j) &\mbox{ if }s_2=0\\
							 0 &\mbox{ if }s_2> 0
\end{cases}
\end{equation}
which is readily checked to be a morphism $\mathcal{S}\twoheadrightarrow \mathcal{A}_2$.
 \end{proof}

\begin{Prop} For every $n\ge 2$, the monoid $2\mathfrak{Cob}^\circ_n$ does not satisfy any nontrivial identity $Z_k\bumpeq W$.
\end{Prop}
\begin{proof} There is a canonical morphism $2\mathfrak{Cob}^\circ_2\to\mathfrak{P}_2$ and the latter monoid is known to have the semigroup $\mathcal{A}_2$ among its divisors; see the proof of \cite[Proposition 2.11]{adv2}. Now Corollary~\ref{cor:a2b2} applies.
\end{proof}
So we formulate the main result of the present subsection.
\begin{Thm} \begin{enumerate}
\item For every $n\ge 1$ the monoids $2\mathfrak{Cob}_n$ and $\overline{2\mathfrak{Cob}_n}$ are non-finitely based.
\item For every $n\ge 2$ the monoids $2\mathfrak{Cob}^\circ_n$ and $\overline{2\mathfrak{Cob}^\circ_n}$ are non-finitely based.
\item For every $n\ge 2$ the monoids $\mathfrak{P}^{\mathrm d}_n$ and $\overline{\mathfrak{P}^{\mathrm d}_n}$ are non-finitely based.
\end{enumerate}
\end{Thm}
For statement (3) one has to take into account that the inverse image of any idempotent under the canonical morphism $\overline{\mathfrak{P}^{\mathrm d}_n}\twoheadrightarrow \mathfrak{P}_n$ is always commutative while the arguments above show that these monoids do not satisfy any nontrivial identity $Z_k\bumpeq W$. Notice that $\mathfrak{P}_1$ and $\overline{\mathfrak{P}_1}$ are commutative hence are finitely based. We are left with the cases $2\mathfrak{Cob}^\circ_1$ and $\overline{2\mathfrak{Cob}^\circ_1}$.
\begin{Thm}\label{Thm:2Cob_1} The monoids $2\mathfrak{Cob}^\circ_1$ and $\overline{2\mathfrak{Cob}^\circ_1}$ are finitely based. The identities
\begin{equation}\label{finite1}
    xt_1yt_2\cdot xy\cdot t_3yt_4x=xt_1yt_2\cdot yx\cdot t_3yt_4x
\end{equation} and
\begin{equation}\label{finite2}
    xt_1yt_2\cdot xy\cdot t_3xt_4y=xt_1yt_2\cdot yx\cdot t_3xt_4y
\end{equation}
form a basis for the monoid identities of $2\mathfrak{Cob}^\circ_1$ as well as of $\overline{2\mathfrak{Cob}^\circ_1}$. A basis of semigroup identities is  obtained by deleting every possible choice of a subset of the set of letters $\{t_1,t_2,t_3,t_4\}$.
\end{Thm}
The (lengthy) proof of this statement will be deferred to Subsection \ref{finite basedness}.

\subsection{The finite basis problem for monoids of annular type}\label{subsec:nfb-annular}
\subsubsection{Affine Temperley--Lieb categories} In the following let $n\ge 2$; we study the canonical morphisms $\Psi\colon\mathfrak{aTL}^{\mathrm e}_n\twoheadrightarrow \mathfrak{Ann}_n$, $\Psi_1\colon\mathfrak{aTL}_n\twoheadrightarrow \mathfrak{Ann}_n$ and $\Psi_2\colon \mathfrak{aTL}^{\mathrm d}_n\twoheadrightarrow \mathfrak{Ann}_n$.

\begin{Prop}\label{aTLe} For every idempotent $\vp\in \mathfrak{Ann}_n$ the semigroup $\vp\Psi^{-1}$ satisfies the identity $axyb\bumpeq ayxb$.
\end{Prop}
\begin{proof} Let $\vp$ be an idempotent of $\mathfrak{Ann}_n$. Assume first that $\vp$ has positive rank. Since $\mathfrak{aTL}^{\mathrm e}_n$ is regular, a well-known fact from the theory of regular semigroups (Lallement's lemma, \cite[Lemma 2.2]{Lall66}) applies, showing the existence of an idempotent $\bm\vp\in \vp\Psi^{-1}$. By Lemma \ref{lem:map atle->ann}, for any $\bm\alpha,\bm\beta\in \vp\Psi^{-1}$ there exist $s,t\in \mathbb{Z}$ such that $\bm\alpha=\bm\lambda_n^s\bm\vp$ and $\bm\beta=\bm\lambda_n^t\bm\vp$. Since by Lemma \ref{lem:action of lambda} $\bm\lambda_n^s$ and $\bm\lambda_n^t$ are central in $\mathfrak{aTL}^{\mathrm e}$, we get $\bm\alpha\bm\beta=\bm\beta\bm\alpha$. Suppose that $\rk(\vp)=0$; then $\vp\Psi^{-1}$ consists of rank $0$ elements of $\mathfrak{aTL}^{\mathrm e}$. The set of all rank $0$ elements of $\mathfrak{aTL}^{\mathrm e}_n$ forms a rectangular band, hence $\vp\Psi^{-1}$ forms a rectangular band. Altogether, $\vp\Psi^{-1}$ is commutative or a rectangular band, so in any case satisfies the identity $axyb\bumpeq ayxb$.
\end{proof}

Next we consider $\Psi_2\colon \mathfrak{aTL}^{\mathrm d}_n\twoheadrightarrow \mathfrak{Ann}_n$.
 \begin{Prop}\label{aTLd}
 Let  $\vp$ be an idempotent of $\mathfrak{Ann}_n$. Then the semigroup $\vp\Psi_2^{-1}$ satisfies  the identity \eqref{eq:Z_3}.
 \end{Prop}
\begin{proof} Let $\vp$ be an idempotent of $\mathfrak{Ann}_n$ and suppose first that $\rk(\vp)>0$; all elements of $\vp\Psi_2^{-1}$ are of the form $(\bm\alpha,0,a)$ ($a$ is the number of $0$-circles) and $\bm\alpha\in \vp\Psi^{-1}$. According to the proof of Proposition \ref{aTLe}, $\bm\alpha=\bm\lambda_n^s\bm\vp$ for some idempotent $\bm\vp\in \vp\Psi^{-1}$ and some $s\in \mathbb Z$ (in fact, the idempotent $\bm\vp$ is unique). For two elements $\bm\lambda_n^s\bm\vp$ and $\bm\lambda_n^t\bm\vp$ of this kind, we note that
$$\mathrm{b}_0(\bm\lambda_n^s\bm\vp,\bm\lambda_n^t\bm\vp)= \mathrm{b}_0(\bm\vp,\bm\vp)=\mathrm{b}_0(\bm\lambda_n^t\bm\vp,\bm\lambda_n^s\bm\vp)$$
since the number $\mathrm{b}_0$ only depends on the in-strings and out-strings of the involved elements, and $\bm\vp$ and $\bm\lambda_n^s\bm\vp$  (for any $s$) have the same in- and out-strings as was shown in the proof of Lemma \ref{lem:action of lambda}. That is, $\mathrm{b}_0(\bm\alpha,\bm\beta)=\mathrm{b}_0(\bm\beta,\bm\alpha)$, so that, in this case, $\vp\Psi^{-1}_2$ is commutative.

We are left with the case $\rk(\vp)=0$. Let $\bm\alpha_1=(\bm\vp_1,a_1,c_2),\dots,\bm\alpha_t=(\bm\vp_t,a_t,c_t)\in\vp\Psi_2^{-1}$ with $\bm\vp_i\in \vp\Psi^{-1}$ and $a_i,c_i\in \mathbb{N}_0$ (the elements $\bm\vp_i$ are rank $0$ elements of $\mathfrak{aTL}^{\mathrm e}_n$). Then for the product of these elements we have
$$(\bm\vp_1,a_1,c_1)\cdots(\bm\vp_t,a_t,c_t)=(\bm\vp_1\bm\vp_t,a,c)$$
where
$$a=\sum_{i=1}^ta_i+\sum_{i=1}^{t-1}\mathrm{b}_0(\bm\vp_i,\bm\vp_{i+1})$$
and
$$c=\sum_{i=1}^tc_i+\sum_{i=1}^{t-1}\mathrm{b}_\omega(\bm\vp_i,\bm\vp_{i+1}).$$
We see that the value of this product only depends
\begin{enumerate}
\item on the first element and the last element  of the sequence $\bm\alpha_1,\dots,\bm\alpha_t$,
\item on the number how often each element $\bm\alpha_i$ occurs in the sequence $\bm\alpha_1,\dots,\bm\alpha_t$,
\item on the number how often each pair $(\bm\alpha_i,\bm\alpha_{i+1})$ occurs in the sequence $(\bm\alpha_1,\bm\alpha_2),\dots,(\bm\alpha_{t-1},\bm\alpha_t)$.
\end{enumerate}
Therefore the images of words $u$ and $v$ satisfying conditions (1)--(3) of Definition~\ref{def:RMAb} under every substitution $X^+\to\vp\Psi_2^{-1}$ coincide. Hence,  the semigroup $\vp\Psi_2^{-1}$ satisfies every identity in the set $\mathrm{RMAb}$, including~\eqref{eq:Z_3}. 
\end{proof}

\begin{Cor}\label{cor:identities subsemigroups} Let $n\ge 2$ and $\vp$ be an idempotent in $\mathfrak{Ann}_n$; then each of the semigroups $\vp\Psi^{-1}\subseteq \mathfrak{aTL}^{\mathrm e}_n$, $\vp\Psi^{-1}_1\subseteq \mathfrak{aTL}_n$ and $\vp\Psi_2^{-1}\subseteq \mathfrak{aTL}^{\mathrm d}_n$ satisfies the identity $Z_3\bumpeq x_1x_3x_1x_2x_1x_2x_1$.
\end{Cor}
\begin{proof} For the first and the third cases this follows from Propositions \ref{aTLe} and \ref{aTLd}, respectively. For $\vp\Psi^{-1}_1$ this follows from the fact that this semigroup is a quotient of $\vp\Psi_2^{-1}$.
\end{proof}

\subsubsection{The categories $\mathfrak{tAnn}$ and $\mathfrak{L}$}\label{sec: tAnn and L}
We start with a `negative' observation. For even $n$, the ideal $\mathfrak{O}_n$ of $\mathfrak{tAnn}_n$ consisting of all t$A$-diagrams of rank $0$ does not satisfy any nontrivial semigroup identity (and hence is finitely based and so is $\mathfrak{tAnn}_n$). Indeed, take the  following t$A$-diagrams from $\mathfrak{O}_n$:
\begin{center}
\begin{tikzpicture}[scale=0.95]
\draw (-0.25,2) node {$\alpha=$};
\draw[dotted](0.625,3)--(4.375,3);
\draw[dotted](6.625,3)--(10.375,3);
\draw[dotted](0.625,1)--(4.375,1);
\draw[dotted](8.625,1)--(10.375,1);
\draw (5.5,1) node{$\mathbf{\centerdot\,\centerdot\,\centerdot}$};
\draw (5.5,3) node{$\mathbf{\centerdot\,\centerdot\,\centerdot}$};
\foreach \x in {1,...,4,7,8,9,10} \foreach \y in {1,3} \filldraw (\x,\y) circle (1.5pt);
\foreach \x in {1,...,4} \draw[above] (\x,3) node{\x};
\foreach \x in {1,...,4} \draw[below] (\x,1) node{\x};
\draw[above] (7,3) node{$n{-}3$};
\draw[above] (8,3) node{$n{-}2$};
\draw[above] (9,3) node{$n{-}1$};
\draw[above] (10,3.1) node{$n$};
\draw[below] (7,1) node{$n{-}3$};
\draw[below] (8,1) node{$n{-}2$};
\draw[below] (9,1) node{$n{-}1$};
\draw[below] (10,.9) node{$n$};
\draw (1,3) to [out=315,in=225] (2,3);
\draw(3,3) to [out=315,in=225] (4,3);
\draw (7,3) to [out=315,in=225] (8,3);
\draw(9,3) to [out=315,in=225] (10,3);
\draw (1,1) to [out=45,in=135] (2,1);
\draw(3,1) to [out=45,in=135] (4,1);
\draw (7,1) to [out=45,in=135] (8,1);
\draw(9,1) to [out=45,in=135] (10,1);
\end{tikzpicture}
\end{center}
and
\begin{center}
\begin{tikzpicture}[scale=0.95]
\draw (-0.25,2) node {$\beta=$};
\draw[dotted](0.625,3)--(4.375,3);
\draw[dotted](6.625,3)--(10.375,3);
\draw[dotted](0.625,1)--(4.375,1);
\draw[dotted](8.625,1)--(10.375,1);
\draw (0.625,2)--(10.375,2);
\draw (5.5,1) node{$\mathbf{\centerdot\,\centerdot\,\centerdot}$};
\draw (5.5,3) node{$\mathbf{\centerdot\,\centerdot\,\centerdot}$};
\foreach \x in {1,...,4,7,8,9,10} \foreach \y in {1,3} \filldraw (\x,\y) circle (1.5pt);
\foreach \x in {1,...,4} \draw[above] (\x,3) node{\x};
\foreach \x in {1,...,4} \draw[below] (\x,1) node{\x};
\draw[above] (7,3) node{$n{-}3$};
\draw[above] (8,3) node{$n{-}2$};
\draw[above] (9,3) node{$n{-}1$};
\draw[above] (10,3.1) node{$n$};
\draw[below] (7,1) node{$n{-}3$};
\draw[below] (8,1) node{$n{-}2$};
\draw[below] (9,1) node{$n{-}1$};
\draw[below] (10,.9) node{$n$};
\draw (1,3) to [out=315,in=225] (2,3);
\draw(3,3) to [out=315,in=225] (4,3);
\draw (7,3) to [out=315,in=225] (8,3);
\draw(9,3) to [out=315,in=225] (10,3);
\draw (1,1) to [out=45,in=135] (2,1);
\draw(3,1) to [out=45,in=135] (4,1);
\draw (7,1) to [out=45,in=135] (8,1);
\draw(9,1) to [out=45,in=135] (10,1);
\end{tikzpicture}
\end{center}
(So, $\beta$ differs from $\alpha$ only by the presence of an $\omega$-circle, shown as the solid line in the middle; the dotted lines indicate the inner and outer circles of the annulus.) The map $\begin{cases} x\mapsto\alpha\\ y\mapsto\beta\end{cases}$ extends to a morphism $\vartheta$ from the free semigroup $\{x,y\}^+$ onto the subsemigroup $\langle\alpha,\beta\rangle$ of $\mathfrak{O}_n$. We prove that every word $w\in\{x,y\}^+$ is uniquely determined by its image $w\vartheta$. Each composition of two t$A$-diagrams from $\langle\alpha,\beta\rangle$ creates $\frac{n}2$ 0-circles but no $\omega$-circles. Hence the number of $\omega$-circles in $w\vartheta$ is equal to the number of occurrences of the letter $y$ in the word $w$; let $k$ be this number. If $k=0$, that is, $w=x^m$ for some $m>1$, the number $c$ of 0-circles in the t$A$-diagram $w\vartheta=\alpha^m$ is $c=\frac{n}2(m-1)$, whence $m=\frac{2c}{n}+1$. If $k>0$ and $c_0+\bigcirc+c_1+\cdots+c_{k-1}+\bigcirc+c_k\in \mathbb{O}\mathop{*}\mathbb{M}$  is the label of the unique transversal region of $w\vartheta$ as in Figure~\ref{fig:rank0label} (see Subsection 3.5.1), then
\begin{itemize}
  \item $\frac{2c_0}{n}$ is the number of occurrences of the letter $x$ preceding the first occurrence of $y$ in $w$;
  \item $\frac{2c_i}{n}-1$ is the number of occurrences of the letter $x$ between the $i$th and the $(i+1)$st occurrences of $y$ in $w$, $1\le i<k$;
  \item $\frac{2c_k}{n}$ is the number of occurrences of the letter $x$ following the last occurrence of $y$ in $w$.
\end{itemize}
We see that the word $w$ can be recovered from its image $w\vartheta$ as
\[
w=\begin{cases} x^{\frac{2c}n+1} &\text{if }\ k=0,\\ 
x^{\frac{2c_0}n}yx^{\frac{2c_1}n-1}y\cdots yx^{\frac{2c_{k-1}}n-1}yx^{\frac{2c_k}n} &\text{if }\ k>0.
\end{cases}
\]
Hence, the morphism $\vartheta$ is injective, and the free semigroup on two generators embeds into the ideal $\mathfrak{O}_n$. The proof also shows that the free semigroup on two generators embeds into the ideal $\mathfrak{O}_n^\flat$ of all rank $0$ elements of $\mathfrak{t^\flat Ann}_n$ (see Subsection~\ref{sec:factored tAnn}).

We have thus proved:
\begin{Thm}\label{thm:tAnn even} For every even $n$, the monoids $\mathfrak{tAnn}_n$ and $\overline{\mathfrak{tAnn}_n}$ as well as $\mathfrak{t^\flat Ann}_n$ and $\overline{\mathfrak{t^\flat Ann}_n}$  do not satisfy any nontrivial identity and hence are finitely based.
\end{Thm}

In order to proceed, recall the representation of the arrows of $\mathfrak{tAnn}$ discussed in Subsection \ref{subsection:combinatorial description} {and how the composition of their arrows work described in Subsection \ref{sec:description of arrow composition}. }
Let $n>0$ and consider the canonical morphism $\Phi\colon\mathfrak{tAnn}_n\twoheadrightarrow \mathfrak{Ann}_n$. Let $\vp$ be an idempotent of $\mathfrak{Ann}_n$ with $\rk(\vp)=r>0$ and let  $\bm\alpha_1,\bm\alpha_2,\dots,\bm\alpha_t\in \vp\Phi^{-1}$ and $\bm\alpha=\bm\alpha_1\cdots\bm\alpha_t$.  We can choose $s_1,\dots,s_t\in \mathbb{Z}$ such that for every $i$, $\bm\lambda_i^{- s_i}\bm\alpha_i=\bm\vp_i$ where $\underline{\bm\vp_i}$ does not depend on $i$ and $\bm\vp:=\underline{\bm\vp_i}$ is idempotent in $\mathfrak{aTL}^{\mathrm e}_n$ and $\bm\alpha_i\Phi=\bm\vp_i\Phi=\bm\vp\Psi=\vp$. We note that the elements $\bm\lambda_n^{s}$ (interpreted in $\mathfrak{tAnn}_n$ as: the diagram of $\mathfrak{aTL}^{\mathrm e}$ with empty label on every region of $\bm\lambda_n^{s}$) are central in $\mathfrak{tAnn}_n$. Hence we have
\begin{equation}\label{labeling the product alpha}
\bm\alpha_1\cdots \bm\alpha_t=\bm\lambda_n^s\bm\vp_1\cdots \bm\vp_t
\end{equation}
for $s=\sum_{i=1}^ts_i$ and we now look at the product $\bm\vp_1\cdots\bm\vp_t$. All elements $\bm\vp_i$ have the same skeleton $\bm\vp$ but they have individual labeling functions $\mathbf{c}_i$. So, $\bm\vp_i=(\bm\vp,\mathbf{c}_i)$ according to the pair representation discussed in Subsection \ref{subsection:combinatorial description}, we are looking at the product
\begin{equation}\label{labeling the product}
\bm\vp_1\cdots \bm\vp_t=(\bm\vp,\mathbf{c}_1)\cdots(\bm\vp,\mathbf{c}_t)=(\bm\vp,\mathbf{c}),
\end{equation}
 and we need to discuss the labeling function $\mathbf c$.

Since $\rk(\bm\vp_1\cdots \bm\vp_t)=\rk(\bm\vp_i)=r$ for all $i$, a first observation is that $\bm\vp_1\cdots \bm\vp_t$ has no more in-regions than $\bm\vp_1$ and no more out-regions than $\bm\vp_t$ and for every such region $R$ we have $\mathbf{c}(R)=\mathbf{c}_1(R)$ for every in-region $R$ and $\mathbf{c}(R)=\mathbf{c}_t(R)$ for every out-region $R$.

For any t$A$-diagram of non-zero rank, the transversal regions are in a bijective correspondence with the transversal strings, for example by assigning to a transversal string the transversal region that is to the right hand side if the transversal string is traversed from the incoming point to the outgoing point. Consider some transversal string in $\bm\vp_1\cdots\bm\vp_t$. The transversal region $T$ to its right, say, is obtained as the union
$$T^{(1)}\cup Q_1\cup T^{(2)} \cup \cdots \cup Q_{t-1}\cup T^{(t)}$$ of the corresponding transversal regions  $T^{(i)}$ of $\bm\vp_i$ and every $Q_i$ is a union of certain out-regions $R^{(i)}_k$ ($k=1,\dots,p)$ of $\bm\vp_i$ and certain in-regions $L^{(i+1)}_\ell$ of $\bm\vp_{i+1}$ ($\ell=1,\dots,q$). Note that these regions are the same for all $i$. For every pair $(i,i+1)$ let $\inc(T;\mathbf{c}_i,\mathbf{c}_{i+1})$ be the \emph{increment of the label of $T$ when $\bm\vp_i$ is composed with $\bm\vp_{i+1}$}. The label $\mathbf{c}(T)$ then can be written as
\begin{equation*}
\mathbf{c}(T)=\sum_{i=1}^t\mathbf{c}_i(T)+\sum_{i=1}^{t-1}\sum_{k=1}^p\mathbf{c}_i(R_k)+\sum_{i=2}^{t}
\sum_{\ell=1}^q\mathbf{c}_i(L_\ell)+\sum_{i=1}^{t-1}\inc(T;\mathbf{c}_i,\mathbf{c}_{i+1}).
\end{equation*}

Since all labels commute with each other, the contribution of $\bm\vp_i$ to the label of $T$ does not depend on its position within the sequence $(\bm\vp_1,\dots,\bm\vp_t)$, but only on how often it occurs; likewise the contribution of the pair $(\bm\vp_i,\bm\vp_{i+1})$ to the label of $T$ (via increments) does not depend on the position of the pair $(\bm\vp_i,\bm\vp_{i+1})$ within the sequence $(\bm\vp_1,\bm\vp_2),\dots,(\bm\vp_{t-1},\bm\vp_t)$, but only on how often it occurs. Altogether we have shown that the labeling function $\mathbf{c}$ in the product (\ref{labeling the product}) only depends
\begin{enumerate}
\item on the first element $\bm\vp_1$ and the last element $\bm\vp_t$ of the sequence $\bm\vp_1,\dots,\bm\vp_t$,
\item on the number of how often every element $\bm\vp_i$ occurs in the sequence $\bm\vp_1,\dots,\bm\vp_t$,
\item on the number of how often every pair $(\bm\vp_i,\bm\vp_{i+1})$ occurs in the sequence $(\bm\vp_1,\bm\vp_2),(\bm\vp_2,\bm\vp_3),\dots,(\bm\vp_{i-1},\bm\vp_t)$.
\end{enumerate}
Moreover, the contribution of $\bm\alpha_i$ to the number $s=\sum s_i$ also does not depend on the position of $\bm\alpha_i$ within the sequence $\bm\alpha_1,\dots,\bm\alpha_t$, but only how often it occurs. Hence we have shown:
\begin{Prop}\label{prop:Z_3} 
The value of the product \eqref{labeling the product alpha} only depends
\begin{enumerate}
\item on the first and the last elements of the sequence $\bm\alpha_1,\dots,\bm\alpha_t$,
\item on the number of how often every element $\bm\alpha_i$ occurs in the sequence $\bm\alpha_1,\dots,\bm\alpha_t$,
\item on the number of how often every pair $(\bm\alpha_i,\bm\alpha_{i+1})$ occurs in the sequence $(\bm\alpha_1,\bm\alpha_2),(\bm\alpha_2,\bm\alpha_3),\dots,(\bm\alpha_{t-1},\bm\alpha_t)$.
\end{enumerate}
\end{Prop}
Proposition~\ref{prop:Z_3}, combined with Lemma~\ref{lem:z3} yields the central result of the present subsection.
\begin{Thm} \label{tAnn}
 Let $n > 0$,  $\Phi\colon \mathfrak{tAnn}_n\twoheadrightarrow \mathfrak{Ann}_n$ be the canonical morphism and let $\vp$ be an idempotent of $\mathfrak{Ann}_n$ which has non-zero rank. Then the semigroup $\vp\Phi^{-1}$ satisfies the identity \eqref{eq:Z_3}.
\end{Thm}
For rank zero idempotents the analogous claim is false; this follows from the discussion at the beginning of the present subsection concerning  the ideal $\mathfrak{Q}_n$ of rank zero elements in $\mathfrak{tAnn}_n$ for even $n$.

If we restrict the morphism $\Phi\colon \mathfrak{tAnn}_{n}\twoheadrightarrow \mathfrak{Ann}_{n}$ to $\mathfrak{L}_n$ we obtain the canonical morphisms $\mathfrak{L}_n\twoheadrightarrow \mathfrak{L}^\flat_n\twoheadrightarrow\mathfrak{J}_n$ where the latter is the \emph{Jones monoid} (or \emph{Temperley--Lieb monoid}) \cite{auingeretal}. In particular, Theorem \ref{tAnn} holds \emph{mutatis mutandis} for $\mathfrak{tAnn}_n$ and $\mathfrak{Ann}_n$ replaced with $\mathfrak{L}_n$ (respectively $\mathfrak{L}^\flat_n$) and $\mathfrak{J}_n$; in this case, the assumption that $\vp$ has non-zero rank can be dropped.
\subsubsection{Summary}
We are ready to formulate the summary of our results concerning the existence of finite equational bases for monoids of annular type.
For even $n$ we set $\mathfrak{tAnn}_n^\circ:=\mathfrak{tAnn}_n/\mathfrak{O}_n$ where $\mathfrak{O}_n$ is the ideal of $\mathfrak{tAnn}_n$ consisting of all rank zero elements; likewise, we set $\mathfrak{t^\flat Ann}^\circ_n:=\mathfrak{t^\flat Ann}_n/\mathfrak{O}^\flat_n$ (with $\mathfrak{O}_n^\flat$ the ideal of rank $0$ elements of $\mathfrak{t^\flat Ann}_n$). Since $\mathfrak{O}_n$ respectively $\mathfrak{O}^\flat_n$ does not satisfy any nontrivial identity, the finite basedness of $\mathfrak{tAnn}_n$ respectively $\mathfrak{t^\flat Ann}_n$ ($n$ even) is essentially caused by this ideal. It seems to be natural to look what happens if we pass to the Rees quotient over that ideal.

The list in the next theorem does not mention the deformed annular monoids $\mathfrak{Ann}^{\mathrm{d}_n}$ since these  are implicitly contained in \cite{auingeretal}: they are non-finitely based if and only if $n\ge 3$.
\begin{Thm}\label{thm:main annular}The following statements hold. In every case, the statement remains true if the monoid $\mathfrak{C}_n$ in question is replaced with its extended regular version $\overline{\mathfrak{C}_n}$.
\begin{enumerate}
\item $\mathfrak{tAnn}_n$ and $\mathfrak{t^\flat Ann}_n$ are non-finitely based for every odd $n\ge 3$.
\item $\mathfrak{tAnn}^\circ_n$  and $\mathfrak{t^\flat Ann}^\circ_n$ are non-finitely based for every even $n\ge 2$.
\item $\mathfrak{L}_n$ is non-finitely based for every $n\ge 2$.
\item $\mathfrak{L}^\flat_n$ is non-finitely based for $n\ge 3$.
\item $\mathfrak{aTL}^{\mathrm d}_n$ is non-finitely based for every $n\ge 2$.
\item $\mathfrak{aTL}_n$ is non-finitely based for every $n\ge 2$.
\item $\mathfrak{aTL}^{\mathrm e}_n$ is non-finitely based for every $n\ge 3$.
\end{enumerate}
\end{Thm}
\begin{proof} We are left with showing that the monoids in question do not satisfy any nontrivial identity $Z_k\bumpeq W$.

Recall the congruence ${\equiv}_5$ used in the definition of the deformed annular category $\mathfrak{Ann}^{\mathrm d}$ (\ref{deformed annular}). We can restrict ${\equiv}_5$ to $\mathfrak{L}$ and get $\mathfrak{TL}=\mathfrak{L}/{\equiv}_5$, the \emph{Temperley--Lieb category} that occurred in Subsection~\ref{sec:factored tAnn}. The morphism $\mathfrak{L}\twoheadrightarrow\mathfrak{TL}$ factors through $\mathfrak{L}^\flat$. Recall that the local monoids of $\mathfrak{TL}$  are the Kauffman monoids $\mathfrak{K}_n$. In \cite[Theorem 4.1]{auingeretal} it was proved that for any $n\ge 3$, the Kauffman monoids $\mathfrak{K}_n$ do not satisfy any nontrivial  semigroup identity $Z_k\bumpeq W$. Moreover, the same proof shows that this result is true for every monoid of the form $\mathfrak{K}_n/J_n$ where $n\ge 4$ is even and $J_n$ denotes the ideal of all rank $0$ elements. From this it follows that the monoids of items (1)--(4) do not satisfy any identity $Z_k\bumpeq W$, except for $\mathfrak{tAnn}^\circ_2$, $\mathfrak{t^\flat Ann}^\circ_2$ and $\mathfrak{L}_2$. The latter three cases are dealt with in Propositions \ref{tAnn2} (taking into account that $\mathfrak{tAnn}^\circ_2\twoheadrightarrow \mathfrak{t^\flat Ann}^\circ_2$) and \ref{L2} below.

Also, for $n\ge 3$, $\mathfrak{Ann}^{\mathrm d}_n$ does not satisfy any nontrivial identity $Z_k\bumpeq W$ (since the Kauffman monoid $\mathfrak{K}_n$ embeds). Since there is a canonical morphism $\mathfrak{aTL}^{\mathrm d}_n\twoheadrightarrow \mathfrak{Ann}^{\mathrm d}_n$ also the former does not satisfy any nontrivial identity $Z_k\bumpeq W$. Hence $\mathfrak{aTL}^{\mathrm d}_n$ appears in the list for every $n\ge 3$. For $n=2$ the claim follows from Proposition \ref{aTL_2} below since $\mathfrak{aTL}^{\mathrm d}_2\twoheadrightarrow \mathfrak{aTL}_2$. Altogether, this completely settles item (5).

For $n\ge 4$, the canonical morphisms $\mathfrak{aTL}_n\twoheadrightarrow \mathfrak{aTL}^{\mathrm e}_n\twoheadrightarrow \mathfrak{Ann}_n$ implies by \cite[Figure 12]{adv2} that the Brandt semigroup $\mathcal{B}_2$ divides all these monoids, and hence, Corollary~\ref{cor:a2b2} applies. The case $\mathfrak{aTL}_2$ is settled by Proposition \ref{aTL_2}. The case $\mathfrak{aTL}_3\cong \mathfrak{aTL}^{\mathrm e}_3$ is treated in Proposition \ref{aTL3}.   Altogether this completely solves items (6) and (7).
\end{proof}

We handle the remaining items. The first statement is not hard to verify. 
\begin{Lemma} Let the group $\mathbb{Z}$ act on $\mathbb{N}_0\times \mathbb{N}_0$ by letting even $n$ fix $\mathbb{N}_0\times \mathbb{N}_0$ and odd $n$ swap the entries of $(a,b)\in \mathbb{N}_0\times \mathbb{N}_0$. Then $\mathfrak{t^\flat Ann}^\circ_2$ is isomorphic to the monoid $(\mathbb{N}_0\times \mathbb{N}_0)\rtimes \mathbb{Z}$ with an extra zero adjoined.
\end{Lemma}

\begin{Prop}\label{tAnn2} The monoid $(\mathbb{N}_0\times \mathbb{N}_0)\rtimes \mathbb{Z}$ does not satisfy any nontrivial identity $Z_k\bumpeq W$.
\end{Prop}
\begin{proof} Suppose by contrast that $(\mathbb{N}_0\times \mathbb{N}_0)\rtimes \mathbb{Z}$ does satisfy such an identity. Then this identity is balanced. We recall that if an identity $u\bumpeq v$ is satisfied by a monoid then so is the identity $u'\bumpeq v'$ where $u'$ and $v'$ are the words obtained from $u$ respectively $v$ by deletion of some fixed letter in both words. Suppose that $(\mathbb{N}_0\times \mathbb{N}_0)\rtimes \mathbb{Z}$ satisfies the nontrivial balanced identity $Z_k\bumpeq W$  and assume that $k$ is minimal with that property. Deletion of the letter $x_1$ in both words leads to the identity $Z_k'\bumpeq W'$ which is also satisfied by $(\mathbb{N}_0\times \mathbb{N}_0)\rtimes \mathbb{Z}$. But $Z_k'$ is, up to renaming of the letters just the Zimin word $Z_{k-1}$ from which (by the minimal assumption on $k$) it follows that $Z_k'$ and $W'$ are identical as words. Since $W$ is obtained from $W'$ by reinserting $2^{k-1}$ copies of the letter $x_1$ and since $W$ is not identical to $Z_k$ it follows that $W$ somewhere contains the factor $x_1^2$.

Now consider the substitution $X\to (\mathbb{N}_0\times \mathbb{N}_0)\rtimes \mathbb{Z}$, 
\[x_i\mapsto\begin{cases}
   ((1,0),1)\mbox{ if }i=1\\
   ((0,1),1)\mbox{ if }i>1.
\end{cases}\]  
Under this substitution, the value of some word $U=y_1\cdots y_n$ is $((a,b),n)$ with $a+b=n$ and the numbers  $a$ and $b$ depend on which letters  are in which positions of $U$ (even or odd). It is readily checked that the value of $Z_k$ 
is just
\[((2^k-1,0),2^k-1).\]
Indeed, the contributions to the  first component of the value of $Z_k$ of the letters $x_i$ with $i>0$ is $(1,0)$ (the swapped version of $(0,1)$) since these letters are all in even positions, while the contribution of the letters $x_1$ is also $(1,0)$ (non-swapped) since all $x_1$ are in odd positions. Since $W$ contains a factor $x_1^2$, at least one occurrence of $x_1$ is in an even position. The contribution of such an occurrence to the first component of the value of $W$ is $(0,1)$ (the swapped version of $(1,0)$). If follows that the value of $W$ is $((a,b),2^k-1)$ for some $b>0$.
Altogether the values of $Z_k$ and $W$ do not agree and the claim is proved.
\end{proof}

\begin{Prop}\label{L2} The monoid $\mathfrak{L}_2$ does not satisfy any nontrivial identity $Z_k\bumpeq W$.
\end{Prop}
\begin{proof} Let $\mathcal{L}$ be the subsemigroup (ideal) of all rank $0$ elements of $\mathfrak{L}_2$. Every element $\alpha$ of $\mathcal{L}$ has one in-string and one out-string, hence one in-region, one out-region and one transversal region. Altogether, $\alpha$ is uniquely determined by the triple $(a,b,c)\in \mathbb{M}\times \mathbb{M}\times \mathbb{M}$ of labels of these regions. We assume that $a$ is the label of the in-region, $b$ the label of the transversal region and $c$ is the label of the out-region. Composition of two diagrams $\alpha\beta$ then retains the label of the in-region of $\alpha$, the label of the out-region of $\beta$. The label of the transversal region is the sum of the two labels of the transversal regions of $\alpha$ and $\beta$ plus the increment which is formed as the encircled sum of the labels of the out-region of $\alpha$ and the in-region of $\beta$. So, $\mathcal{L}$ is isomorphic to the $\mathbb{M}\times \mathbb{M}$ Rees matrix semigroup over $\mathbb{M}$ subject to the matrix $(p_{ij})$ where $p_{ij}=(i+j)$. In other words, the composition of two triples is defined by
\begin{equation*}
(a,b,c)(a',b',c')=(a,b+(c+a')+b',c').
\end{equation*}
Assume by contrast that $\mathcal{L}^1$ does satisfy some nontrivial identity $Z_k\bumpeq W$. Since $(0,0,0)$ is of infinite order, the variety generated by $\mathcal{L}$ is non-periodic. Hence the identity is balanced. As in the proof of Proposition \ref{aTL_2} we may assume that $k$ is minimal and hence $W$ somewhere contains the factor $x_1^2$. We substitute $((0),0,(0))$ for $x_1$ and $(0,0,0)$ for all other letters $x_i$. Computing the value of $Z_k$
\begin{equation*}
((0),0,(0))(0,0,0)((0),0,(0))(0,0,0)((0),0,(0))(0,0,0)\cdots
\end{equation*}
we see that all increments have the form $((0)+0)=((0))=((0)+0)$. Hence
the value of $Z_k$ is $((0),\mu,(0))$ where $\mu$ is a sum of elements of the form $((0))$, more precisely, $\mu=(2^k-2)\cdot((0))$. On the other hand, since $W$ contains the factor $x_1^2$ the value of $W$ is a product containing the factor $((0),0,(0))((0),0,(0))$ which yields the increment $((0)+(0))$. Consequently, the middle component of the value of $W$ contains the summand $((0)+(0))$ which implies that the values of $Z_k$ and $W$ do not coincide.
\end{proof}

\begin{Prop}\label{aTL_2} The monoid $\mathfrak{aTL}_2$ does not satisfy any nontrivial identity $Z_k\bumpeq W$.
\end{Prop}
\begin{proof} Take the partitions $\alpha,\beta\colon \mathbb{Z}\times[2]\leadsto \mathbb{Z}\times[2]$

\begin{center}
\begin{tikzpicture}[scale=1.1]
\draw (-0.75,0.5) node {$\alpha=$};

\draw[fill] (1,1) circle [radius=2pt];
\draw[above] (1,1) node{$(0,1)$};
\draw[fill] (2,1) circle [radius=2pt];
\draw[above] (2,1) node{$(0,2)$};
\draw(3,1) node{$\dots$};
\draw(0.25,0.5) node{$\dots$};
\draw(1,1) .. controls (1.25,0.75) and (1.75,0.75) .. (2,1);

\draw[fill] (3,0) circle [radius=2pt];
\draw[below] (3,0) node{$(1,1)$};
\draw[fill] (2,0) circle [radius=2pt];
\draw[below] (2,0) node{$(0,2)$};
\draw(1,0) node{$\dots$};
\draw(3.75,0.5) node{$\dots$};
\draw(2,0) .. controls (2.25,0.25) and (2.75,0.25) .. (3,0);

\draw(5.25,0.5) node {$\beta=$};

\draw[fill] (7,0) circle [radius=2pt];
\draw[below] (7,0) node{$(0,1)$};
\draw[fill] (8,0) circle [radius=2pt];
\draw[below] (8,0) node{$(0,2)$};
\draw(9,0) node{$\dots$};
\draw(6.25,0.5) node{$\dots$};
\draw(7,0) .. controls (7.25,0.25) and (7.75,0.25) .. (8,0);

\draw[fill] (9,1) circle [radius=2pt];
\draw[above] (9,1) node{$(1,1)$};
\draw[fill] (8,1) circle [radius=2pt];
\draw[above] (8,1) node{$(0,2)$};
\draw(7,1) node{$\dots$};
\draw(9.75,0.5) node{$\dots$};
\draw(8,1) .. controls (8.25,0.75) and (8.75,0.75) .. (9,1);
\end{tikzpicture}
\end{center}

{\noindent and consider them as arrows of $\mathfrak{aTL}$. Then $\alpha\beta$ and $\beta\alpha$ are idempotents (note that $0$-circles are dismissed in $\mathfrak{aTL}$) and every product containing a factor $\alpha^2$ or $\beta^2$ contains an $\omega$-circle. Factoring $\left<\alpha,\beta\right>$ by the ideal consisting of all elements containing at least one $\omega$-circle is isomorphic with the Brandt semigroup $\mathcal{B}_2$. Hence, Corollary~\ref{cor:a2b2} applies.}
\end{proof}

\begin{Prop}\label{aTL3} The monoid $\mathfrak{aTL}^{\mathrm e}_3=\mathfrak{aTL}_3$ does not satisfy any nontrivial identity $Z_k\bumpeq W$.
\end{Prop}
\begin{proof} 
Consider the following rank $1$ elements of the monoid $\mathfrak{aTL}^{\mathrm e}_3$ represented as partitions $\mathbb{Z}\times[3]\leadsto \mathbb{Z}\times[3]$:
\begin{center}
\begin{tikzpicture}[scale=0.6]
\draw (-2,1) node {$\alpha=$};
\draw(-0.5,1) node{$\dots$};
\draw(17.5,1) node{$\dots$};
\foreach \x in {1,4,7,10,13,16} \foreach \y in {0,2} \filldraw (\x,\y) circle (3pt);
\draw[above] (1,2) node{$(-1,3)$};
\draw[above] (4,2) node{$(0,1)$};
\draw[above] (7,2) node{$(0,2)$};
\draw[above] (10,2) node{$(0,3)$};
\draw[above] (13,2) node{$(1,1)$};
\draw[above] (16,2) node{$(1,2)$};
\draw[below] (1,0) node{$(-1,3)$};
\draw[below] (4,0) node{$(0,1)$};
\draw[below] (7,0) node{$(0,2)$};
\draw[below] (10,0) node{$(0,3)$};
\draw[below] (13,0) node{$(1,1)$};
\draw[below] (16,0) node{$(1,2)$};
\draw(1,2)--(7,0);
\draw(4,2) .. controls (4.5,1) and (6.5,1) .. (7,2);
\draw(10,2)--(16,0);
\draw(13,2) .. controls (13.5,1) and (15.5,1) .. (16,2);
\draw(1,0) .. controls (1.5,1) and (3.5,1) .. (4,0);
\draw(10,0) .. controls (10.5,1) and (12.5,1) .. (13,0);
\end{tikzpicture}
\end{center}
and
\begin{center}
\begin{tikzpicture}[scale=0.6]
\draw (-2,1) node {$\beta=$};
\draw(-0.5,1) node{$\dots$};
\draw(17.5,1) node{$\dots$};
\foreach \x in {1,4,7,10,13,16} \foreach \y in {0,2} \filldraw (\x,\y) circle (3pt);
\draw[above] (1,2) node{$(-1,2)$};
\draw[above] (4,2) node{$(-1,3)$};
\draw[above] (7,2) node{$(0,1)$};
\draw[above] (10,2) node{$(0,2)$};
\draw[above] (13,2) node{$(0,3)$};
\draw[above] (16,2) node{$(1,1)$};
\draw[below] (1,0) node{$(-1,2)$};
\draw[below] (4,0) node{$(-1,2)$};
\draw[below] (7,0) node{$(0,1)$};
\draw[below] (10,0) node{$(0,2)$};
\draw[below] (13,0) node{$(0,3)$};
\draw[below] (16,0) node{$(1,1)$};
\draw(1,2)--(7,0);
\draw(4,2) .. controls (4.5,1) and (6.5,1) .. (7,2);
\draw(10,2)--(16,0);
\draw(13,2) .. controls (13.5,1) and (15.5,1) .. (16,2);
\draw(1,0) .. controls (1.5,1) and (3.5,1) .. (4,0);
\draw(10,0) .. controls (10.5,1) and (12.5,1) .. (13,0);
\end{tikzpicture}
\end{center}
(so $\alpha$ is just $\beta$ shifted one step right). It is readily checked that $\alpha$ and $\beta$ are idempotents.  The product
\begin{center}
\begin{tikzpicture}[scale=0.43]
\draw (-3,1) node { $\alpha\beta=$};
\draw(-1.25,1) node{$\dots$};
\draw(23.5,1) node{$\dots$};
\draw(4.5,-.75) node{$\dots$};
\draw(18.5,2.75) node{$\dots$};
\foreach \x in {1,4,7,10,13,16} \filldraw (\x,2) circle (3pt);
\foreach \x in {7,10,13,16,19,22} \filldraw (\x,0) circle (3pt);
\draw[above] (1,2) node{$(-1,3)$};
\draw[above] (4,2) node{$(0,1)$};
\draw[above] (7,2) node{$(0,2)$};
\draw[above] (10,2) node{$(0,3)$};
\draw[above] (13,2) node{$(1,1)$};
\draw[above] (16,2) node{$(1,2)$};
\draw[below] (7,0) node{$(0,2)$};
\draw[below] (10,0) node{$(0,3)$};
\draw[below] (13,0) node{$(1,1)$};
\draw[below] (16,0) node{$(1,2)$};
\draw[below] (19,0) node{$(1,3)$};
\draw[below] (22,0) node{$(2,1)$};
\draw(1,2)--(13,0);
\draw(4,2) .. controls (4.5,1.5) and (6.5,1.5) .. (7,2);
\draw(10,2)--(22,0);
\draw(13,2) .. controls (13.5,1.5) and (15.5,1.5) .. (16,2);
\draw(7,0) .. controls (7.5,0.5) and (9.5,0.5) .. (10,0);
\draw(16,0) .. controls (16.5,0.5) and (18.5,0.5) .. (19,0);
\end{tikzpicture}
\end{center}
is easily seen to have infinite order. Thus, the subsemigroup $\langle\alpha,\beta\rangle$ of $\mathfrak{aTL}^{\mathrm e}_3$ is infinite. 

It is known \cite[Th\'eor\`eme~3]{Benzaken&Mayr:1975} (and easy to verify) that the only (up to isomorphism) infinite semigroup generated by two idempotents is $\mathcal{T}_2:=\langle e,f \mid e^2=e,\ f^2=f\rangle$, the free product of two trivial semigroups. Hence, the monoid $\mathfrak{aTL}^{\mathrm e}_3$ has a submonoid isomorphic to $\mathcal{T}_2^1$. Therefore, our claim follows from the proof of \cite[Theorem 2]{ShnVol:2017} where it was shown that the monoid $\mathcal{T}_2^1$ does not satisfy any nontrivial identity $Z_k\bumpeq W$.
\end{proof}

For $n=0$ or $n=1$ all monoids in items (3--6) of Theorem~\ref{thm:main annular} are finitely based; they are actually commutative. Moreover, $\mathfrak{tAnn}_0\cong \mathbb{O}\mathop{*}\mathbb{M}$ and $\mathfrak{t^\flat Ann}_0^\circ\cong \mathbb{O}\mathop{*}\mathbb{P}$ both of which do not satisfy any nontrivial identity. Altogether we have settled all cases with the exception of the monoids $\mathfrak{L}^\flat_2$ and $\mathfrak{aTL}^{\mathrm{e}}_2$ about which we state:

\begin{Thm} \label{Thm:aTLe_2} The monoids $\mathfrak{L}^\flat_2$ and $\mathfrak{aTL}_2^{\mathrm{e}}$ are finitely based. The identities \eqref{finite1} and \eqref{finite2} form a basis of the monoid identities of $\mathfrak{L}^\flat_2$;  the identities \eqref{finite1} and \eqref{finite2} together with the identity $x^3yx\bumpeq xyx^3$ form a basis of the monoid identities of $\mathfrak{aTL}_2^{\mathrm{e}}$. In both cases, a basis of semigroup identities is  obtained by deleting every possible choice of a subset of the set of letters $\{t_1,t_2,t_3,t_4\}$ in the identities \eqref{finite1} and \eqref{finite2}.
\end{Thm}
The proof of this result is deferred to the next subsection.

\subsection{Three nontrivial cases of finite basedness}\label{finite basedness} The purpose of the present subsection is to prove Theorems \ref{Thm:2Cob_1} and \ref{Thm:aTLe_2}. We start with two auxiliary constructions.

Let $\mathcal{S}$ be either $\langle\mathbb{N}_0;+,0\rangle$ or $\langle\mathbb{Z};+,0\rangle$. We consider the Rees matrix semigroup $\mathcal{M}(\mathcal{S},\mathcal{S},\mathcal{S};P)$ with $\mathcal{S}\times \mathcal{S}$-matrix $P$ defined by $p_{s,t}:=s+t+1$. Let the Cartesian power $\mathcal{S}^3$ act on $\mathcal{M}(\mathcal{S},\mathcal{S},\mathcal{S};P)$ by
\begin{equation}\label{eq:action of S^3}
	\begin{split}
[p,q,r]\cdot(a,b,c):=(q+a,p+b+r,c)\\ (a,b,c)\cdot[p,q,r]:=(a,p+b+r,c+q)
\end{split}
\end{equation}
where the triples in $\mathcal{M}(\mathcal{S},\mathcal{S},\mathcal{S};P)$ are denoted $(a,b,c)$ while those in $\mathcal{S}^3$ are denoted $[p,q,r]$. Define a new monoid $\mathcal{F}(\mathcal{S})$ by setting 
\[
\mathcal{F}(\mathcal{S}):=\mathcal{S}^3\sqcup \mathcal{M}(\mathcal{S},\mathcal{S},\mathcal{S};P),
\] 
retaining the respective operations in the `components' $\mathcal{S}^3$ and $\mathcal{M}(\mathcal{S},\mathcal{S},\mathcal{S};P)$ and letting the operation on `mixed pairs' be defined by~\eqref{eq:action of S^3}.

By the \emph{weight} of a triple $\tau\in\mathcal{F}(\mathcal{S})$ we mean the sum in $\mathcal{S}$ of the entries of $\tau$, regardless of whether $\tau$ comes from $\mathcal{M}(\mathcal{S},\mathcal{S},\mathcal{S};P)$ or from $\mathcal{S}^3$. We need the following observation, which can be derived by a straightforward induction from the rules \eqref{eq:action of S^3} and the multiplication rule of the Rees matrix semigroup $\mathcal{M}(\mathcal{S},\mathcal{S},\mathcal{S};P)$.

\begin{Lemma}\label{lem:fact}
For any $\tau_1,\dots,\tau_n\in\mathcal{F}(\mathcal{S})$ and every $(a,b,c),(f,g,h)\in\mathcal{M}(\mathcal{S},\mathcal{S},\mathcal{S};P)$,
\begin{equation}\label{eq:product}
(a,b,c)\tau_1\cdots \tau_n(f,g,h)=(a,b+c+T+N+1+f+g,h),
\end{equation}
where $N=\sharp\{i\mid \tau_i\in \mathcal{M}(\mathcal{S},\mathcal{S},\mathcal{S};P)\}$ and $T$ is the sum of the weights of the triples $\tau_1,\dots,\tau_n$.
\end{Lemma}

A kind of simplification of the monoid $\mathcal{F}(\mathcal{S})$ can be defined as follows: replace $\mathcal{M}(\mathcal{S},\mathcal{S},\mathcal{S};P)$ by $\mathcal{S}\times \mathcal{S}$, endowed with rectangular band multiplication $(a,b)\cdot(c,d)=(a,d)$, and $\mathcal{S}^3$ by $\mathcal{S}$ (more concretely, by the middle component of $\mathcal{S}^3$). This leads to a monoid $\mathcal{E}(\mathcal{S})$ defined on $\mathcal{S}\sqcup \left.\mathcal{S}\times \mathcal{S}\right.$ where the operation on the components is retained and the operation of `mixed pairs' is defined by
\[s\cdot(a,b):=(s+a,b)\mbox{ and }(a,b)\cdot s:=(a,b+s).\]
Finally, the mapping $\mathcal{F}(\mathcal{S})\to \mathcal{E}(\mathcal{S})$ defined by
\[[p,q,r]\mapsto q,\ (a,b,c)\mapsto (a,c)\]
is a surjective monoid morphism.

\subsubsection{The cases $2\mathfrak{Cob}^\circ_1$, $\overline{2\mathfrak{Cob}^\circ_1}$, $\mathfrak{L}^\flat_2$ and $\overline{\mathfrak{L}^\flat_2}$.}\label{subsubsec:2Cob1}

\begin{Lemma}\label{Lem:model_2Cob^0_1}  $2\mathfrak{Cob}^\circ_1\cong \mathcal{E}(\mathbb{N}_0)$ and $\overline{2\mathfrak{Cob}^\circ_1}\cong \mathcal{E}(\mathbb{Z})$.
\end{Lemma}

\begin{proof} Look at the picture below. The rank $1$ elements of $2\mathfrak{Cob}^\circ_1$ are cylinders carrying numbers from $\mathbb{N}_0$, the rank $0$ elements are pairs of cups/caps (that is, surfaces with exactly one boundary component) both of which carry numbers from $\mathbb{N}_0$. The composition does not create any increments, the numbers are simply added. In case of $\overline{2\mathfrak{Cob}^\circ_1}$ the labels are from $\mathbb{Z}$, otherwise everything is equal to the non-negative case. In addition, components without boundary that may be formed by composition of two rank $0$ elements are neglected. From this it is clear that the mapping depicted below is an isomorphisms $2\mathfrak{Cob}^\circ_1\to \mathcal{E}(\mathbb{N}_0)$ respectively $\overline{2\mathfrak{Cob}^\circ_1}\to\mathcal{E}(\mathbb{Z})$:
\begin{center}
\begin{tikzpicture}[
	tqft,
	cobordism edge/.style={draw},
	view from=incoming,
	cobordism height=1.5cm,
	]
	\begin{scope}[every node/.style={rotate=90}]
		\pic[tqft/cylinder,  at={(-5,0)},
		every incoming boundary component/.style={draw, rotate=90},
		every outgoing lower boundary component/.style={draw,dotted},
		every outgoing upper boundary component/.style={draw},
		];
		\pic[tqft/cap,  at={(0,0)},
		every incoming boundary component/.style={draw, rotate=90},
		every outgoing lower boundary component/.style={draw,dotted},
		every outgoing upper boundary component/.style={draw},
		];
		\pic[tqft/cup,  at={(0,0)},
		every incoming boundary component/.style={draw, rotate=90},
		every outgoing lower boundary component/.style={draw,dotted},
		every outgoing upper boundary component/.style={draw,dotted},
		];
	\end{scope}
	\draw(-4.25,0)node{$s$};
	\draw(0,0)node{$a$};
	\draw(1.5,0)node{$b$};
	\draw(-2.6,0)node{$\longmapsto s$,};
	\draw(2.7,0)node{$\longmapsto(a,b)$.};
\end{tikzpicture}
\end{center}
\end{proof}

\begin{Lemma}\label{lem:model_L^flat2}
	$\mathfrak{L}^\flat_2\cong \mathcal{F}(\mathbb{N}_0)$ and $\overline{\mathfrak{L}^\flat_2}\cong \mathcal{F}(\mathbb{Z})$.
\end{Lemma}
\begin{proof}
The monoids $\mathfrak{L}^\flat_2$ and $\overline{\mathfrak{L}^\flat_2}$ consist of rank $2$ and rank $0$ elements. The rank $2$ elements have two transversal strings and three transversal regions carrying labels $p,q,r\in\mathbb{N}_0$ (respectively $\mathbb{Z}$) --- denote these elements by $p|q|r$. The rank $0$ elements have an in-string and an out-string, an in-region, a transversal region and an out-region; again these regions carry labels from $\mathbb{N}_0$ (respectively $\mathbb{Z}$) --- denote the elements with $a$ on the in-region, $b$ on the transversal region and $c$ on the out-region by $a)b(c$. In the diagrams below it is summarized how the various kinds of diagrams are composed with each other (to do the artwork for the left-right dual of the third diagram, i.e. the diagrammatic representation of the product $\left.a)b(c\right.\circ \left.p|q|r\right.$ is left to the reader). From these diagrams it is clear that the mapping
\[\left.p|q|r\right.\mapsto [p,q,r]\mbox{ and }\left.a)b(c\right.\mapsto (a,b,c)\]
is an isomorphism $\mathfrak{L}^\flat_2\to \mathcal{F}(\mathbb{N}_0)$ (respectively $\overline{\mathfrak{L}^\flat_2}\to \mathcal{F}(\mathbb{Z})$).

\begin{center}
\begin{tikzpicture}[scale=0.5]
	\draw[dotted](0,0)--(2,0)--(2,3)--(0,3)--(0,0);
	\draw(0,1)--(2,1);
	\draw(0,2)--(2,2);
	\draw(1,0)[above]node{$p$};
	\draw(1,1)[above]node{$q$};
	\draw(1,2)[above]node{$r$};
	\draw((2.5,1.5)node{$\circ$};
	\draw[dotted](3,0)--(5,0)--(5,3)--(3,3)--(3,0);
	\draw(3,1)--(5,1);
	\draw(3,2)--(5,2);
	\draw(4,0)[above]node{$u$};
	\draw(4,1)[above]node{$v$};
	\draw(4,2)[above]node{$w$};
	\draw(6,1.5)node{$=$};
	\draw[dotted,shift={(6.9,0)}](0,0)--(2.2,0)--(2.2,3)--(0,3)--(0,0);
	\draw(6.9,1)--(9.1,1);
	\draw(6.9,2)--(9.1,2);
	\draw(8,0)[above]node{$p+u$};
	\draw(8,1)[above]node{$q+v$};
	\draw(8,2)[above]node{$r+w$};
\end{tikzpicture}
\end{center}
\begin{center}
\begin{tikzpicture}[scale=0.5]
	\draw[dotted](0,0)--(3,0)--(3,3)--(0,3)--(0,0);
	\draw(0,1) .. controls (1,1) and (1,2)..(0,2);
	\draw(3,1)..controls (2,1) and (2,2)..(3,2);
	\draw(-0.15,1)[above right]node{$a$};
	\draw(1.5,1)[above]node{$b$};
	\draw((3.1,1)[above left]node{$c$};
	\draw((3.5,1.5)node{$\circ$};
	\draw[dotted](4,0)--(7,0)--(7,3)--(4,3)--(4,0);
	\draw(4,1) .. controls (5,1) and (5,2)..(4,2);
	\draw(7,1)..controls (6,1) and (6,2)..(7,2);
	\draw(3.75,0.95)[above right]node{$d$};
	\draw(5.5,1)[above]node{$e$};
	\draw((7.25,0.9)[above left]node{$f$};
	\draw(8,1.5)node{$=$};
	\draw[dotted](9,0)--(15,0)--(15,3)--(9,3)--(9,0);
	\draw(9,1)..controls(10,1) and (10,2)..(9,2)--(9,2);
	\draw(8.85,1)[above right]node{$a$};
	\draw(15,1) ..controls(14,1) and (14,2).. (15,2);
	\draw(15.25,0.9)[above left]node{$f$};
	\draw(12,1.5)node[rotate=25]{$b+c+1+d+e$};
	\end{tikzpicture}
\end{center}
\begin{center}
\begin{tikzpicture}[scale=0.5]
	\draw[dotted](0,0)--(3,0)--(3,3)--(0,3)--(0,0);
	\draw(0,1)--(3,1);
	\draw(0,2)--(3,2);
	\draw(1.5,0)[above]node{$p$};
	\draw(1.5,1)[above]node{$q$};
	\draw(1.5,2)[above]node{$r$};
	\draw((3.5,1.5)node{$\circ$};
	\draw[dotted](4,0)--(7,0)--(7,3)--(4,3)--(4,0);
	\draw(4,1) .. controls (5,1) and (5,2)..(4,2);
	\draw(7,1)..controls (6,1) and (6,2)..(7,2);
	\draw(3.85,1)[above right]node{$a$};
	\draw(5.5,1)[above]node{$b$};
	\draw((7.1,1)[above left]node{$c$};
	\draw(8,1.5)node{$=$};
	\draw[dotted](9,0)--(14,0)--(14,3)--(9,3)--(9,0);
	\draw(9,1)--(10.5,1)..controls(11.5,1) and (11.5,2)..(10.5,2)--(9,2);
	\draw(8.85,0.9)[above right]node{$q+a$};
	\draw(14,1) ..controls(13,1) and (13,2).. (14,2);
	\draw(14.1,1)[above left]node{$c$};
	\draw(12.25,1.5)node[rotate=45]{$p+b+r$};
	\end{tikzpicture}
\end{center}
\end{proof}
Finally, the diagrams also show what Lemma~\ref{lem:fact} means geometrically when specialized to the monoid $\mathfrak{L}^\flat_2\cong\mathcal{F}(\mathbb{N}_0)$. Suppose we have a product
\begin{equation}\label{eq:product2}
	a)b(c\circ \tau_1\circ\cdots\circ\tau_n\circ d)e(f
\end{equation}
with $\tau_i\in \mathfrak{L}^\flat_2$ (or $\in \overline{\mathfrak{L}^\flat_2}$). The value of this product is obviously equal to $a)T(f$ where $T$ is the sum of all entries of all factors, except the first (which is $a$) and the last (which is $f$),  plus the number of all circles that are formed in the process of composition of the product~\eqref{eq:product2}. This number is equal to $N+1$ where $N$ is the number of elements of rank $0$ --- i.e.~of elements of the form $p)q(r$ --- among the factors $\tau_1,\dots, \tau_n$: indeed, every such element contributes a left half and a right half of a circle, the first factor $a)b(c$ contributes a further left half, while the last factor $d)e(f$ contributes a further right half. Altogether this gives $N+1$ circles since the factors of the form $p|q|r$ do not contribute to the number of circles formed in the product~\eqref{eq:product2}.

Hence, in order to prove Theorem \ref{Thm:2Cob_1} and the first statement of Theorem~\ref{Thm:aTLe_2} we need to find identity bases for the monoids $\mathcal{E}(\mathbb{N}_0)$ and $\mathcal{E}(\mathbb{Z})$ as well as for the monoids $\mathcal{F}(\mathbb{N}_0)$ and $\mathcal{F}(\mathbb{Z})$. In order to do so, we first shall characterize the identities which hold in these monoids and shall see that they are, in fact, the same, that is, the four monoids in question are equationally equivalent.

We need some combinatorial concepts. Let $w$ be a word and $x$ be a letter; we say that $x$ \emph{occurs} in $w$ if $w$ can be factorized as $w=uxv$ for some words $u,v$ (which may be empty). In this situation, the words $u$ and $v$ are referred to as the \emph{left context} and, respectively, the \emph{right context} of the occurrence of $x$. Clearly, it may happen that $x$ has several occurrences in $w$; we order these occurrences according to the lengths of their left contexts so that the \emph{leftmost occurrence} is the one with the shortest left context, and the \emph{rightmost occurrence} is the one with the shortest right context.
We denote by $\occ_x(w)$ the number of occurrences of $x$ in $w$. Recall that an identity $w\bumpeq w'$ is said to be balanced if $\occ_x(w)=\occ_x(w')$ for every letter $x$. It is well known and easy to verify that an identity is balanced if and only if it holds in the monoid $\mathbb{N}_0$. As the four monoids that we consider contain $\mathbb{N}_0$ as a submonoid, every identity holding in these monoids is balanced. In what follows, the latter observation is used without reference.

For a letter $x$ and a word $w$, we denote by $L(w,x)$ the longest prefix of $w$ in which $x$ does not occur. (If $x$ occurs in $w$, this prefix is nothing but the left context of the leftmost occurrence of $x$ in $w$; if $x$ does not occur in $w$, then $L(w,x)=w$.) Dually, $R(x,w)$ stands for the longest suffix of $w$ in which $x$ does not occur. Given an identity $w\bumpeq w'$, its \emph{left} (\emph{right}) \emph{section} at a letter $x$ is the identity $L(w,x)\bumpeq L(w',x)$ (respectively, $R(x,w)\bumpeq R(x,w')$).

\begin{Lemma}
\label{lem:characterization}
For an identity $w\bumpeq w'$, the following are equivalent:

\emph{(i)} All of the left and right sections of $w\bumpeq w'$ are balanced.

\emph{(ii)} The monoid $\mathcal{F}(\mathbb{Z})$ satisfies $w\bumpeq w'$.

\emph{(iii)} The monoid $\mathcal{E}(\mathbb{Z})$ satisfies $w\bumpeq w'$.

\emph{(iv)} The monoid $\mathcal{F}(\mathbb{N}_0)$ satisfies $w\bumpeq w'$.

\emph{(v)} The monoid $\mathcal{E}(\mathbb{N}_0)$ satisfies $w\bumpeq w'$.
\end{Lemma}

\begin{proof}
(i) $\Rightarrow$ (ii). Denote by $X$ the set of all letters involved in either $w$ or $w'$ and consider any substitution $\varphi\colon X\to\mathcal{F}(\mathbb{Z})$. Let 
\[
Y:=\{y\in X\mid y\varphi\in \mathbb{Z}^3\}.
\]
If $Y=X$ then the range of the substitution $\varphi$ is contained in  $\mathbb{Z}^3$, and $w\varphi=w'\varphi$ since the identity $w\bumpeq w'$ is balanced as the left section of itself at any letter beyond $X$ and the commutative monoid $\mathbb{Z}^3$ satisfies every balanced identity.

Suppose that $Y\subsetneq X$, and let $x,z\in X\setminus Y$ be such that $x$ occurs in $w$ first and $z$ occurs in $w$ last amongst all letters from $X\setminus Y$ (notice that $x$ and $z$ may coincide). Let $x\varphi:=(a,b,c)$ and $z\varphi:=(f,g,h)$. All letters that precede the leftmost occurrence of $x$ in $w$ as well as all letters that follow the rightmost occurrence of $z$ in $w$ belong to $Y$. Since the identities $L(w,x)\bumpeq L(w',x)$ and $R(z,w)\bumpeq R(z,w')$ are balanced, the same conclusion holds for all letters that either precede the leftmost occurrence of $x$ in $w'$ or follow the rightmost occurrence of $z$ in $w'$. If $L(w,x)$ and $L(w',x)$ are non-empty then  $L(w,x)\varphi=L(w',x)\varphi=[p,q,r]$ for some $[p,q,r]\in\mathbb{Z}^3$; similarly, if $R(z,w)$ and $R(z,w')$ are non-empty then $R(z,w)\varphi=R(z,w')\varphi=[s,t,u]$ for some $[s,t,u]\in\mathbb{Z}^3$. Let $v=x_1\cdots x_n$, $x_i\in X$, be the (possibly empty) factor of $w$ occurring between the leftmost occurrence of $x$ and the rightmost occurrence of $z$, and $v'$ is the analogous factor of $w'$. By \eqref{eq:product}, we have
\[
(xvy)\varphi=(a,b+c+T+N+1+f+g,h),
\]
where $T$ is the sum of the weights of the triples $x_1\varphi,\dots,x_n\varphi$ and $N=\sharp\{i\mid x_i\in X\setminus Y\}$. Here, the right hand side does not depend of the order in which the letters $x_1,\dots,x_n$ occur in $v$;
it depends only on the multiset $\{\!\{\occ_y(v)\mid y\in\{x_1,\dots,x_n\}\!\}$ of their multiplicities. Since each letter $y\in X$ occurs the same number of times in $w$ and $w'$, in $L(w,x)$ and $L(w',x)$, and in $R(z,w)$ and $R(z,w')$, it also occurs the same  number of times in $v$ and $v'$. Therefore, we conclude that $(xv'y)\varphi=(a,b+c+T+N+1+f+g,h)$ as well. Denoting the triple on the right hand side by $\tau$, we obtain
\[
w\varphi=w'\varphi=\begin{cases}
[p,q,r]\tau[s,t,u] &\text{if $L(w,x)$ and $R(z,w)$ are non-empty},\\
[p,q,r]\tau &\text{if $L(w,x)$ is not empty, but $R(z,w)$ is},\\
\tau[s,t,u] &\text{if $L(w,x)$ is empty, but $R(z,w)$ is not},\\
\tau &\text{if $L(w,x)$ and $R(z,w')$ are empty}.
\end{cases}
\]
Thus, the identity $w\bumpeq w'$  holds in the monoid $\mathcal{F}(\mathbb{Z})$.

\smallskip 

(ii) $\Rightarrow$ (iii) and (ii) $\Rightarrow$ (iv) are clear since $\mathcal{F}(\mathbb{Z})$ has $\mathcal{E}(\mathbb{Z})$ as a morphic image and $\mathcal{F}(\mathbb{N}_0)$ as a submonoid.

\smallskip

(iii) $\Rightarrow$ (v) and (iv) $\Rightarrow$ (v) follow since $\mathcal{E}(\mathbb{N}_0)$ is a submonoid of $\mathcal{E}(\mathbb{Z})$ and a morphic image of $\mathcal{F}(\mathbb{N}_0)$. 

\smallskip

(v) $\Rightarrow$ (i). Suppose that $w\bumpeq w'$ holds in $\mathcal{E}(\mathbb{N}_0)$ and consider its left and right sections at some letter $x$. For an arbitrary letter $y\ne x$, let  $\ell:=\occ_y(L(w,x))$, $\ell':=\occ_y(L(w',x))$, $r:=\occ_y(R(x,w))$, $r':=\occ_y(R(x,w))$. Consider the substitution $\varphi$ that assigns the pair $(0,0)\in\mathbb{N}_0\times\mathbb{N}_0$ to $x$, the number 1 to $y$, and the number 0 to all other letters. Then 
\[
w\varphi=(\ell,r)\quad\text{and}\quad w'\varphi=(\ell',r').
\]
Since  $w\bumpeq w'$ holds in $\mathcal{E}(\mathbb{N}_0)$, we have   $w\varphi=w'\varphi$, whence $\ell=\ell'$ and  $r=r'$. 
Thus, the number of occurrences of an arbitrary letter in $L(w,x)$  equals its number of occurrences in $L(w',x)$, and similarly for  $R(x,w)$ and $R(x,w')$. Hence both  the left section $L(w,x)\bumpeq L(w',x)$ and the right section $R(x,w)\bumpeq R(x,w')$ are balanced.
\end{proof}

Let us now recall the two monoid identities of Theorem \ref{Thm:2Cob_1}. These are
\begin{gather}
\label{eq:b1} xt_1yt_2\cdot xy\cdot t_3yt_4x \bumpeq xt_1yt_2\cdot yx\cdot t_3yt_4x,\\
\label{eq:b2} xt_1yt_2\cdot xy\cdot t_3xt_4y \bumpeq xt_1yt_2\cdot yx\cdot t_3xt_4y.
\end{gather}
In the following, the letter $\mathcal{M}$ stands for any of the four monoids $$\mathcal{E}(\mathbb{N}_0)\mbox{ or }\mathcal{E}(\mathbb{Z})\mbox{ or }\mathcal{F}(\mathbb{N}_0)\mbox{ or }\mathcal{F}(\mathbb{Z}).$$

From Lemma~\ref{lem:characterization}, we readily see that the monoid $\mathcal{M}$ satisfies the identities \eqref{eq:b1} and \eqref{eq:b2}. It remains to show that these two identities imply every identity holding in $\mathcal{M}$. For this, we devise a suitable normal form for words so that the following two properties hold:
\begin{enumerate}
\item[A)] every word can be transformed into a word in normal form via a suitable sequence of applications of the identities \eqref{eq:b1} and \eqref{eq:b2};
\item[B)] if an identity $w\bumpeq w'$ holds in the monoid $\mathcal{M}$, then the normal forms of the words $w$ and $w'$ coincide.
\end{enumerate}

Indeed, if such a normal form does exist, then for each identity $w\bumpeq w'$ that holds in the monoid $\mathcal{M}$, we have the following deduction from \eqref{eq:b1},\eqref{eq:b2}:
\[
w\stackrel{\eqref{eq:b1}, \eqref{eq:b2}}{\mbox{$\to\cdots\to$}}\text{the normal form of $w$} = \text{the normal form of $w'$}\stackrel{\eqref{eq:b1},  \eqref{eq:b2}}{\mbox{$\to\cdots\to$}} w'.
\]

Working towards the definition of our normal form, we need a common name for the leftmost and the rightmost occurrences of a letter in a word. We use the expression `an extreme occurrence' as such a common name. An arbitrary word $w$ in the letters $x_1,\dots,x_k$ can be decomposed as
\begin{equation}
\label{eq:extreme}
w=z_0u_1z_1u_2\cdots u_nz_n,\quad k-1\le n<2k,
\end{equation}
where $z_0,z_1,\dots,z_n$ represent all extreme occurrences of the letters. This means that $\{z_0,z_1,\dots,z_n\}=\{x_1,\dots,x_k\}$ and $z_i=z_j=x\in\{x_1,\dots,x_k\} $ for $i<j$ if and only if $z_i$ and $z_j$ represent the leftmost and respectively the rightmost occurrences  of the letter $x$. (The inequalities $k-1\le n<2k$ hold due to the fact that each of the letters $x_1,\dots,x_k$ has at least one and at most two extreme occurrences.) We refer to the decomposition \eqref{eq:extreme} as to the \emph{extreme representation} of $w$. The word $e(w):=z_0z_1\cdots z_n$ is called the \emph{word of extreme occurrences} of $w$ and the sequence of factors $u_1,\dots,u_n$ (of which some may be empty) is called the \emph{interior block sequence} of $w$.

\begin{Exm}
\label{exm:extreme}
If $w=x^3yxytz^4xyz$, then the extreme representation of $w$ is
\[
w=x\cdot x^2\cdot y\cdot xy\cdot t\cdot 1\cdot z\cdot z^3\cdot x\cdot 1\cdot y\cdot 1\cdot z,
\]
with $e(w)=xytzxyz$ and the interior block sequence $x^2$, $xy$, 1, $z^3$, 1, 1, where 1 stands for the empty word.
\end{Exm}

\begin{Lemma}
\label{lem:key1}
If an identity $w\bumpeq w'$ holds in the monoid $\mathcal{M}$ then the words of extreme occurrences of $w$ and $w'$ coincide.
\end{Lemma}

\begin{proof}
Let $x$ be a letter occurring in the word $w$ and let $w=uxv$ be the factorization that corresponds to an extreme occurrence of $x$ in $w$. If this is the leftmost occurrence of $x$ in $w$, then $u=L(w,x)$. By Lemma~\ref{lem:characterization}, the identity $L(w,x)\bumpeq L(w',x)$ is balanced, in particular, these words have the same length: $|L(w,x)|=|L(w',x)|$. By the definition of the left section, the letter following the prefix $L(w',x)$ in $w'$ must be $x$, and this occurrence of $x$ must be its leftmost occurrence in $w'$. Hence, the leftmost occurrences of $x$ in $w$ and $w'$ happen in exactly the same position from the left, namely, in the position $|L(w,x)|+1$.

Dually, if the occurrence of $x$ in $w$  that corresponds to the factorization $w=uxv$ is the rightmost occurrence of $x$ in $w$, then $v=R(x,w)$. By Lemma~\ref{lem:characterization}, the identity $R(x,w)\bumpeq R(x,w')$ is balanced, in particular, $|R(x,w)|=|R(x,w')|$. By the definition of right section, the letter preceding the suffix $R(x,w')$ in $w'$ must be $x$, and this occurrence of $x$ must be its rightmost occurrence in $w'$. Hence, the  rightmost occurrences of $x$ in $w$ and $w'$ happen in exactly the same position from the right. Since the identity $w\bumpeq w'$ is balanced, $|w|=|w'|$, whence the rightmost occurrences of $x$ in $w$ and $w'$ happen in exactly the same position from the left, namely, in the position $|w|-|R(x,w)|$.

Now imagine that we build the words $e(w)$ and $e(w')$ as follows. We initialize two word variables $e$ and $e'$ with the empty word and scan the words $w$ and $w'$ in parallel from left to right, appending a letter to, respectively, $e$ or $e'$ whenever we encounter an extreme occurrence of this letter in, respectively, $w$ or $w'$. As shown in the two preceding paragraphs, each extreme occurrence appears in $w$ and $w'$ in exactly the same position. Thus, our process appends to $e$ and $e'$ the same letters in the same order. This means that the equality $e=e'$ maintains at each step of the described process. Clearly, at the end of the process we get $e=e(w)$ and $e'=e(w')$, whence $e(w)=e(w')$.
\end{proof}

In view of Lemma~\ref{lem:key1}, if an identity $w\bumpeq w'$ holds in the monoid $\mathcal{M}$ then the interior block sequences of $w$ and $w'$ have the same number of blocks.

\begin{Lemma}
\label{lem:key2}
Suppose that an identity $w\bumpeq w'$ holds in the monoid $\mathcal{M}$ and $u_1,\dots,u_n$ and $u'_1,\dots,u'_n$ are the interior block sequences of $w$ and $w'$, respectively. Then for each $i=1,\dots,n$, the identity $u_i\bumpeq u'_i$ is balanced.
\end{Lemma}

\begin{proof}
Let \eqref{eq:extreme} be the extreme representation of $w$. Then the extreme representation of $w'$ is $w'=z_0u'_1z_1u'_2\cdots u'_nz_n$. We formally define $u_0=u'_0:=1$ and proceed by induction. Thus, let $i>0$ and assume that all identities $u_s\bumpeq u'_s$ with $s<i$ are balanced. Then so is the identity
\begin{equation}
\label{eq:induction}
z_0u_1z_1u_2\cdots u_{i-1}z_{i-1}\bumpeq z_0u'_1z_1u'_2\cdots u'_{i-1}z_{i-1}.
\end{equation}
We aim to show that $u_i\bumpeq u'_i$ is balanced as well. Consider two cases.

\smallskip

\noindent\textbf{\emph{Case 1:}} $z_i$ is the leftmost occurrence of some letter $x$ in $w$. Then we have $L(w,x)=z_0u_1z_1u_2\cdots u_{i-1}z_{i-1}u_i$ and $L(w',x)=z_0u'_1z_1u'_2\cdots u'_{i-1}z_{i-1}u'_i$. Invoking Lemma~\ref{lem:characterization} we get that the identity $L(w,x)\bumpeq L(w',x)$ is balanced. Since~\eqref{eq:induction} is balanced too, we have that for each letter $y$:
\begin{multline*}
\occ_y(u_i)=\occ_y(L(w,x))-\occ_y\left(z_0u_1z_1u_2\cdots u_{i-1}z_{i-1}\right)=\\
\occ_y(L(w',x))-\occ_y\left(z_0u'_1z_1u'_2\cdots u'_{i-1}z_{i-1}\right)=\occ_y(u'_i),
\end{multline*}
whence the identity $u_i\bumpeq u'_i$ is balanced as well.

\smallskip

\noindent\textbf{\emph{Case 2:}} $z_i$ is the rightmost occurrence of some letter $x$ in $w$. Lemma~\ref{lem:characterization} ensures that the identity $R(x,w)\bumpeq R(x,w')$ is balanced and the identity  $w\bumpeq w'$  itself is balanced. Besides, the identity~\eqref{eq:induction} is balanced. Now for each letter $y$, we have
\begin{multline*}
\occ_y(u_i)=\occ_y(w)-occ_y(xR(x,w))-\occ_y\left(z_0u_1z_1u_2\cdots u_{i-1}z_{i-1}\right)=\\
\occ_y(w')-occ_y(xR(x,w'))-\occ_y\left(z_0u'_1z_1u'_2\cdots u'_{i-1}z_{i-1}\right)=\occ_y(u'_i),
\end{multline*}
which again means that the identity $u_i\bumpeq u'_i$ is balanced.
\end{proof}

We take our alphabet $X=\{x_1,x_2,\dots\}$ and order it according to $x_1<x_2<\cdots$. A word $u$ over $X$ is \emph{sorted} if
$u=x_1^{\occ_{x_1}(u)}x_2^{\occ_{x_2}(u)}\cdots$, where some exponents may be $0$ and $x^0$ is the empty word for every letter $x\in X$. In other terms, $u$ is sorted if it is the least (with respect to the lexicographic order) word among words $v$ such that the identity $u\bumpeq v$ is balanced. We say that a word $w$ is \emph{in normal form} if all interior blocks in the extreme representation of $w$ are sorted. For instance, if we assume that $t<x<y<z$, then the word $w=x^3yxytz^4xyz$ from Example~\ref{exm:extreme} is in normal form.

We are in a position to complete the proof of Theorem \ref{Thm:2Cob_1}. In fact, it follows from the next proposition.
\begin{Prop} The identities \eqref{eq:b1} and \eqref{eq:b2} form a basis for the monoid identities of $\mathcal{M}$.
\end{Prop}
\begin{proof}
 As discussed above, the proof amounts to verifying that our normal form enjoys Properties A and B stated after Lemma~\ref{lem:characterization}. Property B readily follows from Lemma~\ref{lem:key2}. Indeed, by the lemma, if an identity $w\bumpeq w'$ holds in $\mathcal{M}$ and $u_1,\dots,u_n$ and $u'_1,\dots,u'_n$ are the interior block sequences of $w$ and respectively $w'$, then for each $i=1,\dots,n$, the identity $u_i\bumpeq u'_i$ is balanced, and therefore, sorting $u_i$ and $u'_i$ leads to the same sorted word. In order to verify Property A, we have to show that the identities \eqref{eq:b1} and \eqref{eq:b2} allow one to sort interior blocks in the extreme representation of any word $w$. For this, it suffices to show that an application of either \eqref{eq:b1} or \eqref{eq:b2} allows one to swap letters in any factor $x_jx_i$ with $i<j$ of an interior block of $w$. By the definition of interior blocks, the occurrences of $x_i$ and $x_j$ in this factor are not extreme in $w$. We locate the extreme occurrences of $x_i$ and $x_j$ in $w$ and consider four cases, depending on the location of these occurrences relative to each other.

\smallskip

\noindent\textbf{\emph{Case 1:}}  The leftmost occurrence of $x_i$ precedes the leftmost occurrence of $x_j$; the rightmost occurrence of $x_i$ precedes the rightmost occurrence of $x_j$.

In this case $w=V_0x_iV_1x_jV_2x_jx_iV_3x_iV_4x_jV_5$ for some (possibly empty) factors $V_0,\dots,V_5$, and one can apply  \eqref{eq:b2} `from right to left' to the factor $x_iV_1x_jV_2x_jx_iV_3x_iV_4x_j$, substituting $x_i$ for $x$, $x_j$ for $y$, $V_k$ for $t_k$, $k=1,2,3,4$.

\smallskip

\noindent\textbf{\emph{Case 2:}}  The leftmost occurrence of $x_i$ precedes the leftmost occurrence of $x_j$; the rightmost occurrence of $x_i$ follows the rightmost occurrence of $x_j$.

Here $w=V_0x_iV_1x_jV_2x_jx_iV_3x_jV_4x_iV_5$, and one can apply  \eqref{eq:b1} `from right to left' to the factor $x_iV_1x_jV_2x_jx_iV_3x_jV_4x_i$, using the same substitution as in Case 1.

\smallskip

\noindent\textbf{\emph{Case 3:}}  The leftmost occurrence of $x_i$ follows the leftmost occurrence of $x_j$; the rightmost occurrence of $x_i$ precedes the rightmost occurrence of $x_j$.

Here $w=V_0x_jV_1x_iV_2x_jx_iV_3x_iV_4x_jV_5$, and one can apply  \eqref{eq:b1} `from left to right' to the factor $x_jV_1x_iV_2x_jx_iV_3x_iV_4x_j$, substituting $x_j$ for $x$, $x_i$ for $y$, $V_k$ for $t_k$, $k=1,2,3,4$.

\smallskip

\noindent\textbf{\emph{Case 4:}}  The leftmost occurrence of $x_i$ follows the leftmost occurrence of $x_j$; the rightmost occurrence of $x_i$ follows the rightmost occurrence of $x_j$.

Here $w=V_0x_jV_1x_iV_2x_jx_iV_3x_jV_4x_iV_5$, and one can apply  \eqref{eq:b2} `from left to right' to the factor $x_jV_1x_iV_2x_jx_iV_3x_jV_4x_i$, using the same substitution as in Case 3.
\end{proof}

\subsubsection{The case $\mathfrak{aTL}^{\mathrm{e}}_2$} Let the group $\mathbb{Z}$ act (on the right as well as on the left) on $[2]=\{1,2\}$ by letting even $n$ fix $[2]$ and odd $n$ swap the two elements of $[2]$. On the set $[2]\times [2]$ we let $\mathbb{Z}$  act by letting every element of $\mathbb{Z}$ act on the adjacent entry of the pair $(x,y)\in [2]\times [2]$. We define the monoid $\mathcal{E}_2(\mathbb{Z}):=\mathbb{Z}\sqcup\, [2]\times [2]$ by retaining the binary operation $+$ of $\mathbb{Z}$, letting $[2]\times [2]$ be endowed with rectangular band multiplication and defining the operation of `mixed pairs' by the action of $\mathbb{Z}$ on $[2]\times [2]$ just mentioned.

The following is readily checked.
\begin{Lemma}\label{lem: aTL_2^e} $\mathfrak{aTL}^{\mathrm{e}}_2\cong \mathcal{E}_2(\mathbb{Z})$.
\end{Lemma}
We start with a characterization of the identities which hold in $\mathcal{E}_2(\mathbb{Z})$. First of all, $\mathcal{E}_2(\mathbb{Z})$ is a quotient of the monoid $\mathcal{E}(\mathbb{Z})$ considered in Subsection~\ref{subsubsec:2Cob1}. Indeed, let us define a relation $\sim$ on $\mathcal{E}(\mathbb{Z})$ by setting
\begin{equation*}p\sim q\Longleftrightarrow
\begin{cases}p,q\in \mathbb{Z}\mbox{ and }p=q,\mbox{ or} \\ p=(a,b),\ q=(c,d)\in \mathbb{Z}\times\mathbb{Z}\mbox{ and }a\equiv c,\ b\equiv d\ (\bmod\  2).\end{cases}
\end{equation*}
Then $\sim$ is a congruence on $\mathcal{E}(\mathbb{Z})$ and the quotient $\mathcal{E}(\mathbb{Z})/{\sim}$ is isomorphic to $\mathcal{E}_2(\mathbb{Z})$. In particular, every identity which holds in $\mathcal{E}(\mathbb{Z})$ also holds in $\mathcal{E}_2(\mathbb{Z})$.

We need some further concepts. For a word $w\in X^+$ we denote by $\con(w)$ the \emph{content} of $w$, that is, the set of all letters of $X$ which occur in $w$. An identity $w\bumpeq w'$ is \emph{balanced $\bmod\ 2$} if (i) $\con(w)=\con(w')$ and (ii) for every letter $x\in X$, $\occ_x(w)\equiv \occ_x(w')\ (\bmod\ 2)$. We recall that every identity which holds in $\mathcal{E}_2(\mathbb{Z})$ must be balanced.
\begin{Lemma}\label{lem:characterization2} A balanced identity holds in $\mathcal{E}_2(\mathbb{Z})$ if an only if all of its left and right sections are balanced $\bmod\ 2$.
\end{Lemma}
\begin{proof} Suppose that $u\bumpeq u'$ holds in $\mathcal{E}_2(\mathbb{Z})$, but its left section at a letter $x$ is not balanced $\bmod\ 2$. Since $u\bumpeq u'$ itself is balanced we may assume that $x$ actually occurs in $u$ (and hence also in $u'$). Assume that $L(u,x)\bumpeq L(u',x)$ is not balanced $\bmod\ 2$. Hence, either there is a letter $y$ which occurs in exactly one of the words $L(u,x)$, $L(u',x)$, e.g., $y\in \con(L(u,x))\setminus \con(L(u',x))$, or $\occ_y(L(u,x))\not\equiv \occ_y(L(u',x))\ (\bmod\ 2)$. In the former case let $v$ and $v'$ be the words obtained from $u$ respectively $u'$ by deleting all letters except $x$ and $y$. The the identity $v\bumpeq v'$ still holds in $\mathcal{E}_2(\mathbb{Z})$. However, $v$ starts with $y$ while $v'$ starts with $x$ which contradicts the fact that $\mathcal{E}_2(\mathbb{Z})$ contains the rectangular band on $[2]\times [2]$ as a subsemigroup so that both sides of every identity holding in $\mathcal{E}_2(\mathbb{Z})$ must start and end with the same letter.

Now assume that  $\occ_y(L(u,x))\not\equiv \occ_y(L(u',x))\ (\bmod\ 2)$. Without loss of generality, we assume that $\occ_y(L(u,x))$ is even and $\occ_y(L(u',x))$ is odd. We substitute $(1,1)\in\mathbb{Z}\times \mathbb{Z}$ for $x$ and $1\in \mathbb{Z}$ for $y$. Then the value of $u$ is $(1,k)$ for some $k$ while the value of $u'$ is $(2,\ell)$ (for some $k,\ell\in [2]$). In any case, these values do not agree which contradicts the assumption that the identity $u\bumpeq u'$ holds in $\mathcal{E}_2(\mathbb{Z})$. The case of right sections follows by symmetry. The proof of the converse  is very similar to the second part of the proof Lemma \ref{lem:characterization}.
\end{proof}
The remainder of the proof of Theorem \ref{Thm:aTLe_2} now follows from arguments which are similar to the earlier arguments of Subsection~\ref{subsubsec:2Cob1} and which we sketch as follows. Suppose the word $w$ is in normal form according to the previous definition: $w=z_0u_1z_1\cdots u_nz_n$. The interior blocks $u_i$ are sorted, that is, $u_i=x_{i_1}^{n_1}x_{i_2}^{n_2}\cdots$ with $i_1<i_2<\cdots$. We say that $w$ is in \emph{canonical form} if in every interior block $u_i$ the exponents $n_j$ are at most $2$ except for those letters $x_j$ which do not occur in an interior block $u_\ell$ with $\ell<i$ (that is, $u_i$ is the leftmost interior block in which the letter in question occurs). It is clear that every word in normal form can be transformed into its canonical form by application of the identity $x^3yx\bumpeq xyx^3$. Finally, from Lemma \ref{lem:characterization2} it follows that for two words $u$ and $u'$ in canonical form, the identity $u\bumpeq u'$ holds in $\mathcal{E}_2(\mathbb{Z})$ if and only of $u$ and $u'$ are equal as words.

\begin{Rmk}\label{rem:clustering}
The monoid identities \eqref{eq:b1} and \eqref{eq:b2} have already occurred in the literature. In \cite{CHL2016}, it was proved that for every integral domain $R$ of characteristic~0, \eqref{eq:b1} and \eqref{eq:b2} constitute an identity basis for the monoid $UT_2(R)$ formed by all upper triangular $2\times2$-matrices over $R$ whose main diagonal entries are 0s and/or 1s. In \cite{CMR23}, the same identities were shown to yield an identity basis for each Baxter monoid of rank greater than one; the same result was independently obtained in \cite{HZ2023}. Baxter monoids are interesting combinatorial objects belonging to another widely studied family of infinite monoids, so-called \emph{plactic-like} monoids. We refer the reader to \cite{Giraudo2012} for the definition, combinatorial properties and applications of Baxter monoids; here we mention only that elements of these monoids are pairs of twin binary search trees and their multiplication is governed by a Robinson–Schensted-like algorithm. Thus, Baxter monoids appear to have nothing in common with the four monoids of Subsection~\ref{subsubsec:2Cob1}, nor with monoids of the form $UT_2(R)$.

The fact that the same identities provide identity bases in such seemingly different situations is yet another instance of a `clustering phenomenon' revealed by the studies of identities of ‘interesting' monoids: monoids arising in various areas of mathematics and differing substantially in both the nature of their elements and their multiplication rules tend to cluster with respect to their identities. See \cite[Section 3.2]{CHKLV20} for a discussion and further examples of this phenomenon. It is not clear at present whether this clustering is merely a kind of pigeonhole effect, caused by `interesting' monoids being more abundant than `interesting' identities, or whether it reflects some deep but as-yet-uncovered connections between objects that appear to be unrelated at the surface level.
\end{Rmk}

\subsection{Adding involutions to the algebraic signature} We met already inverse unary operations on the categories $\mathfrak{P}$, $\mathfrak{Ann}$ and on the regular extensions $\overline{\mathfrak{C}}$ for  all categories $\mathfrak{C}$ in question. The bad news is that our methods are not applicable to semigroups equipped with inverse unary operations $^*$, simply because in this case the nontrivial identity $Z_1\bumpeq x_1x_1^*x_1$ holds.

However, the good news is that the involution $^*$ on $\mathfrak{P}$ and $\mathfrak{Ann}$ gives rise to involutions $^\sigma$ and $^\rho$ on almost all categories in question  which our methods can be applied to. Recall that the involution $^*$ on $\mathfrak{P}$ is defined as follows: for $\alpha\colon[m]\leadsto [n]$ we have that $\alpha^*\colon[n]\leadsto [m]$ is the `same' partition, but with incoming and outgoing vertices swapped. In the pictorial representation in Subsection \ref{partition category} this can be viewed as `reflection along the horizontal axis'. For the same $\alpha\colon[m]\leadsto [n]$ we let $\alpha^\rho\colon[n]\leadsto [m]$ be $\alpha^*$ \emph{reflected along the vertical axis}, that is, after reflecting along the horizontal axis we apply the permutation $1\leftrightarrow n,\ 2\leftrightarrow n-1,\dots$ to the incoming vertices and $1\leftrightarrow m,\ 2\leftrightarrow m-1,\dots$ to the outgoing vertices. In the pictorial representation of Subsection \ref{partition category} this means that the whole diagram is rotated by $180^\circ$ and afterwards the vertices are re-numbered in increasing order from left to right. We call $\alpha\mapsto \alpha^*$ the \emph{reflection} and $\alpha\mapsto \alpha^\rho$ the \emph{rotation} both of which are involutions on the partition category $\mathfrak{P}$ which leave invariant the category $\mathfrak{Ann}$.

\subsubsection{Partition and cobordism categories enriched with involutions} We can extend the reflection and the rotation to $2\mathfrak{Cob}$, $2\mathfrak{Cob}^\circ$ and $\mathfrak{P}^{\mathrm{d}}$ as well as to their regular extensions  $\overline{2\mathfrak{Cob}}$, $\overline{2\mathfrak{Cob}^\circ}$ and $\overline{\mathfrak{P}^{\mathrm{d}}}$ by setting (we use the symbol $^\sigma$ to denote the reflection since the symbol $^*$ has been used for the inverse unary operation on the regular extensions of the respective categories):
\begin{align*}
(\alpha,s)^\sigma:=&(\alpha^*,s)\mbox{ for }(\alpha,s)\in \overline{\mathfrak{P}^{\mathrm d}}\\
(\alpha,g,s)^\sigma:=&(\alpha^*,g^\sigma,s)\mbox{ for }(\alpha,g,s)\in \overline{2\mathfrak{Cob}}\\
(\alpha,g)^\sigma:=&(\alpha^*,g^\sigma)\mbox{ for }(\alpha,g)\in \overline{2\mathfrak{Cob}^\circ}
\end{align*}
where $g^\sigma$ means that $g$ is applied to the reflected blocks ($g^\sigma$ does not coincide with $g^*$ in (\ref{def:gstar})), and the `non-overlined subcategories' are always invariant under this operation. For the reflection the definition is analogous:
\begin{align*}
(\alpha,s)^\rho:=&(\alpha^\rho,s)\mbox{ for }(\alpha,s)\in \overline{\mathfrak{P}^{\mathrm d}}\\
(\alpha,g,s)^\rho:=&(\alpha^\rho,g^\rho,s)\mbox{ for }(\alpha,g,s)\in \overline{2\mathfrak{Cob}}\\
(\alpha,g)^\rho:=&(\alpha^\rho,g^\rho)\mbox{ for }(\alpha,g)\in \overline{2\mathfrak{Cob}^\circ}
\end{align*}
and likewise, $g^\rho$ means that $g$ is applied to the rotated blocks and the operation leaves invariant the non-overlined subcategories. Reflection as well as rotation have an obvious natural geometric interpretation in the categories $2\mathfrak{Cob}$ and $2\mathfrak{Cob}^\circ$. Both operations are respected by the canonical morphisms (\ref{category morphisms}) and (\ref{category morphisms 2}).

In order to apply Theorem \ref{thm:sufficient involutory version} it is sufficient to check which of the monoids do not satisfy any nontrivial involutory semigroup identity $Z_k\bumpeq W$ (which means that $W\in \mathcal{I}(X)$, that is, $W$ may contain `starred' letters); the other requirement follows from the `plain semigroup' case. The main result in this context is
\begin{Thm} The following statements hold. In every case, the statement remains true if the monoid $\mathfrak{C}_n$ in question is replaced with its extended regular version $\overline{\mathfrak{C}_n}$.
\begin{enumerate}
\item The involutory monoid $2\mathfrak{Cob}_n$ is non-finitely based for both involutions $^\sigma$ and $^\rho$ for every $n\ge 1$.
\item The involutory monoid $2\mathfrak{Cob}^\circ_n$ is non-finitely based for both involutions $^\sigma$ and $^\rho$ for every $n\ge 2$.
\item The involutory monoid $\mathfrak{P}^{\mathrm{d}_n}$ is non-finitely based for both involutions $^\sigma$ and $^\rho$ for every $n\ge 2$.
\end{enumerate}
\end{Thm}
\begin{proof} It suffices to show that in any mentioned case no involutory identity $Z_k\bumpeq W$ is satisfied. In case $2\mathfrak{Cob}_1$ we see that the involutions $^\rho$ and $^\sigma$ coincide. Consider the monoid $\mathcal{A}_2^1$ and endow it with the involution $(i,j)^*=(j,i),\ 0^*=0,\ 1^*=1$; on the semigroup $\mathcal{S}$ defined in the proof of Proposition \ref{2Cob1}, the involution ${^\sigma}$ (and likewise $^\rho$) is defined by $(i,s,j)^\sigma=(j,s,i)$ from which it follows that the morphism (\ref{A21 dividing 2Cob1}) respects the involution. It now follows from \cite[Corollaries 2.7, 2.8]{adv} that $2\mathfrak{Cob}_1$ does not satisfy any nontrivial involutory identity $Z_k\bumpeq W$ for both involutions $^\sigma$ and $^\rho$. For every $n\ge 1$ there is an embedding $2\mathfrak{Cob}_1\hookrightarrow 2\mathfrak{Cob}_n$ which respects the involution $^\sigma$ so that the claim is true for the reflection $^\sigma$. Finally, for every $n\ge 2$ the canonical morphism $2\mathfrak{Cob}_n\twoheadrightarrow \mathfrak{P}_n$ respects the involution $^\rho$; it follows from  \cite[Theorem 2.13]{adv2} that $2\mathfrak{Cob}_n$ endowed with rotation does not satisfy any nontrivial identity $Z_k\bumpeq W$. The latter argument also proves the claim for $2\mathfrak{Cob}^\circ_n$ and $\mathfrak{P}^{\mathrm d}_n$ for $n\ge 2$ and for the involution $^\rho$.

So we are left with the latter two cases, but with respect to $^\sigma$ instead of $^\rho$. Consider the following rank $1$ partitions $[2]\leadsto [2]$ (incoming vertices are drawn on the left, outgoing on the right):
\begin{center}
\begin{tikzpicture}[scale=1.3]
\draw (-3+0,0) node {$a_{00}=$};
\draw (-3+0.75,0.5) node {$1$};
\draw (-3+0.75,-0.5) node{$2$};
\draw(-3+2,0.5) node {$1$};
\draw(-3+2,-0.5) node {$2$};

\draw (-3+0,-2) node {$a_{10}=$};
\draw (-3+0.75,-1.5) node {$1$};
\draw (-3+0.75,-2.5) node{$2$};
\draw(-3+2,-1.5) node {$1$};
\draw(-3+2,-2.5) node {$2$};

\draw (2.75,0) node {$=a_{01}$};
\draw (0.75,0.5) node {$1$};
\draw (0.75,-0.5) node{$2$};
\draw(2,0.5) node {$1$};
\draw(2,-0.5) node {$2$};

\draw (2.75,-2) node {$=a_{11}$.};
\draw (0.75,-1.5) node {$1$};
\draw (0.75,-2.5) node{$2$};
\draw(2,-1.5) node {$1$};
\draw(2,-2.5) node {$2$};

\draw[rounded corners=5pt](-2.25,0.7)--(-0.8,0.7)--(-0.8,0.3)--(-2.45,0.3)--
(-2.45,0.7)--(-2.25,0.7);

\draw[rounded corners=5pt](-2.25,-0.3)--(-2.05,-0.3)--(-2.05,-0.7)--(-2.45,-0.7)--(-2.45,-0.3)--(-2.25,-0.3);
\draw[rounded corners=5pt](-1,-0.3)--(-0.8,-0.3)--(-0.8,-0.7)--(-1.2,-0.7)--(-1.2,-0.3)--(-1,-0.3);
\draw[rounded corners=5pt](0.75,-0.3)--(0.55,-0.3)--(0.55,-0.7)--(0.95,-0.7)--(0.95,-0.3)--(0.75,-0.3);
\draw[rounded corners=5pt](-1,-2.3)--(-0.8,-2.3)--(-0.8,-2.7)--(-1.2,-2.7)--(-1.2,-2.3)--(-1,-2.3);

\draw[rounded corners=5pt]
(1,-1.3)--(2.2,-1.3)--(2.2,-1.7).. controls (1.9,-1.7) and (1.9,-2.3 ).. (2.2,-2.3)--(2.2,-2.7)  --(0.55,-2.7)--(0.55,-2.3).. controls (0.85,-2.3) and (0.85,-1.7) .. (0.55,-1.7)--(0.55,-1.3)-- (1,-1.3);

\draw[rounded corners=5pt](-2,-1.3)--(-0.8,-1.3)--(-0.8,-1.7).. controls (-1.4,-1.7) and (-2.05,-2.1) .. (-2.05,-2.7)--(-2.45,-2.7)--(-2.45,-1.3)--(-2,-1.3);


\draw[rounded corners=5pt] (1,0.7)--(2.2,0.7)--(2.2,-0.7)--(1.8,-0.7) .. controls (1.8,-0.1) and (1.15,0.3) .. (0.55,0.3)--(0.55,0.7)--(1,0.7);

\end{tikzpicture}
\end{center}
These elements give rise to elements $\alpha_{ij}:=(a_{ij},0)$ in $\mathfrak{P}^{\mathrm d}_2$ and $\beta_{ij}:=(a_{ij},\mathbf{0})$ in $2\mathfrak{Cob}^\circ_2$ (where in the latter case, $\mathbf{0}$ denotes the labeling function with value constant $0$). In both cases the sets $\{\alpha_{ij}\}$ and $\{\beta_{ij}\}$ are invariant under the reflection $^\sigma$. In the first case, $\alpha_{ij}$ is idempotent if and only if $(i,j)\neq (0,0)$ while in the second case, $\beta_{ij}$ is idempotent if and only if $(i,j)\neq (1,1)$. Let $S$ be the subsemigroup of $\mathfrak{P}^{\mathrm d}_2$ generated by $\{\alpha_{ij}\}$ factored by the ideal consisting of all elements with second component larger than $0$. Then $S$ is isomorphic to $\mathcal{A}_2$ as an involutory semigroup, hence $S^1$ is isomorphic to $\mathcal{A}_2^1$ as involutory monoid (with involution defined above). Likewise, let $T$ be the subsemigroup of $2\mathfrak{Cob}^\circ_2$ generated by $\{\beta_{ij}\}$ factored by the ideal containing a block with positive label (that is, a component with positive genus). Again, $T$ is isomorphic to $\mathcal{A}_2$ as an involutory semigroup, hence $T^1$ is isomorphic to $\mathcal{A}_2^1$ as an involutory monoid with the same involution. It follows that $\mathfrak{P}^{\mathrm d}_2$ as well as $2\mathfrak{Cob}^\circ_2$ do not satisfy any nontrivial involutory identity $Z_k\bumpeq W$ with respect to reflection. Since there are reflection preserving embeddings $\mathfrak{P}^{\mathrm d}_2\hookrightarrow \mathfrak{P}^{\mathrm d}_n$ and $2\mathfrak{Cob}^\circ_2\hookrightarrow 2\mathfrak{Cob}^\circ_n$ for all $n\ge 2$ the claim is proved.
\end{proof}

\subsubsection{Annular categories enriched with involutions} We have already defined the operation $\bm\alpha\mapsto \bm\alpha^*$ for $\bm\alpha\in \mathfrak{aTL}^{\mathrm{e}}$;  in the present context we call this operation reflection.  The definition can be extended to $\mathfrak{aTL}$ and $\mathfrak{aTL}^{\mathrm d}$ by setting
\begin{align*}
(\bm\alpha,k)^\sigma:=&(\bm\alpha^*,k)\mbox{ for }(\bm\alpha,k)\in \mathfrak{aTL},\\
(\bm\alpha,k,\ell)^\sigma:=&(\bm\alpha^*,k,\ell)\mbox{ for }(\bm\alpha,k,\ell)\in \mathfrak{aTL}^{\mathrm{d}}.
\end{align*}
Using the $R$-diagram representation of a t$A$-diagram we can define the reflection in the same manner on $\mathfrak{tAnn}$: reflect the diagram along the horizontal line $\mathbb{R}\times\{\frac{a+b}{2}\}$ to get $\bm\alpha^\sigma$ for $\bm\alpha\in \mathfrak{tAnn}$. The involution $\bm\alpha\mapsto \bm\alpha^\sigma$ leaves invariant the category $\mathfrak{L}$ and hence can be restricted to the latter.

In order to define the rotation for diagrams $\bm\alpha\in\mathfrak{aTL}^{\mathrm{e}}$ recall the embedding of $\left.\mathbb{Z}\times[m]\right.\sqcup\left.\mathbb{Z}\times[n]\right.$ into $R(a,b)$ via (\ref{embedding points in R(a,b)}) and (\ref{embedding points in R(a,b) II}). Then $\bm\alpha^\rho$ is defined to be the diagram obtained by rotating the strip $R(a,b)$ by $180^\circ$ around the point $(\frac{1}{2},\frac{a+b}{2})$ and relabel the boundary points accordingly. This provides an involution $\bm\alpha\mapsto \bm\alpha^\rho$ on the category $\mathfrak{aTL}^{\mathrm{e}}$. Exactly the same procedure leads to the involution $\bm\alpha\mapsto \bm\alpha^\rho$ on $\mathfrak{tAnn}$ which leaves invariant the subcategory $\mathfrak{L}$ (hence can be restricted to the latter). Finally, for the intermediary categories $\mathfrak{aTL}$ and $\mathfrak{aTL}^{\mathrm{d}}$ we proceed as for the reflection:
\begin{align*}
(\bm\alpha,k)^\rho:=&(\bm\alpha^\rho,k)\mbox{ for }(\bm\alpha,k)\in \mathfrak{aTL},\\
(\bm\alpha,k,\ell)^\rho:=&(\bm\alpha^\rho,k,\ell)\mbox{ for }(\bm\alpha,k,\ell)\in \mathfrak{aTL}^{\mathrm{d}}.
\end{align*}

Just as in the plain case at the beginning of Subsection~\ref{sec: tAnn and L}, we start with a `negative' observation analogous to Theorem~\ref{thm:main annular}.
\begin{Thm}
For every even $n$, the monoids $\mathfrak{tAnn}_n$ and $\overline{\mathfrak{tAnn}_n}$, as well as $\mathfrak{t^\flat Ann}$ and $\overline{\mathfrak{t^\flat Ann}_n}$ do not satisfy any nontrivial involutory identity with respect to reflection $\sigma$ as well as rotation $\rho$. Hence these involutory monoids are finitely based.
\end{Thm}
\begin{proof}
It has been shown at the beginning of Subsection~\ref{sec: tAnn and L} that the free monoid on two letters $x,y$ embeds into the ideal $\mathfrak{O}_n$ of all rank $0$ elements of any of the monoids mentioned in the statement of the theorem. It is readily checked that, restricted to the image of this embedding, both involutions $\sigma$ and $\rho$ coincide and correspond to the involution $w\mapsto \mathrm{r}(w)$ where $\mathrm{r}(w)$ denotes the \emph{reversed word} (the word $w$ read from right to left).  This leads to an embedding of the free involutory monoid on countable many generators $x_1,x_2,\dots$ into $\mathfrak{O}_n$ since that monoid embeds into the free monoid on two generators $x,y$ (equipped with word reversion) via
    \[x_i\mapsto xy(yx)^ixy\mbox{ for all }i.\]
\end{proof}
The main result of the present subsection may be formulated as follows.
\begin{Thm}\label{thm:involutory annular}The following statements hold. In every case, the statement remains true if the involutory monoid $\mathfrak{C}_n$ in question is replaced with its extended regular version $\overline{\mathfrak{C}_n}$.
\begin{enumerate}
\item $\mathfrak{tAnn}_n$ and $\mathfrak{t^\flat Ann}_n$ are non-finitely based for every odd $n\ge 3$ with respect to both involutions $^\sigma$ and $^\rho$.
\item $\mathfrak{tAnn}^\circ_n$ and $\mathfrak{t^\flat Ann}^\circ_n$ are non-finitely based for every even $n\ge 2$ with respect to both involutions $^\sigma$ and $^\rho$.
\item $\mathfrak{L}_n$ is non-finitely based for every $n\ge 2$ with respect to both involutions $^\sigma$ and $^\rho$.
\item $\mathfrak{L}^\flat_n$ is non-finitely based for every $n\ge 3$ with respect to both involutions $^\sigma$ and $^\rho$.
\item $\mathfrak{aTL}^{\mathrm d}_n$ is non-finitely based for every $n\ge 2$ with respect to both involutions $^\sigma$ and $^\rho$.
\item $\mathfrak{aTL}_n$ is non-finitely based for every even $n\ge 2$ with respect to both involutions $^\rho$ and $^\sigma$.
\item $\mathfrak{aTL}_n$ and $\mathfrak{aTL}^{\mathrm e}_n$ are non-finitely based for every $n\ge 3$ with respect to the involution $^\rho$.
\end{enumerate}
\end{Thm}
\begin{proof}
In order to apply Theorem \ref{thm:sufficient involutory version} it is sufficient to check that the monoids in question do not satisfy any nontrivial involutory semigroup identity $Z_k\bumpeq W$. Recall the congruence ${\equiv}_5$ used in the definition of the deformed annular category $\mathfrak{Ann}^{\mathrm d}$ (\ref{deformed annular}). As in the case without involutions, we restrict ${\equiv}_5$ to $\mathfrak{L}$ and get $\mathfrak{TL}=\mathfrak{L}/{\equiv}_5$ and the morphism $\mathfrak{L}\twoheadrightarrow \mathfrak{TL}$ factors through $\mathfrak{L}^\flat$. Here $\mathfrak{TL}$ is the Temperley--Lieb category mentioned in Subsection~\ref{sec:factored tAnn}; recall that its local monoids are the Kauffman monoids $\mathfrak{K}_n$. Note that both involutions are well defined on $\mathfrak{L}$ and  the congruence ${\equiv}_5$ respects both of them. In \cite[Theorems 4.1, 4.3]{auingeretal} it was proved that for any $n\ge 3$, the Kauffman monoids $\mathfrak{K}_n$ do not satisfy any nontrivial involutory semigroup identity $Z_k\bumpeq W$ for both $^\sigma$ and $^\rho$. Moreover, for even $n$ let $\mathfrak{K}_n^\circ:=\mathfrak{K}_n/J_n$ where $J_n$ is the ideal of $\mathfrak{K}_n$ consisting of all rank $0$ elements. Similar  arguments as those in the proofs of Theorems 4.1 and 4.3 of \cite{auingeretal} (but more easy ones) also show that $\mathfrak{K}^\circ_n$ does not satisfy any involutory identity $Z_k\bumpeq W$ for every even $n\ge 4$. Finally, recall that there is a canonical morphism $\mathfrak{aTL}^{\mathrm d}\twoheadrightarrow \mathfrak{Ann}^{\mathrm d}$ and an embedding $\mathfrak{TL}\hookrightarrow\mathfrak{Ann}^{\mathrm d}$ (both maps respect both involutions). Altogether these arguments show that all monoids in items (1)--(5), with the exception of $\mathfrak{tAnn}^\circ_2$, $\mathfrak{t^\flat Ann}_2^\circ$, $\mathfrak{L}_2$ and $\mathfrak{aTL}^{\mathrm d}_2$, do not satisfy any nontrivial involutory identity $Z_k\bumpeq W$ for both involutions $^\sigma$ and $^\rho$.

For the first two cases it is sufficient to treat the second one; for $\mathfrak{t^\flat Ann}^\circ_2$ we ignore the $0$ and consider the submonoid of rank $2$ elements in the representation of Proposition \ref{tAnn2} as the semidirect product $(\mathbb{N}_0\times \mathbb{N}_0)\rtimes \mathbb{Z}$. We note that the two involutions $^\sigma$ and $^\rho$ agree and in terms of the semidirect product representation this involution is given by
$$((a,b),k)\mapsto ((a,b),-k).$$
We have already seen that $(\mathbb{N}_0\times\mathbb{N}_0)\rtimes \mathbb{Z}$ does not satisfy any nontrivial semigroup identity $Z_k\bumpeq W$. Assume that $(\mathbb{N}_0\times\mathbb{N}_0)\rtimes \mathbb{Z}$ satisfies some nontrivial \textsl{involutory} identity $Z_k\bumpeq W$. Substituting $((1,0),0)$ for $x_i$ and $((0,0),0)$ for the other letters we see that every letter occurs the same number of times on both sides. Since $Z_k\bumpeq W$ is nontrivial, $W$ must contain a starred letter, say $x_i^*$. We substitute $((0,0),1)$ for $x_i$ and $((0,0),0)$ for all other letters. Then the value of every word is $((0,0),t)$ where $t$ denotes the number of non-starred occurrences of $x_i$ minus the number of starred occurrences of $x_i$. More concretely, the value of $Z_k$ is $((0,0),2^{k-i})$ but the value of $W$ is $((0,0),t)$ for some $t\le 2^{k-i}-2<2^{k-i}$ since among the $2^{k-i}$ occurrences of $x_i$ in $W$ there is at least one starred version.

For the monoid $\mathfrak{L}_2$ we consider the representation of the ideal $\mathcal{L}$ consisting of all rank $0$ elements of $\mathfrak{L}_2$ presented in Proposition \ref{L2}. We have already seen that the  monoid $\mathcal{L}^1$ does not satisfy any nontrivial semigroup identity $Z_k\bumpeq W$. Moreover, as in the previous case, both involutions coincide. In the triple representation of the elements of $\mathcal{L}$, the involution is given by the map $(a,b,c)\mapsto (c,b,a)$. Suppose that the nontrivial involutory identity $Z_k\bumpeq W$ is satisfied by $\mathcal{L}^1$. Fix some letter $x_i$; substituting $(0,0,0)$ for $x_i$ and $1$ for all other letters we get that the value of $Z_k$ is $$(0,0,0)^{2^{k-i}}=(0,(2^{k-i}-1)\cdot (0),0)$$ while the value of $W$ is $(0,0,0)^t=(0,(t-1)\cdot (0),0)$ where $t$ is the number of occurrences of $x_i$ in $W$. We see that these two number must coincide. Let $i$ be such that $x_i^*$ occurs in $W$. Let $Z_k'$ and $W'$ be the words obtained from $Z_k$ and $W$ by deleting all variables except $x_i$. Then $Z_k'$ is the word $x_i^{2^{k-i}}$ while $W'$ is a product consisting of $2^{k-i}$ factors of either $x_i$ or $x_i^*$.  The identity $Z_k'\bumpeq W'$ also holds in $\mathcal{L}^1$. Substituting $(0,0,(0))$ for $x_i$ we see that $W'$ must start and end with $x_i$, otherwise the values of $Z_k'$ and $W'$ would not coincide. Moreover, the value of $Z_k'$ under the latter substitution is $$(0,0,(0))^{2^{k-i}}=(0,(2^{k-i}-1)\cdot((0)),(0)).$$ On the other hand, since $W'$ starts and ends with $x_i$ and contains at least one occurrence of $x_i^*$ it follows that the value of $W'$ contains the factor $(0,0,(0))((0),0,0)$. Consequently, the middle entry of the value of $W'$ contains the summand $((0)+(0))$ which is not the case for the value of $Z_k'$. Altogether, the identity $Z_k\bumpeq W$ fails in $\mathcal{L}^1$ and therefore in $\mathfrak{L}_2$.

Since the canonical morphisms $\mathfrak{aTL}_n\twoheadrightarrow \mathfrak{aTL}^{\mathrm e}_n\twoheadrightarrow \mathfrak{Ann}_n$ respect $^\rho$ it follows from \cite[Theorem 2.13(2)]{adv2} that the monoids in item (6) do not satisfy any nontrivial involutory identity $Z_k\bumpeq W$ with respect to $^\rho$, except for $n=3$. That is, we still have to consider $\mathfrak{aTL}^{\mathrm e}_3=\mathfrak{aTL}_3$. We have seen already in Proposition \ref{aTL3} that this monoid does not satisfy any nontrivial semigroup identity $Z_k\bumpeq W$. Take the idempotents $\alpha,\beta$ defined in the proof of Proposition \ref{aTL3}. Then $\alpha^\rho=\beta$ and $\alpha\beta$ and $\beta\alpha$ are of infinite order. Let $S$ be the involutory semigroup generated by $\alpha$ and let $I$ be the ideal consisting of all elements of infinite order. This ideal is invariant under $^\rho$ whence the involution $^\rho$ is well defined on the Rees quotient $S/I$. The latter involutory semigroup has three elements $\{\alpha,\beta,0\}$ and is a twisted semilattice in the sense of \cite{adpv}. In the same way as in the proof of Theorem 3.1 in \cite{adpv} it now follows that $\mathfrak{aTL}^{\mathrm e}_3$ does not satisfy and nontrivial involutory identity $Z_k\bumpeq W$ with respect to $^\rho$.

So we are left with the cases (i) $\mathfrak{aTL}^{\mathrm d}_2$ and $\mathfrak{aTL}_2$ for both involutions and (ii) $\mathfrak{aTL}_n$ for even $n$ for the reflection $^\sigma$.  Take an even $n\ge 2$ and the partition $\alpha\colon \mathbb{Z}\times[n]\leadsto \mathbb{Z}\times[n]$

\begin{center}
\begin{tikzpicture}[scale=1.4]
\draw (-0.75,0.5) node {$\alpha=$};

\draw[fill] (1,1) circle [radius=2pt];
\draw[above] (1,1) node{$(0,1)$};
\draw[fill] (2,1) circle [radius=2pt];
\draw[above] (2,1) node{$(0,2)$};
\draw(3,1) node{$\dots$};
\draw(0,0.5) node{$\dots$};
\draw(1,1) .. controls (1.25,0.75) and (1.75,0.75) .. (2,1);

\draw[fill] (3,0) circle [radius=2pt];
\draw[below] (3,0) node{$(0,3)$};
\draw[fill] (2,0) circle [radius=2pt];
\draw[below] (2,0) node{$(0,2)$};
\draw(1,0) node{$\dots$};
\draw(4,0) node{$\dots$};
\draw(2,0) .. controls (2.25,0.25) and (2.75,0.25) .. (3,0);

\draw[fill] (5,0) circle [radius=2pt];
\draw[below] (5,0) node{$(0,n)$};
\draw[fill] (6,0) circle [radius=2pt];
\draw[below] (6,0) node{$(1,1)$};
\draw(7,0.5) node{$\dots$};
\draw(5,0) .. controls (5.25,0.25) and (5.75,0.25) .. (6,0);

\draw[fill] (5,1) circle [radius=2pt];
\draw[above] (5,1) node{$(0,n)$};
\draw[fill] (4,1) circle [radius=2pt];
\draw[above] (4,1) node{$(0,n-1)$};
\draw(6,1) node{$\dots$};
\draw(4,1) .. controls (4.25,0.75) and (4.75,0.75) .. (5,1);

\end{tikzpicture}
\end{center}
and consider $\alpha$ as an element of $\mathfrak{aTL}^{\mathrm e}_n$. First of all, $\alpha^\sigma=\alpha^\rho$, that is, $^\rho$ and $^\sigma$ agree on the involutory subsemigroup generated by $\alpha$, hence it is sufficient to treat just one involution, say $^\sigma$. Let $\bm\alpha:=(\alpha,0)\in \mathfrak{aTL}_n$; then also $\bm\alpha^\rho=\bm\alpha^\sigma$ which allows us to treat cases (i) and (ii) simultaneously.

In order to show that no involutory identity is satisfied, we proceed similarly as in the proof of \cite[Theorem 4.3]{auingeretal}. Suppose by contrast that $\mathfrak{aTL}_n$ satisfies some nontrivial involutory identity $Z_k\bumpeq W$. First we observe that each letter $x_i$ for $i=1,\dots,n$ occurs the same number of times in $Z_k$ and $W$. Indeed, take $\vp=\alpha\alpha^\sigma$ and $\bm\psi:=(\vp,1)$. Note that $\vp$ is idempotent in $\mathfrak{aTL}^{\mathrm e}_n$. Then $\bm\psi^t=(\vp,t)$ for all $t\ge 0$. We substitute $\bm\psi$ for $x_i$ and $1$ for all other letters. The value of $Z_k$ under this substitution is $\bm\psi^{2^{k-i}}=(\vp,2^{k-i})$. Since $\bm\psi^\sigma=\bm\psi$, the value of $W$ is $\bm\psi^t$ where $t$ is the number of occurrences of $x_i$ in $W$. Since $Z_k\bumpeq W$ holds in $\mathfrak{aTL}_n$ these values must coincide, hence $t=2^{k-i}$. Similarly one can see that the only letters occurring in $W$ are those of $Z_k$, namely $x_1,\dots,x_k$.

By the plain semigroup case already treated in Proposition \ref{aTL_2} and Theorem \ref{thm:main annular}, we may assume that $W$ contains a `starred' letter $x_i^*$. We substitute $\bm\alpha$ for $x_i$ and $1$ for all other letters. Then the value of $Z_n$ is $\bm\alpha^{2^{k-i}}=(\alpha,2^{k-i}-1)$ since $(\alpha,0)^t=(\alpha,t-1)$ for all $t\ge 1$. Moreover, any product $\beta_1\cdots\beta_t$ in $\mathfrak{aTL}^{\mathrm e}_n$ of elements $\beta_i\in\{\alpha,\alpha^\sigma\}$ coincides with $\beta_1\beta_t$.  It follows that the value of $W$ is a product $p$  containing $2^{k-i}$ factors which are either $\bm\alpha$ or $\bm\alpha^\sigma$ such that the first and last factor is $\bm\alpha$, and which contains at least one occurrence of $\bm\alpha^\sigma$. Also, $(\alpha^\sigma,0)^t=(\alpha^\sigma,t-1)$ for all $t\ge 1$, while $(\alpha,0)(\alpha^\sigma,0)=(\alpha\alpha^\sigma,0)$ and $(\alpha^\sigma,0)(\alpha,0)=(\alpha^\sigma\alpha,0)$. Consequently, that value of $W$ is $(\alpha,2^{k-i}-1-s)$ where $s$ is the number of occurrences in $p$ of factors of the form $\bm\alpha\bm\alpha^\sigma$ and $\bm\alpha^\sigma\bm\alpha$. Since $p$ contains elements $\bm\alpha$ as well as $\bm\alpha^\sigma$ it follows that $s>0$ so that the values of $Z_k$ and $W$ cannot coincide. Altogether $\mathfrak{aTL}_n$ does not satisfy any involutory identity $Z_k\bumpeq W$ for both involutions $^\rho$ and $^\sigma$. Since the canonical morphism $\mathfrak{aTL}^{\mathrm d}_2\twoheadrightarrow \mathfrak{aTL}_2$ respects both involutions the claim follows also for this case.
\end{proof}

Note that the involution $^\sigma$ is an inverse unary operation on every $\mathfrak{aTL}^{\mathrm e}_n$ which coincides with $\mathfrak{aTL}_n$ in case $n$ is odd, hence we cannot apply our theory to these cases. Using Lemma~\ref{lem:characterization2}, it can be verified that $\mathfrak{aTL}^{\mathrm e}_2$ satisfies the semigroup identity $Z_4\bumpeq x_1x_2x_1x_3x_2x_1^2x_4Z_3$, so we cannot apply the theory in the involutory case either.

\subsubsection{Four nontrivial cases of finite basedness in the involutory case}\label{subsec:finite basedness reflection}
We start with the case $\mathfrak{L}_2^\flat$ and $\overline{\mathfrak{L}_2^\flat}$ endowed with rotation. To this end we consider the monoid $\mathcal{F}(\mathcal{S})$ for $\mathcal{S}=\mathbb{N}_0$ or $\mathcal{S}=\mathbb{Z}$ and we let $\mathcal{F}(\mathcal{S})$ be endowed with the involution $^\dag$ defined by
\begin{equation}\label{eq:rotation}
	[p,q,r]^\dag=[r,q,p]\mbox{ and }(a,b,c)^\dag=(c,b,a).
\end{equation}
Denote this involutory monoid $\mathcal{F}(\mathcal{S})^\dag$. We proceed similarly as in Subsection~\ref{finite basedness} but we have to adjust to the involutory case. We need to distinguish between \textsl{letters} and \textsl{variables}. The set of variables is just $X$, while the set of letters is $X\cup X^*$; every variable $x$ may occur in one of its two incarnations, namely as $x$ itself or as its starred version $x^*$. In particular, given words $u,v\in \mathcal{I}(X)$, then the identity $u\bumpeq v$ may be \emph{variable-balanced} or even \textsl{letter-balanced} (with the obvious meaning).  We start by collecting some properties  an identity $u\bumpeq v$ that holds in $\mathcal{F}(\mathcal{S})^\dag$  necessarily enjoys.
\begin{Lemma}\label{lem:letter-balanced}
	Every identity $u\bumpeq v$ that holds in $\mathcal{F}(\mathcal{S})^\dag$ is letter-balanced.
\end{Lemma}
\begin{proof}
	Let $x\in X$ and consider the substitution $\varphi\colon X\to \mathcal{F}(\mathcal{S})^\star$, defined by $x\mapsto [1,0,0]$ and $y\mapsto [0,0,0]$ for all $y\ne x$. Then $u\varphi=[a,0,b]$ where $a$ is the number of occurrences of the letter $x$ while $b$ is the number of occurrences of the letter $x^*$ in $u$. Since $u\varphi=v\varphi$ the claim follows.
\end{proof}

We introduce some notation: for a word $u\in \mathcal{I}(X)$ and $x\in X$ let $\lambda_xu$ be the incarnation of the leftmost occurrence of the variable $x$ in $u$; dually, $u\rho_x$ is the incarnation of the rightmost occurrence of the variable $x$ in $u$. Note that both, $\lambda_xu$ and $u\rho_x$, are letters (either $x$ or $x^*$). The next lemma shows  that these incarnations coincide for words $u$ and $v$ for which the identity $u\bumpeq v$ holds in $\mathcal{F}(\mathcal{S})^\dag$.
\begin{Lemma}\label{lem:first and last occurrence}
	If the identity $u\bumpeq v$ holds in $\mathcal{F}(\mathcal{S})^\dag$ then $\lambda_xu=\lambda_xv$ and $u\rho_x=v\rho_x$ for each $x\in X$.
\end{Lemma}
\begin{proof} Let $x\in X$ and substitute $(1,0,0)$ for $x$ and $[0,0,0]$ for all $y\ne x$. Then, for this substitution $\varphi$,
	\begin{align*}
		&u\varphi = (1,\ldotp,\ldotp)\Leftrightarrow \lambda_xu=x\mbox{ and  }
		u\varphi=(0,\ldotp,\ldotp) \Leftrightarrow \lambda_xu=x^*,\\
		&u\varphi=(\ldotp,\ldotp,1)\Leftrightarrow u\rho_x=x^*\mbox{ and }
		u\varphi=(\ldotp,\ldotp,0)\Leftrightarrow u\rho_x=x.
		\end{align*} 
	The claim follows immediately.
\end{proof}
Let $u\in \mathcal{I}(X)$ and $x\in X$; as in the plain semigroup case, let $L(u,x)$ and $R(x,u)$ be the longest prefix and suffix, respectively, of $u$ not containing the \textsl{variable} $x$ (note that the incarnation of the first/last occurrence of the variable $x$ may be the letter $x^*$).
\begin{Lemma} Let $u,v\in \mathcal{I}(X)$ and $x\in X$.
	If $\mathcal{F}(\mathcal{S})^\dag$ satisfies the identity $u\bumpeq v$ then the left section $L(u,x)\bumpeq L(v,x)$ and the right section $R(x,u)\bumpeq R(x,v)$ are variable-balanced.
\end{Lemma}
\begin{proof}
	Let $x,y\in X$, $x\ne y$. Consider the substitution $\varphi$ defined by $y\mapsto[0,1,0]$, $x\mapsto(0,0,0)$ and $z\mapsto[0,0,0]$ for all $z\ne x,y$. Then $u\varphi=(a,\ldotp,b)$ where $a$ is the number of occurrences of the variable $y$ in $L(u,x)$ and $b$ is the number of the occurrences of the variable $y$ in $R(x,u)$. The claim follows immediately.
\end{proof}
We now proceed similarly as in Subsection~\ref{finite basedness}: for each variable $x$ in a word $u$ consider its \emph{extreme occurrences}; this leads to the \emph{word of extreme occurrences of} $u$, to the \emph{extreme representation of a word} $u$ and to the \emph{interior block sequence of} $u$. The following statements are proved similarly as Lemmas~\ref{lem:key1} and~\ref{lem:key2}:
\begin{Lemma}\label{lem:keyinvolutory}
If an identity $u\bumpeq v$ holds in $\mathcal{F}(\mathcal{S})^\dag$ then
\begin{enumerate}
	\item the words of extreme occurrences of variables of $u$ and $v$ coincide;
	\item if $u_1,\dots,u_n$ and $v_1,\dots,v_n$ denote the interior block sequences of $u$ and $v$, respectively, then the identities $u_i\bumpeq v_i$ are variable-balanced.
\end{enumerate}
\end{Lemma}
\begin{Def}\label{def:ij-identites}
\rm For $i,j=0,1,2,3$ denote by~\eqref{finite1}$_{ij}$ the identity obtained from the identity~\eqref{finite1} by replacing (on both sides) the $i$th $x$ and the $j$th $y$ by $x^*$ and $y^*$, respectively. For $i=0$ or $j=0$ this means that no $x$ and/or $y$ is replaced. Define~\eqref{finite2}$_{ij}$ analogously.  
\end{Def}
Note that~\eqref{finite1}$_{ij}$ denotes a set of sixteen identities and so does~\eqref{finite2}$_{ij}$.
\begin{Lemma}
	The identities~\eqref{finite1}$_{ij}$ and~\eqref{finite2}$_{ij}$ hold in $\mathcal{F}(\mathcal{S})^\dag$ for all $i,j=0,1,2,3$.
\end{Lemma}
\begin{proof}
Let $\varphi$ be a substitution.  If all variables are substituted with elements of $\mathcal{S}^3$ then the claim follows from commutativity of $\mathcal{S}^3$. If some variable is substituted with an element of $\mathcal{M}(\mathcal{S},\mathcal{S},\mathcal{S};P)$ then the result follows from the facts that (i) $z\varphi$ and $z^*\varphi$ have the same weight and hence contribute equally to the middle entry of the result and (ii) that for $z\varphi\in\mathcal{S}^3$ the elements $z\varphi$ and $z^*\varphi$ act identically on elements of $\mathcal{M}(\mathcal{S},\mathcal{S},\mathcal{S};P)$.
\end{proof}
Altogether,~\eqref{finite1}$_{ij}$ and~\eqref{finite2}$_{ij}$ represent 32 identities. Note that the substitutions $x\mapsto x^*$ and $y\mapsto y^*$ lead to 32 further equivalent identities. We need two more identities. These are
\begin{gather}
	xt_1xt_2x^*t_3x\bumpeq xt_1x^*t_2xt_3x, \label{eq:ident_rotat1}\\
	xt_1xt_2x^*t_3x^*\bumpeq xt_1x^*t_2xt_3x^*. \label{eq:ident_rotat2}
\end{gather}
Here the substitution $x\mapsto x^*$ leads to two further equivalent identities.
\begin{Lemma}
	The identities~\eqref{eq:ident_rotat1} and~\eqref{eq:ident_rotat2} hold in $\mathcal{F}(\mathcal{S})^\dag$.
\end{Lemma}
\begin{proof}
	Let $\varphi$ be a substitution. If all variables are substituted with  elements of $\mathcal{S}^3$ then this follows from commutativity of $\mathcal{S}^3$. If at least one variable is substituted with an element of $\mathcal{M}(\mathcal{S},\mathcal{S},\mathcal{S};P)$ then the values of the inner occurrences of $x$ and $x^*$ have the same weight and hence contribute equally to the middle entry of the result.
\end{proof}

For a word $u$ and a variable $x$ let $u[x]$ be the word in $\{x,x^*\}$ obtained from $u$ by deletion of all (occurrences of) variables $y\ne x$. Suppose that $(\lambda_xu,u\rho_x)=(x,x)$ or $(\lambda_xu,u\rho_x)=(x,x^*)$; then $u$ admits a word $u'$, equivalent modulo the identities~\eqref{eq:ident_rotat1} and~\eqref{eq:ident_rotat2}, so that  $u'[x]$ is of one of the following forms:
\begin{equation}\label{eq:canonical_rotat}
x^a\mbox{ or }x^a(x^*)^b\mbox{ or }x^{a-1}(x^*)^bx
\end{equation}
where $a$ is the number of occurrences of (the incarnation of) $x$ and $b$ that of $x^*$ (in the first case $b=0$). A word $u$ in which for \textsl{every} variable $x$ the word $u[x]$ is of one of the forms in~\eqref{eq:canonical_rotat} is said to be in \emph{reduced form}. 
\begin{Lemma}\label{lem:reduced_rotate}
	Every identity $u\bumpeq v$ that holds in $\mathcal{F}(\mathcal{S})^\dag$ admits an identity $\til u\bumpeq \til v$ that is equivalent modulo the identities~\eqref{eq:ident_rotat1} and~\eqref{eq:ident_rotat2} \textup(hence also holds in $\mathcal{F}(\mathcal{S})^\dag$\textup) and such that both words, $\til u$ and $\til v$, are reduced.
\end{Lemma}
\begin{proof} Let $u\bumpeq v$ be given and $x$ be a variable. If necessary, apply the transformation $x\mapsto x^*$ to $u$ so that for the transformed word $u'$ one of the equalities $(\lambda_xu',u'\rho_x)=(x,x)$ or $(\lambda_xu',u'\rho_x)=(x,x^*)$ holds. Since $(\lambda_xu,u\rho_x)=(\lambda_xv,v\rho_x)$ at the beginning, we apply the transformation $x\mapsto x^*$ to $u$ if and only if we apply it to $v$. Doing so for every variable $x$, we end up with an equivalent identity $u'\bumpeq v'$ for which the pairs $(\lambda_xu',u'\rho_x)$ and $(\lambda_xv',v'\rho_x)$ are equal and are of the form $(x,x)$ or $(x,x^*)$. Now we may apply, for every variable $x$, to each of $u'$ and $v'$ transformations induced by the identities~\eqref{eq:ident_rotat1} and~\eqref{eq:ident_rotat2} to get words $u''$ and $v''$ for which $u''[x]$ and $v''[x]$ have the required form. The identity $u''\bumpeq v''$ is equivalent modulo~\eqref{eq:ident_rotat1} and~\eqref{eq:ident_rotat2} to the original $u\bumpeq v$ (and hence also holds in $\mathcal{F}(\mathcal{S})^\dag$).
\end{proof}
\begin{Lemma}\label{lem:interior letter bal}	If for two reduced words $u$ and $v$ the identity $u\bumpeq v$ holds in $\mathcal{F}(\mathcal{S})^\dag$ then all identities $u_i\bumpeq v_i$ in item (2) of Lemma~\ref{lem:keyinvolutory} are letter-balanced.
\end{Lemma}
\begin{proof}
Let $x\in X$. By Lemma~\ref{lem:letter-balanced} $u$ and $v$ contain the same number of occurrences of the letters $x$ and $x^*$, respectively. Since both words are reduced the words $u[x]$ and $v[x]$ coincide. It follows that $u_i[x]=v_i[x]$ for every $i$. Since the identities $u_i\bumpeq v_i$ are variable-balanced by item (2) of  Lemma~\ref{lem:keyinvolutory}, the last assertion implies that these identities are even letter-balanced.	
\end{proof}
 We are ready for the first result in the present subsection.
\begin{Thm}
	The identities~\eqref{finite1}$_{ij}$, \eqref{finite2}$_{ij}$ (for $i,j=0,1,2,3$) together with the identities~\eqref{eq:ident_rotat1} and~\eqref{eq:ident_rotat2} form a basis for the monoid identities of  $\mathfrak{L}_2^\flat$ and $\overline{\mathfrak{L}_2^\flat}$, both equipped with rotation $\rho$. A basis of semigroup identities is obtained by deleting every possible choice of a subset of the set of letters $\{t_1,t_2,t_3,t_4\}$ in \eqref{finite1}$_{ij}$ and \eqref{finite2}$_{ij}$ and every possible choice of a subset of the set of letters $\{t_1,t_2,t_3\}$ in \eqref{eq:ident_rotat1} and \eqref{eq:ident_rotat2}.
\end{Thm}
\begin{proof}
	Let $u\bumpeq v$ be an identity holding in $\mathcal{F}(\mathcal{S})^\dag$. According to Lemma~\ref{lem:reduced_rotate} we may assume that both words $u$ and $v$ are reduced. By Lemma~\ref{lem:keyinvolutory}, $u$ and $v$ have the same words of extreme occurrences, and the identities $u_i\bumpeq v_i$ coming from the interior block sequences are variable-balanced. By Lemma~\ref{lem:interior letter bal} these identities are even letter-balanced. Now, similarly as in the plain case, all these identities $u_i\bumpeq v_i$ are equivalent modulo the identities~\eqref{finite1}$_{ij}$ and~\eqref{finite2}$_{ij}$.
\end{proof}

We proceed with the case of reflection $\sigma$ and 
the monoids $2\mathfrak{Cob}^\circ_1$, $\overline{2\mathfrak{Cob}^\circ_1}$, $\mathfrak{L}^\flat_2$, and $\overline{\mathfrak{L}^\flat_2}$. We will use the representation of these monoids as $\mathcal{E}(\mathbb{N}_0)$, $\mathcal{F}(\mathbb{N}_0)$, $\mathcal{E}(\mathbb{Z})$, and $\mathcal{F}(\mathbb{Z})$ respectively; see Lemmas~\ref{Lem:model_2Cob^0_1} and~\ref{lem:model_L^flat2}.
Let $\mathcal{S}$ be either $\langle\mathbb{N}_0;+,0\rangle$ or $\langle\mathbb{Z};+,0\rangle$. The operation ${}^\star$ on the monoid $\mathcal{E}(\mathcal{S}) = \mathcal{S}\sqcup\, \mathcal{S}\times \mathcal{S}$ is defined by:
\begin{itemize}
  \item $s^\star:=s$ for $s\in\mathcal{S}$,
  \item $(a,b)^\star:=(b,a)$ for $(a,b)\in \mathcal{S}\times \mathcal{S}$.
\end{itemize}
This operation is readily seen to be an involution, and for $\mathcal{E}(\mathbb{N}_0)$ and $\mathcal{E}(\mathbb{Z})$, it corresponds to the involution $\sigma=\rho$ on $2\mathfrak{Cob}^\circ_1$ and, respectively, $\overline{2\mathfrak{Cob}^\circ_1}$. We denote this involutory monoid $\mathcal{E}(\mathcal{S})^\star$.

A similar involution can be defined the monoid $\mathcal{F}(\mathcal{S}) = \mathcal{S}^3\sqcup \mathcal{M}(\mathcal{S},\mathcal{S},\mathcal{S};P)$ by letting:
\begin{itemize}
  \item $[p,q,r]^\star:=[p,q,r]$ for $[p,q,r]\in\mathcal{S}^3$,
  \item $(a,b,c)^\star:=(c,b,a)$ for $(a,b,c)\in\mathcal{M}(\mathcal{S},\mathcal{S},\mathcal{S};P)$.
\end{itemize}
For $\mathcal{F}(\mathbb{N}_0)$ and $\mathcal{F}(\mathbb{Z})$, this involution corresponds to the reflection $\sigma$ on $\mathfrak{L}^\flat_2$ and $\overline{\mathfrak{L}^\flat_2}$, respectively. Notice that $\sigma\ne\rho$ in this case. We denote the corresponding involutory monoid $\mathcal{F}(\mathcal{S})^\star$.

Clearly, the surjective monoid morphism $\mathcal{F}(\mathcal{S})\to \mathcal{E}(\mathcal{S})$ defined by
\[[p,q,r]\mapsto q,\ (a,b,c)\mapsto (a,c)\]
is compatible with the involution ${}^\star$, hence is a surjective morphism $\mathcal{F}(\mathcal{S})^\star\twoheadrightarrow \mathcal{E}(\mathcal{S})^\star$. 
For the equational theories of the four involved involutory monoids there are the following containment relations (containment is from `below to above'). 
\begin{center}
 \begin{tikzpicture}[scale=2,
  every node/.style={circle, inner sep=2pt}]

\node (e1) at (0,1) {$\Id \mathcal{E}(\mathbb{N}_0)^\star$};
\node (e2) at (1,0) {$\Id \mathcal{E}(\mathbb{Z})^\star$};
\node (e3) at (-1,0) {$\Id \mathcal{F}(\mathbb{N}_0)^\star$};
\node (e4) at (0,-1) {$\Id \mathcal{F}(\mathbb{Z})^\star$};

\draw (e1) -- (e2);
\draw (e1) -- (e3);
\draw (e2) -- (e4);
\draw (e3) -- (e4);

\end{tikzpicture}   
\end{center}
In the following we shall derive certain combinatorial invariants of the members of the top $\Id \mathcal{E}(\mathbb{N}_0)^\star$ and then show that all these identities are consequences of a finite set of identities all of which belong to the bottom set $\Id \mathcal{F}(\mathbb{Z})^*$. This gives finite basedness and equational equivalence of the four involutory monoids in question.

We proceed similarly as  in the case of the rotation. We start by listing properties of members of   $\Id\mathcal{E}(\mathcal{S})^\star$.

\begin{Lemma}\label{lem:keyreflection}
    Let $u,v\in \mathcal{I}(X)$; if the identity $u\bumpeq v$ holds in  $\mathcal{E}(\mathcal{S})^\star$ then
    \begin{enumerate}
        \item $u\bumpeq v$ is variable-balanced,
        \item $\lambda_xu=\lambda_xv$ and $u\rho_x=v\rho_x$ for every $x\in X$,
        \item $L(u,x)\bumpeq L(v,x)$ and $R(x,u)\bumpeq R(x,v)$ are variable-balanced for every $x\in X$,
        \item the words of extreme occurrences \textup(of variables\textup) of $u$ and $v$ coincide,
        \item all identities $u_i\bumpeq v_i$ coming form the interior block sequences of $u$ and $v$ are variable-balanced.
     \end{enumerate}
\end{Lemma}
\begin{proof}
    The proofs are very similar to the case of rotation.  For (1) consider, for every fixed $x\in X$, the substitution $x\mapsto 1$ and $y\mapsto 0$ for  $y\ne x$. For (2) consider, for every fixed $x$, the substitution $x\mapsto(1,0)$ and $y\mapsto 0$ for $y\ne x$. For (3) consider the substitution $x\mapsto (0,0)$ and then, for every fixed $y\ne x$, $y\mapsto 1$ and for all $z\ne x,y$, $z\mapsto 0$. (4) and (5) are proved similarly as Lemmas~\ref{lem:key1} and~\ref{lem:key2}. 
\end{proof}
Now we need some more identities.
\begin{Lemma}\label{lem:xx*x=xxx}
The involutory monoid $\mathcal{F}(\mathcal{S})^\star$ satisfies the identities
\begin{gather}
\label{eq:xx*x=xxx}  xt_1x^* t_2x\bumpeq xt_1xt_2x,\\
\label{eq:xx*x*=xxx*} xt_1x^* t_2x^*\bumpeq xt_1xt_2x^*,\\
\label{eq:x*xx*=x*x*x*} x^* t_1xt_2x^*\bumpeq x^* t_1x^* t_2x^*,\\
\label{eq:x*x*x=x*xx} x^* t_1x t_2x\bumpeq x^* t_1x^* t_2x.
\end{gather}
\end{Lemma}

\begin{proof}
Consider any substitution $\varphi\colon\{x,t_1,t_2\}\to\mathcal{F}(\mathcal{S})$. If $x\varphi\in\mathcal{S}^3$, then $x^*\varphi=(x\varphi)^\star=x\varphi$, whence
\[
(xt_1x^* t_2x)\varphi=(xt_1xt_2x)\varphi=(xt_1x^* t_2x^*)\varphi=(xt_1xt_2x^*)\varphi.
\]
If $x\varphi=(a,b,c)\in\mathcal{M}(\mathcal{S},\mathcal{S},\mathcal{S};P)$, then the equalities $(xt_1x^* t_2x)\varphi=(xt_1xt_2x)\varphi$ and $(xt_1x^* t_2x^*)\varphi=(xt_1xt_2x^*)\varphi$ follow from Lemma~\ref{lem:fact} since the weight $a+b+c$ of $x\varphi$ is equal to the weight $c+b+a$ of $x^*\varphi=(x\varphi)^\star=(c,b,a)$.

Substituting $x^*$ for $x$ in \eqref{eq:xx*x=xxx} and \eqref{eq:xx*x*=xxx*} and using the identity $(x^*)^*\bumpeq x$, we obtain the identities \eqref{eq:x*xx*=x*x*x*} and \eqref{eq:x*x*x=x*xx}, respectively.
\end{proof}
The last sentence in the proof essentially says that~\eqref{eq:xx*x=xxx} is equivalent to~\eqref{eq:x*xx*=x*x*x*} and~\eqref{eq:xx*x*=xxx*} is eqivalent to~\eqref{eq:x*x*x=x*xx}.
For a word $w\in\mathcal{I}(X)$ and $x\in X$ let again $w[x]$ be the word obtained from $w$ by deleting all occurrences of all variables $y\ne x$.

It should be clear that the identities \eqref{eq:xx*x=xxx}--\eqref{eq:x*x*x=x*xx} allow one to transform $w$ into a word $w'$ such that $w'[x]$ is a word of exactly one of the following types:
\begin{equation*}\label{eq:4types}
x^k,\ (x^*)^k,\  x^{k-1}x^*,\ (x^*)^{k-1}x,\ \text{ where } k:=\occ_x(w)+\occ_{x^*}(w).
\end{equation*}
(Recall that for a letter $y\in X\cup X^*$, we denote by $\occ_y(w)$ the number of occurrences of $y$ in $w$.) Furthermore, application of the substitution $x\mapsto x^*$ allows one to transform the word $w'$ to a word $w''$ such that $w''[x]=x^k$ or $w''[x]=x^{k-1}x^*$. Call a word $w$ \emph{reduced} (\emph{with respect to reflection}) (\emph{reduced} for short) if for every variable $x$, $w[x]=x^k$ or $w[x]=x^{k-1}x^*$. Observe that in a reduced word every starred letter occurs in an extreme position. In particular, the members $w_i$ of the interior block sequence of $w$ do not contain starred letters. We are ready for the next result. Similarly as in the case of rotation, every identity $u\bumpeq v$ that holds in $\mathcal{F}(\mathcal{S})^\star$ or $\mathcal{E}(\mathcal{S})^\star$ admits an equivalent identity $u'\bumpeq v'$ both sides of which are reduced (with respect to reflection).
\begin{Thm}\label{Thm:2Cob_1 with rho} 
The involutory monoids $2\mathfrak{Cob}^\circ_1$, $\overline{2\mathfrak{Cob}^\circ_1}$, $\mathfrak{L}^\flat_2$, and $\overline{\mathfrak{L}^\flat_2}$ with involution $\sigma$ are finitely based. The identities \eqref{finite1}$_{ij}$, \eqref{finite2}$_{ij}$ for $i,j=0,1,3$, \eqref{eq:xx*x=xxx}, and \eqref{eq:xx*x*=xxx*} form a basis for the involutory monoid identities of each of these monoids. A basis of semigroup identities is obtained by deleting every possible choice of a subset of the set of letters $\{t_1,t_2,t_3,t_4\}$ in \eqref{finite1}$_{ij}$ and \eqref{finite2}$_{ij}$ and every possible choice of a subset of the set of letters $\{t_1,t_2\}$ in \eqref{eq:xx*x=xxx} and \eqref{eq:xx*x*=xxx*}.
\end{Thm}
\begin{proof}
That all identities mentioned hold in $\mathcal{F}(\mathcal{S})^\star$ follows from Lemma~\ref{lem:xx*x=xxx} and otherwise similarly as in the case of the rotation. Now suppose $u\bumpeq v$ is an identity which holds in $\mathcal{E}(\mathcal{S})^\star$.  By Lemma~\ref{lem:keyreflection} item (4), both words have identical words of extreme occurrences; by  item (5), the identities $u_k\bumpeq v_k$ coming from the interior block sequences all are variable-balanced. Similarly as in the case of rotation and by  use of the identities~\eqref{eq:xx*x=xxx} and~\eqref{eq:xx*x*=xxx*} we may assume that both words $u$ and $v$  are reduced. It follows that the members of both interior block sequences $u_1,\dots, u_n$ and  $v_1,\dots, v_n$ do not contain starred letters and hence are even letter-balanced. If follows that all identities $u_k\bumpeq v_k$ are consequences of the identities~\eqref{finite1}$_{ij}$ and~\eqref{finite2}$_{ij}$ for $i,j=0,1,3$ (here the cases $i,j=2$ are not required since the members $u_k$ and $v_k$ of the interior block sequences do not contain starred letters).  
\end{proof}

Finally, let us consider the monoid $\mathfrak{aTL}_2^\mathrm{e}$ endowed with rotation. We look at the model $\mathcal{E}_2(\mathbb{Z})$ introduced in Lemma~\ref{lem: aTL_2^e} and define an involution $^\star$ on $\mathcal{E}_2(\mathbb{Z})$ by setting
\begin{itemize}
    \item $s^\star =s$ for $s\in \mathbb{Z}$,
    \item $^\star$ swaps $(1,1)$ with $(2,2)$ and fixes $(1,2)$ and $(2,1)$.
\end{itemize}
Denote the resulting involutory monoid $\mathcal{E}_2(\mathbb{Z})^\star$. It is readily checked that $\mathcal{E}_2(\mathbb{Z})^\star$ is isomorphic to $\mathfrak{aTL}_2^\mathrm e$ endowed with rotation $\rho$. The following is proved similarly as the non-involutory case in combination with the case treated in Theorem~\ref{Thm:2Cob_1 with rho}.
\begin{Thm} The involutory monoid $\mathfrak{aTL}_2^\mathrm e$ equipped with rotation $^\rho$ is finitely based. 
The identities \eqref{finite1}$_{ij}$, \eqref{finite2}$_{ij}$ for $i,j=0,1,3$, \eqref{eq:xx*x=xxx},  \eqref{eq:xx*x*=xxx*} and $x^3yx\bumpeq xyx^3$ form a basis for the involutory monoid identities of this involutory monoid. A basis of semigroup identities is obtained by deleting every possible choice of a subset of the set of letters $\{t_1,t_2,t_3,t_4\}$ in \eqref{finite1}$_{ij}$ and \eqref{finite2}$_{ij}$ and every possible choice of a subset of the set of letters $\{t_1,t_2\}$ in \eqref{eq:xx*x=xxx} and \eqref{eq:xx*x*=xxx*}.
\end{Thm}

\begin{Problem}
 Are the monoids $\mathfrak{aTL}^{\mathrm e}_n$ ($n\ge 2$) finitely based as involutory semigroups with respect to $^\sigma$? For these monoids, $^\sigma$ is an inverse unary operation, and the method used in the present paper does not apply.
\end{Problem}

\begin{Rmk}
The reader may notice that for all categories $\mathfrak{C}$ we considered, the local monoids of $\mathfrak{C}$ and those of its von-Neumann-regular extension $\overline{\mathfrak{C}}$ behave equally with respect to the property of being finitely or non-finitely based. Moreover, for some $\mathfrak{C}$, the equality $\Id\mathfrak{C}_n=\Id\overline{\mathfrak{C}_n}$ can be established. In this paper, it was explicitly shown only for a few special cases ($2\mathfrak{Cob}^\circ_1$ vs. $\overline{2\mathfrak{Cob}^\circ_1}$, $\mathfrak{L}^\flat_2$ vs.  $\overline{\mathfrak{L}^\flat_2}$, and $\mathfrak{aTL}^{\mathrm e}_2$ vs. $\overline{\mathfrak{aTL}^{\mathrm e}_2}$); in addition, it holds whenever the monoid $\overline{\mathfrak{C}_n}$ can be generated by $\mathfrak{C}_n$ and some central elements\footnote{An element $c$ of a monoid $\mathcal{M}$ is called \emph{central} if $cs=sc$ for all $s\in \mathcal{M}$.} of $\overline{\mathfrak{C}_n}$; see \cite[Lemma 1]{Vo23}. This applies for instance to the pair $\mathfrak{P}^{\mathrm d}_n$, $\overline{\mathfrak{P}^{\mathrm d}_n}$. We do not know, however, whether the  equality $\Id\mathfrak{C}_n=\Id\overline{\mathfrak{C}_n}$ persists for categories with more complicated constructions of regular extension. 
\end{Rmk}

\begin{Rmk} Several (but not all) results of the present subsection can alternatively be obtained from the `plain' semigroup results of Subsections \ref{subsec:nfb-cobord} and \ref{subsec:nfb-annular} by use of Theorem 4 of the paper by Lee \cite{lee}.
\end{Rmk}

\section*{Appendix}

Table 1 below summarizes the status of the  Finite Basis Problem for the monoids we examined. The acronyms NFB and FB stand for `non-finitely based' and `finitely based', respectively; the FB entries are highlighted to make them stand out among the predominant NFB entries. The table omits monoids that are finitely based because they are commutative (such monoids arise in some series when $n=1$, and for $n=0$ in all series).

\begin{table}[p]
\caption{Status of the Finite Basis Problem for cobordism, partition,  topological annular, and Temperley--Lieb monoids}
\centering
\begin{tabular}{|c|c|c|c|}
\hline
Monoid
& Plain
& With involution $\rho$
& With involution $\sigma$ \\
\hline
\rule{0pt}{15pt} $2\mathfrak{Cob}_n$, $\overline{2\mathfrak{Cob}_n}$ & NFB & NFB & NFB \\
\rule[-6pt]{0pt}{6pt} with $n \ge 1$ & (Thm.~5.5(1)) & (Thm.~5.29(1)) & (Thm.~5.29(1))\\
\hline
\rule{0pt}{15pt} $2\mathfrak{Cob}^\circ_n$, $\overline{2\mathfrak{Cob}^\circ_n}$ & NFB & NFB & NFB \\
\rule[-6pt]{0pt}{6pt} with $n \ge 2$ & (Thm.~5.5(2)) & (Thm.~5.29(2)) & (Thm~5.29(2))\\
\hline
\rule{0pt}{15pt} $2\mathfrak{Cob}^\circ_1$, $\overline{2\mathfrak{Cob}^\circ_1}$ & {\blue FB} (Thm. 5.6) &  {\blue FB} (Thm. 5.43) & {\blue FB} (Thm. 5.43) \\[4pt]
\hline
\rule{0pt}{15pt} $\mathfrak{P}^{\mathrm d}_n$, $\overline{\mathfrak{P}^{\mathrm d}_n}$ & NFB & NFB & NFB \\
\rule[-6pt]{0pt}{6pt} with $n \ge 2$ & (Thm.~5.5(3)) & (Thm.~5.29(3)) & (Thm~5.29(3))\\
\hline
\rule{0pt}{15pt} $\mathfrak{tAnn}_n$, $\overline{\mathfrak{tAnn}_n}$ & & & \\
$\mathfrak{t^{\flat}Ann}_n$, $\overline{\mathfrak{t^{\flat}Ann}_n}$ & {\blue FB} (Thm. 5.10) &  {\blue FB} (Thm. 5.30) & {\blue FB} (Thm. 5.30)\\
\rule[-6pt]{0pt}{6pt} with even $n\ge 2$ & & &\\
\hline
\rule{0pt}{15pt} $\mathfrak{tAnn}_n$, $\overline{\mathfrak{tAnn}_n}$ & \raisebox{-.6\normalbaselineskip}{NFB} &  \raisebox{-.6\normalbaselineskip}{NFB} & \raisebox{-.6\normalbaselineskip}{NFB} \\[-2pt]
$\mathfrak{t^{\flat}Ann}_n$, $\overline{\mathfrak{t^{\flat}Ann}_n}$  & \raisebox{-.4\normalbaselineskip}{(Thm.~5.13(1))} & \raisebox{-.4\normalbaselineskip}{(Thm.~5.30(1))} & \raisebox{-.4\normalbaselineskip}{(Thm.~5.30(1))}\\[-4pt]
\rule[-6pt]{0pt}{6pt} with odd $n\ge 3$ &  &  & \\
\hline
\rule{0pt}{15pt} $\mathfrak{tAnn}^\circ_n$, $\overline{\mathfrak{tAnn}^\circ_n}$ & \raisebox{-.6\normalbaselineskip}{NFB} &  \raisebox{-.6\normalbaselineskip}{NFB} & \raisebox{-.6\normalbaselineskip}{NFB} \\[-2pt]
$\mathfrak{t^{\flat}Ann}^\circ_n$, $\overline{\mathfrak{t^{\flat}Ann}^\circ_n}$  & \raisebox{-.4\normalbaselineskip}{(Thm.~5.13(2))} & \raisebox{-.4\normalbaselineskip}{(Thm.~5.30(2))} & \raisebox{-.4\normalbaselineskip}{(Thm.~5.30(2))}\\[-4pt]
\rule[-6pt]{0pt}{6pt} with even $n\ge 2$ &  &  & \\
\hline
\rule{0pt}{15pt} $\mathfrak{L}_n$, $\overline{\mathfrak{L}_n}$ & NFB & NFB & NFB \\
 with $n\ge 2$ & (Thm. 5.13(3))  & (Thm.~5.30(3)) & (Thm.~5.30(3))\\
\hline
\rule{0pt}{15pt} $\mathfrak{L}^\flat_n$, $\overline{\mathfrak{L}^\flat_n}$ & NFB & NFB & NFB \\
\rule[-6pt]{0pt}{6pt} with $n\ge 3$ & (Thm. 5.13(4))  & (Thm.~5.30(4)) & (Thm.~5.30(4))\\
\hline
\rule{0pt}{15pt} $\mathfrak{L}^\flat_2$, $\overline{\mathfrak{L}^\flat_2}$ & {\blue FB} (Thm. 5.19) & {\blue FB} (Thm. 5.40) & {\blue FB} (Thm. 5.43) \\
\hline
\rule{0pt}{15pt} $\mathfrak{aTL}^{\mathrm d}_n$, $\overline{\mathfrak{aTL}^{\mathrm d}_n}$ & NFB & NFB & NFB \\
\rule[-6pt]{0pt}{6pt} with $n\ge 2$ & (Thm. 5.13(5))  & (Thm.~5.30(5) & (Thm.~5.30(5))\\
\hline
\rule{0pt}{15pt} $\mathfrak{aTL}_n$, $\overline{\mathfrak{aTL}_n}$ & NFB & NFB & NFB \\
\rule[-6pt]{0pt}{6pt} with even $n\ge 2$ & (Thm. 5.13(6))  & (Thm.~5.30(6)) & (Thm.~5.30(6))\\
\hline
\rule{0pt}{15pt} $\mathfrak{aTL}_n$, $\overline{\mathfrak{aTL}_n}$ & NFB & NFB & \raisebox{-.5\normalbaselineskip}{\textbf{?}} \\
\rule[-6pt]{0pt}{6pt} with odd $n\ge 3$ & (Thm. 5.13(7))  & (Thm.~5.30(7)) & \\
\hline
\rule{0pt}{15pt} $\mathfrak{aTL}^{\mathrm e}_n$, $\overline{\mathfrak{aTL}^{\mathrm e}_n}$ & NFB & NFB &  \raisebox{-.5\normalbaselineskip}{\textbf{?}} \\[-.1cm]
\rule[-6pt]{0pt}{6pt} with $n\ge 3$ & (Thm.~5.13(7) & (Thm.~5.30(7)) & \\
\hline
\rule[-6pt]{0pt}{21pt} $\mathfrak{aTL}^{\mathrm e}_2$, $\overline{\mathfrak{aTL}^{\mathrm e}_2}$ & {\blue FB} (Thm. 5.19) & {\blue FB} (Thm.~5.44) &  \textbf{?} \\
\hline
\end{tabular}
\end{table}

\subsection*{Acknowledgment} The authors are grateful to Nils Carqueville for introducing them to some secrets of the category of 2-cobordisms and for providing several references.

\end{document}